\documentclass[leqno,11pt]{amsart}
\usepackage{latexsym,amssymb,amscd,amsthm,amsmath}

\theoremstyle{plain}

\theoremstyle{definition}

\theoremstyle{remark}

\numberwithin{equation}{section}
\def\({\left( }
\def\){\right)}

\def\Cal{\mathcal }
\begin{document}
\nonstopmode
\hsize=5 true in
\vsize=7.286 true in
\font\b=cmr6

\title[Riemannian Submanifolds]{Riemannian Submanifolds: A survey}

\author[B.-Y. Chen]{ Bang-Yen Chen}

\thanks{{\it Riemannian Submanifolds}, Handbook of Differential Geometry, vol. I (2000), 187--418.} 

\maketitle

\noindent{\bf Contents}
\vskip.1in
\noindent Chapter 1. Introduction \dotfill 6

\noindent Chapter 2. Nash's embedding theorem and
some related results \dotfill 9

\hskip.0in 2.1. Cartan-Janet's theorem
\dotfill 10

\hskip.0in 2.2. Nash's embedding theorem
\dotfill 11

\hskip.0in 2.3. Isometric immersions with the
smallest possible codimension \dotfill 8

\hskip.0in 2.4. Isometric immersions with
prescribed Gaussian or Gauss-Kronecker

curvature \dotfill 12

\hskip.0in  2.5. Isometric immersions with
prescribed mean curvature\dotfill 13

\noindent Chapter 3. Fundamental theorems, basic
notions and results
\dotfill 14

\hskip.0in 3.1. Fundamental equations \dotfill
14

\hskip.0in 3.2. Fundamental theorems \dotfill
 15

\hskip.0in 3.3. Basic notions \dotfill 16

\hskip.0in 3.4. A general inequality \dotfill
17

\hskip.0in  3.5. Product immersions \dotfill
19

\hskip.0in  3.6. A relationship between
$k$-Ricci tensor and shape
operator\dotfill 20

\hskip.0in 3.7. Completeness of curvature
surfaces\dotfill 22

\noindent   Chapter 4. Rigidity and reduction
theorems
\dotfill 24

\hskip.0in 4.1. Rigidity \dotfill 24

\hskip.0in 4.2. A reduction theorem
 \dotfill 25

\noindent  Chapter 5. Minimal submanifolds \dotfill
26

\hskip.0in 5.1. First and second variational
formulas  \dotfill 27

\hskip.0in 5.2. Jacobi operator, index, nullity
and Killing nullity \dotfill 28

\hskip.0in 5.3. Minimal submanifolds in
Euclidean space \dotfill 28

\hskip.0in 5.4. Minimal submanifolds in
sphere \dotfill 41

\hskip.0in 5.5. Minimal submanifolds in
hyperbolic space \dotfill 49

\hskip.0in 5.6. Gauss map of minimal surfaces
\dotfill 50

\hskip.0in  5.7. Complete minimal
submanifolds in Euclidean space with finite
total 

\hskip .59in curvature\dotfill 54

\hskip.0in  5.8. Complete minimal
surfaces in $E^3$ lying between two
parallel 

\hskip .59in planes\dotfill 67

\hskip.0in 5.9. The geometry of Gauss
image\dotfill 67

\hskip.0in 5.10. Stability and index of
minimal submanifolds\dotfill 68

\noindent  Chapter 6. Submanifolds of finite type
\dotfill 77

\hskip.0in 6.1. Spectral resolution
\dotfill 77

\hskip.0in 6.2. Order and type of immersions
\dotfill 78

\hskip.0in 6.3. Equivariant submanifolds as
minimal submanifolds in their adjoint

\hskip.59in hyperquadrics\dotfill 80

\hskip.0in 6.4. Submanifolds of finite type
\dotfill 80

\noindent Chapter 7. Isometric immersions between
real space forms \dotfill 89

\noindent  Chapter 8. Parallel submanifolds
\dotfill 92

\hskip.0in 8.1. Parallel submanifolds in
Euclidean space \dotfill 92

\hskip.0in 8.2.  Parallel submanifolds in
spheres \dotfill 94

\hskip.0in 8.3.  Parallel submanifolds
in hyperbolic spaces \dotfill 94

\hskip.0in 8.4. Parallel submanifolds in
complex projective and complex
hyperbolic

\hskip .6in spaces \dotfill 94

\hskip.0in 8.5.  Parallel submanifolds in
quaternionic projective spaces \dotfill 95

\hskip.0in 8.6.  Parallel submanifolds in the
Cayley plane \dotfill 95

\noindent  Chapter 9. Standard immersions and
submanifolds with simple geodesics \dotfill 
96

\hskip.0in 9.1. Standard immersions \dotfill
96

\hskip.0in 9.2. Submanifolds with planar
geodesics \dotfill 96

\hskip.0in 9.3. Submanifolds with
pointwise planar normal
sections\dotfill 97

\hskip.0in  9.4. Submanifolds with
geodesic normal sections and helical

\hskip .59in  immersions\dotfill 98

\hskip.0in 9.5. Submanifolds whose
geodesics are generic $W$-curves
\dotfill 99

\hskip.0in 9.6. Symmetric spaces in
Euclidean space with simple geodesics
\dotfill 101

\noindent  Chapter 10. Hypersurfaces of real space
forms \dotfill 103

\hskip.0in 10.1. Einstein hypersurfaces
 \dotfill 103

\hskip.0in 10.2. Homogeneous hypersurfaces
\dotfill 103

\hskip.0in 10.3. Isoparametric hypersurfaces
\dotfill 104

\hskip.0in 10.4. Dupin hypersurfaces \dotfill
109

\hskip.0in 10.5. Hypersurfaces with constant
mean curvature \dotfill 117

\hskip.0in 10.6. Hypersurfaces with constant 
 higher order mean curvature \dotfill 123

\hskip.0in 10.7. Harmonic spaces and
Lichnerowicz conjecture  \dotfill 124

\noindent Chapter 11. Totally geodesic submanifolds
\dotfill 127

\hskip.0in 11.1. Cartan's theorem \dotfill
127

\hskip.0in 11.2. Totally geodesic submanifolds
of symmetric spaces \dotfill 127

\hskip.0in 11.3. Stability of totally geodesic
submanifolds \dotfill 133

\hskip.0in 11.4. Helgason's spheres \dotfill
136

\hskip.0in  11.5. Frankel's theorem\dotfill
1387

\noindent  Chapter 12. Totally umbilical
submanifolds \dotfill 139

\hskip.0in 12.1. Totally umbilical
submanifolds of real space forms \dotfill
139

\hskip.0in 12.2.  Totally umbilical
submanifolds of complex space forms \dotfill
140

\hskip.0in 12.3.  Totally umbilical
submanifolds of quaternionic space forms
\dotfill 140

\hskip.0in 12.4.  Totally umbilical
submanifolds of the Cayley plane \dotfill 
140

\hskip.0in 12.5.  Totally umbilical
submanifolds in complex quadric \dotfill 141

\hskip.0in 12.6.  Totally umbilical
submanifolds of locally symmetric spaces
\dotfill 141

\hskip.0in 12.7.  Extrinsic spheres in locally
symmetric spaces \dotfill 142

\hskip.0in 12.8.  Totally umbilical
hypersurfaces \dotfill 143

\noindent Chapter 13. Conformally flat submanifolds
\dotfill 145

\hskip.0in 13.1. 
Conformally flat hypersurfaces \dotfill 145

\hskip.0in 13.2.  Conformally flat
submanifolds \dotfill 147

\noindent Chapter 14. Submanifolds with parallel
mean curvature vector \dotfill 151

\hskip.0in 14.1. Gauss map and mean
curvature vector\dotfill 151

\hskip.0in 14.2. Riemann sphere with parallel
mean curvature vector \dotfill 151

\hskip.0in 14.3. Surfaces with parallel mean
curvature vector \dotfill 152

\hskip.0in 14.4. Surfaces with parallel
normalized mean curvature vector \dotfill
153

\hskip.0in 14.5. Submanifolds satisfying
additional conditions \dotfill 153

\hskip.0in 14.6. Homogeneous submanifolds with
parallel mean curvature vector \dotfill 154

\noindent Chapter 15. K\"ahler submanifolds of K\"ahler
manifolds \dotfill 156

\hskip.0in 15.1. Basic properties of K\"ahler
submanifolds \dotfill 156

\hskip.0in 15.2. Complex space forms and Chern
classes \dotfill 157

\hskip.0in 15.3. K\"ahler immersions 
of complex space forms in complex space 

\hskip .59in  forms
 \dotfill 158

\hskip.0in 15.4. Einstein-K\"ahler submanifolds
and K\"ahler submanifolds 
satisfying

\hskip.59 in  $Ric(X,Y)={\widetilde
Ric}(X,Y)$ \dotfill 159

\hskip.0in 15.5. Ogiue's conjectures and
curvature pinching \dotfill 160

\hskip.0in 15.6. Segre embedding \dotfill 
162

\hskip.0in 15.7. Parallel K\"ahler
submanifolds \dotfill 162

\hskip.0in 15.8. Symmetric and
homogeneous K\"ahler submanifolds \dotfill
163

\hskip.0in 15.9. Relative nullity of K\"ahler
submanifolds and reduction theorems \dotfill
164

\noindent Chapter 16. Totally real and
Lagrangian submanifolds of K\"ahler 

\hskip .59in  manifolds
\dotfill 165

\hskip.0in 16.1. Basic properties of
Lagrangian submanifolds \dotfill 165

\hskip.0in  16.2. A vanishing theorem and its
applications \dotfill 167

\hskip.0in 16.3.  The Hopf lift of Lagrangian
submanifolds of nonflat complex 

\hskip.59in space forms  \dotfill 168

\hskip.0in 16.4. Totally real minimal
submanifolds of complex space forms \dotfill
170

\hskip.0in 16.5. Lagrangian real space form
in  complex space form \dotfill 171

\hskip.0in  16.6. Inequalities for Lagrangian
submanifolds\dotfill 172

\hskip.0in 16.7.  Riemannian and
topological obstructions to Lagrangian
 
\hskip.6in immersions \dotfill 173

\hskip.0in 16.8. An inequality between
scalar curvature and mean curvature  \dotfill
174

\hskip.0in 16.9.  Characterizations of
parallel Lagrangian submanifolds  \dotfill
177

\hskip.0in 16.10. Lagrangian $H$-umbilical
submanifolds and Lagrangian catenoid  \dotfill
137

\hskip.18in 16.11. Stability of Lagrangian
submanifolds  \dotfill 178

\hskip.18in 16.12.  Lagrangian
immersions and Maslov class \dotfill 179

\noindent Chapter 17. $CR$-submanifolds of K\"ahler
manifolds \dotfill 181

\hskip.0in  17.1. Basic properties of
$CR$-submanifolds of K\"ahler manifolds
\dotfill 181

\hskip.0in  17.2. Totally umbilical
$CR$-submanifolds \dotfill 182

\hskip.0in 17.3. Inequalities for
$CR$-submanifolds \dotfill 183

\hskip.0in  17.4. $CR$-products \dotfill
184

\hskip.0in  17.5. Cyclic parallel
$CR$-submanifolds \dotfill 185

\hskip.0in  17.6. Homogeneous  and mixed
foliate $CR$-submanifolds \dotfill 186

\hskip.0in  17.7. Nullity of
$CR$-submanifolds \dotfill 186

\noindent Chapter 18. Slant submanifolds of K\"ahler
manifolds\dotfill 188

\hskip.0in  18.1. Basic properties of slant
submanifolds \dotfill 188

\hskip.0in  18.2. Equivariant slant
immersions \dotfill 189

\hskip.0in  18.3. Slant surfaces in
complex space forms \dotfill 190

\hskip.0in  18.4. Slant surfaces and
almost complex structures \dotfill 192

\hskip.0in  18.5. Slant spheres in
complex projective spaces \dotfill 193

\noindent Chapter 19. Submanifolds of the nearly
K\"ahler 6-sphere \dotfill 196

\hskip.0in  19.1. Almost complex curves
\dotfill 196

\hskip.0in  19.2. Minimal  surfaces of
constant curvature in the nearly K\"ahler
 
\hskip.6in 6-sphere \dotfill 199

\hskip.0in  19.3. Hopf hypersurfaces
and almost complex curves \dotfill 200

\hskip.0in  19.4. Lagrangian submanifolds in
nearly K\"ahler 6-sphere\dotfill 200

\hskip.0in  19.5.  Further results \dotfill
202

\noindent Chapter 20. Axioms of submanifolds \dotfill
204

\hskip.0in 20.1. Axiom of planes \dotfill
204

\hskip.0in 20.2. Axioms of spheres and 
of totally umbilical submanifolds \dotfill
204

\hskip.0in 20.3. Axiom of
 holomorphic $2k$-planes \dotfill 205

\hskip.0in 20.4. Axiom  of
 antiholomorphic $k$-planes \dotfill 206

\hskip.0in 20.5. Axioms of
 coholomorphic spheres \dotfill 206

\hskip.0in 20.6. Submanifolds
contain many circles \dotfill 207

\noindent Chapter 21. Total absolute curvature \dotfill
209

\hskip.0in 21.1. Rotation index and
total  curvature of a curve \dotfill 209

\hskip.0in 21.2. Total absolute curvature of
Chern and Lashof \dotfill 210

\hskip.0in  21.3. Tight immersions
 \dotfill 211

\hskip.0in 21.4. Taut immersions
\dotfill 214

\noindent Chapter 22. Total mean curvature \dotfill 219

\hskip.0in 22.1. Total mean curvature of
surfaces in Euclidean 3-space \dotfill 219

\hskip.0in 22.2. Willmore's conjecture
\dotfill 220

\hskip.0in 22.3. Further results on total mean
curvature for surfaces in 

\hskip.6in Euclidean space \dotfill 220

\hskip.0in 22.4. Total mean curvature for
arbitrary submanifolds and

\hskip.6in applications \dotfill 2222

\hskip.0in 22.5. Some related results
\dotfill 223

\noindent References \dotfill 227
\vfill\eject

\section{Introduction}

Problems in submanifold theory have been studied
since the invention of calculus and it was
started with differential geometry of plane
curves. Owing to his studies of how to draw
tangents to smooth plane curves, P. Fermat
(1601--1665) is regarded as a pioneer in this
field. Since his time, differential geometry
of plane curves, dealing with curvature,
circles of curvature, evolutes, envelopes,
etc., has been developed as an important part
of calculus. Also, the field has been expanded
to analogous studies of space curves and
surfaces, especially of lines of curvatures,
geodesics on surfaces, and ruled surfaces. 

Some historians date the beginning even before
the invention of calculus.  Already around 
1350, the French bishop Nicole Oresme
(1323--1382) proposed to assign 0-curvature
to the straight lines and curvature
${1\over r}$ to the circles of radius $r$. Along
the lines of previous work by J. Kepler
(1571--1630), R. Descartes (1596--1650) and C.
Huygens (1629-1695) in 1671, I. Newton
(1642--1727) succeeded in defining and
computing the curvature
$\kappa(t)$ at each point of a plane curve
using the ideas of osculating circles and
intersection of neighborhood normals. 

The first major contributor to the subject was
L. Euler (1707--1783). In 1736 Euler introduced
the arc length and the radius of curvature and
so began the study of intrinsic
differential geometry of submanifolds. 

Concerning space curves, G. Monge
(1747--1818) obtained in 1770 the expression
for the  curvature
$\kappa(t)$ of a space curve
$\gamma=\gamma(t)$. The expression for the
torsion $\tau(t)$ was first obtained by
M. A. Lancret in 1806. The works of A. L.
Cauchy (1789--1857) in 1826, F. Frenet
(1816--1900) in 1847 and J. A. Serret
(1819--1885) in 1851, resulted in the
well-known  Frenet-Serret formulas which
give all the successive derivatives of a
curve. The fundamental theorems or
congruence theorems for curves were obtained
by L. S. V. Aoust (1814--1885) in 1876.

C. F. Gauss (1777--1855) established the theory
of surfaces by introducing the concepts of
the geometry of surfaces ({\sl
Disquisitiones circa superficie curvas,}
1827). Since then the subject has
come to occupy a very firm position in
mathematics. The influence of  differential
geometry  of curves and surfaces exerted upon
branches of mathematics, dynamics, physics, and
engineering has been profound. For instance,
the study of geodesics is a topic deeply
related to dynamics, calculus of variations,
and topology; also the study of minimal
surfaces is intimately related to the theory of
functions of a complex variable, calculus of
variations, and topology.  Weierstrass and Schwarz established
its relationship with the theory of functions.
Among others, J. L. Lagrange (1736--1813), K.
Weierstrass (1815--1897), H. A. Schwarz
(1843--1921), J. Douglas (1897--1965), T.
Rad\'o (1897--1965), S. S. Chern (1911-- ),
and R. Osserman (1926-- ) are those who made
major contributions on this subject.

Belgian physicist J. A. Plateau (1801--1883)
showed experimentally that minimal surfaces
can be realized as soap films by dipping wire
in the form of a closed space curve into a
soap solution (around 1850). The Plateau
problem, that is, the problem of proving
mathematically the existence of a minimal
surface with prescribed boundary curve, was
solved by T. Rad\'o (1895--1965) in 1930, and
independently by J. Douglas (1897--1965) in
1931.

 Before Gauss, geometers viewed a
surface as being made of infinitely many
curves, whereas Gauss viewed the surface as
an entity in itself. Influenced by Gauss'
geometry on a surface in Euclidean 3-space, B.
Riemann (1826--1866) introduced in 1854
Riemannian geometry. Riemannian geometry
includes Euclidean and non-Euclidean geometries
as special cases, and it is important for the
great influence it exerted on geometric
and physical ideas of the twentieth century.

Using the concept of the intrinsic Riemannian
structure on the surface, one can compute
the curvature of a surface in two different
ways. One  is to compute the principal
curvatures and the other is done intrinsically
using the induced Riemannian metric on the
surface. The Theorema Egregium of Gauss
 provides a direct relationship between the
intrinsic and the extrinsic geometries of
surfaces. 

Motivated by the theory of mechanics,
G. Darboux (1824--1917)  unified
the theory of curves and surfaces with his
concept of a moving frame. 
This is the beginning of modern
submanifold theory which in turn gave
valuable insight into the field.

Since the celebrated embedding
theorem of J. F. Nash (1928-- ) allows
geometers to view each Riemannian manifold
as a submanifold in a Euclidean space, the
problem of discovering simple sharp
relationships between  intrinsic
 and extrinsic invariants of a
submanifold is one of the most
fundamental problems in submanifold theory.
Many beautiful  results in
submanifold theory, including
the Gauss-Bonnet theorem and isoperimetric
inequalities, are  results in this respect.
In the modern theory of submanifolds, the
study of relations between local and global
properties has also attracted the interest
of many geometers. This view was emphasized
by W. Blaschke (1885--1962), who worked on the
differential geometry of ovals and ovaloids.
The study of rigidity of ovaloids by S.
Cohn-Vossen (1902--1936) belongs in this
category.

An important class of Riemannian manifolds was
discovered by J. A. Schouten (1883--1971), D.
van Dantzig (1900--1959), and E. K\"ahler
(1906-- ) around 1929--1932. This class of
manifolds, called K\"ahler manifolds today,
includes the projective algebraic manifolds.
The study of complex submanifolds of a K\"ahler
manifold from differential geometric points of
view was initiated by E. Calabi (1923-- ) in
the early 1950's. Besides complex submanifolds,
there are some other important classes of
submanifolds of a K\"ahler manifold
determined by the behavior of the tangent
bundle of the submanifold under the action
of the almost complex structure of the
ambient manifold. These classes of
submanifolds have many interesting
properties and many important results have
been discovered in the last quarter of this
century from this standpoint.

Submanifold theory is a very active vast
research field which in turn plays an
important role in the development of modern
differential geometry  in this century. 
This branch of differential geometry
is still so far from being exhausted; only a
small portion of an exceedingly fruitful field
has been cultivated, much more remains to be
discovered in the coming centuries.

\vfill\eject
\section{ Nash's embedding theorem and some related results}

Throughout this article manifolds are assumed
to be connected, of class $C^\infty$, and
without boundary, unless mentioned otherwise.

H. Whitney (1907--1989) proved in 1936 that
an $n$-manifold can always be immersed in the
Euclidean $2n$-space $E^{2n}$, and can always
be embedded in $E^{2n+1}$ as a closed set.

 An immersion $f$ (or an embedding) of a
Riemannian manifold $(M,g)$ into another
Riemannian manifold $(\tilde M,\tilde g)$ is
said to be isometric if it satisfies the
condition $f^*\tilde g=g$.  In this case,
$M$ is called a Riemannian submanifold (or
simply a submanifold) of $\tilde M$.

 We shall identify the image
$f(M)$ with $M$ when there is no danger of
confusion. 

One of fundamental problems in submanifold
theory is the problem of
isometric immersibility. The earliest
publication on isometric embedding appeared
in 1873 by L. Schl\"afli (1814--1895).

The problem of isometric immersion (or
embedding) admits an obvious analytic
interpretation; namely, if $g_{ij}(x),\,
x=(x_1,\ldots,x_n),$ are the components of
the metric tensor $g$ in local coordinates
$x_1,\ldots,x_n$ on a Riemannian
$n$-manifold $M$, and $y=(y_1,\ldots,y_m)$
are the standard Euclidean coordinates in
$E^m$, then the condition for an
isometric immersion in $E^m$ is
$$\sum_{i=1}^n {{\partial y_j}\over
{\partial x_i}} {{\partial y_k}\over
{\partial x_i}}=g_{jk}(x),$$
that is, we have a system of $\frac12
n(n+1)$ nonlinear partial differential
equations in $m$ unknown functions. If
$m=\frac 12 n(n+1)$, then this system is
definite and so we would like to have a
solution. Schl\"afli asserted that any Riemannian
$n$-manifold can be isometrically embedded
in Euclidean space of dimension $\frac12
n(n+1)$. Apparently it is appropriate to
assume that he had in mind of analytic
metrics and local analytic embeddings.
This  was later called
Schl\"afli's conjecture.

\subsection{Cartan-Janet's theorem}

In 1926 M. Janet (1888--1984) published a
proof of Schl\"afli's conjecture which states that a real
analytic Riemannian $n$-manifold $M$ can be
locally isometrically embedded into any real
analytic Riemannian manifold of dimension
${1\over 2}n (n+1)$. In 1927
\'E. Cartan (1869--1951) revised Janet's
paper with the same title; while Janet
wrote the problem in the form of a system
of partial differential equations which he
investigated using rather complicated
methods, Cartan applied his own theory of
Pfaffian systems in involution. Both Janet's
and Cartan's proofs contained obscurities.
In 1931  C. Burstin got rid of them. This
result of Cartan-Janet implies that every
Einstein
$n$-manifold $(n\geq 3)$ can be locally
isometrically embedded in $E^{n(n+1)/2}$.

 The Cartan-Janet theorem
is dimensionwise the best possible, that is,
there exist real analytic Riemannian
$n$-manifolds which do not possess smooth local
isometric embeddings in any Euclidean space of
dimension strictly less than ${1\over 2}n
(n+1)$. Not every Riemannian $n$-manifold
can be isometrically immersed in
$E^m$ with $m\leq \frac12 n(n+1)$. For
instance, not every Riemannian 2-manifold
can be isometrically immersed in $E^3$.

Cartan-Janet's theorem implies that an
analytic Riemannian 3-manifold can be
locally isometrically embedded into $E^6$.
For Riemannian smooth 3-manifolds, R. L.
Bryant, P. A. Griffiths and D. Yang (1983)
proved the following: 

Let $M$ be a smooth Riemannian 3-manifold and
let $x\in M$ be such that the Einstein tensor
at $x$ is neither zero nor a perfect square
$L^2=\sum\ell_i\ell_jdx^idx^j,\,L=\sum\ell_idx^i\in
T^*_xM$. Then there exists a neighborhood of
$x$ which can be smoothly isometrically
embedded into $E^6$. 

G. Nakamura and Y. Maeda (1986) improved the
result to the following:
 Let $M$ be a smooth Riemannian 
3-manifold and let $x\in M$ be a point such
that the curvature tensor $R(x)$ at $x$ does
not vanish, where $R(x)$ is
considered as a symmetric linear operator acting
on the space of 2-forms. Then there exists a local 
smooth isometric embedding of a neighborhood
of $x$ into $E^6$. 

H. Jacobowitz and J. D. Moore
(1973) proved that, for real analytic Riemannian
manifolds $M$ and $\tilde M$ of dimensions $n$
and $N$ respectively with $n\geq 2,\, N\geq
{1\over 2}n(n+1)-1$, if $x_0$ is a point of
$M$, then there exists an open neighborhood
$U$ of $x_0\in M$ which can be conformally
embedded in $\tilde M$.

\subsection{Nash's embedding theorem}

A global isometric embedding theorem was
proved by J. F. Nash (1956) which states as
follows.

\vskip.1in
\noindent {\bf Theorem 2.1} {\it Every compact Riemannian
$n$-manifold can be isometrically embedded  in
any small portion of a Euclidean $N$-space
$ E^N$ with
$N={1\over 2}n (3n+11)$. Every non-compact
Riemannian $n$-manifold can be isometrically
embedded in any small portion of a Euclidean
$m$-space $E^m$ with $m={1\over 2}n
(n+1)(3n+11)$.} 
\vskip.1in   

R. E. Greene (1970) improved Nash's result
and proved that every non-compact Riemannian
$n$-manifold can be isometrically embedded in
Euclidean $N$-space with $N=2(2n+1)(3n+7)$.
Furthermore, Greene (1970) and M. L. Gromov
and V. A. Rokhlin (1970) proved
independently that a local isometric
embedding from Riemannian
$n$-manifold into $E^{{1\over 2}n(n+1)+n}$
always exist.  Gromov and  Rokhlin
(1970) also proved that a compact Riemannian
$n$-manifold of class $C^r$ ($r=\infty$ or
$\omega$) can be isometrically $C^r$-embedded
in  $E^{\frac12 n(n+1)+3n+5}$.

Nash (1954) proved that if a manifold $M$
admits a $C^1$-embedding in $E^m$, where
$m\geq n+2$, then it admits an isometric
$C^1$-embedding in $E^m$. N. H. Kuiper
(1920--1994) improved this result in 1955 by
showing that it is true when $m\geq n+1$.
M. L. Gromov showed that 
every Riemannian $n$-manifold can be $C^1$
isometrically immersed into $E^{2n}$ (cf.
[Gromov 1986]).  

It is known that a Hermitian symmetric
space $G/K$ of compact type can be
equivariantly and isometrically embedded
into $E^m$ with $m=\dim G$  [Lichnerowicz
1958]. This result has been extended to
almost all symmetric spaces of compact type
[Nagano 1965; Kobayashi 1968]. In 1976 J. D. 
Moore proved that every compact
Riemannian homogeneous manifold admits an
equivariant isometric embedding in some
Euclidean space.

\subsection{Isometric immersions with the smallest possible codimension}

According to Nash's embedding theorem, every
Riemannian manifold can be isometrically
embedded in a Euclidean space of sufficiently
large dimension, it is thus natural  to look
for a Euclidean space of smallest possible
dimension in which a Riemannian manifold can be
isometrically embedded. 

D. Hilbert (1862--1943) proved in 1901 that
a complete surface of constant negative
curvature cannot be $C^4$-isometrically
immersed in Euclidean 3-space.

A result of S. S. Chern and N. H. Kuiper
(1952) states that a compact Riemannian
$n$-manifold with non-positive sectional
curvature cannot be isometrically immersed in
$ E^{2n-1}$. A generalization by B. O'Neill
(1924-- )  states that if $N$ is a complete
simply-connected Riemannian $(2n-1)$-manifold
with sectional curvature $K_N\leq 0$, then a
compact Riemannian $n$-manifold $M$ with
sectional curvature $K_M\leq K_N$ cannot be
isometrically immersed in $N$ [O'Neill
1960].

T. Otsuki (1917-- ) proved in 1954 that a
Riemannian $n$-manifold of constant negative
sectional curvature cannot be isometrically
immersed in
$ E^{2n-2}$ even locally. He also proved
that a Riemannian $n$-manifold cannot be
isometrically immersed in a Riemannian
$(2n-2)$-manifold if $K_M<K_N$. 

Using the purely algebraic theory of ``flat
bilinear forms'', J. D. Moore (1977) proved that
if a compact Riemannian $n$-manifold
$M$ of constant curvature 1 admits an isometric
immersion in $ E^N$ with $N\leq {3\over
2}n$, then $M$ is simply-connected, hence
isometric to the unit $n$-sphere. A real
analytic version of this result was obtained
by W. Henke (1976) in the special case where
$n\geq 4$ and $N=n+2$. 

\subsection{Isometric immersions with
prescribed Gaussian or Gauss-Kro\-necker
curvature}

A 1915 problem of H. Weyl (1885--1955) is
that whether a Riemannian 2-manifold of
positive Gaussian curvature that is
diffeomorphic to a sphere can be realized as
a smooth ovaloid in
$E^3$? Weyl himself suggested an incomplete
solution of this problem for analytic
surfaces. In fact, he solved the problem in
the case of analytic metrics sufficiently
close to the metric of a sphere. 

A complete solution of Weyl's problem for
analytic case was given by H. Lewy
(1904--1988) in 1938. L. Nirenberg (1925-- )
proved in 1953 that given  a
$C^\infty$-smooth Riemannian metric $g$ on a
topological 2-sphere $S^2$ with Gaussian
curvature $K>0$, there exists  a
$C^\infty$-smooth global isometric
embedding of $(S^2,g)$ into $E^3$.

A local immersibility for $C^k$-smooth
metric with $K\geq 0$ and $k\geq 10$ in the
form of a $C^{k-6}$-smooth convex surface
was proved in 1985 by C. S. Lin. In 1995 J.
Hong and C. Zuily  extended Nirenberg's
global result to the case $K\geq 0$. J.
Hong (1997) established isometric embedding
in $E^3$ of complete noncompact
nonnegatively curved surfaces. 

C. S. Lin (1986) considered the problem of 
isometric embedding of two-dimen\-sional
metrics of curvature that changes sign and
proved the following: Let the curvature
$K$ of a Riemannian 2-manifold be equal to
zero at the point $P$, but the gradient of
the curvature $\nabla K$ be nonzero. If the
metric of the manifold belongs to the class
$C^s$, $s\geq 6$, then a neighborhood of $P$
admits an isometric embedding of class
$C^{s-3}$ in $E^3$. W. Greub and D.
Socolescu (1994) claimed that the
condition $\nabla K\ne 0$ in Lin's result
can be removed.

N. V. Efimov (1964) proved that a complete
surface with Gaussian curvature $K\leq
-c^2$,  $c$ a positive constant, does not
admit an isometric immersion in $E^3$. C.
Baikoussis and T. Koufogiorgos (1980)
showed that a complete surface with
curvature $-\infty<-a^2\leq K\le 0$ in
$E^3$ is unbounded. 

B. Smyth and F. Xavier (1987) proved  that
if a complete Riemannian $n$-manifold $M$
with negative Ricci curvature is immersed
as a hypersurface in a Euclidean space,
then the upper bound of the Ricci curvature
of $M$ is equal to zero if $n=3$, or if
$n>3$ and the sectional curvatures of $M$
do not take all real values.

For an oriented hypersurface $M$ in $E^{n+1}$,
the determinant of the shape operator of
$M$ with respect to the unit outward normal
is called the Gauss-Kronecker curvature of the
hypersurface.

Given a smooth positive real function $F$ 
on $E^{n+1}$, the problem to find a 
closed convex hypersurface $M$ in $E^{n+1}$
with Gauss-Kronecker curvature
$K=F$ have been studied by various geometers.
For instance, sufficient conditions were found
under which this problem can be solved,
either by topological methods [Oliker 1984;
Caffarelli-Nirenberg-Spruck 1986] or by 
geometric variational approach [Oliker
1986; Tao 1991; X. J. Wang 1996].

\subsection{Isometric immersions with
prescribed mean curvature }

Given a closed $(n-1)$-dimensional
submanifold $\Gamma$ in a Riemannian
manifold $N$, the problem of
finding an oriented
$n$-dimensional submanifold $M$ with a
prescribed mean curvature vector and with
$\Gamma$ as its boundary has been
investigated by  many mathematicians. 

The first necessary conditions for parametric
surfaces were given by E. Heinz in 1969 for
surfaces of constant mean curvature in
Euclidean 3-space with a prescribed
rectifiable boundary.  R. Gulliver gave in
1983 a necessary condition on the magnitude
of the mean curvature vector field for there
to exist an oriented submanifold of a
Riemannian manifold
$M$ having prescribed mean curvature vector
and a given closed submanifold as boundary.

I. J. Bakelman and B. E. Kantor (1974), A.
Treibergs and W. S. Wei (1983), A. Treibergs
(1985), and K. Tso (1989) established the
existence of closed convex hypersurfaces in a
Euclidean space with prescribed mean
curvature.

\vfill\eject

\section{Fundamental theorems, basic notions and results}
 
\subsection{Fundamental equations}

 Let $f:(M,g)\to (\tilde M,\tilde g)$ be an
isometric immersion. Denote by
$\nabla$ and
$\tilde \nabla$ the metric connections of $M$
and $\tilde M$, respectively. For vector fields
$X$ and $Y$ tangent to $M$, the tangential
component of $\tilde\nabla_XY$ is equal to
$\nabla_XY$. 

Let
$$h(X,Y)=\tilde\nabla_XY-\nabla_XY.\leqno(3.1)$$
The $h$ is a normal-bundle-valued symmetric
$(0,2)$ tensor field on $M$, which is called
the second fundamental form of the submanifold
(or of the immersion). Formula (3.1) is known as
the Gauss formula [Gauss 1827].

For a normal vector $\xi$ at a point $x\in M$,
we put
$$g(A_\xi X,Y)=\tilde g(h(X,Y),\xi).\leqno(3.2)$$
Then $A_\xi$ is a symmetric linear
transformation on the tangent space $T_xM$ of
$M$ at $x$, which is called the shape operator
(or the Weingarten map) in the direction of
$\xi$. The eigenvalues of $A_\xi$ are called the
principal curvatures in the direction of $\xi$.

The metric connection on the normal bundle
$T^\perp M$  induced from the metric connection
of $\tilde M$ is called the normal connection of
$M$ (or of $f$).

Let $D$ denote covariant differentiation
with respect to the normal connection. For a
tangent vector field $X$ and a normal vector
field $\xi$ on $M$, we have
$$\tilde\nabla_X\xi=-A_\xi X+D_X\xi,\leqno(3.3)$$
where $-A_\xi X$ is the tangential component of 
$\tilde\nabla_X\xi$. 
 (3.3) is known as
the Weingarten formula, named after
the 1861 paper of J. Weingarten
(1836--1910).

Let $R,\tilde R$ and $R^D$ denote the
Riemannian curvature tensors of
$\nabla,\tilde\nabla$ and $D$, respectively. Then
the integrability condition for (3.1) and (3.3)
implies
$$\aligned \tilde R(X,Y)Z&=R(X,Y)Z+A_{h(X,Z)}Y-
A_{h(Y,Z)}X\\&
+(\bar\nabla_Xh)(Y,Z)-(\bar\nabla_Yh)(X,Z),
\endaligned\leqno(3.4)$$
for tangent vector fields $X,Y,Z$ of $M$, where
$\bar\nabla$ is the covariant differentiation
with respect to the connection in $TM\oplus
T^\perp M$. The tangential and normal
components of (3.4) yield the
following equation  of Gauss
$$
\left<R(X,Y)Z,W\right>=\left<{}\right.\tilde
R(X,Y)Z,W\left.{}\right>+\left<h(X,W),
h(Y,Z)\right>-
\left<h(X,Z),h(Y,W)\right>\leqno(3.5)$$
and the equation of Codazzi
$$ (R(X,Y)Z)^\perp=
(\bar\nabla_Xh)(Y,Z)-(\bar\nabla_Yh)(X,Z),
\leqno(3.6)$$
where $X,Y,Z,W$ are tangent vectors of $M$,
$(R(X,Y)Z)^\perp$ is the normal component of
$R(X,Y)Z$, and $\left<\;,\right>$ is the
inner product.

Similarly, for normal vector fields $\xi$ and
$\eta$, the relation
$$\left<{}\right.\tilde
R(X,Y)\xi,\eta\left.{}\right>=\left<R^D(X,Y)\xi,
\eta\right>
-\left<[A_\xi,A_\eta]X,Y\right>\leqno(3.7)$$
holds, which is called the equation of Ricci.

Equations (3.1), (3.3), (3.5), (3.6) and (3.7)
are called the fundamental equations of the
isometric immersion $f:M\to\tilde M$.

As a special case, suppose the ambient
space $\tilde M$ is a Riemannian manifold of
constant sectional curvature $c$. Then the
equations of Gauss, Codazzi and Ricci reduce
respectively to
$$\left<R(X,Y)Z,W\right>=\left<h(X,W),h(Y,Z)\right>-\left<h(X,Z),h(Y,W)\right>$$
$$+c\{\left<X,W\right>\left<Y,Z\right>-\left<X,Z\right>\left<Y,W\right>\},\leqno(3.8)$$
$$(\bar\nabla_Xh)(Y,Z)=(\bar\nabla_Yh)(X,Z),\leqno(3.9)$$
$$\left<R^D(X,Y)\xi,\eta\right>=\left<[A_\xi,A_\eta]X,Y\right>.\leqno(3.10)$$

Formulas (3.8) and (3.9) for
surfaces in $E^3$ were given in
principal, though not explicitly, in
[Gauss 1827]. The formulas can be found
in a 1860 paper by D. Codazzi
(1824--1875) in his answer to a
``concours'' of the Paris Academy
(printed in the M\'emoires
pr\'esent\'e \`a l'Acad\'emie in 1880;
also in [Codazzi 1868]). These formulas
at that time already published by G.
Mainardi (1800--1879) in [Mainardi
1856]. The fundamental importance of
these formulas was fully recognized by
O.  Bonnet (1819--1892) in [Bonnet
1867]. The  equations of Gauss and
Codazzi for general submanifolds were
first given by A. Voss in 1880. The
equation (3.10) of Ricci was first
given by G. Ricci (1853--1925) in 1888.

\subsection{Fundamental theorems}

The fundamental theorems of submanifolds are
the following.
\vskip.1in

\noindent{\bf Existence Theorem} {\it Let $(M,g)$ be a simply-connected Riemannian $n$-manifold and suppose there is a given $m$-dimensional Riemannian vector bundle
$\nu(M)$ over $M$ with curvature tensor $R^D$ and a $\nu(M)$-valued symmetric $(0,2)$ tensor $h$ on $M$. For a cross section $\xi$ of $\nu(M)$, define $A_\xi$ by $g(A_\xi
X,Y)=\left<h(X,Y),\xi\right>$, where $\left<\;,\;\right>$ is the fiber metric of
$\nu(M)$. If they satisfy $(3.8), (3.9)$ and $(3.10)$, then $M$ can be
isometrically immersed in an
$(n+m)$-dimensional complete simply-connected
Riemannian manifold
$R^{n+m}(c)$ of constant curvature $c$ in such
way that $\nu(M)$ is the normal bundle and $h$ is
the second fundamental form.}

\vskip.1in
\noindent{\bf Uniqueness Theorem}{\it 
Let $f,f':M\to R^{m}(c)$ be two isometric
immersions of a Riemannian $n$-manifold $M$
into a complete simply-connected Riemannian
$m$-manifold of constant curvature $c$
with normal bundles $\nu$ and $\nu'$
equipped with their canonical bundle
metrics, connections and second fundamental
forms, respectively. Suppose there is an
isometry
$\phi:M\to M$ such that
$\phi$ can be covered by a bundle map
$\bar\phi:\nu\to\nu'$ which preserves the
bundle metrics, the connections and the second
fundamental forms. Then there is an isometry
$\Phi$ of $R^{m}$ such that $\Phi\circ
f=f'$.}

The first to give a proof of the fundamental theorems was O. Bonnet 
(cf. [Bonnet 1867]).

\subsection{Basic notions}

Let $M$ be an $n$-dimensional Riemannian
submanifold of a Riemannian manifold $\tilde M$.
A point $x\in M$ is called a geodesic point
if the second fundamental form $h$ vanishes
at $x$. The submanifold is said to be totally
geodesic if every point of $M$ is a geodesic
point. A Riemannian submanifold $M$ is a
totally geodesic submanifold of $\tilde M$
if and only if every geodesic of $M$ is a
geodesic of $\tilde M$.

Let $M$ be a submanifold of $\tilde M$  and let
$e_1,\ldots,e_n$ be an orthonormal basis of
$T_xM$. Then the mean curvature vector
$H$ at
$x$ is defined by $$ H={1\over
n}\sum_{j=1}^n h(e_j, e_j).$$  The length
of $ H$ is called the mean curvature
which is denoted by $H$.
$M$ is called a minimal submanifold of $\tilde
M$ if the mean curvature vector field vanishes
identically. 

A point $x\in M$ is called an umbilical point
if $h=g\otimes  H$ at $x$, that is,
the shape operator $A_\xi$ is proportional to
the identity transformation for all $\xi\in
T^\perp_x M$. The submanifold is said to be
totally umbilical if every point of the
submanifold is an umbilical point. 

A point $x\in M$ is called an isotropic point if
$\,|h(X,X)|/|X|^2\,$ does not depend on
the nonzero vector $X\in T_xM$. A submanifold
$M$ is called an isotropic submanifold if
every point of $M$ is an isotropic point.
The submanifold $M$ is called constant
isotropic if 
$|h(X,X)|/|X|^2$ is also independent of the
point $x\in M$. It is clear that umbilical
points are isotropic points.

The index of relative nullity at  $x\in
M$ of a submanifold
$M$ in $\tilde M$ is defined by $$\mu(x)=\dim
\cap_{\xi\in T^\perp M} \ker A_\xi.$$
If we denote by $N_0(x)$ the null space of the
linear mapping $\xi\to A_\xi$, then the
orthogonal complement of $N_0(x)$ in the normal
space $T^\perp_xM$ is called the first normal
space at $x$. 

A  normal vector field $\xi$ of $M$ in $\tilde
M$ is said to be parallel in the normal bundle
if $D_X\xi=0$ for any vector $X$ tangent to
$M$. A submanifold $M$ is said to have
parallel mean curvature vector if the mean
curvature vector field of $M$ is a parallel
normal vector field.

A submanifold $M$ in a Riemannian manifold is
called a parallel submanifold if its second
fundamental form $h$ is parallel, that is,
$\bar\nabla h=0$, identically.

 A Riemannian submanifold $M$ is said
to have flat normal connection if the curvature
tensor $R^D$ of the normal connection $D$
vanishes at each point $x\in M$.

For a Riemannian $n$-manifold $M$, we denote by
$K(\pi)$ the sectional curvature of a 2-plane
section  $\pi\subset T_xM$. Suppose
$\{e_1,\ldots,e_n\}$ is an orthonormal basis of
$T_xM$. The Ricci curvature $Ric$ and the scalar
curvature $\rho$ of $M$ at $x$ are defined
respectively by $$Ric(X,Y)
 =\sum_{j=1}^n
\left<R(e_j,X)Y,e_j\right>,\leqno(3.11)$$
$$ \rho=\sum_{i\ne
j} K(e_i\wedge e_j),\leqno(3.12)$$
 where $K(e_i\wedge e_j)$ denotes the sectional
curvature of the 2-plane section spanned by $e_i$
and $e_j$.

In general if $L$ is an $r$-plane section in
$T_xM$ and $\{e_1,\ldots,e_r\}$ an orthonormal
basis of $L\subset T_xM$, then the scalar
curvature $\rho(L)$ of $L$ is defined by
$$\rho(L)=\sum_{i\ne j}
K(e_i\wedge e_j),\quad 1\leq i,j\leq r.\leqno(3.13)$$

\subsection{A general inequality}

Let $M$ be a  Riemannian $n$-manifold.
For an integer $k\geq 0$, denote by $\Cal S(n,k)$
the finite set  consisting of $k$-tuples
$(n_1,\ldots,n_k)$ of integers $\geq 2$
satisfying  $n_1< n$ and
$n_1+\cdots+n_k\leq n$. Denote by ${\Cal S}(n)$ the set of
(unordered) $k$-tuples with $k\geq 0$ for a
fixed positive integer $n$. 

The cardinal number $\#\Cal S(n)$ of $\Cal
S(n)$ is equal to $p(n)-1$, where $p(n)$
denotes the number of partition of $n$
which  increases quite rapidly with $n$. 
For instance, for
$$n=2,3,4,5,6,7,8,9,10,\ldots,20,\ldots,
50,\ldots,100,\ldots,200,$$
the cardinal number $\#\Cal S(n)$ are given
 by
$$1,2,4,6,10,14,21,29,41,\ldots,626,
\ldots,
204225,\ldots,
190569291,\ldots,3972999029387,$$
respectively. The asymptotic behavior of
$\#\Cal S(n)$ is given by
$$\#\Cal S(n)\sim {1\over
{4n\sqrt{3}}}\exp\left[\pi\sqrt{{2n}/
3}\,\right]\quad
\hbox{as}\;\;n\to\infty.$$

For each  $(n_1,\ldots,n_k)\in\Cal S(n)$ B.
Y. Chen (1996f,1997d) introduced a Riemannian
invariant
$\delta{(n_1,\ldots,n_k)}$ by
$$\delta{(n_1,\ldots,n_k)}(x)={1\over
2}\left(\rho(x)-\inf\{\rho(L_1)+\cdots+
\rho(L_k)\}\right),
\leqno(3.14)$$where $L_1,\ldots,L_k$ run over
all
$k$ mutually orthogonal subspaces of $T_xM$ such
that $\dim L_j=n_j,\, j=1,\ldots,k$. We put
$\delta(\emptyset)=\frac \rho 2$.

For each $(n_1,\ldots,n_k)\in {\Cal S}(n)$,
let
$$ a(n_1,\ldots,n_k)={1\over2}
{{n(n-1)}}-{1\over2}\sum_{j=1}^k
{{n_j(n_j-1)}}.\leqno(3.15)$$
$$b(n_1,\ldots,n_k)= {{n^2(n+k-1-\sum
n_j)}\over{2(n+k-\sum
n_j)}}.\leqno(3.16)$$
When $k=0$, the left hand sides of (3.15)
and (3.16) are denoted respectively by
$a(\emptyset),\, b(\emptyset)$.

For any $n$-dimensional submanifold $M$ of a
real space form $R^m(c)$ and for any 
$k$-tuple $(n_1,\ldots,n_k)\in
\Cal S(n)$, regardless of dimension and
codimension, there is a sharp
general inequality between the invariant
$\delta(n_1,\ldots,n_k)$ and the squared
mean curvature $H^2$ given by [Chen 1996c,1996f]:
$$\delta(n_1,\ldots,n_k)\leq
b(n_1,\ldots,
n_k)H^2+a(n_1,\ldots,n_k)c.\leqno(3.17)$$
The equality case of inequality
$(3.17)$ holds at a point $x\in M$ if and only
if, there exists an orthonormal basis 
$e_1,\ldots,e_m$ at $x$, such
that  the shape operators of $M$ in $R^m(c)$ at
$x$ take the following form:
$$A_r=\begin{pmatrix} A^r_{1} & \hdots & 0&0 &
\hdots & 0\\
\vdots  & \ddots& \vdots&\vdots  & &\vdots \\
0 &\hdots &A^r_k&0 &\hdots & 0 
\\ 0&\hdots&0&\mu_r&\hdots&0\\
\vdots  & \ddots & \vdots &\vdots &\ddots&\vdots \\
0  &\hdots& 0&0 &\hdots & \mu_r \end{pmatrix}
,\quad  r=n+1,\ldots,m,
\leqno(3.18)$$
where $A^r_j$ are symmetric $n_j\times n_j$ 
submatrices which satisfy 
$$\hbox{\rm trace}\,(A^r_1)=\cdots=\hbox{\rm
trace}\,(A^r_k)=\mu_r.$$

Inequality (3.17) has many interesting
applications. For instance, regardless of
codimension, it implies that the squared
mean curvature of every isometric immersion
of $S^n(1)$ (respectively, of $S^k(1)\times 
E^{n-k}$ or of $S^k(1)\times S^{n-k}(1)$) into
a Euclidean space must satisfy
$$H^2\geq 1\quad \left(\hbox{respectively, }
H^2\geq\left({k\over n}\right)^2\hbox{ or }
H^2\geq \left({k\over
n}\right)^2+\left({{n-k}\over
n}\right)^2\right),\leqno(3.19)$$ with the
equality holding if and only if it is a
standard embedding.

The equality case of (3.17) with $k=0$ holds
at a point $x\in M$ when and only when $x$
is an umbilical point. In general, there
exist ample examples of submanifolds which
satisfy the equality case of (3.17) with
$k>0$. 

Inequality (3.17) also provides the following
sharp estimate of the first nonzero eigenvalue
$\lambda_1$ of the Laplacian $\Delta$ on each
compact irreducible homogeneous Riemannian
$n$-manifold $M$ [Chen 1996f, Chen 1997d]:
$$\lambda_1\geq n\hat\Delta_0,\leqno(3.20)$$
where
$$\hat\Delta_0=\max\left\{{{\delta(n_1,\ldots,n_k)}
\over{b(n_1,\ldots,n_k)}}:(n_1,\ldots,n_k)\in\Cal
S(n)\right\}.\leqno(3.21)$$

Clearly, the invariant $\hat\Delta_0$ is
constant on a homogeneous Riemannian
manifold, since each $\delta$-invariant
$\delta(n_1,\ldots,n_k)$ is constant on such
a space.

The estimate of $\lambda_1$ given in (3.20)
improves a well-known result of T. Nagano
(1930-- ) who proved in 1961 that
$\lambda_1\geq\frac\rho{n-1}$ for each
compact irreducible homogeneous Riemannian
$n$-manifold $M$, with the equality holding if
and only if $M$ is a Riemannian $n$-sphere. We
remark that $\rho/({n-1})$ is nothing but
$\delta(\emptyset)/b(\emptyset)$.

Inequality (3.17) for $\delta(2)$ was
first proved in [Chen 1993]. The equality
case of (3.17) for $\delta(2)$ have been
investigated by D. E. Blair, J. Bolton, B.
Y. Chen, M. Dajczer, F. Defever, R. Deszcz,
F. Dillen, L. A. Florit, C. S. Houh, I.
Mihai, M. Petrovic, C. Scharlach, L.
Verstraelen, L. Vrancken, L. M. Woodward, J.
Yang, and others.

Further applications of (3.17) can be found
in \S5.3.1, \S5.4.1, and \S16.7.

\subsection{Product immersions} 

Suppose that $M_1,\ldots,M_k$ are  Riemannian
manifolds and that $$f:M_1\times\cdots\times
M_k\to E^N$$ is an isometric immersion of
the Riemannian product $M_1\times\cdots\times
M_k$ into Euclidean $N$-space. J. D. Moore
(1971) proved that if the second fundamental
form
$h$ of $f$ has the property that
$h(X,Y)=0$, when $X$ is tangent to $M_i$ and
$Y$ is tangent to $M_j$ for $i\ne j$, then $f$
is a product immersion, that is, there exist
isometric immersions $f_i:M_i\to  E^{m_i},\,
1\leq i\leq k$ such that
$$f(x_1,\ldots,x_k)=(f(x_1),\ldots,f(x_k))$$
when $x_i\in M_i$ for $1\leq i\leq k$.

 Let $M_0,\cdots ,M_k$ be Riemannian
manifolds, $M=M_ 0 \times\cdots
\times M_ k$ their product, and $\pi_i : M
\to M_i$  the canonical
projection. If $\rho_1,\cdots ,\rho_k
: M_0\to \hbox{\bf R}_+$ are positive-valued
functions, then $$\left< X,Y \right>:=
\left<\pi_{0*}X,\pi_{0*}Y
\right> + \sum^k_{i=1} (\rho_i \circ
\pi_0)^2 \left< \pi\sb {i*}X, \pi\sb {i*}Y
\right>$$ defines a Riemannian metric on $M$,
called a warped product metric. $M$
endowed with this metric is denoted by $M_0
\times_{\rho_1} M_1 \times \cdots \times_{\rho\sb
k}M_k$. 

A warped product immersion is
 defined as follows: Let $M_0 \times_{\rho_1}
M_1\times
\cdots\times_{\rho_k}M_k$ be a warped
 product and let $f_i: N_i
\to M_i$, $i=0,\cdots,k$, be isometric
immersions, and define $\sigma_i
:= \rho_i \circ f_0 : N_0
\to \hbox{\bf R}_+$ for $i = 1,\cdots ,k$.
Then the map $$f: N_0 \times_{\sigma_1} N_1
\times \cdots \times_{\sigma_k} N_k
\to M_0 \times_{\rho_1}
M_1\times\cdots \times_{\rho_k} M_k$$ 
given by $f(x_0,\cdots ,x_k)
:= (f_0(x_0),f_1(x_1),\cdots ,
f_k(x_k))$ is an isometric immersion, which
is called a warped product immersion. 

S. N\"olker (1996) extended
 Moore's result to the following.
 
 Let $f: N_0 \times_{\sigma_1} N_1
\times \cdots \times_{\sigma_k} N_k
\to R^N(c)$ be an isometric immersion into
a  Riemannian
manifold of constant curvature
$c$. If $h$ is the second fundamental form of
$f$ and
$h(X_i,X_j)=0$, for all vector fields $X_i$ and
$X\sb j$, tangent to $N_i$ and $N_j$
respectively, with $i\ne j$, then, locally,
$f$ is a warped product immersion. 

\subsection{A relationship between
$k$-Ricci tensor and shape
operator}

Let $M$ be a Riemannian $n$-manifold and
$L^k$ is a $k$-plane section of $T_xM^n,\,
x\in M.$ For each unit vector $X$ in $L^k$, we
choose  an orthonormal basis
$\{e_1,\ldots,e_k\}$ of $L^k$ such that
$e_1=X$. Define the Ricci curvature
$Ric_{L^k}$ of ${L^k}$ at $X$ by
$$Ric_{L^k}(X)=K_{12}+\cdots+K_{1k},\leqno(3.22)$$
where $K_{ij}$ denotes the sectional
curvature of the
$2$-plane section spanned by $e_i, e_j$. 
 We call $Ric_{L^k}(X)$ a
 $k$-Ricci curvature of $M$ at $X$
relative to $L^k$.  Clearly,  the $n$-th
Ricci curvature is nothing but the
 Ricci curvature in the usual sense and second
Ricci curvature coincides with the sectional
curvature.

For each integer $k,\,2\leq k\leq n$, 
let $\theta_k$ denote the Riemannian invariant
defined on $M$ by
$$\theta_k(x)=\left({1\over
{k-1}}\right)\inf_{L^k,X} Ric_{L^k}(X),\quad
X\in T_xM,\leqno(3.23)$$ where $L^k$ runs over
all $k$-plane sections in $T_xM$ and $X$
runs over all unit vectors in $L^k$. 

The following results provide a sharp
relationship between the $k$-Ricci curvature and
the shape operator for an arbitrary submanifold in
a real space form, regardless of codimension
[Chen 1996e,1998c].

 Let $f:M\to R^m(c)$ be an
isometric immersion of a Riemannian
$n$-manifold $M$ into a Riemannian $m$-manifold
$R^m(c)$ of constant sectional curvature $c$.
Then, for any
 integer $k,\,2\leq k\leq n$, and any point
$x\in M^n$, we have  

(1) if $\theta_k(x)\ne c$, then the shape
operator in the direction of the mean curvature
vector satisfies
 $$A_{ H}>{{n-1}\over n}(\theta_k(x)-
c)I\quad
\hbox{at }x,\leqno (3.24)$$ where $I$
denotes the identity map of $T_xM^n$;

Inequality (3.24) means that $A_{
H}-{{n-1}\over n}(\theta_k(x)- c)I$
is positive-definite.

(2) if $\theta_k(x)=c$, then $A_{
H}\geq 0$ at $x$;

(3) a unit vector $X\in T_xM$
satisfies $A_{ H}X={{n-1}\over
n}(\theta_k(x)- c)X$ if and only if
$\theta_k(x)=c$ and $X$ lies in the relative
null space at $x$;

(4) $A_{ H}\equiv {{n-1}\over
n}(\theta_k- c)I$ at $x$ if and only if
$x$ is a  totally geodesic point.
\smallskip

The estimate of the eigenvalues of
$A_{ H}$ given above is sharp.  

In particular, the result implies the
following:

(i) If there is an integer  $k,\,2\leq
k\leq n$, such that $\theta_k(x)>c$
(respectively, $\theta_k(x)\geq c$) for a
Riemannian $n$-manifold $M$ at a point
$x\in M$, then, for any
isometric immersion of $M$ into $R^m(c)$,
 every eigenvalue of the shape operator
$A_{ H}$ is greater than 
${{n-1}\over n}$ (respectively, $\geq 0$),
regardless of codimension.

(ii) If $M$ is a compact hypersurface of $
E^{n+1}$ with $\theta_k \geq 0$
(respectively, with $\theta_k>0$) for a fixed 
$k,\,2\leq k\leq n$, then $M$ is embedded as a
convex (respectively, strictly convex)
hypersurface in $ E^{n+1}$. In particular, if
$M$ has constant scalar curvature, then $M$ is a
hypersphere of $ E^{n+1}$; according to a
result of W. S\"uss (1929) which states that
the only compact convex hypersurfaces with
constant scalar curvature in Euclidean space
are hyperspheres.

Statement (ii) implies 

(ii)$'$ If $M$ is a
compact hypersurface of
$ E^{n+1}$ with nonnegative Ricci curvature
(respectively, with positive Ricci curvature),
then $M$ is embedded as a convex (respectively,
strictly convex) hypersurface in $
E^{n+1}$. In particular, if
$M$ has constant scalar curvature, then $M$ is a
hypersphere of $ E^{n+1}$.

\subsection{Completeness of curvature surfaces}

Let $M$ be a hypersurface in a Euclidean space,  $N_0(x)$ the
null space of the second fundamental form of 
$M$ at $x\in M$, $k$ the minimal value of the
dimensions of the vector spaces $N_0(x)$ on
$M$, and $U$ the open subset of $M$ on
which this minimum occurs. Then $U$ is
generated by $k$-dimensional totally geodesic
submanifolds along which the normal space of
$M$ is constant. Moreover, if $M$ is complete,
then these generating submanifolds of
$U$ are also complete [Chern-Lashof 1957]. This result was generalized
in 1971 by D. Ferus to submanifolds of higher codimension in real
space forms. 

For an arbitrary principal curvature function 
$\lambda$ of an isometric immersion $f:M\to
R^m(c)$  of a Riemannian $n$-manifold $M$ into
a Riemannian $m$-manifold of constant sectional
curvature $c$, a similar result was obtained
by T. Otsuki (1970) and H. Reckziegel
(1976,1979). 

 Let $T^\perp M$ denote
the normal bundle and $T^*_\perp M$ its dual bundle. A 1-form $\mu\in
T^*_\perp M$ at $x\in M$ is called a principal curvature of $f$ at $x$
if the vector space 
$${\Cal E}(\mu)=\{X\in T_xM\,:\, A_\xi
X=\mu(\xi)\cdot X\hbox{ for all }\xi\in
T^\perp_xM\}$$ is at least 1-dimensional.

Suppose that there is given a continuous 
principal curvature function
$\lambda$ of $f$, that is, a continuous section of the bundle $T^*_\perp
M$ with $\dim {\Cal E}(\lambda_x)\geq 1$ for
all $x\in M$, and let $U$ be any open subset
of $M$ on which the function
 $x\mapsto \dim {\Cal E}(\lambda_x)$
is constant, say $\dim {\Cal
E}(\lambda_x)=k$ for all
 $x\in U$.  

H. Reckziegel (1979) obtained the following.

(i) The principal curvature form $\lambda$
is
$C^\infty$-differentiable on $U$;

(ii) The vector spaces ${\Cal
E}(\lambda_x),\, x\in U$, form a vector
subbundle of ${\Cal E}$ of the tangent bundle
$TM|_U$;

(iii) If $L$ denotes the foliation obtained 
by integrating ${\Cal E}$, and
$i:L\to M$ its inclusion, then all leaves of
 $L$ are $k$-dimensional
totally umbilical submanifolds of $M$ and 
$f\circ i:L\to R^m(c)$ is a
totally umbilical immersion into $R^m(c)$;

(iv) If $\lambda$ is covariant constant along 
${\Cal E}$, that is, if
$$(\nabla_X\lambda)(\xi):= X\cdot \lambda(\xi)-
\lambda(D_X\xi)=0$$
for all $X\in \Gamma({\Cal E})$ and $\xi\in
\Gamma(T^\perp M|_U)$, if furthermore
$\gamma:J\to L$ is a geodesic of $L$ with 
$\delta:=\sup J <\infty$, and
if $q:=\lim_{t\to\delta} \gamma(t)$ exists 
in $M$, then also $\dim
{\Cal E}(\lambda_q)=k$;

(v) If, in particular, $U$ is the subset of $M$ on which the function
$x\mapsto\dim {\Cal E}(\lambda_x)$ is
minimal (this subset is open, because
$x\mapsto \dim {\Cal E}(\lambda_x)$ is
upper-semicontinuous), $\lambda$ is
covariant constant along ${\Cal E}$, and $M$
is complete, then all the leaves of $L$ are
also complete spaces.

  If $k\geq	 2$, then $\lambda$ is always
covariant constant along ${\Cal E}$.

The leaves of $L$ are called the curvature surfaces of $f$ in $U$
corresponding to $\lambda$.

There exist submanifolds of codimension $\geq
2$ without any principal curvature. For
instance, the Veronese isometric embedding
of $RP^2\left(\sqrt{3}\right)$ in $S^4(1)$
mentioned in \S5.4.5 has no principal
curvature in the above sense.

\vfill\eject

\section{Rigidity and reduction theorems} 

\subsection{Rigidity}

An isometric immersion $f:M\to \tilde M$ is
called rigid if it is unique up to an isometry
of $\tilde M$, that is, if $f':M\to \tilde M$
is another isometric immersion, then there is
an isometry $\phi$ of $\tilde M$ such that
$f'=\phi\circ f$. If $f:M\to\tilde M$ is rigid,
then every isometry of $M$ can be extended to an
isometry of $\tilde M$.

F. Minding (1806--1885) conjectured in 1839
that a standard sphere a $E^3$ is rigid.
This conjecture was proved by H. Liebmann
(1874--1939) in 1899.   In 1929, S.
Cohn-Vossen (1902--1936) proved that a
closed convex surface in Euclidean 3-space
is rigid.  A. V. Pogorelov (1919--
) proved in 1951 that the requirement of
smoothness can be removed, by proving
that any closed convex surface, that
is, the boundary of a bounded convex
body, is uniquely determined up to a
rigid motion by its metric.

If $f:M\to  E^{n+1}$  is a hypersurface of
 Euclidean $(n+1)$-space, then at each point
$x\in M$, the type number of $f$ at $x$, denoted
by $t(x)$, is defined to be the rank of the
shape operator of $f$ at $x$.

A classical result of R. Beez (1876) states
that if $M$ is an orientable Riemannian
$n$-manifold and $f$ is an isometric
immersion of
$M$ into $ E^{n+1}$ such that the type
number  of $f$ is $\geq 3$ at every point
of $M$, then $f$ is rigid. 

R. Sacksteder (1962) obtained a number of
rigidity theorems for hypersurfaces. Among them
he proved that a complete convex hypersurface in
$ E^{n+1},\, n\geq 3$, is rigid if its type
number is at least 3 at one point.

D. Ferus (1970) proved that if $M$ is a
complete Riemannian $n$-manifold with $n\geq
5$ and if $f:M\to S^{n+1}$ is an isometric
immersion whose type number is everywhere $\geq
2$, then $f$ is rigid. 

C. Harle (1971) proved that if $M$ is a
Riemannian $n$-manifold, $n\geq 4$, with
constant scalar curvature $\rho\ne
 {1\over 2}n(n-1)c$ and $c\ne 0$, then every
isometric immersion of $M$ in a complete
simply-connected Riemannian space form
$R^{n+1}(c)$ is rigid. Y. Matsuyama (1976)
showed that  a hypersurface  with nonzero
constant mean curvature in a real hyperbolic
$(n+1)$-space $H^{n+1}$ with $n\geq 3$ is
rigid.

Given a system $\{A_1,\ldots,A_m\}$ of symmetric
endomorphisms of a vector space $V$ that are
linearly independent. The type number of the
system is defined to be the largest integer $t$
for which there are $t$ vectors $v_1,\ldots,v_t$
in $V$ such that the $mt$ vectors $A_r(v_i)$,
$1\leq  r\leq m,\, 1\leq i\leq t,$ are linearly
independent. When $m=1$, the type number of the
system of one single endomorphism $A$ is just
equal to the rank of $A$.

A generalization of R. Beez's result to higher
codimension was obtained in 1939 by C.
Allendoerfer (1911--1974).  Allendoerfer's
result states as follows: Let $f$ and $f'$
be two isometric embeddings of a Riemannian
$n$-manifold $M$ into $ E^{n+m}$. Assume
that, for a neighborhood $U$ of a point
$x_0\in M$, we have $(1)$ the dimensions of
the first normal space at $x\in U$ for both
$f$ and $f'$ are equal to a constant, say
$k$, and $(2)$ the type number of $f$ is at
least 3 at each point $x\in M$. Then there
exists an isometry $\phi$ of $ E^{n+m}$ 
such that
$f'=\phi\circ f$ on a neighborhood of $x_0$.

The standard $n$-sphere $S^n$ of
constant curvature one is rigid in
$ E^{n+1}$ when
$n\geq2$. However, it is not rigid in 
$ E^{n+2}$, since one can construct an
infinite-dimensional family of compositions
of isometric immersions $S^n\to 
E^{n+1}\to  E^{n+2}$. J. D. Moore (1996)
proved that this is essentially the only way
in which rigidity fails when $n\geq 3$.
Moreover, in this case he proved that any
isometric immersion of $S^n$ into  $
E^{n+2}$ is homotopic through isometric
immersions to a standard embedding into a
hyperplane.

E. Berger, R. Bryant and P. Griffiths (1983)
proved the following: Consider a local
 isometric embedding of a
Riemannian $n$-manifold $M$ into $E^{n+r}$.
Assume the embedding is ``general'' in the sense
that the second fundamental form lies in a
certain Zariski open subset of all such forms.
Then, if $r\leq n$ and $n\geq 8$,
or $r\leq 3$ and $n=4$, or $r\leq 4$ and
$n=5,6$, or $r\leq 6$ and $n=7,8$, then the
embedding is unique up to a rigid motion. 
This rigidity is not a consequence of
algebraic properties of the Gauss
equation, but depends rather on the
properties of the prolonged Gauss equations
involving the higher covariant derivatives
of the curvature tensor. 

By applying inequality (3.17), Chen (1996f)
established some rigidity theorems for
isometric immersions from some homogeneous
Riemannian $n$-manifolds into Euclidean space
$E^{n+k}$, regardless of codimension $k$,
under the assumption that the immersions
have the smallest possible squared mean
curvature.

\subsection{A reduction theorem }

Let $M$ be an $n$-dimensional submanifold of a
 Riemannian $m$-manifold $N$ and $E$  a
subbundle of the normal bundle $T^\perp M$.
Then $E$ is said to be parallel in the normal
bundle if, for each section $\xi$ of $E$
and each vector $X$ tangent to $M$, we have
$D_X\xi\in E$. 

The reduction theorem of J. Erbacher (1971)
states as follows:

Let $M$ be an $n$-dimensional submanifold of a
complete simply-connected Riemannian
$m$-manifold $R^m(c)$ of constant curvature
$c$. If there exists normal subbundle $E$ of
rank $\ell$ which is parallel in the normal
bundle and the first normal space
$N^1_x$, spanned by $\{h(X,Y):X,Y\in
T_xM\}$, is contained in $E_x$ for each $x\in
M$, then $M$ is contained in an
$(n+\ell)$-dimensional totally geodesic 
submanifold of $R^m(c)$.

\vfill\eject

\section{Minimal submanifolds}

The theory of minimal submanifolds is
closely related with the the theory of
calculus of variations. According to
historians, it is not quite certain when L.
Euler (1707--1783) began his study of the
calculus of variations. C. Carath\'eodory
(1873--1950) believed that it occurred
during his period in Basel with John
Bernoulli (1667--1748). Euler considered in
1732 and 1736 problems more or less arising
out the isoperimetric problems of  James
Bernoulli (1654--1705); and even as early
as the end of 1728 or early of 1729. In
effect, Euler in 1744, following John
Bernoulli, examined the equation of
end-curves that cut a family of geodesics
so that they have equal length. In his
famous 1744 book, L. Euler gave the first
systematic treatment of the calculus of
variations for curves. In this book he gave
a general procedure for writing down the
so-called Euler differential equation or
the first necessary condition, and to
discuss the principle of least action.

The history of the theory of minimal
surfaces goes back to J. L.
Lagrange (1736--1813) who  studied minimal
surfaces in Euclidean 3-space. In his famous
memoir, ``Essai d'une nouvelle m\'ethode pour
d\'eterminer les maxima et les minima des
formules int\'egrales 
 ind\'efinies'' which
appeared in 1760--1761, Lagrange developed his
 algorithm for the calculus of variations;
an algorithm which is also applicable in
higher dimensions and which leads to what is
known today as the Euler-Lagrange
differential equation. Lagrange communicated
his method in his first letter, dated 
August 12, 1755 when he was only nineteen, to L.
Euler who applauded his results (cf. [Euler
1755]). The basic idea of Lagrange ushered in a
new epoch in the calculus of variations. After
seeing Lagrange's work, Euler dropped his own
method, espoused that of Lagrange, and renamed
the subject the calculus of variations.

In his  famous memoir Lagrange also discovered
the minimal surface equation:
$$(1+z_y^2)z_{xx}-2z_xz_yz_{xy}+(1+z_x^2)z_{yy}=0.
$$ for a surface  defined by
$z=f(x,y)$ for $(x,y)$ in a domain of $E^2$
as the equation for a critical point of the
area functional. It was J. Meusnier
(1754--1793) in 1776 who gave a geometrical
interpretation of this equation as meaning
that the surface has vanishing mean
curvature function.

Before Lagrange, L. Euler had found in 1744
that a catenoid is a minimal surface; the
earliest nontrivial minimal surface
discovered which remained the only known
nontrivial minimal surface for over twenty
years, until J. Meusnier found in
1766 that a right helicoid is a minimal
surface.  It took almost ninety years until
H. F. Scherk (1798--1885) discovered further
minimal surfaces. Scherk found in 1834 that
the surface defined by $$z=\log(\cos
y)-\log(\cos x)$$ is a minimal surface, which
is  known today as Scherk's surface.  

In 1842 E. Catalan (1814--1894) proved that the
helicoid is the only ruled minimal surface
in
$E^3$. O. Bonnet (1819--1892) proved in 1860
that the catenoid is the only minimal surface
of revolution.

The theory of minimal surfaces experienced a
rapid development throughout of the
nineteenth century. The major achievements
of this period were presented in detail in
the 1903 book of L. Bianchi (1856--1928) and
the 1887 book of J. G. Darboux (1847--1917).
A  detailed account of more recent results
was given in the 1989 book of J. C. C.
Nitsche (1926--1996).

\subsection{First and second variational formulas}

Let $f:M \to \tilde M$ be an immersion
of a compact $n$-dimensional manifold $M$
(with or without boundary $\partial M$) into a
Riemannian $m$-manifold $\tilde M$.
Let $\{f_{t}\}$ be a one-parameter family of
immersions of $M \rightarrow \tilde M$ with
the property that
$f_{0}=f$. Assume the map $F: M\times [0,1]
\rightarrow \tilde M$ defined by $F(p,t)=f_{t}(p)$
is differentiable (we further assume $f_t=0$ on
 $\partial M$ when $\partial M\ne\emptyset$). 
$\{f_{t}\}$ is called a  variation of $f$. 

A variation of $f$ induces a vector field in
$\tilde M$ defined  along the image of $M$
under $f$, called the variational vector field.
We shall denote this field by
$\zeta$ and it is constructed as follows:

Let $\partial/\partial t$ be the standard vector field in
$M\times [0,1]$. We set
$\zeta(p)=F_{*}({\partial\over {\partial t}}
(p,0)).$
Then $\zeta$ gives rise to cross-sections 
$\zeta^T$ and
$\zeta^N$ in $TM$ and $T^{\perp}M$, respectively.
If we have $\zeta^{T}=0$, then $\{f_{t}\}$ is
called a  normal variation of $f$. 

For a given normal vector field $\xi$ on $M$,
exp$\,t\xi\,$ defines a normal variation
$\{f_{t}\}$ induced from $\xi$. We denote by
${V}_t$ the volume of $M$ under $f_t$ 
with respect
to the induced metric. 

The first variational formula is given by:
$$V'(\xi):={{d V_t}\over{d t}}|_{t=0}=-n\int_{M}
\left<\xi,H\right> dV_0.\leqno(5.1)$$ 

Hence, the immersion $f$ is minimal if and
only if ${{d V_t}\over{d t}}|_{t=0}=0$ for all
variations of $f$. Thus, a minimal submanifold
gives an extremal of the volume integral,
though not necessarily of the
least volume.

For a minimal submanifold $M$ of  a
Riemannian manifold $\tilde M$, the second
variational formula is given by
$${ V}''(\xi):={{d^2 V_t}\over{d
t^2}}|_{t=0}=\int_{M} \{||D\xi||^{2}-{\bar
S}(\xi,\xi)- ||A_{\xi}||^{2}\}dV,\leqno(5.2)$$ 
 where ${\bar S}(\xi,\eta)$  is
defined by $${\bar S}(\xi,\eta) =\sum_{i=1}^{n} {\tilde
R}(\xi,e_{i},e_{i},\eta),\leqno(5.3)$$
  $e_{1},\ldots,e_n$ is a local orthonormal frame
of $TM$, and ${\tilde R}$  is
the Riemann curvature tensor of the ambient 
manifold $\tilde M$.

Applying Stokes' theorem to the integral of
the first term of (5.2), we have
$${V}''(\xi)=\int_{M}\left<{\Cal
J}\xi,\xi\right>dV,\leqno(5.4)$$
 in which $\Cal J$ is a self-adjoint 
strongly elliptic linear differential
operator of the second order acting on the
space of sections of the normal bundle,
  given by
$${\Cal J}=-\Delta^{D}-{\hat A}-{\hat
S},\leqno(5.5)$$ where 
$\left<{}\right.{\hat  A}\xi,\left.
\eta\right>=\hbox{trace}\,(A_{\xi}A_{\eta})$, 
$\left<{}\right.{\hat S}\xi,\left.
\eta\right>={\bar S}(\xi,\eta)$ and
$\Delta^D$ is the Laplacian operator
associated with the normal  connection.

\subsection{Jacobi operator, index,
nullity and Killing nullity}

The differential operator ${\Cal J}$ defined
by (5.5) is called the Jacobi operator of
the minimal immersion $f:M\to\tilde M$. The
Jacobi operator  has discrete
eigenvalues $\lambda_{1} <\lambda_{2} <\ldots
\,\nearrow \infty.$ We put 
$$E_{\lambda}=\{\xi\in
\Gamma(T^{\perp}M)\, :\, {\Cal
J}(\xi)=\lambda\xi\,\}.$$ 

 A domain $D$, of a minimal submanifold
$M$, with compact closure is called stable if
the second variation of the induced volume of
$D$ is positive for all variations that
leave the boundary $\partial D$ of $D$ fixed.
The minimal submanifold $f:M\to\tilde M$ is
said to be stable if every such domain $D$
of $M$ is stable.

The number $i(f):=\sum_{\lambda
<0}\dim(E_{\lambda})$ is called the  index
of $f$ which measures how far the minimal
submanifold is from being stable.  

 A vector field $\xi$
in $E_0$ is called a Jacobi field.
The number $n(f):=\dim E_0$ is
called the nullity of $f$. 

Define a subspace $P$ of
$\Gamma(T^{\perp}M)$ by
$$P=\{\xi^N:\,\xi \hbox{ is  a Killing vector
field on}\; \tilde M\},$$ where $\xi^N$
denotes the component of $\xi$ normal to $M$.
Then
$P\subset E_0$. The dimension of
$P$ is called the Killing nullity of $f$,
which is denoted by $n_k(f)$.

\subsection{Minimal submanifolds of
Euclidean space}

An immersion $f:M\to
 E^m$ can be viewed as a $
E^m$-valued function. In this case, Beltrami's
formula relates $ H$ to $f$ by
$$\Delta f=-n H,\leqno(5.6)$$
where $\Delta$ is the Laplacian of $M$ with
respect to the induced Riemannian metric which
is defined by 
$$\Delta f=-\hbox{div}\,(\hbox{grad}\,f).$$

Beltrami's formula implies that each coordinate
function of $f$ is a harmonic function. Hence,
there exists no compact minimal submanifold
without boundary in Euclidean space.
 
This nonexistence result also follows from
the following fact: If $M$ is a  compact 
submanifold (without boundary) in a Euclidean
space, there always exits a point $x\in M$
such that the shape operator of
$M$ at $x$ is positive-definite with respect to
some unit normal vector at $x$. The point $x$
can be chosen to be the farthest point on $M$
from any fixed point in the Euclidean space.

It follows from the equation of Gauss and
Moore's lemma that a minimal isometric
immersion of a Riemannian product into a
Euclidean space is a product of minimal
immersions [Ejiri 1979b].

E. F. Beckenbach (1906--1982) and T. Rad\'o
 proved in 1933 that if the Gaussian
curvature of a surface is
$\leq 0$, then, for any simply-connected
domain $D$ on $S$, one has $L^2\geq 4\pi A$,
where $L$ is the arclength of the boundary
of $D$ and
$A$ is the area of $D$. In particular, for
any immersed simply-connected minimal surface
$M$ in $E^m$ with boundary
$C$, one has the following isoperimetric
inequality: $$L^2\geq 4\pi A,$$
where $L$ is the arclength of $C=\partial M$
and $A$ the area of $M$.

When one drops simple connectivity, the
isoperimetric inequality does not hold for
general domains on surfaces satisfying $K\leq
0$. For instance, on a circular cylinder in
$E^3$, the length of each boundary circle is
$2\pi r$ and the area is $2\pi r h$; thus the
area can be made arbitrarily large. 

The isoperimetric inequality $L^2\geq 4\pi A$
holds for domains $D$ lying on minimal
surfaces in $E^m$ in the following cases
[Osserman-Schiffer 1974; Osserman 1978]:

(1) the boundary of $D$ consists of a single
rectifiable Jordan curve;

(2) $D$ is doubly-connected and is bounded by
two rectifiable curves;

(3) $D$ is bounded by a finite number of
rectifiable curves lying on a sphere centered
at a point of $D$;

(4) $D$ is bounded by a finite number of
rectifiable Jordan curves, and minimizes area
among all surfaces with the same boundary.

L. P. Jorge and F. Xavier (1979) proved that
there are no bounded complete minimal
surfaces in $E^3$ with bounded Gaussian
curvature. N. Nadirashvili (1996) showed
that there exists an  immersed complete
bounded minimal surface in $E^3$ with
negative Gaussian curvature.

\vskip.1in
\noindent{\bf 5.3.1. Obstructions to minimal
isometric immersions}

If $M$ is a minimal surface in $ E^3$
with the induced metric $g$, then the Gauss
curvature $K$ of $M$ is $\leq 0$. Thus,
$\sqrt{-K}g$ defines a new metric on points
where $K\ne 0$. 

G. Ricci (1853--1925) proved in 1894
that a given metric $g$ on a plane
domain $D$ arises locally as the metric
tensor of a minimal surface in $ R^3$
if and only if the Gauss curvature $K$
of
$(D,g)$ is everywhere nonpositive and the
corresponding Gauss curvature $\bar K$ of
$\sqrt{-K}g$ vanishes at each point where
$K\ne 0$. 

 Let $g$ be the metric tensor of a minimal
surface $M$ in $ E^m$. If $g$ satisfies
Ricci's condition, then $g$
corresponds locally to the metric tensor of a
minimal surface $\hat M$ in $ E^3$. H.
Lawson (1971) proved that, in this case, 
 either $M$ lies in
$ E^3$ and belongs to a specific
one-parameter family of surfaces associated to
$\hat M$, or else $M$ lies in $E^6$ and
belongs to a specific two-parameter family of
surfaces obtained from $\hat M$, none of which
lie in any $ E^5$.

 Lawson's result implies
that the Ricci condition is an intrinsic
condition which completely characterizes
minimal surfaces lying in $E^3$ among all
minimal surfaces in $E^4$ or $E^5$;
and also the set of all minimal surfaces in
$ E^m$ isometric to a given minimal
surface in
$ E^3$ consists of a specific
two-parameter family of surfaces lying in
$E^6$.

For a minimal submanifold in Euclidean
space in general,  the equation of Gauss implies
that the Ricci tensor of a minimal submanifold
$M$ of a Euclidean space  satisfies
$$Ric(X,X)=-\sum_{i=1}^n |h(X,e_i)|^2\leq
0,\leqno(5.7)$$ where
$\{e_1,\ldots,e_n\}$ is an orthonormal local
frame field on $M$. Thus, the Ricci
tensor of a   minimal submanifold $M$ of a
Euclidean space is negative semi-definite
and, moreover, the minimal  submanifold is
totally geodesic if and only if its scalar
curvature vanishes identically. 

For minimal hypersurfaces in $ E^{n+1}$,
consider the metric $\hat g=-\rho g$ at the
point where the scalar curvature $\rho$ with
respect to the induced metric $g$ is negative. 
J. L. Barbosa and M. do Carmo (1978) proved
that the scalar curvature $\hat \rho$ of this
new metric $\hat g$ must satisfy the inequality
$\hat \rho \leq 2(n-1)-1$.

Besides the above necessary conditions,
inequality (3.17)  provides many further sharp
necessary conditions for a Riemannian
$n$-manifold to admit a minimal isometric
immersion in a Euclidean space, regardless of
codimension.

 In fact, inequality (3.17) of Chen implies
that, for a given Riemannian $n$-manifold
$M$, if there is a $k$-tuple
$(n_1,\ldots,n_k)\in \Cal S(n)$ such that the
$\delta$-invariant
$$\delta(n_1,\ldots,n_k)>0\quad \hbox{ at some
point } x\in M,\leqno(5.8) $$ then $M$ admits no
isometric minimal immersion in Euclidean
space for any arbitrary codimension. 

In particular, if a Riemannian manifold $M$
satisfies $\delta(2)>0$ at some point $x\in M$
(or equivalently, $\rho(x)>2(\inf K)(x)$ at some
point $x$), then $M$ admits no isometric minimal
immersion in Euclidean space, regardless of
codimension. 

There exist ample examples of Riemannian
manifolds with $Ric\leq 0$ which satisfy 
condition (5.8).

\vskip.1in
\noindent{\bf 5.3.2. Branched minimal
surfaces}

A branch point of a harmonic map $f:M\to
 E^m$ is a point $x\in M$ at which the
differential $(f_*)_x$ is zero. A harmonic
map $f:M\to E^m$ of a Riemann surface
$M$ is called a branched (or generalized)
minimal immersion if it is conformal except
at the branch points, and the image $f(M)$
is called a branched minimal surface.

Branched minimal surfaces have the following
three basic properties:

(1) {\bf Convex hull property}. Every
branched minimal surface with boundary in
$ E^m$ lies in the ``convex hull'' of
its boundary curve, that is, the smallest closed
convex set containing the boundary.

(2) {\bf Minimal principle}. If $M_1$
and $M_2$ are two branched minimal surfaces
in
$ E^3$ such that for a point $x\in
M_1\cap M_2$, the surface $M_1$ lies locally
on one side of $M_2$ near $x$, then $M_1$
and $M_2$ coincide near $x$.

(3) {\bf Reflection principle}. If the
boundary curve of a branched minimal surface
contains a straight line $L$, then the
surface can be analytically continued as a
branched minimal surface by reflection
across $L$. 

Based on reflection principle, H. Lewy
(1904--1988) proved in 1951 the
following: Let
$\Gamma$ be an analytic Jordan curve in
$ E^m$ and
$f:M\to  E^m$ a branched minimal
immersion with boundary $\Gamma$. Then $f$
is analytic up to the boundary, that is,
$f(M)$ is contained in the interior of a
larger branched minimal surface. 

S. Hildebrandt (1969) obtained the smooth
version of H. Lewy's theorem: If $f:M\to
E^m$ is a branched minimal immersion with smooth
boundary curve, then $f$ is smooth up to the
boundary.

\vskip.1in
\noindent{\bf 5.3.3. Plateau's problem}

The famous problem of Plateau states that
given a Jordan curve $\Gamma$ in
$E^3$ (or in $ E^m$) find a surface
of least area which has
$\Gamma$ as its boundary. Plateau's problem
was investigated extensively in the second
half of the nineteenth century by E. Betti
(1823--1892), O. Bonnet (1819--1892),  G.
Darboux (1842--1917), A. Enneper
(1830--1883), \'E. L. Mathieu (1835--1890),
H. Poincar\'e (1854--1912),  B. Riemann
(1826--1866), K. Weierstrass (1815--1897),
and others.  Plateau's problem was finally
solved independently by J. Douglas and T.
Rad\'o around 1930. 

The solution
to Plateau's problem given by Douglas and
Rad\'o is a branched minimal surface.
More precisely, they proved that if $\Gamma$
is a rectifiable Jordan curve in $ E^m$
and $D=\{(x,y)\in  E^2: x^2+y^2<1\}$,
then there exists a continuous map $f:\bar
D\to E^m$ from the closure of $D$
into $ E^m$ such that 

(a) $f|_{\partial D}$ maps
homeomorphically onto $\Gamma$, 

(b) $f|_D$ is a harmonic map and almost
conformal, that is $\left<f_x,f_y\right>=0$
and
$|f_x|=|f_y|$ in $D$ with $|df|>0$ except at
isolated branch points, and 

(c) the induced area of $f$ is the least
among the family of piecewise smooth surfaces
with $\Gamma$ as their boundary.

The map $f$ given above is called the
classical solution or the Douglas-Rad\'o
solution to Plateau's problem for $\Gamma$.
The resulting surface $M$ is a branched
minimal disk.

A branched minimal disk $M$ bounded by a
smooth curve $\Gamma$ in $ E^3$
satisfies the following formula of
Gauss-Bonnet-Sasaki-Nitsche:
$$1+\sum (m_\alpha-1)+\sum M_\beta
+{1\over{2\pi}}\int_M |K|dA\leq
{1\over{2\pi}}\kappa(\Gamma),$$
where $m_\alpha-1$ denote the orders of the
interior branch points, $2M_\beta$ the
orders of the boundary branch points which
must be even, $K$ the Gaussian curvature,
and $\kappa(\Gamma)$ the total curvature of
$\Gamma$.

R. Osserman (1970) proved that every
classical solution to Plateau's problem in
$E^3$ is free of branch points in its
interior. Thus, it is a regular immersion.

If  $\Gamma$ is a curve in $ E^3$ which
is real analytic or a  smooth curve with 
total curvature $\kappa(\Gamma)<4\pi$, then
the minimal disk of least area with $\Gamma$
as its boundary has no boundary branch
points [Gulliver-Leslie 1973].

A  Douglas-Rad\'o solution is not
necessarily an embedding. In fact, if
$\Gamma$ is knotted in $ E^3$, then
every solution must have self-intersections.
However, immersed minimal disks of least
area in $ E^3$ which can self-intersect
only in their interiors are embeddings.
Also, there exists an unknotted Jordan curve
which never bounds an embedded minimal disk.
 Meeks and Yau (1982) proved
that if the Jordan curve lying entirely on
the boundary of a convex body, then a
Douglas-Rad\'o solution is embedded (see also
[Almgren-Simon 1979; Gr\"uter-Jost 1986]).

In general a Douglas-Rad\'o
solution to Plateau's problem does not have
the uniqueness property. However, Rad\'o
(1930) provided a sufficient condition on the
boundary curve for which the solution is
always unique. More precisely, he showed that
if a Jordan curve $\Gamma$ in $ E^m$ admits a
one-to-one orthogonal projection onto a
convex curve in a plane $E^2$ in $E^m$, then
the classical  solution for $\Gamma$ is free
of branch points and can be expressed as the
graph over this plane. Furthermore, when
$n=3$, the solution is unique.

J. C. C. Nitsche (1989) proved that if
$\Gamma$ is an analytic Jordan curve in
$E^3$ with total curvature
$\kappa(\Gamma)\leq 4\pi$ or a smooth curve
with $\kappa(\Gamma)<4\pi$, then $\Gamma$
bounds a unique immersed minimal disk.

 A. J. Tromba (1977) proved the following: (1)
Any rectifiable Jordan curve $\Gamma$
``sufficiently close'' to a plane curve has
a unique simply-connected minimal surface
spanning it, and (2) If $\Cal F$ denotes the
set of embeddings $f$ of $S^1$ into $E^3$
with the property that for the curve
$f(S\sp 1)$ every simply connected
minimal surface spanning
$f(S^1)$ is free from branch
points, then there is an $\Cal F$ of
embeddings for which the number of minimal
surfaces spanning the image is finite. 

F. Tomi (1986) showed that an analytic Jordan
curve in $E^3$ bounds only finitely many
minimal disks and Nitsche proved that a
smooth Jordan curve $\Gamma$ with total
curvature $\leq 6\pi$ also bounds only
finitely many minimal disks. On the other
hand, P. L\'evy and R. Courant constructed an
example of rectifiable Jordan curve that is
smooth with the exception of one point and it
bounds uncountably many minimal discs.

F. Morgan (1978) and A. J. Tromba
(1977) proved that, generically, there are at
most finitely many minimal surfaces with a
given boundary. More precisely, in the space
$\Cal A$ of all smooth Jordan curves in
$ E^m$ with suitable topology, there
exists an open and dense subset $\Cal B$
such that for any
$\Gamma$ in $\Cal B$, there exists a unique
area-minimizing minimal disk.

H. Iseri (1996) studied the Plateau's
problem in which the boundary curves  may
have self-intersections and provided a
condition which guarantees that the 
minimal surfaces they generated will not be
degenerate.

The method was carried further in 1948  by C.
B. Morrey (1907--1984) for  Plateau's
problem in a complete Riemannian manifold
which is metrically well behaved at
infinity, and includes the class of compact
or homogeneous Riemannian manifolds.

Morrey's setting of the generalized
Plateau's problem is as follows: Suppose a
homotopically trivial rectifiable Jordan
curve $\Gamma$ is given in a Riemannian
$m$-manifold $N$. Let $D$ denote a disk in
$ E^2$. Find a mapping $f:\bar D\to N$
such that (i) $f$ maps ${\partial D}$
homeomorphically onto $\Gamma$ and (ii) the
induced area of $f$ is the least among the
class of piecewise smooth surfaces in $N$
bounded by $\Gamma$ satisfying (i).

Morrey gave a solution under the assumption
that $N$ is homogeneously regular, that is,
there exist $0<k<K$ such that, for any point
$y\in N$, there is a local coordinate system
$(U,\Psi)$ around $y$ for which 
$\Psi(U)=\{x\in E^m: ||x||<1\}$ and the
Riemannian metric $g=\sum g_{ij}dx^i dx^j$
satisfies $$k\sum v_i^2\leq \sum
g_{ij}v_iv_j \leq K\sum v_i^2$$ for any
$x$ and $v=(v_1,\ldots,v_m)\in E^m$. 

C. B. Morrey (1907--1984) proved in 1948 that
if $N$ is a homogeneously regular Riemannian
manifold and if $\Gamma$ is a homotopically
trivial rectifiable Jordan curve in $N$,
then there exists a branched minimal
immersion $f:D\to M$ with least area bounded
by $\Gamma$ such that  $f|_{\partial D}$ maps
homeomorphically onto $\Gamma$. If $\Gamma$
is smooth, then so is the solution $f$ up to
the boundary. Furthermore, when $m=3$, the
solution $f$ is an immersion in its
interior. If $N$ and $\Gamma$ are real
analytic and $m=3$, then the solution $f$ is
an immersion up to the boundary.

M. Ji and G. Y. Wang (1993) investigated 
disk type minimal surfaces spanned by a given
Jordan curve $\Gamma$ in $N$ and proved
the following:  

(1) Each smooth Jordan curve $\Gamma $ in
$S\sp {n}$ bounds at least two minimal
surfaces, sometimes infinitely many ones; and

(2) Let $N$ be a compact oriented Riemannian
manifold embedded in a Euclidean
$m$-dimensional space. Suppose that the
compact manifold $N$ admits no minimal
sphere. If there are two strictly stable
minimal disks bounded by a Jordan curve
$\Gamma $, then there exists another
minimal surface bounded by $\Gamma$.
\vskip.1in

Given two Jordan curves $\Gamma_1,\Gamma_2$
in $E^3$, does $\Gamma=\Gamma_1\cup\Gamma_2$
bound a minimal annulus? This is called the
Douglas-Plateau problem which is a
generalization of the original Plateau
problem [Douglas 1931b]. 

In many cases the answers to the
Douglas-Plateau problem are negative. One
example is that of two coaxial unit circles
$C_1$ and $C_2$. If the distance between
their centers is large, then $C_1\cup C_2$
cannot bound a minimal annulus. When
$\Gamma_1$ and $\Gamma_2$ are smooth convex
planar Jordan curves lying in parallel but
different planes, the Douglas-Plateau problem
has a satisfactory answer. 

In fact,  [Meeks-White 1991] proved the
following: If $\Gamma_1$ and $\Gamma_2$ are
smooth convex planar Jordan curves lying in
parallel but different planes, then exactly
one of the following three cases occurs:

(a) There are exactly two minimal annuli
bounded by  $\Gamma=\Gamma_1\cup\Gamma_2$,
one is stable and one is unstable;

(b) There is a unique minimal annulus $M$
bounded by $\Gamma$; it is almost stable in
the sense that the first eigenvalue of the
Jacobi operator of $M$ is zero;

(c) There are no minimal annuli bounded by
$\Gamma$. 

If $M$ is a minimal annulus bounded
by $\Gamma=\Gamma_1\cup\Gamma_2$, then the
symmetry group of $M$ is the same as the
symmetry group of $\Gamma$.

R. Hardt and L. Simon (1979) proved that if
$\Gamma$ is the union of any finite
collection of disjoint smooth Jordan curves
in $E^3$, then there exists a compact
embedded minimal surface with boundary
$\Gamma$ which is smooth up to the boundary.

In 1985 J. Jost proved the existence
of a minimal surface of finite prescribed
genus and connectivity spanning a
configuration of oriented Jordan curves in
a homogeneously regular manifold in the sense
of Morrey, under the condition that
the area infimum over combinations of
surfaces of this topological type is
strictly less than the one over combinations
of surfaces of lower genus or connectivity. 

Although the same problem as Jost's was
treated in classical papers by J. Douglas, R.
Courant, and M. Shiffman in the 1930s, there
was some doubt about the validity of their
proofs. Namely, in order to get the
compactness of a minimizing sequence, Douglas
compactified the moduli space of surfaces of
the topological type considered, but he did
not show that the boundary of this
compactification consists of surfaces of
lower genus or connectivity. Courant
provided a complete proof of the case of
higher connectivity but genus zero, his
considerations about the case of higher genus
were pointed out by A. J. Tromba to be too
vague and not detailed enough to be accepted
as a correct proof. Shiffman assumed a priori
a condition which is equivalent to the
compactness of a minimizing sequence and
therefore could only prove a weaker
statement. 

In 1993 F. Bernatzki treated the same
Plateau-Douglas problem as Jost's for
nonorientable surfaces, using a method
similar to Jost's. 

\vskip.1in
\noindent{\bf 5.3.4. Weierstrass' representation formula}

In 1866 Weierstrass gave
a general formula to express a
simply-connected minimal surface in terms of a
complex analytic function $f$ and a
meromorphic function $g$ with certain
properties. His formula allows one to
construct a great variety of minimal surfaces
by choosing these functions. 

 Weierstrass' representation formula states as
follows: Every simply-connected
(branched) minimal surface in $ E^3$ is
represented in the form:
$$x_k(z)=\hbox{Re}\,\left\{ \int_0^z
\phi_k(\zeta)d\zeta\right\}+c_k,\quad
k=1,2,3,\leqno(5.9)$$
where $\phi_1={1\over 2}f(1-g^2),\,
\phi_2={1\over 2}\sqrt{-1}f( 1+g^2),\,
\phi_3=fg$, and $c_k$ are constants. Here
$g(z)$ is a meromorphic function on 
$D$ (= the unit disc or
the entire  complex plane), $f(z)$ is an
analytic function on $D$
satisfying the
property that at each point $z$, where g(z)
has a pole of order $m$ and $f(z)$ has a zero
of order $2m$. 

For instance, Enneper's surface with
coordinates 
$$(x_1,x_2,x_3)=(u-{{u^3}\over
3}+uv^2,v-{{v^3} \over 3}+vu^2,u^2-v^2),\quad
(u,v)\in E^2\leqno(5.10)$$
 is obtained from the Weierstrass formula by
setting $f=1$ and $g(z)=z$. The Richmond
surface, the catenoid, and helicoid are the
minimal surfaces obtained from the
Weierstrass formula by setting $(f,g)=
(z^2,z^{-2}), ({1\over 2}z^{-2},z)$, and
$(-ie^{-z},e^z)$, respectively. 

A minimal surface described by $(f,g)$
via Weierstrass's representation formula
has an associated one-parameter family of
minimal surfaces given respectively by
$(e^{it}f,g)$. Two surfaces of the family
described by $t_0$ and $t_1$ are called 
adjoint by O. Bonnet (1853a,1853b) if
$t_1-t_0={\pi\over 2}$. The catenoid and
the helicoid are a pair of adjoint minimal
surfaces for a suitable choice of
constants. 

All the surfaces of an associated
family are locally isometric. Conversely,
H. A. Schwarz (1890) proved that if a
simply-connected minimal surface $S_1$ is
isometric to a simply-connected
minimal surface $S$, then $S_1$ is congruent
to an associate minimal surface of $S$.

\vskip.1in
\noindent{\bf  5.3.5. Bernstein's problem}

 Let $x_1,\ldots,x_n,z$ be standard
coordinates in $ E^{n+1}$. Consider
minimal hypersurfaces which can be
represented by an equation of the form
$$z=z(x_1,\ldots,x_n)$$ for all $x_i$, that is,
which have a one-to-one projection onto a
hyperplane. 

 The answer to
Bernstein's problem is  known to be
affirmative in the following cases: 

(a) $n=2$ [Bernstein 1914];

(b) $n=3$ [de Giorgi 1965];

(c) $n=4$   [Almgren 1966]; and 

(d) $n=5,6,7$ [Simons 1968]. 

\noindent And the answer is negative for
$n\geq 8$ [Bombieri-de Giorgi-Giusti 1969].

In general, when a bounded domain $D$ in
$ E^n$ and a continuous function $\phi$
on its boundary $\partial D$ are given, the
problem of finding a minimal hypersurface
$M$ defined by the graph of a real-valued
function $f$ on $\bar D$, the closure of $D$,
with $f|_{\partial D}=\phi$ gives rise to a
typical Dirichlet problem. The basic
questions are those of existence, uniqueness,
and regularity of solutions. These problems
were studied by T. Rad\'o for $n=2$ and later
by L. Bers, R. Finn, H.
Jenkins, J. Serrin, R. Osserman, and others.

\vskip.1in
\noindent{\bf 5.3.6. Periodic minimal
surfaces and minimal surfaces with
many symmetries}

 A minimal surface in $ E^3$ is called
periodic if it is invariant under a group
$G$ of  isometries that acts freely on $
E^3$.

The Gauss-Bonnet theorem implies that if $M_g$
is a  minimal surface of genus $g$ in a 3-tours
$T^3$, then its Gauss map $G:M_g\to S^2$
represents $M_g$ as a $(g-1)$ conformal
branched covering of $S^2$. Thus, a surface
of genus 2 is never periodic, and a minimal
surface of genus $g$ in a 3-torus $T^3$ has
$4(g-1)$ zeros of Gaussian curvature, counted
with multiplicities (cf. [Meeks 1993]).

In 1867 H. A. Schwarz established a procedure
for generating periodic minimal surfaces in
$ E^3$ using octahedral or tetrahedral
symmetry. A. H. Schoen (1970) constructed
infinitely periodic minimal surfaces in $
E^3$ without self-intersections. Among them is
a surface containing no straight lines built
out of an infinite number of congruent
curvilinear hexagons whose sides form a family
of curves which are almost, but not exactly,
circular helices. A. H. Schoen (1970) also described
some triply periodic minimal
surfaces in $E^3$. H. Karcher (1989) 
gave a geometric description of the
construction of Schoen's surfaces, by
solving a conjugate Plateau problem
for a polygonal contour. Karcher
showed that the solution of an
analogous Plateau problem in
$S\sp 3$ provides a deformation of
the simpler of these minimal surfaces
into constant mean curvature
surfaces. The spherical polygon can
be obtained from the Euclidean
contour of the conjugate minimal
surface by considering the straight
lines in $E^3$ as integral
curves of parallel vector fields and
the contour in $S\sp 3$ as formed by
integral curves of right-invariant
vector fields. 

In 1978 T. Nagano and B. Smyth gave a
construction procedure to construct periodic
minimal surfaces in $ E^m$ or $m$-tori with
symmetry corresponding to a Weyl group of
rank $m$.

A surface is said to have finite topology  if
it is homeomorphic to a closed surface with a 
finite number of points removed.
W. H. Meeks and H. Rosenberg (1993) proved that
a properly embedded minimal surface in a
3-torus $T^3$ has finite total curvature if
and only if it has finite topology.  They
also proved that the plane and the helicoid
are the only properly embedded
simply-connected minimal surfaces in
$ E^3$ with infinite symmetry group.

D. Hoffman, F. Wei and H.
Karcher (1993) constructed a complete
embedded singly periodic minimal
surface in $E^3$ that is asymptotic to
the helicoid, has infinite genus and
whose quotient by translations has
genus one. The quotient of the
helicoid by translations has genus
zero and the helicoid itself is
simply-connected. Using  the techniques of N.
Kapouleas, M. Traizet (1996)  constructed simply
periodic minimal surfaces in Euclidean 3-space
by glueing together Scherk surfaces.
 F. J. L\'opez, M. Ritor\'e and F. Wei
(1997) found all the properly immersed
minimal tori with two parallel embedded
planar ends.

W. Fischer and E. Koch (1996)
classified triply periodic minimal surfaces
containing straight lines, by using their
associated crystallographic groups or, more
precisely, by group-subgroup pairs with
index 2. C. Frohman (1990) proved that if
$F$ and $F'$ are triply periodic
minimal surfaces in $E^3$,
then there is a homeomorphism
$h: E^3\to E^3$ such that $h(F)=h(F')$. 

M. Callahan, D. Hoffman and W. Meeks (1990)
proved that a properly
embedded minimal surface with more than
one end and with infinite symmetry group
 is either the catenoid or has an
infinite number of flat ends and is
invariant under a screw motion. 
They also established the existence of a
family of complete embedded minimal surfaces
$M_{k,\theta}$ invariant under a rotation of
order $k+1$ and a screw motion of angle
$2\theta$ about the same axis, where $k>0$ is
any integer and $\theta$ is any angle with
$\vert\theta\vert <\pi/(k+1)$. In 1993
Callahan, Hoffman and Karcher gave an
explicit construction of these surfaces using
generalized Weierstrass
representation; generalized in the sense of
using the logarithm derivative of the Gauss
map rather than the Gauss map itself as in
the usual Weierstrass representation.

Further results on periodic minimal surfaces
were  obtained by  Nagano and
Smyth (1975,1976,1978,1980), Meeks and
Rosenberg (1989,1993), Meeks (1990,1993), F.
 Wei (1992),  J. Hass, J. T. Pitts and J.
H. Rubinstein (1993), and others.

\vskip.1in
\noindent{\bf 5.3.7. Ruled minimal
submanifolds}

In 1835 P. Scherk tried unsuccessfully to
determine all ruled minimal surfaces in $E^3$,
that is, those minimal surfaces which contain
a straight line through each point of the
surfaces. The problem was finally solved by
E. Catalan in 1842, who proved that the
helicoid is the only nonplanar ruled minimal
surface in $E^3$.

An $n$-dimensional submanifold $M^n$ of a
Riemannian manifold $\tilde M$ is called 
ruled if $M^n$ is foliated by 
$(n-1)$-dimensional totally geodesic
submanifolds of
$\tilde M$.

\"U. Lumiste (1958) showed that an
$n$-dimensional minimal ruled submanifold of
Euclidean space is either 

(a) generated by an $(n-1)$-dimensional affine
subspace $P$ under a screw motion in $
E^{2n+1}$ such that the axis cuts $P$
orthogonally, or

(b) generated  by an $(n-1)$-dimensional affine
subspace $P$ under a rotation in $ E^{2n}$
around a point in $P$, or

(c) a cylinder on a submanifold of the type 
(i) or (ii).

Analytically any ruled minimal submanifold
therefore can be given by
\vskip.03in
$$\aligned
X(s,t_1,\ldots,t_{n-1})=(t_1\cos(a_1s),t_1
\sin(a_1s),\ldots,\\ t_k\cos(a_ks),t_k
\sin(a_ks),t_{k+1},\ldots,t_{n-1},bs)
 \endaligned$$
where $a_1,\ldots,a_k$ and $b$ are real numbers. 

A submanifold with this kind of parameterization is called a
 generalized helicoid. 

If $b\ne0$ (respectively,
$b=0$), then this gives a cylinder on a submanifold of the
type (i) (respectively, of the type (ii)). 
A ruled submanifold of the type (ii) is a cone on a minimal
ruled submanifold of some hypersphere of the Euclidean space.

 J. M. Barbosa, M. Dacjzer and L. P. Jorge (1984)
proved that any minimal ruled submanifold is generated by an
affine subspace $P$ under a one-parameter
subgroup $A$ of rigid motions of the
Euclidean space such that $P$ is orthogonal
to the orbits of $A$. Then they showed that
the resulting submanifold (at least if it is
minimal) has the same parameterization. They
also extended their result to ruled
submanifolds of real space forms. 

Complete ruled minimal hypersurfaces of Euclidean space were
classified by D. E. Blair and J. R. Vanstone
(1978); and the classification of general ruled
minimal hypersurfaces of Euclidean space were
done by G. Aumann (1981).

\vskip.1in
\noindent{\bf  5.3.8. Minimal immersions of K\"ahler manifolds}

The theory of minimal surfaces in Euclidean
space profits substantially from the study of
the underlying complex structure. Thus, it is
natural to study minimal immersions
$f:M^{2n}\to E^{2n+p}$ of K\"ahler manifolds
into Euclidean space. 

K\"ahler manifolds which
are isometrically immersed into Euclidean space
as real hypersurfaces are called real
K\"ahler hypersurfaces. Such hypersurfaces
have been studied by T. Takahashi (1972),
P. Ryan (1973), K. Abe (1974), M. Dajczer
and D. Gromoll (1985),  M. Dajczer and L.
Rodriguez (1986,1991), H. Furuhata (1994),
and others. For instance, K. Abe (1974)
proved that if
$M^{2n}$ is a complete real K\"ahler
hypersurface of Euclidean $(2n+1)$-space, then
$M^{2n}$ is a product of a Riemann surface and
$ C^{n-1}= E^{2n-2}$, provided either 

(a) $M^{2n}$ has nonnegative scalar curvature,
or

(b) $M^{2n}$ has strictly negative scalar
curvature, or 

(c) the immersion is real analytic.

Minimal real K\"ahler hypersurfaces are abundant
and have been classified by M. Dajczer and D.
Gromoll (1985). It turns out that none of them
is complete unless $M^{2n}=M^2\times 
E^{2n-2}$ and $f=f_1\times \hbox{id}$ splits,
where $f_1:M^2\to  E^3$ is a complete
minimal surface.

An isometric immersion $f:M^{2n}\to 
E^{2n+2}$ is said to be complex ruled if
$M^{2n}$ is a K\"ahler manifold and admits a
continuous codimension two foliation such that
any leaf is a K\"ahler submanifold of $M^{2n}$
and whose image under $f$ is an affine subspace
of $ E^{2n+2}$. $f$ is called completely
complex ruled if in addition the leaves are all
complete Euclidean $ E^{2n-2}$ spaces. 

We say that the scalar curvature $\rho$ of a
complete manifold has subquadratic growth 
along geodesics if its growth along any
geodesic is less than any quadratic
polynomial in the parameter. 

M. Dajczer and L. Rodriguez (1991) investigated
minimal immersions of complete K\"ahler
manifolds of codimension two in Euclidean
space and obtained the following:

Let $f:M^{2n}\to E^{2n+2},\,n\geq 2$, be a
minimal immersion of a complete K\"ahler
manifold. Then one of the following occurs:

(i) $f$ is a holomorphic;

(ii) $f$ is completely complex ruled; 

(iii) $M^{2n}=M^4\times  E^{2n-4}$ and
$f=f_1\times\hbox{id}$.

Moreover,

(a) if the scalar curvature of $M^{2n}$
has subquadratic growth along geodesics, then
$f$ is of type (i) or (iii), and

(b) if the index of relative nullity of $f$ is
$\geq 2n-4$ everywhere, then  $f$ is of type
(i) or (ii).

M. Dajczer and L. Rodriguez (1986)
also investigated minimal immersions of 
K\"ahler manifolds into Euclidean
space of higher codimension and obtained the
following: 

(1) Let $f:M^{2n}\to  E^{2n+p}$ be
an isometric immersion of a K\"ahler manifold.
If the type number of $f$ is $\geq 3$
everywhere, then $f$ is holomorphic.

(2) Let $f:M^{2n}\to  C^{n+1}$ be
a holomorphic isometric immersion of a K\"ahler
manifold. If $g: M^{2n}\to E^{2n+p}$ is a
minimal isometric immersion, then $g$ is
congruent to $f$ in $ E^{2n+p}$.

(3) Let $f:M^{2n}\to C^{n+q}$ be
a full (or substantial) isometric immersion
of a K\"ahler manifold with type number $\geq
3$ everywhere. If $g: M^{2n}\to E^{2n+p}$ is
a minimal isometric immersion, then $g$ is
congruent to $f$ in $E^{2n+p}$.

(4)  Let $f:M^{2n}\to R^{2n+p}(c)$ be
a minimal isometric immersion of a K\"ahler
manifold into a Riemannian manifold of constant
curvature $c$. 

(4.1) If $c<0$, then $n=1$.

(4.2) If $c=0$, then $f$ is circular, that is,
the second fundamental form of $f$ satisfies
$h(X,JY)=h(JX,Y)$ for all $X,Y\in TM$.

(4.3) If $c>0$, then the Ricci curvature
$Ric_M\leq nc$, with equality implying that the
second fundamental form is parallel.

M. Dajczer and D. Gromoll (1995) proved that if
$f:M^{2n}\to  E^{2n+2},\, n\geq 3$ is a
minimal isometric immersion of a complete
K\"ahler manifold and if $f$ is irreducible and
not holomorphic, then $M^{2n}$ contains an open
dense subset $M^*$ on which $f$ is completely
holomorphic ruled. Furthermore, along any
holomorphic section, $f$ has a ``Weierstrass
type representation''.

 Dajczer and Gromoll (1985) also
proved that if an isometric minimal immersion
$f$ of a simply-connected K\"ahler manifold into
a Euclidean space  is not holomorphic, then
there is a one-parameter family,
called the associated family, of non-congruent
isometric minimal immersions with the same
Gauss map. Moreover, the immersion can always
be made the real part of a holomorphic
isometric immersion,  called the
holomorphic representative of $f$.

H. Furuhata (1994) gave a parametrization of
the set of isometric minimal immersions of a
simply-connected K\"ahler manifold into a
Euclidean space by a set of certain complex
matrices, which is described in terms of a full
isometric holomorphic immersion of the K\"ahler
manifold into a complex Euclidean space.
Furuhata's result is an extension of a result
of E. Calabi (1968) on minimal surfaces.

\subsection{Minimal submanifolds of
spheres} 

A submanifold $M$ of a Euclidean
$m$-space is contained in a hypersphere as
a minimal submanifold if and only if it is a
pseudo-umbilical submanifold
with nonzero parallel mean curvature vector
[Yano-Chen 1971]. A submanifold of
codimension two in a Euclidean space is
contained in a hypersphere as a minimal
submanifold if and only if it is a
pseudo-umbilical submanifold with nonzero
constant mean curvature [Chen 1971]. 

By a pseudo-umbilical submanifold we mean a
submanifold of a Riemannian manifold whose
shape operator in the direction of the mean
curvature vector is proportional to the
identity transformation.

The geometry of minimal submanifolds of a
sphere $S^n$ takes a quite different course
than in the Euclidean case, because
there do exist many compact  minimal
submanifolds in spheres.

\vskip.1in
\noindent{\bf 5.4.1. Necessary conditions}

 Similar to the Euclidean
case, there are some necessary conditions
for a Riemannian
$n$-manifold to admit an isometric minimal
immersion in the unit
$m$-sphere $S^m$. In fact,
the equation of Gauss implies that the Ricci
tensor of a minimal submanifold in $S^m$
satisfies $Ric\leq (n-1)g$.

 Inequality (3.17) of Chen provides many
further necessary conditions. In fact, (3.17)
implies that, regardless of codimension, if a
Riemannian $n$-manifold  admits an isometric
minimal immersion in a unit sphere, it must
satisfies 
$$\delta(n_1,\ldots,n_k)\leq
{1\over 2}n(n-1)-{1\over 2}\sum_{j=1}^k
n_j(n_j-1)\leqno(5.11)$$ for any $k$-tuple
$(n_1,\ldots,n_k)\in \Cal S(n)$.

Since the center of gravity of a compact
minimal submanifold of $S^m$ is exactly the
center of the $S^m$, where
$S^m$ is viewed as an ordinary hypersphere in
$ E^{m+1}$, there exists no compact minimal
submanifold in $S^m$ which is contained in an
open hemisphere of $S^m$. 

\vskip.1in
\noindent{\bf  5.4.2. Takahashi's theorem}

 A fundamental result of T. Takahashi
(1933-- ) obtained in 1966 states that  an
isometric immersion $f$ of a Riemannian
$n$-manifold $M$ in the unit $m$-sphere
$S^m$, viewed as a vector-valued function
in $ E^{m+1}$, is minimal if and only if
$\Delta f=nf$. 

An immediate
application of Takahashi's theorem is that any
compact $n$-dimen\-sional homogeneous Riemannian
manifold whose linear isotropy group is
irreducible can be minimally isometrically
immersed into the $m$-sphere of curvature
$\lambda/n$, corresponding to any nonzero
eigenvalue $\lambda$, where $m+1$ is the
dimension of the corresponding eigenspace.

\vskip.1in
\noindent{\bf 5.4.3. Minimal isometric immersions of spheres into spheres}

M. do Carmo and N. R. Wallach (1971) proved
that if an $n$-sphere of constant curvature
$c$ is minimally isometrically immersed into
the unit $m$-sphere, but not in any great
hypersphere, then, for each non-negative
integer $k$, we have
$$c={n\over{k(n+k-1)}},\quad m\leq
(n+2k-1)\left({{n+k-2)!}\over{k!(n-1)!}}
\right)-1.$$ 
The immersion is rigid if  $n=2$ [E.
Calabi, 1967] or $k\leq 3$. 

For general $n$,  do Carmo and Wallach showed
that the space of minimal isometric immersions
from $S^n(1)$ into
$S^m(r)$ can be parametrized by a compact 
convex body in some finite-dimensional vector
space. The immersions corresponding to interior
points of this convex body all have images that
are embedded spheres or embedded real projective
spaces. 

A spherical space form is a compact
manifold of positive constant sectional 
curvature.  D. DeTurck and W. Ziller (1992)
proved that every homogeneous spherical space
form admits a minimal isometric embedding into
some sphere. C. M. Escher (1996) gave a
necessary condition for the existence of a
minimal embedding of nonhomogeneous
3-dimensional spherical space forms. In
particular, she showed that the lens space
$L(5,2)$ cannot be minimally embedded into
any sphere. 

G. Toth (1997) provided a general method that
associates to set of spherical minimal
immersions from $S^n$ a spherical minimal
immersion from $S^{n+1}$. In particular, he
proved the following: Let $n\geq 3$ and $p\geq
4$. Given full spherical minimal immersions
$f_i:S^n\to S^{m_i}$, $i=1,\ldots, p$, there
exists a full spherical minimal immersion
$\tilde f:S^{n+1}\to S^N$, where
$N=\sum_{i=1}^p (m_i+1)$.

\vskip.1in
\noindent{\bf 5.4.4. Minimal surfaces in spheres}

In contrast to Euclidean case, there exist
many compact minimal submanifolds in spheres.
In fact, H. B. Lawson (1970) proved the
following:

 (1) any compact surface of any
genus, except the real projective plane, can
be minimally immersed in $S^3$; 

(2) there exist minimal immersions of every
surface of negative Euler characteristic
into $S^5$ such that none of the images
lies in a totally geodesic $S^4$; and 

(3) there is a countable
family of minimal immersions of the torus
into $S^4$ where none of the images lies in
a totally geodesic $S^3$. 

Many further examples of minimal surfaces
in $S^3$ have also been constructed in 
[Karcher-Pinkall-Sterling 1988,
Pitts-Rubinstein 1988].

F. J. Almgren (1933--1977) proved in 1966
that the only minimal immersion from
$S^2$ into $S^3$ is the totally geodesic
one. H. I. Choi and R. Schoen (1985) showed
that the space of embedded compact minimal
surfaces of any fixed genus in $S^3$ is
compact in the $C^k$ topology, $k\geq 1$.

E. Calabi (1968) described in principle
minimal immersions of $S^2$ into $S^n$. He
also proved in 1967 that, for a minimal
full immersion of $S^2$ into a $(m-1)$-sphere
$S^{m-1}(r)$ of radius $r$, $m$ is an odd
integer and the area of the immersed $S^2$ is
an integral multiple of $2\pi r^2$, at least
${1\over 2}\pi r^2(m^2-1)$. S. S.
Chern (1970) provided a general construction
method of minimal spheres in the unit
4-sphere using their directrix curves. 

A holomorphic curve $\Xi:S^2\to CP^{2m}$ is called totally isotropic
if  any of its local
representations $\xi$ in homogeneous coordinates
 satisfies $(\xi ,\xi )=(\xi',\xi')=\cdots
=(\xi^{(m-1)},\xi^{(m-1)})=0$, where  the upper
indices stand for derivatives and
$(\; ,\; )$ denotes the canonical symmetric
product in $C^{2m+1}$. J. L. Barbosa (1975)
established a one-to-one correspondence
between the set of all full generalized
minimal immersions
 $f: S^2\rightarrow S^{2m}(1)$ and
the set of all linearly full totally isotropic
curves $\Xi: S^2\to CP^{2m}$, where $CP^{2m}$
denotes the $2m$-dimensional complex
 projective space of constant
holomorphic curvature 4, with such immersions corresponding to 
their directrices. It
is then natural to define the degree of the minimal
 immersion $f$ as
the degree of its directrix curve. By considering a very particular
local expression for the directrix curve, Barbosa obtained a set of
minimal immersions such that, for any
multiple of $4\pi $ greater than or equal to
$2\pi m(m+1)$, there is one having that value
as its area. He also showed that the group
${SO}(2m+1,\hbox{\bf C})$ of all complex
matrices $A$ satisfying det$(A)=1$ and
$AA^T=I$ acts on the space of totally isotropic
curves.  Identifying the minimal immersions that are
isometric, he found that each orbit of this
action is diffeomorphic
 to ${SO}(2m+1,\hbox{\bf C})/{
SO}(2m+1,\hbox{\bf R})$. Barbosa also showed
that the space of
 minimal immersion of degree $2m$
consists of exactly one such orbit. X. X. Li (1995)
extended Barbosa's result to the case of
minimal immersions of degree $2m+2$. By solving the
 totally 
isotropic conditions, he concludes that the
set of full minimal immersions $f:S^2\to
S^{2m}(1)$ of degree $2m+2$ is, modulo
isometries, diffeomorphic to a disjoint union of
$m-1$ such orbits. 

N. Ejiri (1986a) investigated equivariant
minimal immersions of $S^2$ into $S^{2m}$
and proved the following:

(a) There are no full minimal immersions
from the real projective 2-plane $RP^2$ into
$S^{2(2n-1)}$;

(b) The minimal cone of a full minimal
immersion of $S^2$ into $S^{2m}$ is stable; 

(c) If $f:S^2\to S^{2m}$ is a full minimal
immersion whose minimal cone has a parallel
calibration, then $m=3$ and $f$ is
holomorphic in the near K\"ahler $S^6$; and

(d) Circle bundles of $S^2$ of positive even
Chern number $(\leq 4)$ can be minimally
immersed in the near K\"ahler $S^6$. 

There do exist many minimal isometric
immersions from $E^2$ into $S^m$. A
description of such minimal immersions has
been obtained by K. Kenmotsu in 1976.

A. Ros (1995) proved that  any 2-equator in a
3-sphere divides each embedded compact minimal
surface into two connected pieces, and closed
regions in the sphere with mean convex
boundary containing a null-homologous great
circle are the intersection of two closed
half-spheres. As an application of these
results Ros proved that the normal surface of
an embedded minimal torus is also embedded.
He also showed that the Clifford torus is
the only embedded minimal torus in $S^3$
that is symmetric with respect to four
pairwise orthogonal hyperplanes in $E^4$. 
In his 1997 doctoral thesis at Universidad
Federal do Cear\'a, F. A. Amaral claimed
that the Clifford torus is the only embedded
minimal torus in $S^3$.

Y. Kitagawa (1995) proved that if $f: M\to
S^3$ is an isometric embedding of a flat
torus $M$ into $S^3$, then the image $f(M)$
is invariant under the antipodal map of $S^3$.
He also showed that there exist an immersed
flat torus $M$ in $S^3$ whose image is not
invariant under the antipodal map of $S^3$. 
H. Hashimoto and K. Sekigawa (1995) showed
that  a complete minimal surface in $S^4$
with nonnegative Gaussian curvature is either
superminimal or congruent to the Clifford
torus. 

H. Gauchman (1986) showed that if the second
fundamental form $h$ of a compact minimal
submanifold in a unit sphere satisfies
$|h(u,u)|^2\leq 1/3$ for any unit
tangent vector at any point, then it is
totally geodesic. P. F. Leung (1993) proved
that if an $n$-dimensional ($n\geq 3$) compact
oriented submanifold (not necessarily
minimal) in a unit sphere satisfies the same
condition as  Gauchman's, then $M$ is
homeomorphic to a sphere when $n> 3$ and $M$
is homotopic to a sphere when $n=3$.  In
1997 C. Y. Xia extended Leung's result
to the following: 

Let $M$ be an $n$-dimensional $(n\geq 4)$
compact simply-connected submanifold
isometrically immersed in a $\delta$-pinched
$(\delta>\frac14)$ Riemannian manifold. If
the second fundamental form of $M$ satisfies
$|h(u,u)|^2<\frac49 (\delta-\frac 14)$ for
any unit tangent vector at any point, then
$M$ is homeomorphic to an $n$-sphere.

\vskip.1in
\noindent{\bf 5.4.5. Simons' theorem}

J. Simons (1968) proved  that if the squared
length of the second fundamental form, denoted
by $S$, of a compact $n$-dimensional minimal
submanifold $M$ of the unit
$(n+p)$-sphere $S^{n+p}$ satisfies $$S\leq
{n\over {2-{1\over p}}},\leqno(5.12)$$ then either
$M$ is totally geodesic or $S={n\over
{2-{1\over p}}}$. 
If the second case occurs, then either 

(a) $M$ is a  generalized Clifford torus:
$$S^k\left(\sqrt{k\over n}\right)\times
S^{n-k}\left(\sqrt{{n-k}\over n}\right),$$
which is the standard product embedding of
the product of two spheres of radius
$\sqrt{k/n}$ and
$\sqrt{(n-k)/n}$, respectively, or 

(b) $M$ is a Veronese surface in $S^4$. 

 A. M. Li and J. M. Li (1992) showed that if
$M$ is a compact $n$-dimensional  minimal
submanifold of $S^{n+p}$ with
$p\geq 2$ and $S\leq 2n/3$, then $M$ is
either a totally geodesic submanifold or a
Veronese surface in $S^4$. 
 
The Veronese surface in $S^4$ is
defined as follows: Let $(x,y,z)$ be the
natural coordinate system of $ E^3$ and
$(u_1,\ldots,u_5)$ that of $E^5$. The
mapping defined by
$$u_1={{yz}\over\sqrt{3}},\;\;
u_2={{xz}\over\sqrt{3}},\;\;
u_3={{xy}\over\sqrt{3}},\;\;
u_4={{x^2-y^2}\over{2\sqrt{3}}},\;\;
u_5={1\over 6}(x^2+y^2-2z^2)\leqno(5.13)$$
gives rise to an isometric immersion of
$S^2(\sqrt{3})$ into $S^4$. Two points
$(x,y,z)$ and $(-x,-y,-z)$ of $S^2(\sqrt{3})$
are mapped into the same point. Thus, the
mapping defines an embedding of the real
projective plane $RP^2$ into $S^4$. This
embedding of $RP^2(\sqrt{3})$ into $S^4$ is
called the Veronese surface.

\vskip.1in
\noindent{\bf 5.4.6. Chern-do Carmo-Kobayashi's
theorem and related results}

S. S. Chern, M. do Carmo and S. Kobayashi
(1970) proved that the open pieces of the
generalized Clifford torus and the Veronese
surface are the only minimal submanifolds of
$S^{n+p}$ with  $S={n\over {2-{1\over
p}}}$.

A large number of examples of minimal
hypersurfaces in $S^{n+1}$ were constructed by
W. Y. Hsiang (1967), using Lie group methods.
For example, he showed that, for each $n\geq
4$, there exist infinitely many, mutually
incongruent minimal embeddings of $S^1\times
S^{n-2}$ (respectively, $S^2\times S^{n-3}$)
into $S^n(1)$. 

 Hsiang also considered the problem of finding
algebraic minimal cones, obtained by setting a
homogeneous polynomial equal to zero. For
quadratic polynomials, they are
$$
q(x_1^2+\cdots+x_{p+1}^2)-p(x^2_{p+1}+
\cdots+x^2_{p+q+2})=0,$$ $$
p\geq 1,\;\; q\geq 1,\;\; p+q+2=n,
\leqno (5.14)$$ whose intersection with $S^{n+1}$
are the generalized Clifford tori. Hsiang showed
that these are in fact the only algebraic
minimal cones of degree 2.

S. S. Chern conjectured that for a compact
 minimal hypersurface with constant scalar
curvature in $S^{n+1}$ the values $S$ are
discrete. C. K. Peng and C. L.
Terng (1983) proved that if $M$ is
a compact minimal hypersurface of $S^{n+1}$
with constant scalar curvature, then there
exists a constant $\epsilon(n)>1/(12n)$ such
that if $n\leq S\leq n+\epsilon(n)$, then
$S=n$, so that $M$ is a generalized Clifford
torus. Furthermore, they showed that if $n =
3$ and $S > 3$, then $S \geq 6$; this bound is
sharp, since the principal curvatures of the
Cartan minimal isoparametric hypersurface
$SO(3)/(\hbox{\bf Z}_2\times \hbox{\bf Z}_2)$
in $S^4$ are given by $\sqrt{3},0,-\sqrt{3}$.
Peng and Terng's result still holds if the
3-dimensional minimal submanifold is assumed
to be complete [Cheng 1990]. Peng and Terng
conjectured that the third value of $S$
should be $2n$, since there  exist Cartan's
isoparametric minimal hypersurfaces in
$S^{n+1}$ satisfying $S=2n$.

H. C. Yang and Q. M. Cheng
(1997) proved that, for a compact minimal
hypersurface $M$ with constant scalar
curvature in $S^{n+1}$, if $S>n>3$, then
$S>n+\frac13n$. In particular, if the shape
operator $A_\xi$ of $M$ in $S^{n+1}$
with respect to a unit normal vector $\xi$
satisfying trace$\,(A_\xi^3)=
\hbox{constant}$, then $S\geq n+\frac23n$. 

Q. M. Wang (1988) constructed examples of
compact noncongruent minimal hypersurfaces
in odd-dimensional spheres which have the
same constant scalar curvature. Thus, the
compact minimal hypersurfaces 
with given constant scalar curvature in a
sphere are not necessary unique.

It is still an open problem to determine
whether $S\geq 2n$ for a compact minimal
hypersurface $M$ with constant scalar
curvature in $S^{n+1}$ with $S>n>3$. 

For an $n$-dimensional compact minimal
manifold  $M$  in $S^{n+p}$ with $p\geq 2$,
C. Xia (1991) proved the following:

(1) If $n$ is even and $S\leq
n(3n-2)/(5n-4)$, then $M$ is either totally
geodesic or a Veronese surface in $S^4$;

(2) If $n$ is odd and
$S\leq n(3n-5)/(5n-9)$, then 

(2-i) when $n>5$, $M$ is totally geodesic in
$S\sp {n+p}$; 

(2-ii) when $n=5$, $M$ is either totally
geodesic or homeomorphic to $S\sp 5$ and
$S=25/8$ on $M$; and

(2-iii) when $n=3$, $S$ is identically equal
to 0 or 2; in the latter case $M$ is
diffeomorphic to $S^3$ or $RP^3$. 

T. Itoh (1978) proved that if $f:M\to
S^{n+p}$ is a minimal full
isometric immersion of a compact orientable
Riemannian $n$-manifold into  $S^{n+p}$ and
the sectional curvature $K$ of $M$
satisfies $K\geq n/2(n+1)$,
then either $M$ is totally geodesic or $M$
is of constant sectional curvature
$n/2(n+1)$ and the immersion is given by the
second standard immersion of an  $n$-sphere 
of sectional curvature $n/2(n+1)$. 

N. Ejiri (1979a) showed that if the Ricci
tensor of  an
$n$-dimensional $(n\geq 4)$ compact minimal
submanifold of
$S^{n+p}$ satisfies 
$Ric\geq(n-2)g$, then  $M$ is 
totally geodesic, or $n=2m$ and
$M$ is $S^m(\sqrt{1/2})\times
S^m(\sqrt{1/2})\subset S^{n+1}\subset
S^{n+p}$ embedded in a standard way, or $M$
is a 2-dimensional complex projective space
$CP^2$ of constant holomorphic sectional
curvature $\frac 43$ which is isometrically
immersed in a totally geodesic $S^7$ via
Hermitian harmonic functions of degree one.

G. Chen and X. Zou (1995) showed that  if the
sectional curvature is 
$\geq {1\over 2}-{1\over {3p}}$, then either
$M$ is totally geodesic or the Veronese
surface in $S^4$.

\vskip.1in
\noindent{\bf 5.4.7. Otsuki's theorem
and Otsuki's equation}

In 1970, T. Otsuki proved the following.

 Let $M$ be a complete minimal hypersurface of
$S^{n+1}$ with two principal curvatures. If
their multiplicities  $k$ and $n-k$ are
$\geq 2$, then $M$ is the generalized
Clifford torus $S^k(\sqrt{k/n})\times
S^{n-k}(\sqrt{(n-k)/n})$. If one of the
 multiplicities is one, then $M$ is a
hypersurface of $S^{n+1}$ in $
E^{n+2}= E^n\times  E^2$ whose
orthogonal projection into $ E^2$ is a
curve of which the support function $x(t)$
is a solution of the following nonlinear
differential equation:
$$nx(1-x^2)x''(t)+x'(t)^2+(1-x^2)(nx^2-1)=0.
\leqno(5.15)$$

Furthermore, there are countably many compact
minimal hypersurfaces immersed but not
embedded in $S^{n+1}$. Only
$S^{n-1}(\sqrt{(n-1)/n})\times
S^1(\sqrt{1/n})$ is minimally embedded in
$S^{n+1}$, which corresponds to the trivial
solution $x(t)=1/\sqrt{n}$ of  Otsuki's
equation (5.15).

Applying Otsuki's result Q. M. Cheng (1996)
proved that if $M$ is a compact minimal
hypersurface of $S^{n+1}(1)$ with two
distinct principal curvatures such that
$$n\leq S\leq n+
{{2n^2(n+4)}\over{3(n+2)^2}},\leqno(5.16)$$ then
$S=n$ and hence
$M$ is a generalized Clifford torus.

Otsuki's result was extended by L. P. Jorge and
F. Mercuri (1984) to submanifolds of higher
codimension: If $f\:M^
n\to S^{n+p}\, (n\geq3, p >1)$
is a full minimal immersion such that the
shape operator $A_\xi$ in any normal
direction
$\xi$ has at most two distinct eigenvalues, then
$M^n$ is an open subset of a projective space
over the complex, quaternion or Cayley numbers,
and $f$ is a standard embedding with parallel
second fundamental form. 

\subsection{Minimal submanifolds in hyperbolic space}

N. Ejiri (1979b) proved that every minimal
submanifold in a hyperbolic space is
 irreducible as a Riemannian manifold.  Chen
(1972) showed that there exists no minimal
surface of constant Gaussian curvature in
$H^3$ except the totally geodesic one.

 On the other hand,  M. do
Carmo and M. Dajczer (1983) constructed many
minimal rotation hypersurfaces in hyperbolic
space, in particular, in $H^3$. They also proved
that 
 complete minimal rotation surfaces of
$H^3$ are embedded.

X. Li-Jost (1994) studied Plateau type
problem in hyperbolic space and proved that  if 
$\Gamma$ is a closed Jordan curve of class $C\sp
{3,\alpha}$ in $H^3$ with total curvature $\leq
4\pi$, then there exists precisely one minimal
surface of disk type,  free of branch
points, spanning $\Gamma$.

Let $\gamma$ be a geodesic in $H^3$,
$\{\psi_t\}$ the translation along
$\gamma$, and  $\{\varphi_t\}$ the
one-parameter subgroup of isometries of $H\sp
3$ whose orbits are circles centered on
$\gamma$. Given any $\alpha\in
\bold R$,  $\lambda=\{\lambda_t\}=\{\psi_t\circ
\varphi_{\alpha t}\}$ is a one-parameter
subgroup of isometries of $H\sp 3$, which is
called a helicoidal group of isometries with
angular pitch $\alpha $. Any surface in
$H\sp 3$ which is
$\lambda$-invariant is called a helicoidal surface. 

J. B. Ripoll (1989) proved the following:
 Let $\alpha \in \hbox{\bf R}$, $\vert
\alpha
\vert <1$. Then there exists a
one-parameter family
$\Sigma$ of complete simply-connected minimal
helicoidal surfaces in $H^3$ with angular
pitch $\alpha $ which foliates $H^3$.
Furthermore, any complete helicoidal minimal
surface in $H^3$ with angular pitch
$\vert \alpha \vert <1$ is congruent to an element of $\Sigma$. 

G. de Oliveira Filho (1993)
considered complete minimal immersions in
hyperbolic space and proved the following.

(1) If $M^n\to H^m$ is a  complete minimal
immersion and $\int_M S^{n/2}dV <\infty$, then $M$
is properly immersed and is diffeomorphic to the
interior of a compact manifold $\overline M$ with
boundary. Furthermore, the immersion $M^n\to H^m$
extends to a continuous map $\overline M\to
\overline H^n$, where $\overline H^n$
 is the compactification of $H^n$. 
 
(2) If $M^2\to H^m$ is a  complete minimal
immersion  with $\int\sb
M SdV<\infty$, then $M$ is conformally
equivalent to a compact surface $\overline M$
with a finite number of disks removed and the
index of the Jacobi operator is finite.
 Furthermore, the asymptotic boundary
$\partial_\infty M$ is a Lipschitz curve.

K. Polthier (1991)
 constructed complete embedded
minimal surfaces in  $H^3$ having the symmetry of
a regular tessellation by Coxeter orthoschemes
and proved that there exist complete
minimal surfaces in
$H\sp 3$ with the symmetry of tessellations given
by (a) all compact and noncompact Platonic
polyhedra; (b) all Coxeter orthoschemes $(p,q,r)$
with $q\in\{3,4, \cdots, 1000\}$ and small $p$ and
$r$; (c) all ``self-dual'' Coxeter
orthoschemes
$(p,q,r)$ with $p=r$. 

Recently, M. Kokubu (1997) established the
Weierstrass type representation for minimal
surfaces in hyperbolic space.

For stable minimal submanifolds in hyperbolic
space,  do Carmo and  Dajczer (1983) proved that
there exists an infinite family of simply
connected stable complete minimal surfaces in
hyperbolic space
$H^3$ that are not totally geodesic.
Furthermore, let $f:M\to H^3$ be an isometric
immersion of
$M$ into $H^3$. If  $D$ is a simply-connected
domain in  $M$ with compact closure $\bar D$
and piecewise smooth boundary $\partial D$,
Barbosa and do Carmo (1981) proved that if
$$\int_{\bar D}\left(|K|+{2\over
3}\right)dV<2\pi,$$ then $D$ is stable.
Also, Ripoll's result implies that  any
complete helicoidal minimal surface in
$H^3$ with angular pitch $\vert \alpha \vert
<1$ is globally stable. 

 Let $M^{p-1}\to S^{n-1}(\infty)$, $p=n-1\leq 6$,
be an immersed compact submanifold in the
$(n-1)$-sphere at infinity of
$H^n$. M. T. Anderson (1982)  proved that 
there exists a complete embedded absolutely
area-minimizing submanifold asymptotic to
$M^{p-1}$ at infinity. In particular, there are
lots of  embedded complete minimal
submanifolds in case $p=n - 1\leq 6$.

\subsection{Gauss map of minimal surfaces}

The Gauss map $G:M\to S^2$ of a surface
$f:M\to  E^3$ is a map from the surface $M$
to the unit sphere $S^2$ given by
$G(x)=\xi(x)$, where $\xi(x)$ is the unit
normal of $M$ at $x$. Since $\xi(x)$ is a
unit vector in $E^3$, one may represent it as
a point in $S^2$.

O. Bonnet (1860) proved that the Gauss map of
a minimal surface in $E^3$ is conformal.
Conversely, E. B. Christoffel proved in 1867
that if the Gauss map of a surface in $E^3$
is conformal, then it is either a minimal
surface or a round sphere.

For a surface $f:M\to 
E^m,\, m\geq 3$, the Gauss map $G$ is defined
to be the map which assigns to each point
$x\in M$ the oriented tangent space
$f_*(T_xM)\subset E^m$. The Gauss map
$G$ can be considered as a map from $M$ into
the Grassmann manifold
$G^R(2,m-2)=SO(m)/SO(2)\times SO(m-2)$ of
oriented 2-planes in $ E^m$, which in
turn can be identified with the complex
quadric $Q_{m-2}(C)$:
$$Q_{m-2}(C)=\{(z_1,z_2,\ldots,z_m)\in
CP^{m-1}:\sum z_j^2=0\}\leqno (5.17)$$
in the complex projective space $CP^{m-1}$
in a natural way. The Gauss map of an
$n$-dimensional submanifold in $ E^m$ is
a map from $M$ into $G^R(n,m-n)$ defined in a
similar way.

The complex projective space
admits a unique K\"ahler metric with
constant holomorphic sectional curvature 2.
The induced metric on $Q_{m-2}$ defines a
metric $\hat g$ on the Grassmannian
$G^R(2,m-2)$ under the identification,
satisfying
$$G^*(\hat g)=-Kg\leqno (5.18)$$
for any minimal surface $M$ in $ E^m$,
where $g$ is the metric on $M$ and $K$ the
Gaussian curvature of $M$. Thus, the Gauss
map $G$ is conformal for a minimal surface
$f:M\to E^m$.
S. S. Chern (1965) showed that an immersion
$f:M\to E^m$ is minimal if and only if the
Gauss map $G$ of
$f$ is antiholomorphic.

Since the area $\hat A(G(M))$ of the Gauss
image $G(M)$ is related with the total
curvature of $M$ by
$$\hat A(G(M))= -\int_M KdA,\leqno (5.19)$$
for the minimal surface, one is able to
translate statements about the total
curvature of a minimal surfaces in $
E^m$ into corresponding statements about
the area of holomorphic curves in $CP^{m-1}$.

\vskip.1in
\noindent{\bf 5.6.1. Chern-Osserman's theorem}

 S. S. Chern and R. Osserman (1967) proved  the following fundamental results: 

If $f:M\to E^m$ is a complete orientable
minimal surface with finite total curvature
$\int_M KdA=-\pi C <\infty,$ then

(1) $M$ is conformally a compact Riemann
surface $\bar M$ with finite number, say $r$,
of points deleted;

(2) $C$ is an even integer and satisfies
$$C\geq 2(r-\chi)=4g+4r-4,$$ where $\chi$
is the Euler characteristic  and $g$ is the
genus of $M$ (= the genus of $\bar M$);

(3) if $f(M)$ does not lie in any proper
affine subspace of $ E^m$, then $$C\geq
4g+r+m-3\geq 4g+m-2\geq m-2;$$

(4) if $f(M)$ is simply-connected and
nondegenerate, that is, $G(M)$ does not lie
in a hyperplane of $CP^{m-1}$, then $C\geq 2n-2$
and this inequality is sharp;

(5) when $m=3$, $C$ is a multiple of 4, with
the minimum value 4 attained only by
Enneper's surface and the catenoid;

(6) the Gauss map $G$ of $f$ extends to a
map of $\bar M$ whose Gauss image $G(\bar M)$
is an algebraic curve in $CP^{m-1}$ lying
in $Q_{m-2}$; the total curvature of $f(M)$
is equal to the area of
$G(\bar M)$ in absolute value, counting
multiplicity;

(7) $G(M)$ intersects a fixed number of
times, say $n$ (counting multiplicity),
every hyperplane in $CP^{m-1}$ except
for those hyperplanes containing any of the
finite number of points of $G(\bar M-M)$;
the total curvature of $f(M)$ equals
$-2n\pi$.

(8) Enneper's surface and the catenoid are
the only two complete minimal surfaces in
$ E^3$ whose Gauss map is one-to-one.

For a complete oriented (not necessary
minimal) surface
$M$ in $E^m$, B. White (1987) proved that if
$\int_M SdA$ is finite, $S$ the squared
length of the second fundamental form, then
the total curvature,
$\int_M KdA$, is an integral multiple of
$2\pi$, or of $4\pi$ in case $m=3$.

\vskip.1in
\noindent{\bf 5.6.2. Value distribution of  Gauss map of complete minimal surfaces}

The Gauss map of Scherk's surface in
Euclidean 3-space omits exactly 4 points of
$S^2$. F. Xavier (1981) proved that the
Gauss map of any complete nonflat minimal
surface in $ E^3$ can omit at most 6
points of $S^2$. F. L\'opez and A. Ros (1987,
unpublished) gave a 1-point improvement by
showing  that the Gauss map of any complete
nonflat minimal surface in $E^3$ can omit at
most 5 points of $S^2$. Finally, H. Fujimoto
(1988) proved that the Gauss map of any
complete nonflat minimal surface in $ E^3$
can omit at most 4 points of $S^2$. Clearly,
Fujimoto's estimate is sharp.

For an orientable complete minimal
surface $M$ in $E^3$ with finite total
curvature, a theorem of A. Huber
(1957) implies that $M$ is conformally
equivalent to a compact Riemann surface
punctured at a finite number of points; thus
there is a closed Riemann surface $M_k$ of
genus $k$ and a finite number of points
$Q_1,\ldots,Q_r$ on
$M_k$ such that $M$ is conformally 
$M=M_k-\{Q_1,\ldots,Q_r\}$ 
[Osserman 1969b]. 
R. Osserman (1964) extended this result to
complete surfaces of finite total curvature in
$E^3$ with nonpositive Gaussian curvature.

For a complete minimal surface $f:M\to E^3$
of finite total curvature, the Gauss map $G$
of $f$ can be extended to a meromorphic
function $G:M_k\to S^2$.

The total curvature of the catenoid is
$-4\pi$ and its Gauss map misses 2 values. R.
Osserman (1961) proved that if the Gauss map
of a complete minimal surface of finite total
curvature in $E^3$ omits more than 3 values,
then it is a plane. One important consequence
of this is a sharpening of Fujimoto's result:
If the Gauss map of a complete nonplanar
minimal surface in $E^3$ omits 4 points on
$S^2$, then every other point of $S^2$ must be
covered infinitely often; and hence the total
curvature of the minimal surface must be
infinite. 

There is no known example of a complete
minimal surface of finite total curvature
whose Gauss map misses 3 values. Osserman
(1964) proved that if the Gauss map of a
complete minimal surface of finite total
curvature in $E^3$ misses 3 values, then the
genus of the minimal surface is at least one
and the total curvature is less than or equal to
$-12\pi$;  A.
Weitsman and F. Xavier proved in 1987 that
the total curvature is less than or equal to
$-16\pi$; and Y. Fang proved in 1993 that the
total curvature must be at most $-20\pi$, and
the degree of the Gauss map is at least five.

In embedded case, the Gauss map a complete
minimal surface in $E^3$ with finite total
curvature cannot omit more than 2 values,
since the limit normal direction at each end
belongs to a certain pair of antipodal points
[Jorge-Meeks 1983]. In particular, if the
minimal surface is embedded or the minimal
surface has parallel embedded ends, then it
has at least two catenoid type ends [Fang
1993].

In 1990 X. Mo and R. Osserman showed that if
the Gauss map of a complete minimal surface
in $E^3$  takes on 5 distinct values only a
finite number of times, then the minimal
surface has finite total curvature. Mo and
Osserman's result is sharp, since there is an
embedded complete minimal surface, due to
Scherk, in $E^3$ whose Gauss map misses four
points and takes any other points infinitely
many times. Mo and Osserman (1990) also
proved that the Gauss map of a nonplanar
complete minimal surface in $E^3$ of infinite
total curvature takes on every value
infinitely often, with the possible exception
of four points. 

Since the complex quadric surface
$Q_2$ is holomorphically isometric to the
product of two spheres of radii
$1/\sqrt{2}$, the Gauss map of
a surface $M$ in $ E^4$ is thus described
by a pair of maps $G_j:M\to S_j$, $j=1,2$.

M. Pinl (1953) showed that for a given
minimal surface in $ E^4$, the maps
$G_1$ and $G_2$ defined above are both
conformal.

W. Blaschke (1949) proved the following:
Let $M$ be a compact surface immersed in
$E^4$ and let $\chi$ be its Euler
characteristic. Denote by $A_j$
the algebraic area of the image of $M$
under the map $G_j,\, j=1,2$. Then
$$A_1+A_2=4\pi \chi.\leqno(5.20)$$

S. S. Chern (1965) proved that if $M$ is a
complete minimal surface in $ E^4$ and
if the image of $M$ under each of the maps
$G_1,G_2$ omits a neighborhood of some
point, then it is a plane.

 X. Mo and R. Osserman (1990) proved that
 if each of the factors
$G_j$ of the Gauss map of  a complete nonflat
minimal surface in $E^4$ omits 4
distinct points, then each of the $G_j$ must
cover every other point infinitely often. If
one of the $G_j$ is constant, then the other
must cover every point infinitely often with at
most 3 exceptions. 

For a complete minimal surface $f:M\to 
E^m$ with $m\geq 3$,
 H. Fujimoto (1990) proved that $G$ can omit 
at most $m(m+1)/2$
hyperplanes in general position if the Gauss map $G$ of $f$ is
nondegenerate, that is, $G(M)$ is not contained
in any hyperplane in $CP^{m-1}$. For arbitrary
odd number $m$, the number $m(m+1)/2$ is sharp
(cf. [Fujimoto 1993, \S5.5]).

M. Ru (1991)  improved Fujimoto's result to the
following: If the Gauss map $G$ of $f$ omits
more than $m(m+1)/2$ hyperplanes in
$CP^{m-1}$, located in general position, then
the minimal surface must be a plane. 

Recently, R. Osserman and M. Ru (1997)
extended the above result to the following:
Let $f:M\to  E^m$ be a minimal surface
immersed in $ E^m$. Suppose that its
Gauss map $G$ omits more than
$m(m+1)/2$ hyperplanes in $CP^{m-1}$,
located in general position. Then there exists
a constant $C$, depending on the set of
omitted hyperplanes, but not on the surface,
such that $|K(x)|^{1\over 2}d(x)\leq C$,
where $K(x)$ is the Gaussian curvature of $M$
at $x$ and $d(x)$ is the geodesic distance
from $x$ to the boundary of $M$.

Related to the above results are some results
for minimal surfaces defined on the complex
plane {\bf C} which are given by P. Hall.
Consider a minimal surface
$x:\hbox{\bf C}\to E^m$ and the reduced
representation $F=(f_1,\ldots,f_m)$ of its
Gauss map $G:\hbox{\bf C}\to CP^{m-1}$. A
direction $v=(v_1,\ldots,v_m)\in E^m$ is called
a normal to $M$ at $p\in M$ if it is orthogonal
to $T_pM$, that is, $\sum_{i=1}^mv_if_i(p)=0$.

P. Hall (1989,1991) proved the following:

(1) If the normals to a minimal surface
$x:\hbox{\bf C}\to E^m$ omits $m$ directions in
general position, then $x:\hbox{\bf C}\to E^m$
has a holomorphic factor, namely, there is an
orthogonal decomposition $E^m=E^2\oplus E^{m-2}$
such that the projection of $x$ into the first
factor is holomorphic or antiholomorphic with
respect to an orthogonal almost complex
structure on $E^2$;

(2) If  the normals to a minimal
surface $x:\hbox{\bf C}\to E^4$ omit four
directions in general position, then $x$ is
holomorphic in some orthogonal almost complex
structure. Moreover, it they omit five
directions in general position, then $x$ is a
plane.

(3) If the normals to a minimal surface
$x:\hbox{\bf C}\to E^m$ omits $n+k\,(k\geq 0)$
directions in general position, then the
dimension $d$ of the linear subspace of
$CP^{m-1}$ generated by the image of the Gauss
map and the dimension $a$ of the affine
subspace of $E^m$ generated by the image of $x$
satisfy $$ 1\leq d\leq {{m-3}\over{k+1}}, \quad
d+3\leq a\leq \min(n-kd,2d+2).$$

\subsection{Complete minimal submanifolds
in Euclidean space with finite total curvature}

Let $M_k$ be a compact surface of
genus  $k$,  $Q_1,\cdots,Q_r$ be $r$ points of
$M_k$, and  $$f: M= M_k-
\{Q_1,\cdots,Q_r\}\to  E^n$$
be a complete minimal immersion in Euclidean
$n$-space $ E^n$. If $D_j\subset M_k$ is
a topological disk centered at $Q_j$,
$j=1,\cdots,r$, $Q_i\not\in D_j$,
$i\not=j$, then $E_j=f(D_j\cap M)$ is an
end of the immersion $f$, and $f$ is called a
complete minimal immersion in $ E^n$, of
genus $k$ and with $r$ ends. 

A surface is said to have finite topology if it
is homeomorphic to a compact Riemann surface
from which a finite number of points have been
removed. A complete immersed minimal surface of
finite total curvature in $E^3$ has finite
topology. In fact, it is conformal to a
punctured compact Riemann surface. 

The helicoid is a complete embedded
simply-connected minimal surface of genus
zero and with one end. Since it is periodic
and nonflat, its total curvature is infinite.
This example shows that finite topology does
not imply finite total curvature.

T. Klotz and L. Sario (1965) proved that
there exist complete minimal surfaces in
$E^3$ of arbitrary genus with any finite
number of ends.

The Gauss map $G:M\to S^2$ of a complete
minimal surface $M$ of finite total curvature
in $E^3$ can be extended to $M_k$ such that
the extension $\tilde G:M_k\to S^2$ is a
holomorphic function. Moreover, the total
curvature of $M$ is $-4\pi \deg G$, where
$\deg G$ is the degree of $G$ [Osserman
1969b].

\vskip.1in
\noindent{\bf 5.7.1. Jorge-Meeks' formula and
its generalization}

Let $f:M\to E^m$ be a complete
minimal surface with finite total
curvature. Assume $M$ is conformally
$M=M_k-\{Q_1,\ldots,Q_r\}$,
$n\geq 1$, where $M_k$ is a closed Riemann
surface of genus $k$. Each $Q_j$
corresponds to an end $E_j$ of $M$. For
each end $E_j$ of $M$ the immersed circles
$\Gamma^j_t={1\over t}(E_j\cap S^{m-1}(t))$
converge smoothly to closed geodesics
$\gamma^j$ on $S^{m-1}(1)$ with multiplicity
$I_j$, where $S^{m-1}(t)$ is the sphere
centered at$(0,0,0)$ with radius $t$.

L. P. Jorge and W. H. Meeks (1983) proved
that if $M$ is a complete minimal surface
of finite total curvature in $E^m$ with $r$
ends, then
$$\int_MK\,dA=2\pi\Big({\chi
}(M)-\sum_{j=1}^r I_j
\Big)\leq 2\pi({\chi}(M)-r),\leqno(5.21)$$ where
$\chi(M)=2(1-k)-r$ is the Euler
characteristic of $M$. Furthermore, if $m$
equals 3, then 
$$\int_M K\,dA=  2\pi(\chi(M)-r)\leqno (5.22)$$if
and only if all the ends of $M$ are embedded.

For a branched complete minimal surface $M$
in $E^3$  Y. Fang (1996a) extended
Jorge-Meeks' formula to the following:
 $$\int_MK\,dA=2\pi\Big(\chi(M)-
\sum^r_{i=1}(J_ i-1)+\sum^N_
{i=1}K_j\Big),\leqno (5.23)$$
where $J_i\, (1\leq i\leq r)$ and $K_j\,(
1\leq j\leq N)$ denote respectively the
order of the ends $E_i$ and of the branch
points $q_j$ of $M$. 

\vskip.1in
\noindent{\bf 5.7.2. Topology of complete minimal surfaces}

\vskip.1in
\noindent{\bf 5.7.2.1. Ends of complete
minimal surfaces in $E^3$}

Let $M$ be a complete minimal surface of
finite total curvature in $E^3$. Suppose the
ends of $M$ are embedded. Then after a
suitable rotation of the coordinates, each end
of $M$ can be written as
$$z=a\log(x^2+y^2)+b+r^{-2}(cx+dy)
+O\left( r^{-2}
\right)\leqno 5.24$$ for suitable constants
$a,b,c$ and
$d$, where $r^2=x^2+y^2$. An end $E_j$ of
$M$ is called a flat (or planar) end, if
$a=0$ at $E_j$. Otherwise, $E_j$ is called a
catenoid end. An end $E_j$ of a complete
minimal surface $M$ in $E^3$ is said to be of
Enneper type if its multiplicity (or its
winding number) is 3.
The complete minimal surface
$M$ is said to be of flat type if all of its
ends are flat ends.

An end of a complete minimal surface $M$ in
$E^3$ is called annular, if it is
homeomorphic to a punctured disk.

Geometrically, all the topological ends of a
complete minimal surface of finite total
curvature in $E^3$ are conformally equivalent
to a punctured disk, and there is a
well-defined limit tangent plane at each end.
Outside of a sufficiently large compact set,
such an end is multisheeted graph over the
limit tangent plane, and if the end is
embedded, it is asymptotic to either a plane
or a half-catenoid [Schoen 1983].

Y. Fang (1996b) proved that if all of the ends
of a complete minimal surface of finite total
curvature in $E^3$ are embedded, then either
it is of flat type  (that is, no catenoid type
ends) or it has at least two catenoid type
ends.

An end of a complete embedded minimal surface
in $E^3$ is called a Nitsche end if it is
fibered by embedded Jordan curves in parallel
planes. Meeks and Rosenberg (1993b) proved that
if a complete embedded minimal surface of
finite topology has more than one end, then any
end of infinite total curvature is a Nitsche
end. 

L. P. Jorge and W. H. Meeks (1983) have
shown  the following topological uniqueness
result: Suppose $M_1$ and $M_2$ are complete
embedded minimal surfaces in $E^3$ with
finite total curvature. If $M_1$ and $M_2$
are diffeomorphic with two topological ends,
then there is an orientation preserving
diffeomorphism $\phi: E^3\to E^3$ with $\phi
(M_1)=M_2$.

For an embedded complete minimal surface in
$E^3$ with finite total curvature and with
$k$ ends, Jorge and Meeks (1983) proved the
following:

(1) If $M$ has an odd number of ends, then
$M$ disconnects $E^3$ into two regions
diffeomorphic to the interior of a solid
$(g+(k-1)/2)$-holed torus where $g$
is the genus of the associated compact
surface $\hat M$; and

(2) If $M$ has an even number of ends and
$\hat M$ has genus $g$, then $M$ disconnects
$E^3$ into $N_1$ diffeomorphic to the
interior of a $(g+k/2)$-holed solid torus
and into $N_2$ diffeomorphic to the interior
of a $(g+(k-2)/2)$-holed solid torus. 

\vskip.1in
\noindent{\bf 5.7.3. Properly embedded complete
minimal surfaces}

A mapping $f:M\to N$ between two topological
spaces is called proper if, for any compact
set $C\subset N$, $f^{-1}(C)$ is also
compact.  R. Osserman proved that every
complete minimal immersion $f:M\to E^3$ of
finite type is proper.

Until 1982 the only known examples
of properly embedded minimal surfaces in $E^3$
of finite total curvature were the plane, the
catenoid, and the helicoid. The total
curvature of the catenoid is $-4\pi$ whose
Gauss map has degree one. The catenoid has
genus zero and two embedded ends. The
helicoid is locally isometric to the catenoid
and it has infinite total curvature. The
Gauss map of the helicoid has an essential
singularity.

 Hoffman and Meeks (1989) proved that on a
properly embedded minimal surface in
$E^3$, at most two distinct annular ends can
have infinite total curvature. Thus, all
other ends have finite total curvature and
are therefore geometrically well behaved,
that is, asymptotic to the plane or a
half-catenoid. Fang and Meeks (1991) showed
that if a properly embedded minimal surface in
$E^3$ with two annular ends having infinite
total curvature, then these ends lie in
disjoint closed halfspaces and all other
annular ends are flat ends parallel to the
boundary of the halfspaces.

P. Collin (1997) proved that if a properly
embedded minimal surface in $E^3$ has at
least two ends, then it has finite topology
if and only if it has finite total curvature. 

\vskip.1in
\noindent{\bf 5.7.4. Half-space theorem}

D. A. Hoffman and W. H. Meeks (1990b) showed
that a nonplanar proper, possible branched,
complete minimal surface in $E^3$ is not
contained in a half-space. This result implies
that every nonplanar embedded complete
minimal surface of finite total curvature has
at least two annular ends. 

\vskip.1in
\noindent{\bf 5.7.5. Complete minimal
surfaces of genus zero}

F. J. L\'opez and A. Ros (1991) showed that
the plane and the catenoid are the only
embedded complete minimal surfaces of finite
total curvature and genus zero in $E^3$. L.
P. Jorge and W. H. Meeks (1983) proved that
the only complete finite total curvature
minimal embedding of $S^2-\{Q_1,\cdots,Q_r\}
\to  E^3$ for $1\leq r\leq 5$ are the plane
$(r=1)$ and the catenoid $(r=2)$. The cases
$r=3, 4$ or 5 do not occur.  

K. Yang (1994) showed that, for any finite
subset $\Sigma$ of $S^2$, one can conformally
immerse $S^2-\Sigma$ into $E^3$ as a complete
minimal surface of finite total curvature.

R. Miyaoka and K. Sato (1994) classified all
complete minimal surfaces in $E^3$ of
genus zero and two ends.

C. J. Costa (1993) constructed
several examples of complete immersed minimal
surfaces of finite topology and infinite
total curvature. The first of these
examples is a one-parameter family
perturbation of the catenoid. Each of the
perturbation surfaces has a catenoid-type
end and an end of infinite total curvature.
The second is a one-parameter family of
complete minimal surfaces of genus one and
three ends. Among the three ends, one is
flat, one is catenoid-type, and the third
is an end of infinite total curvature. This
second family is a perturbation of Costa's
embedded minimal surface of genus one and
three ends. The surfaces in these two
families are not embedded. 

W. Rossman (1995) classified all genus zero
catenoid-ended complete minimal surfaces with
at most $2n+1$ ends and high symmetry.

\vskip.1in
\noindent{\bf 5.7.6. Complete minimal
surfaces with one or two ends}

Planes are the only complete minimal surface of
finite total curvature in $E^3$ with one end.
R. Schoen (1983) proved that 
catenoids are the only complete embedded
minimal surfaces in $E^3$ of finite total
curvatures with two ends. C. Costa (1989)
classified complete minimal surfaces of finite
total curvature in $E^3$ with genus one and
three ends. Hoffman and Meeks (1989,1990a)
classified  complete minimal surfaces of finite
total curvature in $E^3$ with 
three ends and a symmetry group of order at
least $4(k+1)$.

K. Sato (1996) showed the existence of complete
immersed minimal surfaces of higher genus in
$E^3$ with finite total curvature and one
Enneper-type end. 

D. Hoffman, F. Wei and H. Karcher (1993)
constructed a properly embedded complete
minimal surface of infinite total curvature
with genus one and one end that is asymptotic
to the end of a helicoid; a genus one
helicoid. The surface contains two lines, one
vertical and corresponding to the axis of the
helicoid, the other horizontal crossing the
axial line. Rotation about these lines
generates the full isometry group, which is
isomorphic to $\hbox{\bf Z}_2\oplus \hbox{\bf
Z}_2$.

\vskip.1in
\noindent{\bf 5.7.7. Complete minimal
surfaces of higher genus}

R. Miyaoka and K. Sato (1994) constructed
examples  of complete minimal surfaces in
$E^3$ of genus $k$ with $r$ ends for
$k=0,\,r=3$ and for $k=1,\,r\geq 3$, via the
method of generating higher genus algebraic
curves through taking branched  coverings of
the Riemann sphere. Using Weierstrass $\frak
p$ functions on $M_1-\{\hbox{4 points}\,\}$,
they have constructed two series of examples. 
As a consequence Miyaoka and Sato have shown
that there exist complete minimal surfaces 
of finite total curvature in $E^3$, missed 2
values, for $M_k-\{r\hbox{ points}\,\}$ with
(1) $r\geq 2$ when $k=0$,  (2)  $r\geq 3$
when $k=1$, or  (3) $r\geq 4$ when $k\geq 2$.

For every positive integer $k$, D. Hoffman and
W.  Meeks constructed in 1990 an infinite
family of examples of properly  embedded
minimal surfaces of  genus $k$  with
three ends in  $E^3$. The total curvature is
$-4\pi(k+2)$. E. C. Thayer (1995) discovered a
family of complete minimal surfaces with
arbitrary even genus. Recently, For every
$k\geq 2$, Hoffman and Meeks discovered a one
parameter family, $M_{k,x},\,x\geq 1$, of
embedded minimal surfaces of genus $k-1$ and
finite total curvature. The surfaces
$M_{k,x},\, x>1$ have all three ends of
catenoid type and a symmetry group generated by
$k$ vertical planes of reflectional symmetry.

In 1995 W. Rossman constructed examples
of complete minimal surfaces in $E^3$ of
finite total curvature with catenoid-type
ends, of genus zero; and also of higher genus.
Rossman's examples include minimal surfaces
with symmetry group $D_n\times
\hbox{\bf Z}_2$ (dihedral symmetry) and
Platonic symmetry, where
$D_n$ is the dihedral group. 

For a closed Riemann surface $M_k$ of genus
$k$, a positive integer $r$ is called a
puncture number for $M_k$ if $M_k$ can be
conformally immersed in $E^3$ as a complete
finite total curvature minimal surface with
exactly $r$ punctures. The set of all puncture
numbers for $M_k$ is denoted by $P(M_k)$. K.
Yang (1994) proved that given any $M_k$ its
puncture set $P(M_k)$ always contains the set
$\{r\in \hbox{\bf Z}:r\geq 4k\}$.

 J. P\'erez and A.
Ros (1996) showed that the moduli space of
nondegenerate, properly embedded minimal
surfaces in $E^3$ with finite total
curvature is a real analytic
$(r+3)$-dimensional manifold if the fixed
number of ends is $r$.

\vskip.1in
\noindent{\bf 5.7.8.  Minimal annuli of finite total curvature}

The catenoid is topologically an annulus,
that is, it is homeomorphic to a
punctured disk. It follows from
Jorge-Meeks' formula that the catenoid is the
only embedded complete minimal annulus in
$E^3$ with finite total curvature. 

P. Collin (1997) proved that a
properly embedded complete minimal annulus in
$E^3$ with at least two ends has finite total
curvature.

\vskip.1in
\noindent{\bf 5.7.9. Riemann's minimal surfaces}

The catenoid is a rotational surface, hence 
is foliated by circles in parallel planes.
In 1867 B. Riemann found a one-parameter
family of complete embedded singly-periodic
minimal surfaces foliated by circles and
lines in parallel planes.  Each minimal annulus in this
one-parameter family is contained in a slab
and foliated by circles, and its boundary
is a pair of parallel straight lines.
Rotating repeatly about these boundary
straight lines gives a one-parameter family
of singly periodic minimal surfaces. These
surfaces known today as Riemann's minimal
surfaces. Riemann's minimal surfaces were
characterized by Riemann (1892) as the
only minimal surfaces fibered by circles in
parallel planes besides the catenoid. 

A. Enneper (1869) proved that a
minimal surface fibered by circular arcs
was an open part of a Riemann's minimal
surface or an open part of the catenoid.

 M. Shiffman (1956) proved that a minimal
annulus spanning two circles in parallel
planes was foliated by circles in parallel
planes and hence a part of Riemann's
examples or a part of the catenoid.
Hoffman, Karcher and Rosenberg (1991)
showed that an embedded minimal annulus
with boundary of two parallel lines on
parallel planes and lying between the
planes extended by Schwarz reflection to a
Riemann's minimal surface. 

\'E. Toubiana (1992) characterized Riemann's
minimal surfaces as the only properly embedded
minimal annuli between a pair of parallel
planes bounded by any pair of lines.  He
also generalized  Riemann's examples to
produce a countable family of immersed
minimal annuli between a pair of parallel
planes bounded by a pair of parallel lines.
These surfaces are  then extended, via the
reflection principle, to produce complete
immersed minimal surfaces.

 In 1993 P. Romon proved that a properly
embedded annulus with one flat end, lying
between two parallel planes and bounded by
two parallel lines in the planes, is a part
of a Riemann example.

 J. P\'erez (1995) proved
that a properly embedded minimal torus in
$E^3/T$ ($T$ is the group generated
by a nontrivial translation in $E^3$) with
two planar type ends is a Riemann's minimal
surface provided that it is symmetric with
respect to a plane. A. Douady and R. Douady
(1995) showed that Riemann examples are the
only singly-periodic with translational
symmetries minimal surfaces of genus one
with planar ends and a symmetry with
respect to a plane. In 1994 Y. Fang proved
that a properly embedded minimal annulus in
a slab with boundary consisting of two
circles or planes must be part of a
Riemann's minimal surface. Y. Fang and F.
Wei (1998) showed that  a properly 
embedded  minimal annulus with a planar end
and boundary consisting of circles or lines
in parallel planes is a part of a Riemann
example. F. J. L\'opez, M. Ritor\'e and F.
Wei (1997) characterized Riemann's minimal
surfaces as the only properly embedded
minimal tori with two planar ends in
$E^3/T$, where $T$ is the group generated
by a nontrivial translation in $E^3$. Using
numerical methods, F. Wei (1995)
constructed a properly embedded minimal
surface of genus two and two planar ends in
$E^3/T$ by adding handles to the Riemann
examples.

\vskip.1in
\noindent{\bf 5.7.10. Examples and classification
of complete minimal surfaces of finite
total curvature in $E^3$ }

Clearly, planes in $E^3$ are embedded
complete minimal surfaces with zero total
curvature. There are only two complete
minimal surfaces in
$E^3$ whose total curvature is $-4\pi$. These
are the catenoid and the Enneper surface; the
only embedded one is the catenoid. Also, it is
known that the only complete embedded minimal
surfaces with total curvature $\geq -8\pi$ in
$E^3$ are the plane and the catenoid with total
curvature 0 and $-4\pi$ respectively.

In 1981 W. H. Meeks  showed
that if $M$ is diffeomorphic to a real
projective plane minus two points, then it
does not admit a complete minimal immersion
into $E^3$ with total curvature $-6\pi$. 
A complete M\"obius strip in $E^3$ with total
curvature $-6\pi$ was constructed by Meeks
(1975).
M. Barbosa and A. G. Colares (1986) showed
that, up to rigid motions of $E^3$, there
exists a unique complete  minimal immersion
of the M\"obius strip into $E^3$ with total
curvature $-6\pi$.

Osserman, Jorge and Meeks proved that if $M$
is a complete minimal surface in $E^3$ with
total curvature greater than $-8\pi$, then,
up to a projective transformation of $E^3$,
$M$ is the plane, the catenoid, the Enneper
surface, or Meeks' minimal M\"obius strip.

By adjoining a handle
on Enneper's surface,  C. C. Chen and F.
Gackstatter (1982) constructed a complete
minimal surface of total curvature $-8\pi$ in
$E^3$; which was characterized by D. Bloss
(1989) and F. J. L\'opez (1992) as the only
complete minimal once punctured torus in $
E^3$ with total curvature $-8\pi$. 

It follows from a formula of L. Jorge and W. Meeks
(1983) that when the total curvature of $M$ is
$-8\pi$, the genus of the underlying Riemann
surface has to be either 0 or 1. Moreover, if the
genus is 1, the number of punctures (or ends)
has to be 1; and if the genus is 0, the
number of punctures can be 1, 2 or 3.
 The genus zero surfaces were classified rather
easily using the Weierstrass representation.

In 1992 F. J. L\'opez classified orientable
complete minimal surfaces in $ E^3$ with
total curvature $-8\pi$. In 1993 he gave an
example of a once-punctured minimal Klein 
bottle with  total curvature  $-8\pi$, and
proved in 1996 that this minimal Klein bottle 
is the only complete nonorientable minimal
surface in $ E^3$ with total
curvature $-8\pi$. 

M. E. G. G. Oliveira (1984) constructed an
example of a nonorientable complete minimal
surface of genus one with two ends and total
curvature $-10\pi$ in $E^3$. S. P. Zhang
(1989) observed that there is only
one minimal two-punctured projective plane in
$E^3$ of total curvature $-10\pi$ and such
that the branch number of the Gauss map at
the ends is greater than or equal to three. 

C. C. Chen and F. Gackstatter (1982)
constructed a complete minimal surface of
genus two with total curvature $-12\pi$ and
one end in $E^3$. A  complete
minimal surface of genus one with three ends
was discovered by C. J. Costa in 1984 which
satisfies the following two properties: (a)
the total curvature is
$-12\pi$, and (b) the ends are embedded.
Hoffman and Meeks (1985) showed that Costa's
minimal surface is properly embedded. They
also showed that  it contains two straight
lines meeting at right angles, it is composed
of eight congruent pieces in different
octants, each of which is a function graph,
and the entire surface is invariant under a
dihedral group of 3-space rigid motions.

C. J. Costa (1991) classified  orientable
complete minimal surfaces in $E^ 3$ with
total curvature $-12\pi$, assuming that they
are embedded. F. J. L\'opez,
F. Martin and D. Rodriguez (1997) proved that
the genus two Chen-Gackstatter example is the
unique complete orientable minimal surface of
genus two in $ E^3$ with total curvature
$-12\pi$ and eight symmetries. 

N. Do Espirito-Santo (1994)
showed the existence of a complete minimal
surface of genus $3$ with total curvature
$-16\pi$ and one Enneper-type end. F. F.
Abi-Khuzam (1995) constructed a one-parameter
family of complete minimal surfaces of genus
one with total curvature $-16\pi$ and having
four embedded planar ends.  

Costa also constructed an example of complete
minimal surface of genus one with two ends
and total curvature $-20\pi$ in $E^3$.

Complete M\"obius strips in $E^3$ with total
curvature $-2\pi n$, for any odd integer
$n\geq 5$, were constructed  by Oliveira
in 1984. In particular, this implies that
there exist complete M\"obius strips with
total curvature $-10\pi,-14\pi$ or $-18\pi$ in
$E^3$. In 1993 de Oliveira and Toubiana 
constructed, for any integer $n\geq 3$, an
example of complete minimal Klein bottles in
$E^3$ with total curvature $-2\pi(2n+3)$. 

In 1989 H. Karcher obtained a
generalization of Chen-Gackstatter surface by
increasing the genus and the order of the
symmetry group.  For each $k\geq 1$ he proved
that there exists a complete orientable
minimal surface of genus $k$ with one end,
total curvature $-4\pi (2k-2)$, and $4k+4$
symmetries. R. Kusner (1987) constructed a
family of immersed projective planes with
$k\, (k\geq 3)$ embedded flat ends and total
curvature $-4\pi(2k-1)$. In 1996 A. Ros
 proved that if $M$ is an  embedded
complete minimal surface of genus $k>0$
with finite total curvature, then the
symmetry group of
$M$ has at most $4(k+1)$ elements, and it
has
$4k+4$ elements if and only if $M$ is the
Hoffman-Meeks surface $M\sb k$ (1990). 

F. Martin (1995,1997) discovered a family of
complete non-orientable highly symmetrical
complete minimal surfaces with arbitrary
topology and one end and provided
characterizations of such minimal surfaces.
F. Martin and D. Rodriguez (1997) classified
complete minimal surfaces of total curvature
$-4\pi(3k-3)$ with $4k$ symmetries and one
end in $E^3$, for $k$ not a multiple of $3$. 

Jorge and Meeks constructed in 1983 complete
minimal surfaces of genus zero in $E^3$ with
total curvature $-4\pi r$ with $r$ embedded
ends. In 1993 \'E. Toubiana proved  that there
exist nonorientable minimal surfaces of genus
$k$ with two ends and total curvature
$-10(k+1)\pi$. 

The new examples of complete embedded minimal
surfaces of finite total curvature were
discovered by using the global version of the
Enneper-Riemann-Weierst\-rass representation,
which is essentially due to Osserman; The 
method involves knowledge of the compact Riemann
surface structure of the minimal surface as
well as its Gauss map and other
geometric-analytic data.

\vskip.1in
\noindent{\bf 5.7.11. Maximum principle at infinity}

The maximum principle at infinity for 
minimal surfaces in $E^{3}$ was
first studied by R.  Langevin and H.
Rosenberg (1988), who proved that the
distance between two disjoint embedded
complete minimal surfaces in $E^3$ with
finite total curvature and compact
boundaries must be greater than zero, that
is, the surfaces cannot touch each other at
infinity. W. Meeks and Rosenberg (1990) 
extended their result to the following: 

Let $M_1$ and $M_2$ be disjoint, properly
immersed minimal surfaces with nonempty
compact boundaries in a complete flat
3-manifold. Then $${\hbox{\rm
dist}}(M_1,M_2)=\min(\hbox{\rm dist}
(\partial M_1,M_2), \hbox{\rm
dist}(M_1,\partial M_2)).\leqno (5.25)$$

M. Soret (1995)  studied the maximum principle
at infinity for minimal surfaces with
noncompact boundaries and proved that if
$M_1$ and $M_2$ are disjoint properly
embedded minimal surfaces with bounded
curvature in a complete flat 3-manifold and
one of the surfaces is of parabolic type, then
$$\hbox{\rm dist}(M_1,M_2)=\min(\hbox{\rm dist} (\partial M_1,M_2),
\hbox{\rm dist} (M_1,\partial M_2)).\leqno (5.26)$$ Consequently, if $M_1$ and $M_2$ are
disjoint, properly embedded stable minimal
surfaces with noncompact boundaries in a
complete flat 3-manifold, then $$\hbox{\rm
dist}(M_1,M_2)=\min(\hbox{\rm dist} (\partial M_1,M_2),\hbox{\rm dist}(M_1,\partial M_2)).\leqno
(5.27)$$   In particular, if the boundary of one
surface, say $M_1$, is empty then
$$\hbox{\rm dist}(M_1,M_2)=\hbox{\rm dist}(M_1,\partial M_2).\leqno (5.28)$$

\vskip.1in
\noindent{\bf 5.7.12. Further results on complete
minimal surfaces in $E^3$ with finite total
curvature}

H. I. Choi, W. H. Meeks,  and B. White (1990)
proved that any intrinsic local symmetry of
the minimal surface in $E^3$ with finite
total curvature can be extended to a rigid
motion of $E^3$. Y. Xu (1995) observed that 
this property yields the identity of the
intrinsic and exterior symmetry groups for
the minimal surfaces with embedded catenoid
ends. As a consequence he proved
that, for any closed subgroup $G\subset SO(3)$
different from $SO(2)$, there exists a
genus zero complete minimal surface whose
symmetry group is $G$. The proof relies on
the fact that if $G$ is the symmetry group of
the minimal surface, then there exists an
appropriate M\"obius transformation
$\gamma$ which is conjugate to $G$ by the
Weierstrass representation. To
construct the corresponding examples, Xu
described all $\gamma$-invariant polynomials
which generate the Gauss map of symmetric
minimal surfaces.

\vskip.1in
\noindent{\bf 5.7.13. Complete minimal surfaces in $E^m,\, m\geq 4$ with finite total curvature}

A complete minimal surface in $E^m$ is said
to have quadratic area growth if 
$$\hbox{\rm Area}\,(M\cap B(R))\leq C_0
R^2\leqno (5.29)$$ for all $R>0$, where $C_0$ is a
constant and $B(R)$ is a ball of radius $R$ in
$E^m$ centered at 0. According to the
fundamental result of Chern-Osserman (1967)
if a complete minimal surface in $E^m$ has
finite total curvature, it is of quadratic
area growth and has finite topological type.
Conversely, Q. Chen (1997) proved that if a
complete minimal surface in $E^m$ has finite 
topological type and is of quadratic area
growth, then it has finite total curvature;
the result is false if one drops the
assumption of finite topological type, since
the surface $\sin z=\sinh x\sinh y$, a Scherk
surface in $E^3$, has infinite genus and
quadratic area growth.

For a complete oriented minimal surface $M$ of
finite type in $E^4$, S. Nayatani (1990a)
showed that if $M$ has finite total curvature
and degenerate Gauss map, then $M$ is of
finite total curvature or a holomorphic
curve with respect to some orthogonal almost
complex structure on $E^4$.

For complete minimal surfaces of $E^m$ with
$m\geq 4$, C. C. Chen (1979) proved that if a
complete minimal surface in
$ E^m$ has total curvature $-2\pi$, then it
lies in an affine 4-space $E^4\subset E^m$,
and with respect to a suitable complex
structure on $ E^4$, $M$ is a
holomorphic curve in $C^2$. C. C. Chen (1980)
also proved that if a complete minimal
surface in $ E^m$ has total curvature
$-4\pi$, it must be either simply-connected
or doubly-connected. In the former case, it
lies in some affine 6-space $
E^6\subset E^m$, and in the latter
case, in some $ E^5\subset  E^m$. 

D. Hoffman and R. Osserman (1980) gave
complete description of complete minimal
surfaces in Euclidean space with total
curvature $-4\pi$. In particular, they showed
that the dimensions 5 and 6 given by Chen are
sharp. It turns out that doubly-connected
surfaces are all a kind of ``skew catenoid''
generated by a one-parameter family of
ellipses.

 The Chern-Osserman
theorem implies that the total curvature of
a complete orientable minimal surface $M$ in
$ E^m$ is a negative integer multiple of
$4\pi$. Osserman showed that if the total
curvature is $-4\pi$, then $M$ must be either
Enneper's surface or the catenoid. 

\vskip.1in
\noindent{\bf 5.7.14. Complete minimal
submanifolds with finite total scalar
curvature}

Let $f:M\to E^m$ be a minimally immersed
submanifold of $E^m$. The total scalar
curvature of $f$ is defined to be $\int_M
S^{n/2} dV$, where $S$ is the squared length
of the second fundamental form. This integral
is called total scalar curvature because, for
minimal submanifolds in a Euclidean space,
the scalar curvature is equal to $-S$.

M. Anderson (1984, 1985) studied
$n$-dimensional complete minimal submanifolds
of dimension $n\geq 2$ with finite total
scalar curvature in a Euclidean space and
proved the following:

(1) A complete minimal
submanifold $M$ of $E^m$ with finite total
scalar curvature is conformally diffeomorphic
to a compact Riemannian manifold minus a
finite number of points, thus $M$ has only
finitely many ends. Moreover, each of finite
topological type; 

(2) Let $M$ be an $n$-dimensional complete
minimal submanifold
of $E^m$. If $n\geq 3$ and $M$ has finite
total scalar curvature and one end, then $M$
is an $n$-plane;

(3) If a complete minimal submanifold $M$ of
$E^m$ has finite total scalar curvature, then
each end of $M$ has a unique $n$-plane as its
tangent cone at infinity; 

(4) If a complete minimal submanifold $M$ of
$E^m$ has finite total scalar curvature, then
$M$ is properly immersed, that is, the
inverse image of any compact set is compact.

H. Moore (1996) also investigated complete
 minimal submanifolds of dimension $\geq
3$ with finite total scalar curvature. She
obtained the following results:

(5)  Let $M$ be a complete
minimal hypersurface of $E^{n+1}$ with $n\geq
3$. If $M$ has finite total scalar curvature,
then $M$ lies between two parallel $n$-planes
in $E^{n+1}$;

(6) Let $M$ be an $n$-dimensional complete
minimal submanifold of $E^m$ with $n\geq 3$.
If $M$ has finite total scalar curvature and
it has two ends, then either
$M$ is the union of two $n$-planes or $M$ is
connected and embedded; 

(7)  Let $M$ be an $n$-dimensional complete
minimal submanifold of $E^m$ with $n\geq 3$
and $n>m/2$. If $M$ has finite total scalar
curvature and it has two ends, then $M$
lies between two parallel
$n$-planes in some affine $(n+1)$-subspace
$E^{n+1}\subset E^m$; and 

(8)  Let $M$ be an $n$-dimensional complete
nonplanar minimal submanifold of $E^m$ with
$n\geq 3$ and $n>m/2$. If $M$ has
finite total scalar curvature and it has two
ends, then $M$ is a catenoid.

J. Tysk (1989) proved that a complete minimal
hypersurface $M$ in $E^{n+1}$ has
finite index if and only if $M$ has finite
total scalar curvature for
$n=3,4,5,6$, provided that the volume growth
of $M$ is bounded by a constant times $r^n$,
where $r$ is the Euclidean distance function.
Tysk also showed that the result is not valid
in $E^9$ and in higher-dimensional
Euclidean spaces.

\subsection{Complete minimal surfaces in $E^3$ lying between two parallel planes}

In 1980 Jorge and Xavier exhibited
a nontrivial example of a complete minimal
surface which lies between two parallel
planes in $E^3$. Rosenberg and Toubiana 
constructed in 1987 a complete minimal
surface  of the topological type of a
cylinder in $E^3$ which lies between two
parallel planes; this surface intersects
every parallel plane transversally.  Hoffman
and Meeks (1990b) proved that there does not
exist a properly immersed minimal surface in
$E^3$ that is contained between two parallel
planes; this follows from their result that a
nonplanar proper minimal surface $M$ in
$E^3$ is not contained in a half-space.
In 1992 F. F. de Brito  constructed a large
family of complete minimal surfaces which lie
between two parallel planes in $E^3$.  

For each positive integer $k$ and each
integer $N$, $1\leq N\leq 4$, C. J. Costa and
P. A. Q. Simoes constructed in 1996 an
example of complete minimal surface of
genus $k$ and $N$ ends in a slab of $E^3$.
More precisely, they showed that there is a
complete minimal immersion $f_{k,N}:
M_{k,N}\to E^3$ with infinite
total curvature such that: 

(a) $M_{k,1}$ and $M_{k,2}$ are respectively a
compact Riemann surface of genus $k$ minus
one disk and two disks, 

(b) $M_{k,j+2}$, $j=1,2$ are respectively
$M_{k,j}$ punctured at two points, and 

(c) $f_{k,N}(M_{k,N})$ lies between two
parallel planes of $E^3$ and $f_{k,3}$,
$f_{k,4}$ have two embedded planar ends.

\subsection{The geometry of Gauss image}

For a minimal surface $M$ in $ E^m$
let $\hat M=G(M)$ denote the Gauss image of
$M$ under its Gauss map. At all nonsingular
points of $\hat M$, we have a well-defined
Gaussian curvature $\hat K$. It follows from
the normalization of the metric on $CP^{m-1}$
that $\hat K\leq 2$.

D. Hoffman and R. Osserman (1980) proved
that  the Gaussian curvature $\hat K$ of the
Gauss image of a minimal surface $M$ in
$ E^m$ is equal to 2 everywhere if and only if
$M$ lies in some affine $ E^4\subset  E^m$
and is a holomorphic curve in $ C^2$ with
respect to a suitable orthogonal almost
complex structure on $ E^4$. Moreover, the
Gaussian curvature of the Gauss image of a
minimal surface $M$ in $ E^m$ is equal to 1
everywhere if and only if $M$ is locally
isometric to a minimal surface in $ E^3$.

 B. Y. Chen and S. Yamaguchi (1983) proved
that a submanifold $M$ of a Euclidean
$m$-space has totally geodesic Gauss image
if and only if the second fundamental form
$h$ of $M$ in $ E^m$ satisfies
$$(\bar\nabla_Xh)(Y,Z)=h(\nabla^G_XY,Z)-
h(\nabla_XY,Z)\leqno (5.30)$$
for vector fields $X,Y,Z$ tangent to $M$,
where $\nabla^G$ denotes the Levi-Civita
connection of the Gauss image with respect
to the metric induced from the Gauss map $G$.
By applying this necessary and sufficient
condition, Chen and Yamaguchi (1983) proved
that the Gauss image of a minimal surface in
$ E^m$ is totally geodesic in
$G^R(2,m-2)$ if and only if either $M$ lies
in an affine $ E^3\subset
 E^m$ or $M$ is a complex curve lying
fully in $ C^2$, where $ C^2$ is an
affine $ E^4\subset E^m$ endowed
with some orthogonal almost complex structure.
If the second case occurs, the Gaussian
curvature $\hat K$ of the Gauss image is 2.

 Chen and Yamaguchi (1983) also
completely classified surfaces in Euclidean
space with totally geodesic Gauss image:

Let $M$ be a surface in $ E^m$ whose
Gauss image is regular. If the Gauss image
$G(M)$ of
$M$ is totally geodesic in
$G^R(2,m-2)$, then $M$ is one of the
following surfaces:

(1) a surface  in an affine $
E^3\subset E^m$;

(2) a surface  in $ E^m$ with
parallel second fundamental form, that is,
 a parallel surface;

(3) a surface  in an affine 4-space
$ E^4\subset  E^m$ which is locally
the Riemannian product of two plane curves
of nonzero curvature;

(4) a complex curve lying fully in $
C^2$, where $C^2$ denotes an affine
$ E^4\subset  E^m$ endowed with
some orthogonal almost complex structure.

Conversely, surfaces of type (1), (2),
(3) and (4) have totally geodesic Gauss
image. 

Yu. A. Nikolaevskii (1993)  extended
Chen-Yamaguchi's result to the following
(see, also [Chen-Yamaguchi 1984]): 

Let $M$ be an $n$-dimensional submanifold in
$ E^m$ whose Gauss image is regular.
Then the Gauss image $G(M)$ of
$M$ is totally geodesic in
$G^R(n,m-n)$ if and only if $M$ is the
product of submanifolds, each of the factors
is either 

(a) a real hypersurface, or

(b) a submanifold with parallel second
fundamental form, or

(c) a complex hypersurface.

Chen and Yamaguchi (1984)
proved that a submanifold $M$ of $ E^m$ is
locally the product of real hypersurfaces if
and only if the Gauss image is totally
geodesic and the normal connection is flat.

\subsection{Stability and index of minimal submanifolds}.

\vskip.1in
\noindent{\bf 5.10.1. Stability and $\lambda_1$}

 If $f:M\to  E^m$ is a
minimal submanifold and $\xi$ is a normal
vector field on $M$,
 then $f+t\xi$ gives rise to a normal variation
$F:(-\epsilon,\epsilon)\times M\to 
E^m$ for some sufficiently small
$\epsilon>0$. 

The second variational formula
for the volume functional is given by
$${{d^2V_t}\over{dt^2}}|_{t=0}=\int_M
\left(||D\xi||^2-||A_\xi||^2\right)dV_0.
\leqno(5.31)$$
In particular, if $M$ is a minimal surface
in
$ E^3$, (5.25) reduces to
$${{d^2V_t}\over{dt^2}}|_{t=0}=\int_M
\left(|\nabla
\phi|^2+2K\phi^2\right)dA,\leqno(5.32)$$ where
$\xi=\phi e_3$ and $e_3$ is a unit normal
vector field of $M$.  A minimal surface $M$
in $ E^3$ is called stable if the second
variation is positive for all variations on
any bounded domain $D$ in $M$.

Therefore, a
minimal surface $M$ in $ E^3$ is stable
if and only if 
$$\int_D  \left(|\nabla
\phi|^2+2K\phi^2\right)dA >0\leqno (5.33)$$
 for any smooth function $\phi$ with
compact support on $M$. 

It is convenient to rewrite (5.33) using a
new metric $\hat g=-Kg$, where $g$ is the 
metric of $M$. Then we have
$$d\hat A=-KdA\leqno (5.34)$$ 
and 
$$|\nabla \phi|^2=-K|\hat\nabla
\phi|^2,\leqno (5.35)$$
where $\hat\nabla$ denotes the gradient in the
new metric. We can rewrite (5.33) as
$$\int_D|\hat\nabla\phi|^2d\hat A> 2
\int_D \phi^2d\hat A.\leqno (5.36)$$

The ratio $$Q(\phi)={{\int_D|\nabla\phi|^2d\hat
A}\over{\int_D \phi^2dA}}\leqno (5.37)$$ is
called the Rayleigh quotient of $D$, and
the quantity
$$\lambda_1(D)=\inf Q(\phi)\leqno(5.38)$$
represents the first eigenvalue of the problem
$$ \begin{cases}\Delta\phi+\lambda
\phi=0\quad\hbox{in}\;\; D,\\\phi
=0\quad\quad\quad\quad\hbox{on}\;\; \partial
D.\end{cases} \leqno(5.39)$$

The ``inf'' in (5.38) may be taken over all
piecewise smooth functions in $\bar D$ that
vanish on the boundary, where $\Delta$ in 
(5.39) is the Laplacian with respect to a
given metric on $D$. If $D$ has reasonably
smooth boundary, then (5.39) has a solution
$\phi_1$ corresponding to the eigenvalue
$\lambda_1$, and the ``inf'' in (5.38) is
actually attended when $\phi=\phi_1$.

From these it follows that the stability
condition  (5.27) is simply the condition:
$$\lambda_1(D)>2.\leqno (5.40)$$

Since for a minimal surface in $E^3$ the
metric
$\hat g$ is nothing but the pullback under the
Gauss map of the metric on the unit sphere
$S^2$, thus we have the following [Barbosa-do
Carmo 1976]: 

Let $D$ be a relatively compact
domain on a minimal surface $M$ in $ E^3$.
Suppose that the Gauss map $G$ of the minimal
surface maps
$D$ one-to-one onto a domain $\hat D$ on the unit
sphere. If
$\lambda_1(\hat D)<2$, then $D$ cannot be
area-minimizing with respect to its boundary.

Since $\lambda_1(D_1)=2$ for a
hemisphere $D_1$ on the unit sphere, this result
implies in particular a well-known result of
H. A. Schwarz: 

If the Gauss map $G$ of a minimal surface $M$
in $ E^3$ maps a relatively compact domain
$D$ of a minimal surface $M$ in
$ E^3$ one-to-one onto a domain containing a
hemisphere, then $D$ cannot be area-minimizing.

H. A. Schwarz also obtained in 1885 a
 sufficient condition for a domain $D$ in a
minimal surface to be stable; namely, suppose a
minimal surface $M$ in $ E^3$ has one-to-one
Gauss map $G:M\to S^2$, then a relatively
compact domain
$D\subset M$ is stable if $G(D)$ is contained in
a hemisphere of $S^2$. 

Schwarz's result was
generalized by J. L. Barbosa and M. do Carmo
(1976) to the following: If the area
$A(G(D))$ of the Gauss image $G(D)$ is less
than $2\pi$, then $D$ is stable. 

For a minimal surface $M$ in 
Euclidean $m$-space, J. L. Barbosa and M. do
Carmo (1980a) proved that if $D\subset M$ is
simply-connected and that $\int_M |K|dV
<{4\over 3}\pi$, then $D$ is stable.

J. Peetre (1959) obtained the following: Let
$D$ be a domain on the unit sphere $S^2$ and
$\tilde D$ a geodesic disc on the sphere
having the same area as $D$. Then
$\lambda_1(D)\geq\lambda_1(\tilde D)$.

As an analogue to Bernstein's theorem, M. do
Carmo and C. K. Peng (1979), and
independently by Fischer-Colbrie and
Schoen (1980), proved that planes are the
only stable complete minimal surfaces
in $E^3$. 

H. Mori (1977) studied minimal surfaces in
3-sphere and proved the following:

Let $D$ be a relatively compact domain on a
minimal surface $M$ of a unit 3-sphere
$S^3$. Suppose that $\sup_D K=K_0<1$ and 
$$\int_D(1-K)dA<{1\over{54\pi}}\cdot{{1-K_0}
\over{2-K_0}}.\leqno (5.41)$$ Then $D$ is
stable.

 Barbosa and do Carmo
(1980a) studied the stability of minimal
surfaces in 3-sphere and in hyperbolic
3-space and  improved Mori's result to the
following:

(1) Let $f:M\to S^3$ be
a minimal immersion of a surface $M$ into the
unit 3-sphere. Assume that $D\subset M$ is
simply-connected and that $$\int_D
(2-K)dV<2\pi,\leqno (5.42)$$ then $D$ is stable.

Furthermore, the result is sharp in the
following sense: given $\delta>0$ there
exists a minimal immersion $f:M\to S^3$ and
an unstable domain $D_\delta\subset M$ such
that
$$\int_{D_\delta}(2-K)dV=2\pi+\delta,\leqno(5.43)$$ and 

(2) Let $f:M\to H^3$ be
a minimal immersion of a surface $M$ into the
unit hyperbolic 3-space with constant curvature
$-1$. Assume that
$D\subset M$ is simply-connected and that
$$\int_D |K|dV<2\pi,\leqno (5.44)$$ then $D$ is
stable.

Barbosa and do Carmo (1980b) also considered
stability for minimal immersions in higher
dimensional real space form and obtained the
following:

Given a minimal surface $M$ of a hypersphere of
radius $r$ in $ E^m$, let $D$ be a
simply-connected relatively compact domain in
$M$. If $$\int_D\left({2\over
{r^2}}-K\right)
dA<{{2n-6}\over{2n-7}}\pi,\leqno (5.45)$$ then
$D$ is stable.

D. Hoffman and R. Osserman (1982) were able to
prove the stability of $D$ under a weaker
condition:$$\int_D\left({2\over
{r^2}}-K\right)dA<{4\over 3}\pi.\leqno (5.46)$$

\vskip.1in
\noindent{\bf 5.10.2. Indices of minimal submanifolds}

The index of every compact minimal surface in
a Riemannian manifold is always finite. For
complete minimal surfaces in $E^3$, D.
Fischer-Colbrie (1985) obtained a direct
relationship between index and total
curvature. She proved that the index of a
complete minimal surface of finite total
curvature in $E^3$ is equal to the index of its
Gauss map; thus the index of a complete minimal
surface in $E^3$ is finite if and only if its
total curvature is finite. 

Since a nonplanar periodic minimal surface in
$E^3$ has infinite total curvature, they
have infinite index; hence, the index of a
complete Scherk surface is infinite.
The helicoid have infinite index as well
[Tuzhilin 1992].

J. Tysk (1987) showed that for complete
minimal surfaces $M$ in $E^3$ one has
$i_M\leq 7.68183d$, where $d$ is the degree
of the Gauss map of $M$. The number 7.68183
is not optimal, since a catenoid has index
one and $d=1$. It is not known whether the
optimal value is 1.

S. Nayatani (1993) related the upper and lower
bounds for the index with the degree of the
Gauss map and the genus of the minimal
surface.

The indices of the catenoid and the Enneper
surface  are both equal to one. This follows
immediately from the fact that the extended
Gauss map of these genus one  surfaces is a
conformal diffeomorphism to the sphere.
Osserman (1964) proved that the catenoid and
the Enneper surface were the only complete
minimal surfaces satisfying this property.
S. Montiel and A. Ros (1990) showed that the
catenoid and the Enneper surface are the only
complete minimal surface in $E^3$ with index
one.

S. Y. Cheng and J. Tysk
(1988) showed that if $M$ is a complete 
orientable minimal surface in $E^3$ with
embedded ends which is neither a plane nor a
catenoid, then the index of $M$ is at least
2.  F. J. L\'opez and A. Ros (1989) showed
that the catenoid and the Enneper surface
are in fact the only complete orientable
minimal surfaces in $E^3$ with index one. 

M. Ritor\'e and A. Ros (1996) studied
the structure of the space of compact
index one minimal surfaces embedded
in flat $3$-tori and obtained the following: 

Let $M$ be a complete noncompact orientable
index one minimal surface properly embedded in
the quotient of $E^3$ by a discrete subgroup
$\Gamma$ of translations. Then one of
the following must occurs:

(i) $M$ is a catenoid in $E^3$; 

(ii) $M$ is a Scherk surface with genus zero
and four ends in $E^3/\Gamma$;

(iii) $M$ is a Scherk surface with genus zero
and four ends in $T^2\times \hbox{\bf R}$; 

(iv) $M$ is a helicoid with total
curvature $-4\pi$ in $E^3/\Gamma$.

If $M$ is a complete oriented minimal
surface of genus zero in $E^3$ which is not
the plane, the catenoid, or the Enneper's
surface, the index of $M$ is at least 3
[Nayatani 1990b].

Montiel and Ros (1990) and N. Ejiri and M.
Kotani (1993) proved that a generic complete
orientable finitely branched minimal surface
of genus zero in $E^3$ with finite total
curvature $4d\pi$ has index $2d-1$ and nullity
$3$. Here ``generic'' means that the Gauss map
of $M$ belongs to the complement of an analytic
subvariety of the space of such maps.

Ejiri and Kotani (1993) also defined the
notion of (nonembedded) flat ends 
for finitely branched complete minimal surface
in $E^3$ and proved the following:

(a) A complete orientable
finitely branched minimal surface in $E^3$
with finite total curvature has nullity $\geq
4$ if and only if its Gauss map can be the
Gauss map of a complete finitely
branched flat-ended minimal surface in $E^3$;

(b) The index and the nullity of a
complete orientable finitely branched minimal
surface  of genus zero in $E^3$ with total
curvature $-8\pi$ both equal to 3.

Let $M$ be a complete  submanifold of
arbitrary codimension in a Riemannian
manifold $N$, and $\varphi$ a smooth vector
field on $N$. The horizon of
$M$ with respect to $\varphi$, denoted by
$H(M;\varphi)$, is the set of all points of
$M$ at which $\varphi$ is a tangent vector
of $M$. A connected subset $D$ of $M$
is called visible with respect to $\varphi$
if $D$ is disjoint from $H(M; \varphi)$. The
number of components of $M- H(M;\varphi)$
is called the vision number of $M$ with
respect to $\varphi$ and is denoted by
$\nu(M;\varphi)$.

Let $\varphi\sb n,\varphi\sb l$ and $\varphi
\sb p$, respectively, be the variation
vector fields in $E^3$
associated with a 1-parameter family of
translations $\tau^n_t$ in the
direction of a unit vector $n$, a
1-parameter family of rotations $p^l_t$
around a straight line $l$, and a
1-parameter family of homothetic expansions
$\mu^p_t$ with center $p$. 

J. Choe (1990) proved the following: 

(a) For any unit vector $n$ in $E^3$ and any
minimal surface $M$ in $E^3$ of finite
total curvature, orientable
or nonorientable, $i_M\geq
\nu(M;\varphi_n)-1$; 

(b) Let $M$ be a complete minimal surface in
$E^3$ of finite total curvature. If each
end of $M$ is embedded and the normal
vectors at the points of $M$ at infinity are
all parallel to a line $l$, then
$i_M\geq\nu(M;\varphi)-1$; and

(c) Let $M$ be a complete minimal
submanifold in a real space form $R^m(c)$
and $\varphi$ a Killing vector field on
$R^m(c)$. If
$M$ is compact, then
$i_M\geq \nu(M;\varphi)-1$,
and otherwise $i_M\geq
\overline\nu(M;\varphi)$, where
$\overline\nu(M;\varphi)$ is the number of
bounded components of $M- H(M;\varphi)$.

Using these Choe showed the following:

(i) The index of the Jorge-Meeks minimal
surface with $k$ ends is at least $2k-3$;

(ii) The index of the Hoffman-Meeks minimal
surface of genus $g$ [Hoffman-Meeks 1990a] is
at least $2g+1$;

(iii) The index of Lawson's minimal surface
$\xi_{m,k}$ of genus $mk$ in $S^3$ [Lawson
1970] is at least $\max(2m+1,2k+1)$;

(iv) The index of the minimal hypersurface
$S^p(\sqrt{p/(p+q)})\times
S^q(\sqrt{q/(p+q)})$ in $S^{p+q+1}$ is at
least 3; 

(v) The index of any complete immersed
nonorientable minimal surface in $E^3$ of finite
total curvature which is conformally equivalent
to a projective plane with finite punctures is
at least 2; 

(vi) The plane, Enneper's surface, and the
catenoid are the only three complete
immersed orientable minimal surfaces of
genus zero and index less than three in $E^3$;

(vii) The index of Chen-Gackstatter surface is
3. 

S. Nayatani (1993) showed that the index of a
Hoffman-Meeks minimal surface in $E^3$ with
3-ends of genus $k$ is $2k+3$ for $k\leq 37$.
He also proved that if a complete oriented
minimal surface in $E^m$ has finite total
curvature, then it has finite index [Nayatani
1990a]. 

Let $M$ be an $n$-dimensional complete
minimal submanifold of
$E^m$. Then there exists a constant
$C_{n,m}$ depending only on the dimensions
$n,m$ such that the index of $M$ is less
than or equal to $C_{n,m}$
times the total scalar
curvature, that is $$i_M\leq C_{n,m}\int_M
S^{n/2}dV,\leqno (5.47)$$ which was proved by P.
B\'erard and G. Besson (1990), S. Y. Cheng
and J. Tysk (1988) and N. Ejiri (1991). 

In 1994 Cheng and Tysk proved that there 
exist constants $C_m$ depending only on $m$
such that the index of a branched complete
minimal surface in $E^m$ satisfies $$i_M\leq
C_m\int_M(-K)dA,\leqno (5.48)$$ where $K$ is the
Gaussian curvature of $M$.

For complete oriented minimal surfaces in
$E^4$, S. Y. Cheng and J. Tysk (1988) proved
that the index is less than or equal to
$12.72\left({1\over{2\pi}}\int(-KdA)\right)$.

For a minimal hypersurface $M$ of $E^{n+1}$
with $n\geq 3$, J. Tysk (1989) showed that
the index of $M$ satisfies
$$i_M\leq \omega_n^{-1}\left(
{{\sqrt{e}(n-1)2^{
2n+3}}\over{n-2}}\right)^n\int_M
S^{n/2}dV,\leqno(5.49)$$ where $\omega_n$ is the
volume of the unit ball in $E^n$. 

The index of a great 2-sphere in $S^3$ is
one. For a compact orientable minimal
surface $M$ in $S^3$ which is not totally
geodesic, F. Urbano (1990) proved that the
index of $M$ is at least 5, and equal to
5 if and only if $M$ is the Clifford torus. 

N. Ejiri (1983) showed that if $M$ is a closed
minimal surface of genus zero, fully immersed
in $S^{2n}$ with $n\geq 2$, then the index
of $M$ is greater than or equal to
$2(n(n+2)-3)$. 

In 1994 S. Y. Cheng and J. Tysk proved that if
$M$ is a bounded immersed complete minimal
surface in $S^m$, possibly with boundary, then 
$$i_M\leq C_m\left(2\;\hbox{Area}\,(M)-
\int_M KdA\right),\leqno(5.50)$$ where $C_m$ is
a constant depending only on $m$.

A. A. Tuzhilin (1992) proved that the
indices for all the hyperbolic catenoids and
the parabolic catenoids in the hyperbolic
3-space $H^3$ are zero. He also showed that
the indices for the spherical catenoids in
$H^3$ do not exceed one. A result of do
Carmo and Dajczer (1983) implies the
existence of unstable catenoids in the
one-parameter family of spherical
catenoids. Tuzhilin showed that the indices
for these catenoids equal one.

\vskip.1in
\noindent{\bf 5.10.3. Stability of minimal submanifolds }

M. Ross (1992) proved that every complete
nonorientable minimal surface in $E^3$ of
finite total curvature is unstable. By a refined analysis of moduli and
Teichmuller spaces, J. Jost and M. Struwe
(1990) succeed in applying saddle-point
methods to prove the existence of unstable
complete minimal surfaces of prescribed
genus. 

Let $L\subset E^3$ be a 
discrete lattice with  rank $L=1$ or $2$ and
$f: M\to E^3/L$ be a complete and
connected minimal immersion. M. Ross and C.
Schoen (1994) proved that if
$f$ is stable and $M$ has finite genus,
then $f(M)$ is a quotient of the plane, the
helicoid or a Scherk surface. 

For complete stable minimal surfaces in
a given Riemannian manifold,  D.
Fischer-Colbrie and R. Schoen (1980)
proved the following:

(1)  Stable complete orientable minimal
surfaces in complete orientable
$3$-mani\-folds with non-negative Ricci
curvature are totally geodesic;  
 
(2) Let $N$ be a complete Riemannian
3-manifold with nonnegative scalar curvature
$\rho$ and let $M$ be a complete stable
orientable minimal surface in $N$. Then 

(2.1) if $M$ is compact, then either
$M$ is conformally equivalent to the Riemann
sphere $S^2$ or else it is a totally
geodesic flat torus. Furthermore, if
$\rho>0$, the latter case cannot occur, and

(2.2) if $M$ is not compact, then $M$ is
conformally equivalent to the complex plane
or a cylinder.

H. Lawson and W. Y. Hsiang gave a complete
classification of equivariant stable cones of
codimension one in $E^{n+1}$.

B. Palmer (1991) proved that if $M$
is a complete orientable minimal
hypersurface of $E^{n+1}$ and if there exists
a codimension one cycle $C$ in
$M$ which does not separate $M$, then
$M$ is stable. 

M. do Carmo and C. K. Peng (1980) showed that
if an oriented stable complete minimal
hypersurface $M$ in $E^{n+1}$ satisfies
$\int_M SdV<\infty$, then $M$ is a
hyperplane. 
Using Simons' result, P. B\'erard (1991) 
proved that if an oriented stable complete
minimal hypersurface $M$ in $E^{n+1}$ has
finite total scalar curvature, that is,
$\int_M S^{n/2}dV<\infty$, then $M$ is a
hyperplane for $n\leq 5$. The same result was
proved recently  by  Y. B. Shen and  X. H.
Zhu (1998) for $n>5$, using the interior 
curvature estimate and Gromov's compactness
theorem.

There are  results on stability of higher
dimensional minimal submanifolds in real space
forms obtained by various geometers. For
instance, B. H. Lawson and J. Simons (1973)
proved that there does not exist stable
compact minimal submanifold in $S^m$. 

S. P. Wang and S. W. Wei (1983) showed that,
for any dimension $n\geq 2$, there is a
non-totally geodesic complete absolute
area-minimizing hypersurface in a
hyperbolic $(n+1)$-space $H^{n+1}$.

Further results on indices and stabilities of
minimal submanifolds can be found in \S11.3,
\S11.4 and \S16.11.

\vfill\eject

\section{Submanifolds of finite type}

The study of order and submanifolds of finite type
began in the 1970s through Chen's
attempts to find the best possible estimate of
the total mean curvature of an isometric
immersion of a compact manifold in Euclidean
space and to find a notion of ``degree'' for
submanifolds in Euclidean space.

The main objects  in algebraic geometry are
algebraic varieties. One can define the
degree of an algebraic variety by  its
algebraic structure.  On the other hand,
although every Riemannian manifold can be
realized as a submanifold in Euclidean space
according to Nash's embedding theorem, one
lacks the notion of the degree of a
submanifold in  Euclidean space.  Inspired by
this observation, the notions of order and  
submanifolds of finite type were 
introduced in [Chen 1979a,1984b]. 

\subsection{Spectral resolution}

 Let $(M,g)$ be a compact Riemannian 
$n$-manifold. Then 
the eigenvalues of the Laplacian $\Delta$ form
a discrete infinite sequence:
$0=\lambda_0<\lambda_1<\lambda_2<\ldots
\nearrow\infty.$
Let $$V_k=\{f\in C^{\infty}(M):\Delta
f=\lambda_k f\}$$ be the  eigenspace of
$\Delta$ associated with eigenvalue
$\lambda_k$. Then each
$V_k$ is finite-dimensional. Define an inner
product $(\;,\,)$ on $C^{\infty}(M)$
by $(f,h)=\int_M fh\, dV$. 
Then $\sum_{k=0}^{\infty} V_k$ is dense in
$C^{\infty}(M)$ (in
$L^2$--sense). If we denote by $\hat{\oplus} V_k$ the
completion of $\sum V_k$, we have
$C^{\infty}(M)=\hat\oplus_k V_k.$

For each function $f\in C^{\infty}(M)$, let
$f_t$ denote the projection of $f$ onto the
subspace $V_t$. We
have the spectral resolution (or
decomposition):
$f=\sum_{t=0}^{\infty} f_t\; (\hbox{in
$L^2$-sense}).$

Because $V_0$ is 1-dimensional,  there is a
positive integer $p\geq 1$ such that
$f_p\not=0$ and
$f-f_0=\sum_{t\geq p}f_t,$
where $f_0\in V_0$ is a constant. If there are infinite
many $f_t$'s which are nonzero, put
$q=+\infty$; otherwise, there is an integer
 $q\geq p$ such that
$f_q\not=0$ and $f-f_0=\sum_{t=p}^q
f_t.$ 

If $x: M\to  E^m$ is an isometric
immersion  of a compact Riemannian
$n$-manifold $M$ into $ E^m$ (or, more
generally, into a pseudo-Euclidean space), for
each coordinate function
$x_A$ we
have $x_A=(x_A)_0+\sum_{t=p_A}^{q_A}(x_A)_t.$
We put $$p=\inf_{A} \{p_A \}\quad\hbox{and}\quad
q=\sup_A\{q_A\},\leqno(6.1)$$ where $A$ ranges
over all $A$ such $x_A-(x_A)_0\not=0.$ 

Both $p$ and $q$ are well-defined geometric
invariants such that $p$ is a positive integer and $q$ is
either $+\infty$ or an  integer $\geq p$. 
Consequently, we have the  spectral
decomposition of $x$ in vector form:
$$x=x_0+\sum_{t=p}^q x_t,\leqno(6.2)$$
which is called the spectral
resolution (or decomposition) of the immersion
$x$.

\subsection{Order and type of immersions}

For a compact manifold $M$, the set  $T(f)=\{ t\in  \hbox{\bf Z}: f_t\not=\hbox{constant}\}$ of
a function $f$ on $M$ is called the order of $f$. The smallest element in  $T(f)$ is called the lower order of $f$ and the supremum of $T(f)$ is called the  upper order of $f$. 

A function $f$ is said to be of  finite type if $T(f)$ is a finite set,  that is, if its  spectral resolution contains only finitely many non-zero terms. Otherwise $f$ is said to be  of infinite type.

Let $x: M\to  E^m$ be an isometric
immersion  of a compact Riemannian
$n$-manifold $M$ into $ E^m$ (or, more
generally,  a pseudo-Euclidean space). Put
$$T(x)=\{t\in \hbox{\bf Z}:
x_t\not=\hbox{constant map}\}.\leqno(6.3)$$ The
immersion
$x$ or the submanifold $M$ is said to be of
k-type if $T(x)$ contains exactly
$k$ elements. Similarly one can define the
lower order  and the upper order of the
immersion. 
The immersion $x$ is said to be  of finite
type if  its upper order $q$ is finite; and
the immersion is said to be of infinite type
if its upper order  is $+\infty$.
The constant vector $x_0$ in the spectral
resolution is the center of mass of
$M$ in $ E^m$.

 One cannot make the spectral resolution of
a function  on a non-compact Riemannian
manifold in general. However, it remains possible to define
the notion of a  function or an immersion of
finite type and the related notions of
order and type. 

For example, a function $f$ is said to be  of
finite type if it is a finite sum of
eigenfunctions of the Laplacian and an
immersion $x$ of a
non-compact manifold is said to be of finite
type if it admits a finite spectral
resolution 
$x=\sum_{t=p}^q  x_t$ for some natural
numbers $p$ and
$q$; otherwise, the immersion is said to be
of infinite type. A $k$-type immersion is
said to be of null $k$-type if the component
$x_0$ in the spectral resolution is
non-constant.

A result of [Takahashi 1966] can now be
rephrased by saying that 1-type submanifolds
of $ E^m$ are precisely those which are
minimal in $E^m$ or minimal in some
hypersphere of $E^m$. In that regard,
submanifolds of finite type provide vast
generalization of minimal submanifolds. 

Let $x: M\to E^m$ be a $k$-type isometric 
immersion whose spectral
resolution is  given by
 $$x =c+x_{1}+\ldots +x_{k},\hskip.0in 
\Delta x_{i}=\lambda_{i}x_{i},\hskip.2in 
\lambda_{1}<\ldots<\lambda_{k},$$ where
$c$ is  a constant vector in $E^m$
and $x_{1},\ldots,x_{k}$ are non-constant
eigenmaps of the Laplacian. 
For each $i\in\{1,\ldots,k\}$ we put
$$E_{i}=\hbox{\rm Span}\{x_{i}(u):u\in
M\}. $$
Then each $E_i$ is a linear subspace of 
$E^m$. The immersion  is said to be 
linearly independent if the $k$ subspaces
$E_{1},\ldots,E_{k}$ are 
linearly independent, that is,
the dimension of the subspace spanned by vectors in
$E_{1}\cup\ldots \cup E_{k}$ is equal to
$\dim E_{1}+\ldots+\dim E_k$. The
immersion is said to be  orthogonal if
the $k$ subspaces $E_{1},\ldots,E_{k}$
defined above are mutually orthogonal in
$E^m$ [Chen 1991a]. 

Let $x:M\to E^m$ be an isometric
immersion of finite type. B. Y. Chen and
M. Petrovic (1991) proved the following:

(a)  The immersion $x$ is linearly
independent if and only if it satisfies
Dillen-Pas-Verstraelen's condition, that
is, it satisfies $\Delta x=Ax+B$ for some
$m\times m$ matrix $A$ and some vector
$B\in E^m$;

(b) The immersion $x$ is orthogonal if and
only if it satisfies $\Delta x=Ax+B$ for
some symmetric $m\times m$ matrix $A$ and
some vector $B\in E^m$.

Linearly independent submanifolds,
equivalently submanifolds satisfying
condition $\Delta x=Ax+B$, have also been
studied by C. Baikoussis, D. E. Blair, F.
Defever, F. Dillen, O. J. Garay, T. Hasanis,
J. Pas, M. Petrovic, T. Vlachos, L.
Verstraelen, and others.

If $x: M\to E^m$ is an isometric immersion
of  null $k$-type whose spectral resolution
satisfies $\Delta x_{j}=0$, then the
immersion is called weakly linearly
independent if the $k-1$ subspaces
$E_{1},\ldots,E_{j-1},E_{j+1}, \ldots,E_{k}$
are linearly independent; and the immersion
is called weakly orthogonal if 
$E_{1},\ldots,E_{j-1},E_{j+1},\ldots,E_{k}$ 
are mutually orthogonal. 

The mean curvature vector of a submanifold  of
non-null finite type in $E^m$ satisfies 
$\Delta   H =A  H$ for
some $m\times m $ matrix $A$ if and only if
$M$ is linearly independent. On the other
hand, if $M$ is a submanifold of null finite
type  
$E^m$, then $M$ satisfies
$\Delta   H =A  H$ for
some $m\times m $ matrix $A$ if and only if
$M$ is weakly linearly independent.

For surfaces in $E^3$ Chen proved in 1994
the following results:

(1) Minimal surfaces and open parts of 
circular cylinders are the only ruled 
surfaces satisfying
$\Delta  H =A  H$ for
some $3\times 3$  matrix $A$;

(2) Minimal surfaces and open
parts of circular cylinders are the only 
finite type surfaces satisfying the condition
$\Delta  H =A  H$ for
some $3\times 3$ singular matrix $A$;

(3) Open parts of circular cylinders are the
only tubes satisfying the condition $\Delta
  H=A  H$  for some
$3\times 3$ matrix A.

See [Chen 1996d] for the details and for 
related results.

\subsection{Equivariant submanifolds
as minimal submanifolds in their adjoint hyperquadrics}

Let $f:M\to E^m$ be a nonminimal
linearly independent isometric immersion
and let $A$ denote the $m\times m$ matrix 
 associated with the immersion $f$
defined in \S6.2. Then, for any point
$u\in M$, the equation
 $$\left<Au,u\right>:=\sum_{i,j}^{m}
a_{ij}u_{i}u_{j}=c_{u},\,\,\,\,\,
$$ with $c_{u}=\left<Ax,x\right>(u)$
defines a quadric
$Q_u$ of $E^m$, where  $u=(u_{1},
\ldots,u_{m})$ be a Euclidean coordinate
system on $E^m$. The hyperquadric defined
above is called the  adjoint
hyperquadric at $u$.
If $f(M)$ is contained in
an adjoint hyperquadric $Q_u$ for some point $u \in
M$, then all of the adjoint hyperquadrics $\{ Q_{u}: u\in
M\}$ coincide, which give rise to a common adjoint
hyperquadric, denoted by $Q$. This common
hyperquadric $Q$ is called the 
adjoint hyperquadric of the immersion.

Suppose $f:M\to E^m$ is a linearly
independent isometric immersion  of a
compact Riemannian manifold into $E^m$.
Then $M$ is immersed into its adjoint
hyperquadric by the immersion $f$ if and
only if the immersion  is spherical, that
is, $f(M)$ is contained in a hypersphere
of $E^m$.

A nonminimal linearly independent
isometric immersion $f:M\to E^m$ of a
Riemannian manifold is orthogonal if and
only if $M$ is immersed  as a minimal
submanifold of its adjoint hyperquadric
by the immersion $f$.

Although an equivariant isometric 
immersion of a compact homogeneous
Riemannian manifold into Euclidean
$m$-space is of finite type, it is not
necessary a minimal submanifold of any
hypersphere of the Euclidean $m$-space in
general. However, we have the following
general result of Chen (1991a):
 
If $f: M\to E^{m}$
is an equivariant isometric immersion of a
compact homogeneous  Riemannian manifold
into Euclidean $m$-space, then $M$ is
isometrically immersed as a minimal
submanifold in its adjoint hyperquadric by
the immersion $f$.

\subsection{ Submanifolds of finite type}

 Although the class of
submanifolds of finite type is huge, it
consists of ``nice'' submanifolds of
Euclidean spaces. For example, all  minimal
submanifolds of Euclidean space and all
minimal submanifolds of hyperspheres are of
1-type and vice versa. Also, all parallel
submanifolds of Euclidean space and all
compact homogeneous Riemannian manifolds
equivariantly immersed in Euclidean space
are of finite type.

Given a  natural number
$k$, there do exist infinitely many 
non-equivalent $k$-type
submanifolds of codimension 2 in Euclidean 
space. The simplest examples of such codimension
two $k$-type submanifolds of Euclidean space 
are the Riemannian products of the
$(n-1)$-dimensional Euclidean space $
E^{n-1}$ with  any $(k-1)$-type closed curves
in $ E^3$. 

Also, according to a result of C. Baikoussis, F. Defever, T.
Koufogiorgos and L. Verstraelen (1995),
for any natural number $k$, there exist $k$-type isometric
immersions of flat tori in
$ E^6$ which are not product immersions.

\vskip.1in
\noindent{\bf 6.4.1. Minimal polynomial criterion}
 
Compact finite type submanifolds are 
 characterized by the minimal polynomial
criterion which establishes the existence of a
polynomial $P$ of the least degree for
 which $P(\Delta) H=0$, where
$ H$ is the mean curvature vector of
the submanifold  and deg $ P = k$ for a
$k$-type immersion [Chen 1984b].

For general submanifolds Chen and M.
Petrovic (1991) proved the following: Let
$f:M\to  E^m$ be an isometric immersion.
Then $f$ is of finite type if and only if
there exists a vector $c\in  E^m$ and
a polynomial $P(t)$ with simple roots such
that $P(\Delta)(x-c)=0$. Furthermore, in this
case, the type number of $f$ is $\leq \deg P$.

\vskip.1in
\noindent{\bf 6.4.2. A variational minimal principle}

Just like minimal submanifolds, finite type
submanifolds  are characterized by a
variational minimal principle in a natural
way; namely as critical points of directional
deformations 
[Chen-Dillen-Verstraelen-Vrancken 1993]. 

Let $f:M\to  E^m$  be an isometric 
immersion of a compact
Riemannian manifold $M$ into $ E^m$.
Associated with each $E^m$-valued
vector field $\xi$ defined on $M$, there 
is a variation $\phi_t$, defined by
$$\phi_t(p):=f(p)+t \xi(p),\quad p\in M,\quad 
t\in (-\epsilon,\epsilon),\leqno(6.4)$$ where $\epsilon$ is a
 sufficiently small positive number. 

 Let $\Cal D$ denote the class of all
 variations acting on
the submanifold $M$ and let  $\Cal{E}$ denote a nonempty
subclass of
$\Cal D$. A compact submanifold $M$ of 
$E^m$ is said to
satisfy the variational minimal principle 
in the class $\Cal{E}$
if $M$ is a critical point of  the volume 
functional for all
variations in $\Cal{E}$. 

Directional deformations were introduced 
by K. Voss in 1956. 
Directional deformations are defined as 
follows : let $c$ be a fixed
vector in $ E^m$ and let $\phi$ be a
smooth  function defined on
the submanifold $M$.  Then we have a
variation given by
$$\phi_t^{\phi c}(p):=f(p)+t \phi(p) c,\quad p\in
M\quad t\in (-\epsilon,
\epsilon).\leqno(6.5)$$
Such a variation is called a directional 
deformation in the direction
$c$.  For each natural number $q\in
\hbox{\bf N}$, define
$\Cal{C}_q$ to be the class of all
directional deformations given by smooth 
functions $\phi$ in $\sum_{i\geq q}V_i$. 

Chen, Dillen, Verstraelen and Vrancken
(1993) proved the following:

(1) There are no compact submanifolds 
in $E^m$ which satisfy the variational
minimal principle in the classes
$\Cal{C}_0$ and $\Cal{C}_1$.

(2) A compact submanifold  $M$ of  $ E^m$
is of finite type if and only if it satisfies
the variational minimal principle in the class
$\Cal{C}_q$ for some $q\geq2$.

\vskip.1in
\noindent{\bf 6.4.3. Diagonal immersions}

 Let  $y_i:M\to  E^{n_i}$,
$i=1,\ldots,k,$ be $k$ isometric immersions
of a Riemannian manifold $M$ into
$ E^{n_i}$, respectively. For any $k$ real
numbers
$c_1,\ldots,c_k$
 with $c_1^2+\ldots+c_k^2=1$, the immersion
$$f=(c_1y_1,\ldots,c_ky_k):M\to
E^{n_1+\cdots+n_k}\leqno(6.6)$$ is also an
isometric immersion, which is called a diagonal
immersion of $y_1,\ldots,y_k$. If 
$y_1,\ldots,y_k$ are of finite type, then
each diagonal immersion of $y_1,\ldots,y_k$
is also of finite type.

\vskip.1in
\noindent{\bf 6.4.4. Curves of finite type}

 A closed curve in $E^m$ is of finite type if
and only if the Fourier series expansion of
each coordinate function of the curve has
only finite nonzero terms. 

The only curves  of finite type in $ E^2$
are  open portions of  circles or  lines, hence
plane curves of finite type are of 1-type. In
contrast with plane curves,  there exist infinitely
many non-equivalent curves of $k$-type in
$E^3$ for each $k\in\{2,3,4,\cdots\}$. 

 Closed curves of finite type in a
Euclidean space are rational curves.
Furthermore, a closed curve of finite type
in $ E^3$ is of 1-type  if and only
if it lies in a 2-sphere
[Chen-Deprez-Dillen-Verstraelen-Vrancken
1990]. 

\vskip.1in
\noindent{\bf 6.4.5. Finite type submanifolds 
in Euclidean space}

B. Y. Chen, J. Deprez and P. Verheyen
(1987) proved that an isometric immersion
of a symmetric space $M$ of compact type
into a Euclidean space is of finite type
if and only if the immersion maps all
geodesics of $M$ into curves of finite type.

 In 1988, Chen proved that a surface in $
E^3$ is of null 2-type if and only if it
is an open portion of a circular
cylinder. Also, Chen and H. S. Lue (1988)
proved the following:

(a) a $2$-type submanifold $M$ in a Euclidean
$m$-space with parallel mean curvature vector
is either spherical or null; 

(b) every 2-type hypersurface of constant mean
curvature in Euclidean space is of null
2-type;

(b) every compact 2-type hypersurface of a
Euclidean space has non-constant mean
curvature.

Null 2-type hypersurfaces and open
portions of  hyperspheres are the only
hypersurfaces of Euclidean space with
nonzero constant mean curvature and
constant scalar curvature  (cf. [Chen
1996d]). T. Hasanis and T. Vlachos
(1995a) proved that null 2-type
hypersurfaces of $E^4$ have nonzero mean
curvature and constant scalar curvature. 
They also showed that a 3-type surface in
$E^3$ has non-constant mean curvature
[Hasanis-Vlachos 1995a].

Chen proved in 1987 that a tube in $E^3$
is of finite type if and only if it is an
open portion of a circular cylinder.
 O. J. Garay (1988b) showed that
open portions of hyperplanes are the only
cones of finite type in $ E^{n+1}$.

A ruled surface in $E^3$ is of finite
type if and only if it is open portion of
a plane, a circular cylinder or a
helicoid. In particular, a flat surface
in  $ E^3$ is of finite type if and only
if it is an open portion of a plane or a
circular cylinder
[Chen-Dillen-Verstraelen-Vrancken 1990]. 
F. Dillen (1992) considered finite type
ruled submanifolds of Euclidean space
and  proved that a ruled submanifold of
Euclidean space is of finite type if and
only if it is a part of a cylinder on a
curve of finite type or  an open portion
of a generalized helicoid.

Chen and Dillen (1990a) proved that 
open portions of spheres and circular
cylinders are the only quadrics of finite type
in $ E^3$.  Further, Chen, Dillen and H. Z.
Song (1992) showed that a  quadric
hypersurface
$M$ of $ E^{n+1}$ is of finite type if and
only if it is one of the following
hypersurfaces:

(1) a hypersphere;

(2) a minimal algebraic cone $C_{k,n-k-1}, \,
0<k<n-1,$ over an $(n-1)$-dimen\-sional
generalized Clifford torus (defined
by (5.14));

(3) a spherical hypercylinder $ E^k\times
S^{n-k},\, 0<k<n$;

(4) the standard product embedding of the
product of a linear subspace
$ E^\ell$ and one of the algebraic cones
$C_{k,n-\ell-k-1}$ with $0<k<n-\ell-1$.

F. Defever, R. Deszcz and L. Verstraelen
(1993/4) proved that all compact and noncompact
cyclides of Dupin are of infinite type.

A hypersurface $M$ of $ E^{n+1}$ is called
a  translation hypersurface if it is a
non-parame\-tric hypersurface of the form:
$$x_{n+1}=P_1(x_1)+\cdots+P_n(x_n),$$
where each $P_i$ is a function of one variable. If each
function $P_i$ is a polynomial, the hypersurface is
called a  polynomial translation hypersurface.

F. Dillen, L. Verstraelen, L. Vrancken and G. 
Zafindratafa (1995)
proved that a polynomial translation
hypersurface  of a Euclidean space is of finite
type if and only if it is a hyperplane.

A surface in $ E^3$ is called a  surface of 
revolution if it is generated by a curve $C$ on a
plane $\pi$ when $\pi$ is rotated around a
straight line $L$ in
$\pi$. By choosing $\pi$ to be the 
$xz$-plane and   line $L$ to be the $z$-axis, the surface of
revolution can be parameterized by
$$x(u,v)=(f(u)\cos v,f(u)\sin v,g(u))$$

A surface of revolution  is said to be
of  polynomial kind if $f(u)$ and $g(u)$ are 
 polynomial functions in $u$;
and it  is said to be  of rational kind if  
$g$ is  a rational function in $f$, that is $g$ is
the quotient of two polynomial functions in $f$. 

For finite type surfaces of revolution, Chen and
S. Ishikawa (1993) proved the following.

(1)  A  surface of revolution of
 polynomial kind is of finite type if
and only if either it is an open portion of a plane or it
is an open portion of a circular cylinder; 

(2) A  surface of revolution of
rational kind is of finite type if and only if 
it is an open portion of a plane.

 T. Hasanis and T. Vlachos (1993) proved that 
 a  surface of revolution with  constant mean
curvature in $ E^3$  is  of finite type if
and only if it is an open  portion of a plane,
of a sphere, or of a circular cylinder. J.
Arroyo, O. J. Garay and J. J. Menc\'ia (1998)
proved that the only finite type surfaces in
$E^3$ obtained by revolving an ellipse around
a suitable axis are the round spheres.

A  spiral surface  is a surface in $E^3$
generated by rotating a plane curve $C$
about an axis $A$ contained in the plane of
the curve
$C$ and simultaneously transforming $C$
homothetically relative to a point of $A$.

 C. Baikoussis and L. Verstraelen (1995) proved
that a spiral surface is of finite type if and
only if it is a minimal surface.

It was conjectured that round spheres are the
only compact surfaces of finite type in $
E^3$ [Chen 1987]. All of the results mentioned
above support the conjecture.

In 1988, Chen proved that if $f:M\to E^m$
is an isometric immersion of a Riemannian
$n$-manifold into $E^m$, then the mean
curvature vector $ H$ of $f$ is
an eigenvector of the Laplacian on $M$,
that is, $\Delta H=\lambda
 H$ for some
$\lambda\in\hbox{\bf R}$, if and only if
$M$ is one of the following submanifolds:

(a) 1-type submanifold;

(b) a null 2-type submanifold;

(c) a biharmonic submanifold, that is, a
submanifold satisfies $\Delta
H=0$.

I. Dimitri\'c proved in his doctoral thesis
(Michigan State University 1989) that
open portions of straight lines are the
only biharmonic curves of a Euclidean
space. Chen  proved in 1985 that
minimal surfaces are the only biharmonic
surfaces in $E^3$. I. Dimitri\'c
(1989) extended Chen's result to the
following: Minimal hypersurfaces are the
only biharmonic hypersurfaces of a
Euclidean space with at most two distinct
principal curvatures. 

T. Hasanis and T. Vlachos (1995b) showed
that minimal hypersurfaces are the only
biharmonic hypersurfaces of Euclidean
$4$-space. An alternative proof was given in
[Defever 1996].

It was conjectured by Chen that minimal
submanifolds are the only biharmonic
submanifolds in Euclidean spaces. 

The conjecture was also proved to be true
if the biharmonic submanifold is one of
the following submanifolds:

(1) a spherical submanifold [Chen 1991];

(2) a submanifold of finite type [Dimitri\'c
1989];

(3) a pseudo-umbilical submanifold of
dimension $\ne 4$ [Dimitri\'c 1989].

The conjecture is false if
the ambient space is replaced by a
pseudo-Euclidean $m$-space
with $m\geq 4$ [Chen-Ishikawa 1991].

\vskip.1in
\noindent{\bf 6.4.5. Finite type submanifolds in sphere}

Standard 2-spheres in $S^3$ and products of
plane circles are the only finite type
compact surfaces with constant Gauss
curvature in $S^3$  [Chen-Dillen 1990b].

Every hypersurface  with constant
mean curvature and constant scalar curvature
in $S^{n+1}$ is either totally umbilical or
of 2-type. Consequently, every isoparametric
hypersurface in $S^{n+1}$ is either of 1-type
or of 2-type [Chen 1984b].  Furthermore, every
spherical 2-type hypersurface has constant
mean curvature and constant scalar curvature
 [Hasanis-Vlachos 1991]. 

2-type surfaces in $S^3$ are open portions of
the standard product embedding of the product
of two circles [Hasanis-Vlachos 1991].  

Chen (1996d) proved that a compact
hypersurface of  $S^4(1)\subset
 E^5$ is of 2-type if and only if it is
 one of the following hypersurfaces:

(1) $S^1(a)\times S^2(b)\subset S^4(1)\subset
 E^5$ with $a^2+b^2=1$ and
$(a,b)\not= (\sqrt{{1\over 3}},\sqrt{2\over 3}\,)$ embedded in
 the standard way;

(2) a tubular hypersurface of constant radius
$r\not={{\pi}\over 2}$ about the Veronese surface of constant
curvature ${1\over 3}$ in $S^4(1)$.

For higher dimensional spherical hypersurfaces,
Chen (1991b,1984b) proved  the following:

(a) Let $M$ be a hypersurface of a
unit hypersphere $S^{n+1}(1)$ of $E^{n+2}$ with
at most two distinct principal curvatures. Then
$M$ is of 2-type if and only if $M$ is an open
portion of the product of two spheres
$S^k(a)\times S^{n-k}(b),\,1\leq k\leq n-1,$ such that
 $a^2+b^2=1,$
$(a,b)\not=\left(\sqrt{k/n},\sqrt{(n-k)/n}\,
\right)$.
 
(b)  Let $M$ be a
conformally flat hypersurface of a unit
hypersphere $S^{n+1}(1)$ of
$E^{n+2}$. Then
$M$ is of 2-type if and only if
$M$ is an open portion of the product of two
 spheres $S^k(a)\times S^{n-k}(b)$ for some $k$,
$1\leq k\leq n-1,$ such that 
$a^2+b^2=1$,
$(a,b)\not=\left(\sqrt{k/n},\sqrt{(n-k)/n}
\,\right)$.
 
(c) If $M$ is a 2-type Dupin
hypersurface of a hypersphere $S^{n+1}$ with 
three principal curvatures, then $M$ is an
isoparametric hypersurface.

Chen and S. J. Li (1991) proved that every
3-type hypersurface of a sphere has nonconstant
mean curvature.

A compact spherical submanifold $M\subset
S^{m-1}\subset  E^m$ is called
mass-symmetric if the center of the mass of
$M$ coincides with the center of the
hypersphere $S^{m-1}$. Similarly, a
non-compact submanifold $f:M\to
S^{m-1}_c(r)$ is called mass-symmetric if
its position function admits a spectral
resolution of the
form:
$$x=x_{t_p}+\cdots+x_{t_q},\quad\Delta
x_{t_j}=\lambda_{t_j}x_{t_j}.$$ 

Regardless of codimension, every
mass-symmetric 2-type submanifold of a
hypersphere has constant squared mean
curvature, which is determined by the order
of the submanifolds [Chen 1984b]. 

 Although every compact  2-type surface in
$S^3$ is mass-symmetric [Barros-Garay 1987],
there do not exist mass-symmetric 2-type
surfaces which lie fully in
$S^4$ [Barros-Chen 1987a].

 M. Kotani (1990) studied mass-symmetric
2-type immersions of a topological 2-sphere
into a hypersphere of $E^m$. She proved
that such an immersion is the diagonal
immersion of two 1-type immersions. Y.
Miyata (1988) classified mass-symmetric
2-type spherical immersions of surfaces of
constant curvature (see, also [Garay
1988a]).

\vskip.1in
\noindent{\bf 6.4.6. Finite type submanifolds
in hyperbolic space}

Let $ E^{m}_{s}$ denote the $m$-dimensional
pseudo-Euclidean space with index $s$ endowed with
the standard flat metric given by
$$g=-\sum_{i=1}^{s}\, dx_{i}^{2}+\sum_{j=s+1}^{m}\,
dx_{j}^{2},\leqno(6.7)$$
where $(x_{1},\ldots,x_{m})$ is a standard
coordinate system of  $ E^{m}_{s}$. 

For a  number $r>0$, we denote by
$S_{s}^{m-1}(r)$
 the pseudo-Riemannian sphere and by 
$H_{s-1}^{m-1}(-r)$ the pseudo-hyper\-bolic space
defined respectively by
$$S_{s}^{m-1}(r)=\{ u\in
E_{s}^{m}:\left<u,u\right>=r^{2}\},\leqno(6.8)$$ $$
H_{s-1}^{m-1}(-r)=\{u\in
E^{m}_{s}:\left<u,u\right>=-r^{2}\},\leqno(6.9)$$
where $\left<\; ,\, \right>$ denotes the
indefinite inner product on the
pseudo-Euclidean space. 

Denote by $H^{m-1}$ the 
hyperbolic space which is embedded
standardly in the Minkowski space-time $
E^{m}_1$ by
$$H^{m-1} = \{u \in L^m:\,\left<u,u\right>
=-1\;\;\hbox{\rm and}\;\; t>0\},\leqno(6.10)$$
where $L^m= E^{m}_{1}$  and $t=x_1$ is the
first coordinate of  $ E_{1}^m$.

Let $f: M\to  E^{m}_{s}$ be an isometric
immersion of a pseudo-Riemannian manifold $M$ into
the pseudo-Euclidean $m$-space with index $s$.
Then $M$ is of 1-type if and only if either 

(1) $M$ is a minimal submanifold of $
E^{m}_s$; or

(2)  up to  translations, $M$ is a minimal
submanifold of a pseudo-Riemannian sphere
$S^{m-1}_{s}(r),\; r>0$; or

(3) up to translations, $M$ is
 a minimal submanifold of a pseudo-hyperbolic space
$H^{m-1}_{s-1}(-r),\; r>0$.
\medskip

For 2- and  3-type hypersurfaces in hyperbolic
space, the following are known (cf. [Chen
1996d]).

If $M$ is a hypersurface of the hyperbolic space
$H^{n+1}$, embedded standardly in $
E^{n+2}_1$, then

(a) if $M$ has constant mean curvature and
constant scalar curvature,  $M$ is either of
1-type or of 2-type;

(b)  every 2-type hypersurface of the
hyperbolic space $H^{n+1}$ has constant mean
curvature and constant scalar curvature.

Furthermore, we have

(c) there do not exist compact 2-type
hypersurfaces in the hyperbolic space;

(d) there do not exist null 2-type hypersurfaces
in the hyperbolic space;

(e) there do not exist 3-type hypersurfaces
with constant mean curvature in the
hyperbolic space;

(f) if $M$ is a 2-type hypersurface with at most
two distinct principal curvatures in the 
hyperbolic space, then  up to rigid motions of
$H^{n+1}$, $M$ is an open portion of
$M^{n}_{k,r}$ for  positive integer $k, 2\leq
k \leq n$ and for some  $r>0$, where
$M^n_{k,r}$ is defined by
 $$\aligned M^{n}_{k,r}=\{(t,x_{2},&\ldots,
x_{n+2}):t^{2}-x_{2}^{2}-\ldots -
x_{k}^{2}=1+r^{2},\\ &x_{k+1}^{2}+\ldots
+x^{2}_{n+2}=r^{2}\}.\endaligned\leqno(6.11)$$

(g) if $M$ is a 2-type surface in the hyperbolic
3-space $H^{3}$, then it is a flat surface
and, up to rigid motions of $H^3$, M is an
open portion of $M^{2}_{2,r}$ for some $r>0$.

\vskip.1in
\noindent{\bf 6.4.7. Finite type immersions of
irreducible homogeneous spaces}

Every finite type isometric immersion of a
compact irreducible homogeneous Riemannian
manifold is a screw diagonal immersion, that
is, it is the composition of a diagonal
immersion followed by a linear map [Chen
1996d, Deprez 1988].

 J. Deprez (1988) and T. Takahashi (1988)
proved that every equivariant isometric
immersion of a compact  irreducible
homogeneous Riemannian manifold in
Euclidean space is the diagonal immersion
of some standard 1-type isometric immersions.

 \vfill\eject

\section{Isometric immersions between real space forms}

Historically the study of surfaces of
negative constant curvature in $E^3$ was
closely related with the problem of
interpretation of non-Euclidean geometry.
During 1839-1840 F. Minding investigated
properties of surfaces of constant negative
Gaussian curvature in $E^3$, he discovered
the so called helical surfaces of constant
curvature. Minding showed that surfaces of
revolutions of constant negative curvature
in $E^3$ can be divided into three types;
one of them is known as the pseudo-sphere;
the surface obtained by rotating the
asymptote the so-called curve
of pursuit. In 1868 E. Beltrami
(1835--1900) discovered a close
connection between hyperbolic geometry
and the pseudo-sphere. 

In 1901 D. Hilbert proved that any
complete immersed surface in $E^3$  with
constant positive Gaussian curvature is
a round sphere. The analytic case of
this result was given by H. Liebmann in
1900. A result of D. Hilbert and E.
Holmgren states that there do not 
exist  complete immersed surfaces in
$E^3$ with constant negative Gaussian
curvature in $E^3$. A. W. Pogorelov
(1956a,1956b) and P. Hartman and L.
Nirenberg (1959) showed that a complete
flat surface immersed in $E^3$ is a
generalized cylinder.

Isometric immersions of a Riemannian
$n$-manifold $R^n(c)$ with constant curvature $c$
into another Riemannian $(n+p)$-manifold
$R^{n+p}(\bar c)$ of constant curvature
$\bar c$ were studied by \'E. Cartan (1919). 
He proved, for example, that if $c<\bar c$,
then the existence of an isometric immersion
implies $p\geq n-1$.

Cartan's study was confined to
local phenomena of isometric immersions, which
are generally complicated. However, sometimes,
definite results can be obtained under a global
assumption of completeness or compactness. 
\smallskip
\subsection{Case: $c=\bar c$}

The following global results for $c=\bar c$ and
$p=1$ have been known.

\noindent {\bf 7.1.1. Euclidean case: ($c=\bar
c=0$)}

 If a complete flat Riemannian
$n$-manifold $R^n(0)$ is isometrically immersed
in a Euclidean $(n+1)$-space $ E^{n+1}$,
then $R^n(0)$ is a hypercylinder over a plane
curve, that is, $R^n(0)= E^{n-1}\times C$,
where $ E^{n-1}$ is a Euclidean
$(n-1)$-subspace of $ E^{n+1}$ and $C$ is a
curve lying in a Euclidean plane perpendicular
to $ E^{n-1}$. This result was proved in
[Hartman-Nirenberg 1959].
\smallskip
\noindent {\bf 7.1.2. Spherical case: ($c=\bar
c=1$)}

 If the unit $n$-sphere $S^n$ is
isometrically immersed in $S^{n+1}$, then $S^n$
is embedded as a great $n$-sphere in $S^{n+1}$.
This result is a special case of a theorem in
[O'Neill-Stiel 1963].
\smallskip

\noindent {\bf 7.1.3. Hyperbolic case: ($c=\bar
c=-1$)}

 In this case the situation is more
complicated. K. Nomizu (1924-- ) exhibited
three distinct types of isometric
immersions of hyperbolic plane $H^2(-1)$
into the hyperbolic 3-space $H^3(-1)$ which
are not totally geodesic [Nomizu 1973a.

If an isometric immersion
$f:H^n(-1)\to H^{n+1}(-1)$ has no umbilical
points, the second fundamental form has  nullity
$n-1$ at every point. In this case, the relative
nullity distribution is integrable and each leaf
is an $(n-1)$-dimensional complete totally
geodesic submanifold. 

In general, a $k$-dimensional foliation on a
Riemannian manifold is called totally
geodesic if each leaf is a $k$-dimensional
complete totally geodesic submanifold. D.
Ferus (1973) has showed how to obtain all
$(n-1)$-dimensional totally geodesic
foliations on the hyperbolic space $H^n(-1)$
and that every totally geodesic foliation on
$H^n(-1)$ is the nullity distribution of a
suitable isometric immersion of 
 $H^n(-1)$ into $H^{n+1}(-1)$ with umbilical
points. 

Recently, K. Abe, H. Mori and H.
Takahashi (1997) parametrized the space of
isometric immersions of 
$H^n(-1)$ into $H^{n+1}(-1)$ by a family of
properly chosen (at most) countable $n$-tuples
of real-valued functions. 

Concerning the problem of describing 
isometric immersions of the hyperbolic plane
$H\sp 2(-1)$ into the hyperbolic 3-space $H\sp
3(-1),$ by using the stereographic projection
of the upper sheet of the hyperboloid $x\sp
2+y\sp 2-z\sp 2 = -1$ from the origin on the
plane $z = 1,$ it is possible to parameterize
$H\sp 2(-1)$ by points $ \xi = (\xi\sb 1, \xi\sb
2) \in D $ of the open unit disc. In this way,
the problem of determining isometric immersions
$f\colon H\sp 2(-1)\to H\sp 3(-1) $ is reduced
to solving a degenerate Monge-Ampere equation on
the unit disc: $\det(\partial\sp
2u/\partial\xi_i \partial\xi_j)=0,\ \xi\in
D$, where each isometric immersion corresponds
to a solution $u(\xi\sb 1, \xi\sb 2)$ of such an
equation. 

Using a special family of solutions of the
Monge-Ampere equation, Z. J. Hu and G. S. Zhao
(1997a)  gave a way of constructing
many incongruent examples of immersions
$f\colon H\sp 2(-1) \to H\sp 3(-1)$ with or
without an umbilic set and with bounded or
unbounded principal curvature. 

\smallskip
\subsection{Case: $c\ne\bar c$}

For $c\ne \bar c$ a classical result of 
K. Liebmann and D. Hilbert states that a
complete Riemannian 2-manifold of constant
negative curvature cannot be isometrically
immersed in Euclidean 3-space. L. Bianchi 
proved in 1896 that there exist infinitely
many isometric immersions from a complete
flat surface into $S^3$ (cf. [Bianchi
1903]). Ju. A. Volkov and S. M. Vladimirova
in 1971 and S. Sasaki (1912--1987) in 1973
showed that an isometric immersion of a
complete flat surface into a complete
simply-connected
$H^3(-1)$ is either a horosphere or a set of
points at a fixed distance from a geodesic. 

An isometric immersion of a Riemannian
$n$-manifold
$R^n(c)$ with $n\geq 3$ into a Riemannian
$(n+1)$-manifold $R^{n+1}(\bar c)$ with $c>\bar
c$ is totally umbilical. 
If $p\leq n-1$, Ferus (1975) proved that every
isometric immersion of a complete Riemannian
$n$-manifold $R^n(c)$ into a complete Riemannian
$(n+p)$-manifold $R^{n+p}(c)$ with the same
constant curvature is totally geodesic.

D. Blanusa (1955) proved that a hyperbolic
$n$-space $H^n(-1)$ can be isometrically
embedded into $E^{6n-5}$. On the other hand,
J. D. Moore (1972) proved that if $p\leq
n-1$, then there do not exist  isometric
immersions from a complete Riemannian
$n$-manifold $R^n(c)$ of constant curvature
$c$ into a complete simply-connected
Riemannian manifold
$R^{n+p}(\bar c)$ with $\bar c>c>0$.

 By applying Morse theory, Moore (1977)
proved that if a compact Riemannian
$n$-manifold $R^n(1)$ of constant curvature
1 admits an isometric immersion in $
E^N$ with $N\leq {3\over 2}n$, then $R^n(1)$
is simply-connected, hence isometric to
$S^n(1)$. 

Flat surfaces in $ E^4$ with
flat normal connection were classified in
[Dajczer-Tojeiro 1995a].

D. Ferus and F. Pedit (1996) gave a method for
finding local isometric immersions between
real space forms by integrable systems
techniques. The idea is that, given an isometric
immersion between real space forms of nonzero
different curvatures, the structural equations
can be rewritten as a zero curvature equation
involving an auxiliary (spectral) parameter,
that is, as the flatness equation for a 1-form
with values in a loop algebra. The isometric
immersion thus generates a
one-parameter family of isometric
 immersions with flat normal bundle. A large
class of solutions to the flatness equation can
then be found by integrating certain commuting
vector fields on a loop algebra. 

The isometric
immersions so constructed are real-analytic and
depend on the same number of functions as
predicted by Cartan-K\"ahler theory, though not
all isometric immersions are real-analytic.

\vfill\eject

\section{Parallel submanifolds}

The first fundamental form, that is, the
metric tensor, of a submanifold of a
Riemannian submanifold is automatically
parallel, thus, $\nabla g\equiv 0$
with respect to the Riemannian connection
$\nabla$ on the tangent bundle $TM$.  A
Riemannian submanifold is said to be
parallel if its second fundamental form
$h$ is parallel, that is $\bar\nabla h\equiv
0$  with respect to the connection
$\bar\nabla$ on $TM\oplus T^\perp M$. 

\subsection{Parallel submanifolds in Euclidean space}

The first result on parallel submanifolds 
was given by V. F. Kagan in 1948  who
showed that the class of parallel surfaces in
$E^3$ consists of open parts of planes, round
spheres, and circular cylinders $S^1\times
E^1$. U. Simon and A. Weinstein (1969)
determined parallel hypersurfaces of
Euclidean $(n+1)$-space. A general
classification theorem of parallel
submanifolds in Euclidean space was obtained
by D. Ferus in 1974.

An affine subspace of $E^m$ or a
symmetric $R$-space $M\subset  E^m$,
which is minimally embedded in a hypersphere
of $ E^m$ as described in
[Takeuchi-Kobayashi 1965] is a parallel
submanifold of $ E^m$. The class of
symmetric $R$-spaces includes:

(a) all Hermitian symmetric spaces of compact type,

(b) Grassmann manifolds $O(p+q)/O(p)\times O(q), Sp(p+q)/Sp(p)\times Sp(q),$

(c) the classical groups $SO(m),\,U(m),\,Sp(m)$,

(d) $U(2m)/Sp(m),\, U(m)/O(m)$,

(e) $(SO(p+1)\times  SO(q+1))/S(O(p)\times O(q))$, where $S(O(p)\times O(q))$ is 
the subgroup of $SO(p+1)\times SO(q+1)$
consisting of matrices of the form
$$\begin{pmatrix} \epsilon& 0 & & \\ 0 & A&&\\ &&\epsilon& 0\\
&&0&B
\end{pmatrix},\quad \epsilon=\pm1,\quad A\in
O(p),\quad B\in O(q),$$

(f) the Cayley projective plane ${\Cal
O}P^2$, and 

(g) the three exceptional spaces
$E_6/Spin(10)\times T, E_7/E_6\times
T,$ and $E_6/F_4.$

\vskip.1in
D. Ferus (1974) proved that essentially these
submanifolds exhaust all parallel
submanifolds of $ E^m$ in the following
sense: A complete full parallel submanifold
of the Euclidean $m$-space $ E^m$ is
congruent to 

``(1)" $M=
E^{m_0}\times M_1\times \cdots\times
M_s\subset E^{m_0}\times 
E^{m_1}\times\cdots\times 
E^{m_s}= E^m$, $s\geq 0$, 
or to

``(2)" $M= M_1\times \cdots\times
M_s\subset 
E^{m_1}\times\cdots\times  
E^{m_s}= E^m$, $s\geq 1$,

\noindent where each $M_i\subset E^{m_i}$ is an
irreducible symmetric $R$-space.

Notice that in case (1) $M$ is not contained
in any hypersphere of $ E^m$, but in case
(2) $M$ is contained in a hypersphere of
$E^m$. 

For an $n$-dimensional submanifold $f:M\to
E^m$, for each point $x\in M$ and each unit
tangent vector $X$ at $x$, the vector
$f_*(X)$ and the normal space $T_x^\perp$
determine an $(m-n+1)$-dimensional subspace
$E(x,X)$ of $ E^m$. The intersection of
$f(M)$ and $E(x,X)$ defines a curve $\gamma$
in a neighborhood of $f(x)$, which is called
the normal section of $f$ at $x$ in the
direction $X$. A point $p$ on a plane curve
 is called a vertex if its curvature function
$\kappa(s)$ has a critical point at $p$. 

Parallel submanifolds of $E^m$ are 
characterized by the following simple geometric
property: normal sections of $M$ at each 
point $x\in M$ are plane curves with $x$ as
one of its vertices [Chen 1981a]. 

A submanifold $f:M\to  E^m$ is said to be
extrinsic symmetric if, for each
$x\in M$, there is an isometry $\phi$ of $M$
into itself such that $\phi(x)=x$ and $f\circ
\phi=\sigma_x\circ f$, where $\sigma_x$
denotes the reflection at the normal space
$T_x^\perp M$ at $x$, that is the motion of
$ E^m$ which fixes the space through
$f(x)$ normal to $f_*(T_xM)$ and reflects
$f(x)+f_*(T_xM)$ at $f(x)$. The submanifold
$f:M\to  E^m$ is said to be extrinsic
locally symmetric, if each point $x\in M$ has
a neighborhood $U$ and an isometry $\phi$ of
$U$ into itself, such that $\phi(x)=x$ and $f\circ
\phi=\sigma_x\circ f$ on $U$. In other
words, a submanifold $M$ of $E^m$ is
extrinsic locally symmetric if each point
$x\in M$ has a neighborhood which is
invariant under the reflection of $E^m$ with
respect to the normal space at $x$.

D. Ferus (1980) proved that extrinsic locally
symmetric submanifolds of $ E^m$
have parallel second fundamental form and
vice versa.

A canonical connection on a Riemannian
manifold $(M,g)$ is defined as any metric
connection $\nabla^c$ on $M$ such that
the difference tensor $\hat D$ between
$\nabla\sp c$ and the Levi-Civita connection
$\nabla$ is $\nabla\sp c$-parallel. 

An embedded submanifold $M$ of $ E^m$ is
said to be an extrinsic homogeneous
submanifold with constant principal
curvatures if, for any given
$x,y\in M$ and  a given piecewise
differentiable curve $\gamma$ from $x$ to
$y$, there exists an isometry $\varphi$ of
$ E^m$ satisfying (1) $\phi(M)=M$, (2)
$\phi(x)=y$, and (3)
$\phi_*{}_x{}_{|T^\perp_xM}:T^\perp_x
M\to T^\perp_yM$ coincides with
$\hat D$-parallel transport along $\gamma$.

C. Olmos and C. S\'anchez (1991)
extended Ferus' result and obtained the
following: Let
$M$ be a connected compact Riemannian
submanifold fully in  $ E^m$, and let
$h$ be its second fundamental form. Then the
following three statements are equivalent:

(i) $M$ admits a canonical connection
$\nabla^c$ such that $\nabla^c h=0$, 

(ii) $M$ is an extrinsic homogeneous
submanifold with constant principal
curvatures, 

(iii) $M$ is an orbit of
an $s$-representation, that is, of an
isotropy representation of a semisimple
Riemannian symmetric space. 

The notion of extrinsic $k$-symmetric
submanifold of $ E^m$ was introduced and
classified for odd $k$ in [S\'anchez 1985].
Furthermore, S\'anchez (1992) proved that the
extrinsic $k$-symmetric submanifolds are
essentially characterized by the property of
having parallel second fundamental form with
respect to the  canonical connection of
$k$-symmetric space. Thus, the above result
implies that every
extrinsic $k$-symmetric submanifold of a
Euclidean space is an orbit of an
$s$-representation. 

\subsection{Parallel submanifolds in spheres}
 
Regarding the unit
$(m-1)$-sphere $S^{m-1}$ as an ordinary
hypersphere of $ E^m$, a
submanifold $M\subset S^{m-1}$ is parallel if
and only if $M\subset S^{m-1}\subset 
E^{m}$ is a parallel submanifold of $
E^m$. 

Consequently, Ferus' result implies that
$M$ is a parallel submanifold of $S^{m-1}$ if
and only if $M$ is obtained by a submanifold of
type (2).

\subsection{Parallel submanifolds in hyperbolic spaces}
 
Parallel submanifolds of hyperbolic spaces
were classified in 1981 by M. Takeuchi
(1921-- ) which is given as follows:
For each $c<0$,  let $H^m(c)$ denote the
hyperbolic $m$-space defined by
$$H^m(c)=\{(x_0,\ldots,x_{m})\in
 E^{m+1}:-x_0^2+x_1^2+\cdots+x_m^2=1/c,
x_0>0\}.$$ Assume $M$
is a parallel submanifold of $H^m(\bar
c),\,\bar c<0$. Then

(1) if $M$ is not
contained in any complete totally geodesic
hypersurface of
$H^m(\bar c)$, then
$M$ is congruent to the product $$
H^{m_0}(c_0)\times M_1\times \cdots\times
M_s\subset H^{m_0}(c_0)\times
S^{m-m_0-1}(c')\subset H^{m_0}(\bar c)$$ with
$c_0<0,\, c'>0, 1/c_0+1/c'=1/\bar c,\, s\geq
0$, where  $M_1\times \cdots\times
M_s\subset$ $ S^{m-m_0-1}(c')$ is a parallel
submanifold as described   in Ferus'
result; and 

(2) if $M$ is contained in a
complete totally geodesic hypersurface $N$ of
$H^m(\bar c)$, then $N$ is either isometric to
an  $(m-1)$-sphere, or to a Euclidean
$(m-1)$-space, or to a hyperbolic
$(m-1)$-space. Hence, such parallel
submanifolds reduce to the parallel
submanifolds  described before.

\subsection{Parallel submanifolds in complex projective and complex hyperbolic spaces}

A parallel submanifold $M$ of a
Riemannian manifold $\tilde M$ is
curvature-invariant, that is, for each point
$x\in M$ and $X,Y\in T_xM$, we have $$\tilde
R(X,Y)T_xM\subset T_xM,$$ where
$\tilde R$ is the curvature
tensor of $\tilde M$. Thus,
according to a result of [Chen-Ogiue 1974b],
parallel submanifolds of complex projective
and complex hyperbolic spaces are either
parallel K\"ahler submanifolds or parallel
totally real submanifolds.  

Complete parallel K\"ahler
submanifolds of complex projective spaces and
of complex hyperbolic spaces have been
completely classified in [Nakagawa-Takagi
1976] and in [Kon 1974], respectively  (see
\S15.9 for details).

[Naitoh 1981] showed that the classification
of complete totally real parallel
submanifolds in complex projective spaces is
reduced to that of certain cubic forms of
$n$-variables and [Naitoh-Takeuchi 1982]
classified these submanifolds by the theory
of symmetric bounded domains of tube type. 

The complete classification of parallel
submanifolds in complex projective spaces and
in complex hyperbolic spaces was given in
[Naitoh 1983].

\subsection{Parallel submanifolds in quaternionic projective spaces}

Parallel submanifolds of a quaternionic
projective $m$-space or its non-compact dual
were classified in [Tsukada 1985b]. 
Such submanifolds are parallel totally real
submanifolds
 in a totally real totally
geodesic submanifold $RP^m$, or  parallel
totally real submanifolds  in a
totally complex totally geodesic submanifold
$CP^m$, or parallel complex submanifolds
in a totally complex totally geodesic
submanifold $CP^m$, or parallel totally
complex submanifolds in a totally geodesic
quaternionic submanifold $HP^k$ whose
dimension is twice the dimension of the
parallel submanifold. 

\subsection{Parallel submanifolds in the Cayley plane}

 Parallel submanifolds of the Cayley plane
$\Cal OP^2$ are
contained either in a totally geodesic
quaternion projective 
plane $HP^2$ as 
parallel submanifolds or in a totally
geodesic
$8$-sphere as parallel submanifolds [Tsukada
1985c].
\smallskip

All parallel submanifolds in $ E^m$ are
 of finite type. Furthermore, if a
compact symmetric space $N$ of rank one is
regarded as a submanifold of a Euclidean
space via its first standard embedding, then
a parallel submanifold of $N$ is also of
finite type [Chen 1996d].

\vfill\eject

\section{Standard immersions and submanifolds with simple geodesics}

\subsection{Standard immersions} 

Let $M=G/K$ be  a compact  irreducible
homogeneous Riemannian manifold. For each
positive eigenvalue $\lambda$ of the Laplacian
on $M$, we denote by $m_\lambda$ the
multiplicity of the eigenvalue
$\lambda$. 
Let $\phi_1,\ldots,\phi_{m_\lambda}$ be an
orthonormal basis of the eigenspace of the 
Laplacian with eigenvalue
$\lambda$. Define a map $f_\lambda:M\to
E^{m_\lambda}$ by
$$f_\lambda(u)={{c_\lambda}\over{m_\lambda^2}}
(\phi_1(u),\ldots,\phi_{m_\lambda} (u)),\leqno(7.1)$$ where
$c_\lambda$ is a positive number. The map $f_\lambda$ defines
an isometric minimal immersion of
$M$ into $S^{m_\lambda-1}_0(1)$ for some suitable constant
$c_\lambda >0$. 

According to a result
of T. Takahashi (1966) each such $f_\lambda$ is
an isometric minimal  immersion of $M$ into a
hypersphere of $E^{m_\lambda}$. 

If $\lambda_i$ is the $i$-th positive
eigenvalue of Laplacian of $M$, then the
immersion
$\psi_i=f_{\lambda_i}$ is called the  i-th
standard immersion of $M=G/K$.

Every full isometric minimal immersion of
a Riemannian $n$-sphere into a hypersphere of
a Euclidean space is a standard immersion if
either $n=2$ or $n\geq 3$ and the order of the
immersion is either $\{1\}, \{2\}$ or $\{3\}$.

Not every full isometric minimal  immersion
of a Riemannian $n$-sphere into a
hypersphere is a standard immersion. For
instance, N. Ejiri (1981) constructed a
full minimal isometric immersion of
$S^3\left({1\over {16}}\right)$ into
$S^6(1)$ of order
$\{6\}$, which is not a standard immersion.
An explicit construction was given by 
K. Mashimo (1985) and by  F. Dillen, L.
Verstraelen and L. Vrancken (1990), who also
showed that the immersion is a 24-fold cover
onto its image. The image of this minimal
immersion in
$S^6(1)$ was  identified in [DeTurck-Ziller
1992] as
$S\sp 3/{T\sp*}$, where ${T\sp*}$ is the binary
tetrahedral group of order 24.

According to a result of Moore (1972), the
minimum number $m$ for which $S^3(c)$ can
admit a non-totally geodesic isometric minimal
immersion into $S^m$ is $6$.  D. DeTurck and
W. Ziller (1992) showed that every non-totally
geodesic  ${SU}(2)$-equiva\-riant minimal
isometric immersion of $S^3({1\over {16}})$ into
$S^6(1)$  is congruent to the
immersion mentioned above. 

\subsection{Submanifolds with planar geodesics} 

A surface in $ E^3$ whose geodesics are
all planar curves is open portion of a
plane or sphere. S. L. Hong (1973) was the first to ask for all
submanifolds of Euclidean space whose
geodesics are plane curves. He showed
that if $f:M\to  E^m$ is an isometric
immersion which is not totally geodesic
and such that for each geodesic
$\gamma$ in $M$, $f\circ \gamma$ is a
plane curve in $ E^m$, then  $f\circ
\gamma$ is a plane circle.

Let $f:M\to R^m(c)$ be an isometric immersion
of a Riemannian manifold into a complete
simply-connected real space form of constant
curvature $c$. If the image of each geodesic
of $M$ is contained in a 2-dimensional
totally geodesic submanifold of $R^m(c)$,
then  $f$ is either a totally
geodesic immersion, a totally umbilical
immersion or a minimal immersion of a compact
symmetric space of rank one by harmonic
functions of degree 2.

 The later case occurs only when $c>0$ and in
this case the immersions are the first standard
embeddings of the real, complex and
quaternionic projective spaces or the
Cayley plane [Hong 1973, Little 1976, Sakamoto
1977].

\subsection{Submanifolds with pointwise planar normal sections} 

Let $M$ be an $n$-dimensional
submanifold of a Euclidean $m$-space
$E^m$. For a point $x$ in $M$ and a
unit vector $X$ tangent to $M$ at
$x$, the vector $X$ and the normal
space $T_x^\perp M$ to $M$ at $x$
determine an $(m-n+1)$-dimensional
affine subspace $E(x,X)$ of $E^m$
through $x$. The intersection of $M$
and $E(x,X)$ gives rise to a curve
$\gamma$ in a neighborhood of $x$
which is called the normal section
of $M$ at $x$ in the direction $X$. 

A submanifold $M$ of $E^m$ is said
to have planar normal sections if
each normal section $\gamma_X(s)$ at
$x$ of $M$ in $E^m$ is a planar curve
where $s$ is an arclength
parametrization of
$\gamma_X$; thus the first three
derivatives
$\{\gamma_X'(s),\gamma{}''_X(s),
\gamma{}'''_X(s)\}$ of
$\gamma(s)$ are linearly dependent as
vectors in $E^m$. Hypersurface of 
Euclidean space and
Euclidean submanifolds with planar
geodesics are examples of submanifolds
with planar normal sections.
Conversely, B. Y. Chen (1983a) proved that
if a surface
$M$ of $E^m$ has planar normal
sections, then either it lies
locally in an affine 3-space $E^3$ of
$E^m$ or it has planar geodesics. If
the later case occurs, $M$ is an open
portion of a Veronese surface in an
affine 5-space $E^5$ of $E^m$. 

A  submanifold $M$ of $E^m$ is said to have
pointwise planar normal sections if, for each
normal section $\gamma$ at $x$, $x\in M$, the
three vectors $\{\gamma'(0),\gamma{}''(0),
\gamma{}'''(0)\}$ at $x$ are
linearly dependent. Clearly, every
hypersurface of $E^{n+1}$ has planar
normal sections, and hence
has pointwise planar normal sections.
A submanifold $M$ of $E^m$ is called
spherical if $M$ is contained in a
hypersphere of $E^m$.

In 1982 B. Y. Chen  proved that a spherical
submanifold of a Euclidean space has 
pointwise planar normal sections if 
and only if it has parallel second
fundamental form. K. Arslan and A.
West (1996) showed that if an
$n$-dimensional submanifold $M$ of a
Euclidean $m$-space $E^m$ has 
pointwise planar normal sections and
does not have parallel second
fundamental form, then locally it must
lies in an affine $(n+1)$-space
$E^{n+1}$ of $E^m$ as a
hypersurface, that is, for each
point $x\in M$, there exists a
neighborhood $U$ of $x$ such that
$U$ is contained in an affine
$(n+1)$-space $E^{n+1}$ of $E^m$.

W. Dal Lago, A. Garc\'ia, and C. U.
S\'anchez (1994) studied the set ${\Cal X}
[M]$ of pointwise planar normal sections on
the natural embedding of a flag manifold $M$
and proved that it is a real algebraic
submanifold of the real projective
$n$-space $RP^n$ with $n=\dim M$. They also
computed the Euler characteristic of ${\chi}
[M]$ and its complexification $\chi_c[M]$ and
showed that the Euler characteristics of
$\chi [M]$ and of $\chi_c[M]$ depend only on
the dimension of $M$ and not on the
nature of $M$ itself. 

\subsection{ Submanifolds with geodesic normal sections and helical immersions} 

A submanifold
$f:M\to  E^m$ is said to have geodesic
normal sections if, for each point $x\in M$
and each unit tangent vector $X$ at $x$, the
image of the geodesic $\gamma_X$ with
$\gamma_X'(0)=X$ is the normal section of $f$
at $x$ in the direction $X$. 

Submanifolds in
Euclidean space with planar geodesics also have
geodesic normal sections. 

 Chen and Verheyen (1981) asked for all
submanifolds of Euclidean space with geodesic
normal sections. They proved that a
submanifold $M$ of $ E^m$ has geodesic normal
sections if and only if all normal sections
of $M$, considered as curves in $ E^m$, have
the same constant first curvature $\kappa_1$; 
also if and only if every curve $\gamma$ of $M$
which is a normal section of $M$ at $\gamma(0)$
in the direction $\gamma'(0)$ remains a normal
section of $M$ at $\gamma(s)$ in the
direction $\gamma'(s)$, for all
$s$ in the domain of $\gamma$. 

In particular, Chen and Verheyen's result
implies that if a compact symmetric space,
isometrically immersed in Euclidean space, has
geodesic normal sections, then it is of rank
one.

 Chen and Verheyen also proved the
following results:

(1) If  $M$ is a   submanifold of $ E^m$
all of whose geodesics are 3-planar, that
is, each geodesic lies in some 3-plane, then
$M$ has geodesic normal sections if and only if
it is isotropic; and 

(2) if $M$ is a   submanifold of Euclidean
space all of whose geodesics are 4-planar,
then $M$ has geodesic normal sections if and
only if it is constant isotropic.

 Recall that a submanifold $M$ of
a Riemannian manifold $\tilde M$ is called
isotropic if, for each point $x\in M$ and each
unit vector $X\in T_xM$, the length
$|h(X,X)|$ of $h(X,X)$ depends only on $x$
and not on the unit vector
$X$. In other words, each geodesic of
$M$ emanating from $x$, considered as a curve
in $\tilde M$, has the same first curvature
$\kappa_1$ at $x$. In particular, if the
length of $h(X,X)$ is also independent of the
point $x\in M$, then $M$ is called
constant isotropic. 

Chen and Verheyen (1981)
showed that  submanifolds in Euclidean space
with geodesic normal sections are constant
isotropic.

Let $M$ be a compact Riemannian manifold. It has a unique
kernel of the heat equation:
$K:M\times M\times \hbox{\bf R}^*_+\to
\hbox{\bf R}.$ If there exists a function
$\Psi:\hbox{\bf R}_+\times
\hbox{\bf R}^*_+\to \hbox{\bf R}$ such that
$K(u,v,t)=\Psi(d(u,v),t)$ for every $u,v\in
M$ and $r\in \hbox{\bf R}^*_+$, then $M$ is
called  strongly harmonic. 

Related to submanifolds with planar geodesics
and to submanifolds with geodesic normal
sections is the notion of helical immersions.
An isometric immersion
$f:M\to E^m$ is called a  helical
immersion if  each geodesic $\gamma$ of $M$
is mapped to a curve with constant Frenet
curvatures, that is, to a
$W$-curve, which are independent of the chosen geodesic. 

A. Besse (1978) constructed helical immersions
of strongly harmonic manifolds into a unit
sphere.
Conversely, K. Sakamoto (1982) proved that if a 
 complete Riemannian manifold  admits a helical
minimal immersion into a hypersphere of $
E^m$, then $M$ is a strongly harmonic manifold.

Y. Hong (1986) proved that  every
helical immersion of a compact homogeneous
Riemannian manifold into  Euclidean space is
spherical. Hong also proved that every helical
immersion of  a compact
rank one symmetric space into a Euclidean space
is a diagonal immersion of some 1-type
standard isometric immersions.

Chen and Verheyen (1984) showed that a
helical submanifold of Euclidean space is a
submanifold with geodesic normal
sections. Conversely, P. Verheyen (1985) proved
that   every submanifold of Euclidean space with
geodesic normal sections is a helical
submanifold.

Helical submanifolds 
were further investigated by K. Sakamoto,
B. Y. Chen, P. Verheyen, Y. Hong, C.-S. Houh,
K. Mashimo, K. Tsukada, H. Nakagawa, and
others. 

\subsection{Submanifolds whose geodesics are generic $W$-curves}

A $W$-curve $\gamma:\hbox{\bf R}\supset
I\to  E^N$ is said to be of rank $r$,
if for all
$t\in I$ the derivatives
$\gamma'(t),\ldots,\gamma^{(r)}(t)$ are linearly
independent and the derivatives 
$\gamma'(t),\ldots,\gamma^{(r+1)}(t)$ are
linearly dependent. Let $\gamma:\hbox{\bf
R}\supset I\to  E^N$ be a $W$-curve of
infinite length, parametrized by arc
length. If the image
$\gamma(I)$ is bounded, then the rank of
$\gamma$ is even, say $r=2k$. There are positive
constants $a_1,\ldots,a_k$, unique up to order,
corresponding positive constants
$r_1,\ldots,r_k$ and orthonormal vectors
$e_1,\ldots,e_{2k}$ in $ E^N$ such that 
$$\gamma(t)=c+\sum_{i=1}^k r_i(e_{2i-1}\sin
a_it+ e_{2i}\cos a_it),$$ where $c$ is a
constant vector. 

The rank of unbounded
$W$-curves is odd and the expression of
$\gamma(t)$ contains an additional term linear
in $t$.  A $W$-curve $\gamma$ is called a
generic
$W$-curve if the $a_i$ are independent over the
rationals, that is, if the closure of
$\gamma(\hbox{\bf R})$ is a standard torus
$S^1(r_1)\times\cdots\times S^1(r_k)$ up to a
motion.

D. Ferus and S. Schirrmacher (1982) proved that
if $f:M\to  E^m$ is an isometric
immersion of a compact Riemannian manifold into
$ E^m$, then $f$ is extrinsic symmetric if
and only if, for almost every  geodesic $\gamma$
in $M$, the image $f(\gamma)$ is a generic
$W$-curve. 

The above result is false if the
condition on $f(\gamma)$ were replaced by the
condition: for each geodesic $\gamma$ in $M$,
the image $f(\gamma)$ is a $W$-curve.

 For a compact Riemannian 2-manifold $M$,
D. Ferus and S. Schirrmacher (1982) proved
that if
$f:M\to E^4$ is an isometric immersion
such that, for every geodesic $\gamma$ in $M$,
the image
$f(\gamma)$ is a $W$-curve, then either one
of the following holds:

(1) if $M$ contains a non-periodic geodesic,
$f$ covers (up to a motion) a standard torus
$S^1(r_1)\times S^2(r_2)\subset  E^4$,
or 

(2) if all geodesics in $M$ are
periodic, $f$ is (up to a motion) an isometry
onto a Euclidean 2-sphere
$S^2(r)\subset  E^3\subset  E^4$.

Y. H. Kim and E. K. Lee (1993) proved the
following: Let $M$ be a complete surface
in $ E^4$. If there is a point $x\in
M$ such that every geodesic through $x$,
considered as a curve in $ E^4$, is a
$W$-curve, then $M$ is an affine 2-space, 
a round sphere  or a
circular cylinder in an affine 3-space, a
product of two plane circles, or a Blaschke
surface at a point $o\in  E^4$,
diffeomorphic to a real projective plane,
and up to a motion, it  is immersed in
$ E^4$ by
$$({1\over\kappa}\sin\kappa s\cos\theta,
{1\over\kappa}\sin\kappa s\sin\theta,
{1\over\kappa}(1-\cos\kappa s)\cos 2\theta,
{1\over\kappa}(1-\cos\kappa s)\sin
2\theta),$$ where $\kappa$ is the Frenet
curvature of geodesic through $o$. The
converse is also true.

\subsection{Symmetric spaces in Euclidean space with simple  geodesics}

Submanifolds in Euclidean space with finite type
geodesics were studied by Chen, Deprez
and Verheyen (1987). They proved the following
results:

(a) An isometric immersion of a compact
symmetric space
$M$ into a Euclidean space is of finite type
if and only if each geodesic of $M$ is mapped
into curves of finite type; 

(b) if $f:S^n\to  E^m$ is an isometric
immersion of a unit
$n$-sphere into $ E^m$, then the immersion
maps all geodesics of $S^n$ into 1-type or
2-type curves if and only if  the immersion is
of finite type with order
$\{1\},\{2\},\{3\},\{1,2\},\{1,3\}$ or $\{2,4\}$;

(c) if $f:FP^n\to
E^m$ is an isometric immersion of a
real, complex, or quaternionic projective
space, or the Cayley plane into
$E^m$, then the immersion maps all
geodesics of
$FP^n$ into 1-type or 2-type curves if and only
if the immersion   is of finite type with order
$\{1\},\{2\}$ or $\{1,2\}$;

(d) a finite type isometric immersion of a
unit $n$-sphere in
$ E^m$ of order $\{1,2\}$ is a diagonal
immersion of the first and the second standard
immersions of the $n$-sphere; 

(e)  a finite type isometric immersion of a
unit $n$-sphere in
$ E^m$ of order $\{1,3\}$ for which all
geodesics are mapped to $W$-curves is a
diagonal immersion of the first and the third
standard immersions of the
$n$-sphere; and  

(f) there exist
finite type isometric immersions of the unit
2-sphere of order $\{1,3\}$ or of order
$\{2,4\}$ which are not diagonal immersions.

Chen, Deprez and Verheyen (1987) also
studied an isometric immersion
$f:M\to E^m$ which satisfies the condition:
there is a point
$x_0\in M$ such that every geodesic through
$x_0$ is mapped to a circle. 

They proved the following results: 

(g)  Let
$f:S^n\to E^m$ be an isometric embedding.
If there exists a point
$x_0\in S^n$ such that $f$ maps all
geodesics of
$S^n$ through $x_0$ to 1-type curves, then the
embedding is  the first standard
embedding of $S^n$ into a totally geodesic
$ E^{n+1}$; 

(h) Let $f:FP^n\to
 E^m$ be  a finite type isometric
immersion. If there exists a point $x_0\in
FP^n$ such that the immersion maps all
geodesics of
$FP^n$ through
$x_0$ to 1-type curves, then the immersion is
 the first standard embedding of
$FP^n$; and 

(i) up to motions, the set
of isometric immersions of a projective
$n$-space $FP^n\,(F$ $=R, C$ or
$ H)$ into Euclidean $m$-space  which
map all geodesics through a fixed point
$x_o\in FP^n$ to circles is in one-to-one
correspondence with the set of isometric
immersions of
$FP^{n-1}$ into $S^{m-dn-1}$, where $d$ is $1,\,2$ or $4$,
according to the field $F$ is real, complex or quaternion.

\vfill\eject

\section{Hypersurfaces of real space forms}
Complete simply-connected Riemannian
$n$-manifolds of constant curvature are frame
homogeneous, that is, for any pair of points
$x$ and $y$ and any orthonormal frames $u$ at
$x$ and $v$ at $y$ there is an isometry $\phi$
such that $\phi(x)=y$ and $\phi_*$ maps $u$
onto $v$. Such Riemannian manifolds are
Euclidean $n$-space $ E^n$, 
Riemannian $n$-spheres and  real hyperbolic
$n$-spaces. 

Consider an isometrically immersed
orientable hypersurface $M$ in a complete
simply-connected real space form $R^{n+1}(c)$
of constant curvature $c$ with a unit normal
vector field $\xi$. We simply denote the shape
operator $A_\xi$ at $\xi$ by $A$.

 Let $$\kappa_1\leq \kappa_2\leq\cdots\leq
\kappa_n$$ denote the $n$ eigenvalues of $A$ at
each point $x$ of $M$. Then each $\kappa_i,
\,(1\leq i\leq n)$, is a continuous function on
$M$, and is called a principal curvature of
$M$.

 For each $x\in M$ and each $\kappa\in
\{\kappa_1,\ldots,\kappa_n\}$, we define a
subspace $\Cal D_x(\kappa)$ of $T_xM$ by $$\Cal
D_x(\kappa)=\{X\in T_xM:
AX=\kappa(x)X\}.\leqno(10.1)$$ Let
$\Cal D(\kappa)$ assign each point $x\in M$
the subspace $\Cal D_x(\kappa)$. 

The following basic results are well-known: If $\dim \Cal D_x(\kappa)$ is constant on $M$, say $m$,
then 

(1)  $\kappa$ is a
differentiable function on $M$;

(2) $\Cal D(\kappa)$ is a
differentiable distribution on $M$;

(3)  $\Cal D(\kappa)$ is completely
integrable, called the principal foliation;

(4) If $m\geq 2$, then $\kappa$ is
constant along each leaf of  $\Cal
D(\kappa)$;

(5) If $\kappa$ is constant along a leaf $L$
of  $\Cal D(\kappa)$, then $L$ is locally an
$m$-sphere of $R^{n+1}(c)$, where an
$m$-sphere means a hypersphere of an
$(m+1)$-dimensional totally geodesic
submanifold of $R^{n+1}(c)$. 

Suppose that a continuous principal
curvature $\kappa$ has constant multiplicity
$m$ on an open subset $U \subset M$. Then
$\kappa$ and its principal foliation
$\Cal D(\kappa)$ are smooth on
$U$. The leaves of this principal
foliation are the curvature surfaces 
corresponding to $\kappa$ on $U$.  The
principal curvature
$\kappa$ is constant along each of its
curvature surfaces in $U$ if and only if
these curvature surfaces are open subsets of
$m$-dimensional Euclidean spheres or
planes.  The focal map $f_{\kappa}$ corresponding to $\kappa$ is the map which
maps
$x\in M$ to the focal point $f_{\kappa}(x)$ corresponding to $\kappa$, that is,
$$f_{\mu}(x) = f(x) + \frac{1}{\mu(x)}
\xi(x),$$
where $\xi$ is a unit normal vector field.
The principal curvature $\kappa$ is constant
along each of its curvature surfaces in $U$
if and only if the focal map $f_{\mu}$
factors through an immersion of  the
$(n-1-m)$-dimensional space of 
leaves $M/\Cal D(\kappa)$ into $E^n$.

\subsection{Einstein 
hypersurfaces }

A Riemannian manifold is said to
be Einstein if the Ricci tensor is a constant
multiple of the metric tensor, that is
$Ric=\rho g$, where $\rho$ is a constant. 

A. Fialkow (1938) classified Einstein
hypersurfaces in real space forms.

Let $M^n,\,n>2$, be an Einstein hypersurface in
$R^{n+1}(\bar c)$. Then 

(a) if $\rho>(n-1)\bar c$, then $M$ is totally umbilical, hence, $M$  is also a real space form;

(b) if $\rho=(n-1)\bar c$, then the type number, $t(x)=$ rank$\,A_x$, is $\leq 1$
for all $x\in M$ and $M$ is of constant curvature $\bar c$; and

(c) if $\rho<(n-1)\bar c$,
then $\bar c>0$, $\rho=(n-2)\bar c$, and $M$
is locally a standard product embedding of
$S^p\left(\left({{n-2}\over {p-1}}\right)\bar c\right)\times S^{n-p}\left(\left({{n-2}\over {n-p-1}}\right)\bar c\right)$, where $1<p<n-1$. 

In particular,  complete
Einstein hypersurfaces of $ E^{n+1}$ are
hyperspheres, hypercylinders over complete
plane curves and hyperplanes; and complete
Einstein hypersurfaces in $S^{n+1}$ are the
small hyperspheres, the great hyperspheres and
certain standard product embedding of products
of spheres. 

P. Ryan (1971) studied hypersurfaces of real
space forms with parallel Ricci tensor and
proved that if a hypersurface $M$ of
dimension $>2$ in a real space form of
constant curvature $c$ is not of constant
curvature $c$ and if it has parallel Ricci
tensor, then $M$ has at most two distinct
principal curvatures. Furthermore, if $c\ne
0$, then both principal curvatures are
constant.

\subsection{ Homogeneous hypersurfaces}

Let $M$ be a homogeneous Riemannian
$n$-manifold isometrically immersed into an
$(n+1)$-dimensional complete simply-connected
real space form $R^{n+1}(c)$. Then

 (1) if $c=0$, then
$M$ is isometric to the hypercylinder $S^k\times
 E^{n-k}$ [Kobayashi
1958, Nagano-Takahashi 1960];

 (2) if $c>0$, then $M$ is represented as an
orbit of a linear isotropy group of a
Riemannian symmetric space of rank 2; and
$M$ is isometric to $ E^2$ or else is
given as an orbit of a subgroup of the
isometry group
 of $R^{n+1}(c)$ [Hsiang-Lawson 1971]; and 

(3) if $c<0$, then $M$ is isometric to  a
standard product embedding of $ E^n$,
$S^k\times H^{n-k}$,  or a 3-dimension group
manifold 
$$B=\begin{pmatrix} e^t & 0 & x\\
0&e^{-t} &y\\ 0&0&1\end{pmatrix}:
x,y,t\in  \hbox{\bf R}$$
with the metric
$ds^2=e^{-2t}dx^2+e^{2t}dy^2+dt^2$ [Takahashi
1971]. 

Each of the hypersurfaces above except
$ E^2$ in (2) and $B$ in (3) is given as
an orbit of a certain subgroup of the isometry
group of $R^{n+1}(c)$.

\subsection{Isoparametric
hypersurfaces} 

The history of isoparametric hypersurfaces
can be traced back to 1918 of the work of
E. Laura (1918) and of C. Somigliana (1918)
on geometric optics.  T. Levi-Civita
(1873--1941) and B. Segre  (1903--1977)
studied such hypersurfaces of Euclidean
space during the period of 1924--1938. A
major progress on isoparametric
hypersurfaces was made by \'E. Cartan during
the period of 1938--1940.

A hypersurface $M$ of a Riemannian
manifold $\tilde M$ is called isoparametric if
$M$ is locally defined as the level set of a
function $\phi$ on (an open set of) $\tilde M$
with property: $d\phi\wedge d||d\phi||^2=0$
and $d\phi\wedge d(\Delta\phi)=0.$

A hypersurface $M$ of a complete
simply-connected Riemannian manifold
$R^{n+1}(c)$ of constant curvature $c$ is
isoparametric if and only if $M$ has
constant principal curvatures. Each
isoparametric hypersurface of $R^{n+1}(c)$
determines a unique complete embedded
isoparametric hypersurface in $R^{n+1}(c)$.
 Thus, every
open piece of an isoparametric hypersurface
can be extended to a unique complete
isoparametric hypersurface. 

\vskip.1in

\noindent{\bf 10.3.1. Isoparametric family of
hypersurfaces}

 Let $f:M\to R^{n+1}(c)$ be a
hypersurface and $\xi$  is a unit normal
vector field of $f$. For each
$t>0$, let $f_t(x),\,x\in M$, be the point of
$R^{n+1}(c)$ on the geodesic from
$f(x)$ starting in the direction $\xi$ at $x$
which has geodesic distance $t$ from $f(x)$.
In the Euclidean case $(c=0)$, we have
$f_t(x)=f(x)+t\xi_x$. In the spherical case
$(c=1)$, we have $f_t(x)=(\cos t)f(x)+(\sin
t)\xi_x$ by considering $S^{n+1}$ as the unit
sphere in $E^{n+2}$. In all cases, $f_t$
is an immersion of $M$ in $R^{n+1}(c)$ for
sufficiently small $t$.

In $R^{n+1}(c)$ an isoparametric family of
hypersurfaces is a family of hypersurfaces
$f_t:M\to R^{n+1}(c)$ obtained from a
hypersurface $f:M\to R^{n+1}(c)$ with constant
principal curvatures.

\vskip.1in

\noindent{\bf 10.3.2. Cartan's basic identity}

The starting point of Cartan's work on
isoparametric hypersurfaces is the
following basic identity concerning all the
distinct constant principal curvatures
$a_1,\ldots,a_g$ with their respective
multiplicities $\nu_1,\ldots,\nu_g$.

For $g\geq2$, Cartan's basic identity is the
following: For each $i,\,1\leq i\leq g$, we
have $$\sum_{j\ne i}
\nu_j{{c+a_ia_j}\over{a_i-a_j}}=0,\quad
1\leq j\leq g.\leqno(10.2)$$
Cartan's basic identity holds for all $c$,
positive, negative and zero.

K. Nomizu (1975)
observed that the focal set of an
isoparametric hypersurface in a real space
form admits a submanifold structure locally,
and that the mean curvature vector of this
submanifold at $f_t(x)$ can always be
expressed as a multiple of $\xi_x$. Moreover,
if $f_t(x)$ is a focal point corresponding 
to a principal curvature $\kappa$, then
the multiple is the left hand side of
(10.2). Thus, Cartan's identity is
equivalent to the minimality of the focal
varieties.

\vskip.1in
\noindent{\bf 10.3.3. Isoparametric
hypersurfaces in Euclidean space}

In the Euclidean case $(c=0)$, Cartan's
identity implies $g\leq 2$. If $g=2$, one of
the principal curvatures must be 0.

Isoparametric hypersurfaces of a Euclidean
$(n+1)$-space $ E^{n+1}$ are locally
hyperspheres,  hyperplanes or a standard
product embedding of  $S^k\times E^{n-k}$.

This result was proved in 1937 by T.
Levi-Civita for $n=2$ and in 1938 by B. Segre
for arbitrary $n$.

\noindent{\bf 10.3.4. Isoparametric
hypersurfaces in hyperbolic space}

Cartan's basic identity also yields $g\leq 2$
in the case $c<0$. In the hyperbolic space
$H^{n+1}(-1)$, \'E. Cartan proved that the
number of the distinct principal curvatures
of an isoparametric hypersurface is one or
two.  An isoparametric hypersurface in
hyperbolic space is totally umbilical or
locally a standard  product embedding of
$S^k\times H^{n-k}$.

\noindent{\bf 10.3.5. Isoparametric
hypersurfaces in spheres}

Cartan's basic identity does not restrict
$g$ if $c>0$. In fact, let $\theta$ be any
number such that $\sin\theta\ne 0$ and let $g$
be any positive integer. Set
$$\lambda_k=\cot\left(\theta+{{k-1}\over
g}\pi\right)$$ for $k=1,\ldots,g$. Then such
a collection of
$\lambda_k$ with equal multiplicities
satisfies (10.2).

One major part
of Cartan's work on isoparametric
hypersurfaces is to give an algebraic method
of finding an isoparametric family of
hypersurfaces in $S^{n+1}$ with $g$ distinct
principal curvatures with the same
multiplicity. His result is that
such a family is defined by $$M_t=\{x\in
S^{n+1}:\phi (x)=\cos gt\},$$ where $\phi$ is
a harmonic homogeneous polynomial of degree
$g$ on $E^{n+2}$ satisfying
$$||d\phi||^2=g^2(x_1^2+\cdots+x^2_{n+2})^{g-1},
$$ where $(x_1,\cdots,x_{n+2})$ are the
Euclidean coordinates on $E^{n+2}$.

In $S^{n+1}$ an isoparametric hypersurface
with $g=2$ is locally a standard product
embedding of the product of two spheres with
appropriate radii.

Cartan (1939a) proved that if $g=3$,
then $n=3,6,12,$ or 24; and the
multiplicities are the same in
each case. Moreover, the isoparametric
hypersurface must be a tube of constant
radius over a standard Veronese embedding of
a projective plane
$FP^2$ into $S^{3m+1}$, where $F$ is the
division algebra of reals, complex
numbers, quaternions, or Cayley
numbers for the (common)
multiplicity $m=1,2,4,8$, respectively. 
Thus, up to congruence, there is only one
such family for each value of $m$.  These
isoparametric hypersurfaces with three
principal curvatures are known  as  Cartan
hypersurfaces. If $g=4$ and if the
multiplicities are equal, then $n=4$ or 8.

For each of the above cases, Cartan gave an
example of a hypersurface with constant
principal curvatures. All the compact isoparametric hypersurfaces in
$S^{n+1}$ given by Cartan are homogeneous;
indeed, each of them is the orbit of a certain
point by an appropriate closed subgroup of the
isometry group of $S^{n+1}$.

H. F.  M\"unzner(1980,1981) showed that a
parallel family of isoparametric
hypersurfaces in $S^n$ always consists of the
level sets in $S^n$ of a homogeneous
polynomial defined on $E^{n+1}$.

R. Takagi and T. Takahashi (1972)
determined all the orbit hypersurfaces in
$S^{n+1}$. They showed that the number of
distinct principal curvatures of a
homogeneous isoparametric hypersurfaces in
a sphere is 1, 2, 3, 4, or 6, and they
listed the possible multiplicities. The
list of homogeneous isoparametric
hypersurfaces of Takagi and Takahashi
contains  five different classes of orbit
hypersurfaces in $S^{n+1}$ each with four
distinct principal curvatures not of the
same multiplicity. It also contains an
example with $g=6$ and $m_1=m_2=1$, and
one with
$g=6$ and $m_1=m_2=2$. H. Ozeki and M.
Takeuchi (1975,1976) produced two infinite
series of isoparametric families which are
not on the list of Takagi and Takahashi
and are, therefore, inhomogeneous. These
examples all have four distinct principal
curvatures.

H. F. M\"unzner (1980,1981) proved that
isoparametric hypersurfaces in $S^{n+1}$ with
$g$ distinct principal curvatures exist only
when $g=1, 2 , 3, 4,$ or 6.
Through a geometric study of the focal
submanifolds and their second fundamental
forms, M\"unzner showed that if $M$ is an
isoparametric hypersurface with principal
curvatures $\cot\theta_i,\, 0<\theta_1<\cdots
<\theta_g<\pi$, with multiplicities $m_i$,
then $$\theta_k=\theta_1+{{k-1}\over
g}\pi,\quad 1\leq k\leq g,$$ and the
multiplicities satisfy $$m_i=m_{i+2}\quad
\hbox{(subscripts mod}\;g).$$  As a
consequence, if $g$ is odd, then all of
the multiplicities must be equal; and if
$g$ is even, then $m_1=m_3=\cdots=m_{g-1}$
and $m_2=m_4=\cdots=m_g$.

Moreover, H. F. M\"unzner (1981) proved
that if a hypersurface $M$ in  $S^{n+1}$
splits $S^{n+1}$ into two disks bundles
$D_1$ and $D_2$ over compact manifolds with
fibers of dimensions $m_1+1$ and $m_2+1$
respectively, then $\dim H^*(M;\hbox{\bf
Z}_2)=2h$, where $h$ is $1, 2, 3,4$ or 6.
In the case that $M$ is an isoparametric
hypersurface, $M$ splits $S^{n+1}$ into
two disk bundles such that the numbers
$m_1$ and $m_2$  coincide with the
multiplicities of the principal
curvatures of $M$.

Although isoparametric hypersurfaces with four
principal curvatures have not been
completely classified, there is a large
class of examples due to D. Ferus, H.
Karcher and M\"unzner (1981). In fact, 
Ferus,  Karcher and  M\"unzner
constructed isoparametric
hypersurfaces with four distinct principal
curvatures using representations of
Clifford algebras, which include all known
examples, except two. They are able to
show geometrically that many of their
examples are not homogeneous.  The isoparametric
hypersurfaces which belong to  the
Clifford series discovered by Ferus,
Karcher and M\"unzner are the regular
level sets of an isoparametric function on
$S^{2k-1}$ determined by an orthogonal
representation of the Clifford algebra
$C_{m-1}$ on $E^k$. 

J. Dorfmeister and E. Neher (1983)
extended the work of Ferus,  Karcher and 
M\"unzner to a more algebraic setting
involving triple systems. The topology of
isoparametric hypersurfaces of the
Clifford examples have been studied by Q.
M. Wang (1988). Among others he showed
that there exist noncongruent hypersurfaces
in two different isoparametric families
which are diffeomorphic. 

S. Stolz (1997) proved that the only
possible triples $(g,m_1,m_2)$ with $g=4$
are exactly those that appear either in
the homogeneous examples or those appear
in the Clifford examples of Ferus, Karcher
and M\"unzner.

 M\"unzner (1981) proved that $m_1=m_2$
always hold in the case $g=6$. U. Abresch
(1983) showed that the common multiplicity
of $m_1$ and $m_2$ must be 1 or 2.
M\"unzner's computation also yields that
the focal submanifolds must be minimal.

J. Dorfmeister and E. Neher (1983) proved
that every compact isoparametric
hypersurface in the sphere with
$g=6$ and $m_1=m_2=1$ must be homogeneous.
In particular, an isoparametric hypersurface
in $S^7$ with $g=6$ is an orbit of the
isotropy action of the symmetric space
$G_2/SO(4)$, where $G_2$ is the automorphism
group of the Cayley algebra.  A geometric
study of such isoparametric hypersurfaces in
$S^7$ was made in [Miyaoka 1993]. In fact, 
Miyaoka  has shown that a
homogeneous isoparametric hypersurface $M$
in $S^7$ can be obtained as the
inverse image under the Hopf fibration
$\psi:S^7 \to S^4$ of an isoparametric
hypersurface with three principal curvatures
of multiplicity one in $S^4$.  She also
showed that the two focal submanifolds of
$M$ are not congruent, even though they are
lifts under $\psi^{-1}$ of congruent Veronese
surfaces in $S^4$. Thus, these focal
submanifolds are two  noncongruent minimal
taut homogeneous embeddings of $RP^2 \times
S^3$ in $S^7$.   Fang (1995) studied the
topology of isoparametric hypersurfaces with
six principal curvatures. Peng and Hou (1989)
gave explicit forms for the isoparametric
polynomials of degree six for the
homogeneous isoparametric hypersurfaces with
$g=6$.  Recently, R. Miyaoka (1998) proves that
all isoparametric hypersurfaces in the sphere
with $g=6$ are homogeneous.

A smooth immersion $f:M\to E^m$ from a
compact manifold $M$ into a Euclidean
$m$-space is called taut if every
nondegenerate Euclidean squared distance
function has the minimum number of critical
points. T. E. Cecil and P. Ryan (1979a)
showed that all isoparametric hypersurfaces
and their focal submanifolds in the spheres
are taut, and  every isoparametric
hypersurface in a sphere is totally focal,
that is, every squared distance function of
the hypersurfaces is either nondegenerate or
has only degenerate critical points. 
S. Carter and A. West (1982) proved  that
every totally focal hypersurface in the
sphere is isoparametric. 

R. Miyaoka (1982) proved that if a complete
hypersurface in $S^{n+1}$ has constant mean
curvature and three non-simple principal
curvatures, then it is isoparametric.   B.
Y. Chen (1984b) proved that an
isoparametric hypersurface in the sphere
is either of 1-type or of 2-type. S. Chang
(1993) showed that a compact  hypersurface
in $S^4$ with constant scalar curvature
and constant mean curvature is
isoparametric.

M. Kotani (1985) proved that if $M$ is an
$n$-dimensional compact homogeneous minimal
hypersurface in a unit sphere with $g$
distinct principal curvatures, then the
first nonzero eigenvalue $\lambda_1$ of the
Laplacian on $M$ is $n$ unless $g=4$.
 B. Solomon (1990a,1990b) determined
 the spectrum of the Laplacian operator on
isoparametric minimal hypersurfaces of
spheres with $g=3$. All such hypersurfaces
are algebraic and homogeneous. Solomon
(1992) also studied the spectrum of the
Laplacian on quartic isoparametric
hypersurfaces in the unit sphere. These are
hypersurfaces with $g=4$. J. H. Eschenburg
and V. Schroeder (1991) investigated the
behavior of the Tits metric on isoparametric
hypersurfaces.  B. Wu (1994) showed that
for each $n$ there are only finitely many
diffeomorphism classes of compact
isoparametric hypersurfaces of $S^{n+1}$
with four distinct principal curvatures.
J. K. Peng and Z. Z. Tang (1996) derived an
explicit formula for the Brouwer degree of
the gradient map of an isoparametric
function $f$ in terms of the multiplicities
of the principal curvatures of the
isoparametric hypersurface defined by $f$.
They applied this formula to determine the
Brouwer degree in a variety of examples. 

\subsection{Dupin hypersurfaces}.

\vskip.1in

\noindent{\bf 10.4.1. Cyclide of Dupin}

C. Dupin (1784--1873) defined in 1922 a
cyclide to be a surface $M$ in $ E^3$ which
is the envelope of the family of spheres
tangent to three fixed spheres. This was
shown to be equivalent to requiring that
both sheets of the focal set degenerated
into curves. The  cyclides are equivalently
characterized by requiring that the lines of
curvatures in both families be arcs of
circles or straight lines. Thus, one can
obtain three obvious examples: a torus of
revolution, a circular cylinder and a
circular cone. It turns out that all
cyclides can be obtained from these three by
inversions in a sphere in $E^3$.

The cyclides were studied by many prominent
mathematicians of the last century
including A. Cayley (1821--1895), J. G.
Darboux (1842--1917), F. Klein (1849--1925),
J. Liouville (1809--1882) and J. C. Maxwell
(1931--1879).  A detailed treatment of the
cyclides can be found in the book
[Fladt-Baur 1975]. A recent survey on Dupin
hypersurfaces was given in [Cecil 1997].

\vskip.1in
\noindent{\bf 10.4.2. Proper Dupin
hypersurfaces}

The study of Dupin hypersurfaces was
initiated by T. E. Cecil and P. J. Ryan in
1978.  Let
$M$ be a hypersurface in a complete
simply-connected real space form
$R^{n+1}(c)$. A submanifold $S$ of $M$ is
called a curvature surface if, at each point
$x\in S$, the tangent space $T_xS$ is  a
principal space, that is, it is an eigenspace
of the shape operator.

If a principal curvature $\kappa$ has constant
multiplicity $m$ on some open set $U\subset
M$, then the corresponding distribution of
principal spaces is a foliation of rank $m$,
and the leaves of this principal foliation are
curvature surfaces. Furthermore, if the
multiplicity $m$ of $\kappa$ is greater than
one, then $\kappa$ is constant along each of
these curvature surfaces.

 A hypersurface $M$ is called
a Dupin hypersurface if, along each  curvature
surface, the corresponding principal curvature
is constant. A Dupin hypersurface is called
proper if each principal curvature has constant
multiplicity on $M$, that is, the number of
distinct principal curvatures is constant.

 In $ E^3$ the only proper Dupin
hypersurfaces are spheres, planes, and the
cyclides of Dupin.
There exist many examples of Dupin
hypersurfaces which are not proper, for
instance, a tube $M$ of sufficiently small
constant radius $r$ in $ E^4$ over a
torus of revolution $T^2\subset
 E^3\subset E^4$, since there are only two
distinct principal curvatures on the set
$T^2 \times \pm \{r\}$ but three
distinct principal curvatures elsewhere on
$M$.

G. Thorbergsson (1983a) proved that
the possible number $g$ of distinct principal
curvatures of a compact embedded proper Dupin
hypersurface in $E^{n+1}$ is 1, 2, 3, 4 or 6;
the same as for an isoparametric
hypersurface. Thorbergsson's result implies
that compact Dupin hypersurfaces in a
sphere satisfy the same periodicity
$m_i=m_{i+2}$ (subscripts mod $g$) for the
multiplicities of the principal curvatures,
just like isoparametric hypersurfaces.

S. Stolz (1997) proved that the possible
multiplicities $m_1,m_2$ of compact proper
Dupin hypersurfaces are exactly the same as
in the isoparametric case.
\vskip.1in

\noindent{\bf 10.4.3. Local construction
of Dupin hypersurfaces}

U. Pinkall (1985a) gave four local
constructions for obtaining a proper Dupin
hypersurface with $g+1$ distinct principal
curvatures from a lower dimensional proper
Dupin hypersurface with $g$ distinct principal
curvatures. Using these Pinkall proved that
there exists a proper Dupin hypersurface in
Euclidean space with an arbitrary number of
distinct principal curvatures with any given
multiplicities. 

Pinkall's construction is done by using the
following basic constructions:

Start with a Dupin hypersurface $W^{n-1}$ 
in $E^n$ and then consider $E^n$ as the
linear subspace $E^n\times\{0\}$ in
$E^{n+1}$. Then the following constructions
yield a Dupin hypersurface $M$ in $E^{n+1}$:

(1) Let $M$ be the cylinder
$W^{n-1}\times E^1$ in $E^{n+1}$;

(2) Let $M$ be the hypersurface in $E^{n+1}$
obtained by rotating $W^{n-1}$ around an
axis $E^{n-1}\subset E^n$;

(3) Project $W^{n-1}$ stereographically onto
a hypersurface $V^{n-1}\subset S^n\subset
E^{n+1}$. Let $M$ be the cone over $V^{n-1}$
in $E^{n+1}$;

(4) Let $M$ be a tube in $E^{n+1}$ around
$W^{n-1}$.

These constructions give rise to a new
principal curvature of multiplicity one
which is constant along its lines of
curvature. The other principal curvatures
are determined by the principal curvatures
of $W^{n-1}$, and the Dupin property is
preserved for these principal curvatures.
These construction can be easily be
generalized to produce a new principal
curvature of multiplicity $m$ by considering
$E^n$ as a subset of $E^n\times E^m$ rather
than $E^n\times E^1$. These constructions
only yield a compact proper Dupin
hypersurface if the original manifold
$W^{n-1}$ is itself a sphere. Otherwise, the
number of distinct principal curvatures is
not constant on a compact manifold obtained
in this way.

A Dupin hypersurface which is
obtained as the result of one of the four
constructions is said to be reducible.
A proper Dupin hypersurface which does not
contain any reducible open subset is called
locally irreducible.

\vskip.1in

\noindent{\bf 10.4.4. Dupin hypersurfaces with 2, 3 or 4
distinct principal curvatures}

Let $M$ be a complete proper Dupin
hypersurface in $ E^{n+1}$ with two
distinct principal curvatures. If one of the
principal curvatures is identically zero,
$M$ is a standard product embedding of
$S^k(r)\times  E^{n-k}$, where
$S^k(r)$ is a round sphere in a Euclidean
subspace
$ E^{k+1}$ orthogonal to $E^{n-k}$.
Otherwise, T. E. Cecil and P. Ryan (1985)
proved that a compact cyclide $M$ of
characteristic $(k,n-k)$ embedded in
$S^{n+1}$ must be M\"obius equivalent to 
 a standard product embedding of
two spheres $S^k(r)\times S^{n-k}(s)\subset
S^{n+1}(t),\,r^2+s^2=t^2$. The proof of
Cecil and Ryan used the compactness
assumption in an essential way, whereas the
classification of Dupin surfaces in $E^3$
obtained in the nineteenth century does not
need such an assumption. Using Lie sphere
geometry of S. Lie (1842--1899), Pinkall
(1985a) proved that every cyclide of Dupin is
contained in a unique compact connected
cyclide, and any two cyclides of the same
characteristic are locally Lie equivalent.

If the pole of the
projection does not lie on
$S^k(r)\times S^{n-k}(s)$, $M$ is called a ring
cyclide. Otherwise, $M$ is noncompact and is
called a parabolic ring cyclide. In both
cases, the two sheets of the focal set are a
pair of focal conics.

Complete proper Dupin
hypersurfaces embedded in $ E^{n+1}$ with
two distinct principal curvatures were
completely  classified by Cecil and Ryan
(1985). They proved that if such a
hypersurface is compact, it is a ring cyclide;
if such a hypersurface is noncompact and one
of the two distinct principal curvatures is
zero identically, then it is a standard product
embedding of $S^k\times E^{n-k}$ or a
parabolic ring cyclide.

 Pinkall (1985b) gave local classification
of Dupin hypersurfaces with three distinct
principal curvatures in $E^4$ up to
Lie-equivalence, by using the method of
moving frames.
Niebergall (1991) proved that every proper
Dupin hypersurface in $E^5$ with three
distinct principal curvatures at each
point is reducible. 

T. E. Cecil and G. R. Jensen (1998) proved
that if a proper Dupin hypersurface in $E^n$
with three distinct principal curvatures
does not contain a reducible open subset,
then it is equivalent by a Lie sphere
transformation to an isoparametric
hypersurface in a sphere $S^n$.

A necessary condition on a Dupin
hypersurface with at least 4 distinct
principal curvature to be Lie equivalent
to a piece of an isoparametric
hypersurface is the constancy of the Lie
curvatures. Niebergall (1992)  gave a
local classification of proper Dupin
hypersurfaces in $E^5$ with four
distinct principal curvatures which are
Lie equivalent to isoparametric
hypersurfaces. He showed that if the Lie
curvature of the hypersurface $M$ is
constant and a condition on certain
half-invariants are satisfied, then $M$ is
Lie equivalent to an isoparametric
hypersurface. Cecil and Jensen (1997) proved
that the condition on the half-invariants
can be removed.  In 1989, R.  Miyaoka  gave
necessary and sufficient conditions for a
compact embedded Dupin hypersurface with
four or six principal curvatures to be Lie
equivalent to an isoparametric
hypersurface. She showed that a compact
proper Dupin hypersurface embedded in
$E^{n+1}$ is Lie equivalent to an
isoparametric hypersurface if it has
constant Lie curvatures and it satisfies
certain global conditions regarding the
intersections of leaves of its various
principal foliations. 

\vskip.1in

\noindent{\bf 10.4.5. Dupin hypersurfaces and Lie
sphere transformations}

The classes of Dupin and proper
Dupin hypersurfaces in
$S^{n+1}$ are invariant under conformal
transformations of $S^{n+1}$ and under
stereographic projection from $S^{n+1}$ to
$E^{n+1}$ (cf. [Cecil-Ryan 1985]). U. Pinkall
(1985a) proved that they are invariant
under parallel transformations; and thus
under the group of Lie sphere
transformations. Hence, Dupin hypersurfaces
are most naturally studied in the Lie
geometric framework.

A Lie sphere transformation is a projective
transformation of $RP^{n+2}$ which takes the
Lie hyperquadric $Q^{n+1}$ into itself, where
$Q^{n+1}$ in $RP^{n+2}$ is defined by
$$-x_0^2+x_1^2+\cdots+x_{n+1}^2-x_{n+2}^2=0$$
in terms of homogeneous coordinates.

 In terms of the geometry of $E^n$ a Lie
sphere transformation preserves the family of
oriented spheres and planes. The group of
Lie sphere transformations is isomorphic to
$O(n+1,2)/\{\pm I\}$, where $O(n+1,2)$ is
the group of orthogonal transformations of
$E^{n+3}_2$. The group of Lie sphere
transformations contains the group of
conformal transformations of
$S^{n+1}$ as a proper subgroup.

A proper Dupin hypersurface with one principal
 curvature at each point is, of course, totally
umbilical and must, therefore, be an open
subset of a great or small sphere.

Cecil and Ryan (1978) showed
that a compact proper Dupin hypersurface in
$S^{n+1}$ with two distinct principal
curvatures must be a ring cyclide, that is, the
image under a conformal transformation of $S\sp
{n+1}$ of a standard product of two spheres.
Pinkall (1985) was able to drop the
assumption of compactness and showed that any
proper Dupin hypersurface in $S\sp
{n+1}$ with two distinct principal curvatures
is the image under a Lie sphere
transformation of an open subset of a
standard product embedding of the product of
two spheres.

It was conjectured by T. E. Cecil and P.
Ryan that every compact proper Dupin
hypersurface in $S^{n+1}$ is Lie
equivalent to an isoparametric
hypersurface. As noted above, this is true
for $g$ equal to 1 or 2.

R. Miyaoka (1984a) proved
that a compact embedded proper Dupin
hypersurface in a real space form with three
principal curvatures is Lie equivalent to an
isoparametric hypersurface.

Cecil and Ryan's conjecture has  been
shown to be false  by
U. Pinkall and G. Thorbergsson (1989a) and
also independently by R. Miyaoka and T.
Ozawa around 1988 who have produced different
counterexamples to the conjecture.

For a proper Dupin hypersurface with $g=4$,
one can order the principal curvatures so that
$\nu_1<\nu_2<\nu_3<\nu_4$, the Lie curvature
$\Psi$ is defined to be the cross ratio:
$$\Psi={{(\nu_4-\nu_3)(\nu_1-\nu_2)}\over
{(\nu_4-\nu_2)(\nu_1-\nu_3)}}.\leqno10.3$$

The examples of Miyaoka and
Ozawa involve the Hopf fibration of $S^7$ over
$S^4$. Let $ E^8= H\times H$,
where $ H$ is the division ring of
quaternions. The Hopf fibering of the unit
sphere $S^7$ in $ E^8$ over the unit
sphere $S^4$ in $ E^5= H\times E^1$
is given by
$$\psi(u,v)=(2u\bar v,|u|^2-|v|^2),\quad
u,v\in  H.\leqno(10.4)$$

Miyaoka and Ozawa showed that if $M^3$ is a
compact proper Dupin hypersurface embedded in
$S^4$, then $\psi^{-1}(M^3)$ is  a proper
Dupin hypersurface embedded in $S^7$.
Furthermore, if $M^3$ has $g$ distinct
principal curvatures, then $\psi^{-1}(M^3)$
has $2g$ distinct principal curvatures. If
$M^3$ is not isoparametric, then the Lie
curvature of $\psi^{-1}(M^3)$ is not constant,
and so $\psi^{-1}(M^3)$ is not Lie
equivalent to an isoparametric hypersurface.

The examples of  Pinkall and Thorbergsson
are given as follows: 

Let  $E^{2n+2} = E^{n+1} \times E^{n+1}$ and
let $S^{2n+1}$ denote the unit sphere in
$E^{2n+2}$. The Stiefel manifold $V$
of orthogonal 2-frames in $E^{n+1}$ of
length $1/ \sqrt{2}$ is given by 
$$V = \{ (u,v) \in E^{2n+2} | u \cdot
v = 0, |u| = |v| = 1/ \sqrt{2} \}.$$
The submanifold $V$ lies in $S^{2n+1}$ with 
codimension 2, so $V$ has dimension $2n-1$.
Let $\alpha$ and
$\beta$ be positive real numbers satisfying $\alpha ^2 + \beta ^2 = 1$,
and let $T_{\alpha ,\beta }$ be the linear transformation of
$E^{2n+2}$ defined by
$$T_{\alpha ,\beta }(u,v) = \sqrt{2} (\alpha
u, \beta v).$$
Then the image $W^{\alpha ,\beta} = T_{\alpha
,\beta}V$ is contained in $S^{2n+1}$ and it
is proper Dupin with $g=4$. On can show
that the Lie curvature is not constant on
$W^{\alpha ,\beta}$ if $\alpha \neq 1/
\sqrt{2}$. Thus $W^{\alpha ,\beta}$ is  not
Lie equivalent to an isoparametric
hypersurface if $\alpha \neq 1/\sqrt{2}$.
\vskip.1in

If $M^3$ is an isoparametric hypersurface
with $k$ principal curvatures, the inverse
image $\psi^{-1}(M^3)$ is an isoparametric
hypersurface with $2k$ principal curvatures
in $S^7$. When $k=3$, $M^3$ must be a tube
over a Veronese surface, and so the unique
family of isoparametric hypersurfaces with
$g=6$ in $S^7$ has a precise geometric
characterization in terms of $\psi$. Miyaoka
(1993) showed that in this case
$\psi^{-1}(M^3)$ is homeomorphic to $M^3\times
S^3$ and it has a foliation whose leaves are
isoparametric hypersurfaces with three
principal curvatures of multiplicity one.
The two focal submanifolds of $\psi^{-1}(M^3)$
are obtained from the two focal submanifolds
of $M^3$ via $\psi^{-1}$. However, although
the two focal submanifolds of $M^3$ are
Veronese surfaces which are congruent in
$S\sp 4$, the two focal submanifolds of
$\psi^{-1}(M^3)$ are not congruent in $S^7$.
Thus, they are two noncongruent minimal taut
homogeneous embeddings of $RP^2\times S^3$
into $S^7$.

Pinkall and Thorbergsson (1989) introduced 
the  M\"obius curvature which can
distinguish among the Lie equivalent
parallel hypersurfaces in a family of
isoparametric hypersurfaces. 
C. P. Wang (1992) applied the method of
moving frames to determine a complete set of
M\"obius invariants for surfaces in $E^3$
without umbilic points and for
hypersurfaces in $E^4$ with three
distinct principal curvatures at each
point.   He then applied this result to
derive a local classification of Dupin
hypersurfaces in $E^4$ with three
principal curvatures up to M\"obius
transformation. 

Ferapontov (1995a,1995b) studied the
relationship between Dupin and 
isoparametric hypersurfaces and Hamiltonian
systems of hydrodynamic type.

\vskip.1in
\noindent{\bf 10.4.6. Tubes as Dupin
hypersurfaces}

M. Takeuchi (1991) studied Dupin
hypersurfaces in real space forms which are
tubes around symmetric submanifolds in
$R^{n+1}(c)$ and proved the following:
Let $M$  be a non-totally geodesic
symmetric submanifold of $R^{n+1}(c)$ of
codimension $>1$. Then  the
$\varepsilon$-tube $T_\varepsilon(M)$
around $M$ is a proper Dupin hypersurface if
and only if either
\smallskip

(i) $M$ is a complete extrinsic sphere of
$R^{n+1}(c)$ of codimension $>1$; or

(ii) $M$ is one of the following symmetric
submanifolds of the
$n$-dimensional sphere $S\sp n$:

(ii-a) the projective plane $
FP^2\subset S^{3d+1}$,
$d=\dim_{R} F$, over
 $ F= R, C$, the quaternions
$ H$, or octonions $\Cal O$;

(ii-b) the complex quadric $Q_3\subset S^9$;

(ii-c) the Lie quadric $Q^{m+1}\subset
S^{2m+1}$, $m\geq 2$;

(ii-d) the unitary symplectic group
${Sp}(2)\subset S^{15}$.
\smallskip

In case (i),  $T_\varepsilon(M)$ is a Dupin
cyclide, that is, a proper Dupin
hypersurface with two distinct principal
curvatures, but it is not isoparametric. In
case (ii), $T_\varepsilon(M)$ is a
homogeneous isoparametric hypersurface with
three or four distinct principal
curvatures, and it is irreducible in the
sense of Pinkall.

\vskip.1in
\noindent{\bf 10.4.7. Topology of Dupin hypersurfaces}

G. Thorbergsson (1983a) proved the following:

(1) a compact embedded proper Dupin hypersurface $M$ in $R^{n+1}(c)$ satisfies $\dim
H^*(M;\hbox{\bf Z}_2)=2g$, where $g$ is the number of distinct principal curvatures; and

(2) a compact embedded proper Dupin hypersurface $M$ divides $S^{n+1}$ into two ball bundles.

 By using (2),  Thorbergsson showed that $g=1,2,3,4$ or 6.
\smallskip

The integral homology, fundamental group and rational homotopy type of a compact Dupin
hypersurface  in the sphere were determined in [Grove-Halperin 1987].

In particular, K. Grove and S. Halperin proved that if $M$ is a compact Dupin hypersurface of $S^{n+1}$, then

(1) there are two integers $k,\ell$ (possible equal) such that each principal
curvature has multiplicity $k$ or $\ell$;

(2) the integral homology of $M$
determines $k,\ell$, and the number $g$ of
principal curvatures. Conversely, $g,k$
and $\ell$ determine the fundamental group,
integral homology, and rational homotopy
type of $M$;

(3) the integers $g,k,\ell$ satisfy the following restrictions:

(3-a) if $k\ne\ell$, then $g=2$ or
$4$, and $k$ and $\ell$ are each the
multiplicity of $g/2$ principal curvatures;

(3-b) if $g=3$, then $k=1,2,4$ or
$8$,

(3-c) if $g=4$ and $k=\ell$, then
$k=1$ or $2$; furthermore, if $g=4$ and
$k>\ell\geq 2$, then $k+\ell$ is odd,

(3-d) if $g=6$, then $k=1$ or $2$.

It is immediate from the above results that
$n={1\over 2}(k+\ell)g$. Thus, $g=1$ if and
only if $k+\ell>n$;  $g=2$ if and
only if $k+\ell=n$; and  $g=3,4$ or 6 if and
only if $k+\ell<n$.

\vskip.1in
\noindent{\bf 10.4.8. Dupin hypersurfaces of
$T_1S^{n+1}$}

Let $T_1S^{n+1}$ denote the unit tangent
bundle of $S^{n+1}$. Consider $T_1S^{n+1}$ as
the  $(2n+1)$-dimensional submanifold of 
$S^{n+1} \times S^{n+1} \subset E^{n+2}
\times E^{n+2}$ given by
$$T_1S^n = \{(x,\xi): \;|x| =1, \, |\xi| =1,
\,\left<x, \xi\right> =0 \}.$$ Then
$T_1S^{n+1}$ admits a canonical contact
1-form $\omega$ induced from the
canonical almost complex structure  on
$C^{n+2}= E^{n+2}
\times E^{n+2}$; thus,
$\omega\wedge(d\omega)^n\ne 0$ everywhere.
The contact structure gives rise to a
codimension one distribution on
$T_1S^{n+1}$ which has integrable
submanifolds of dimension $n$, but none of
higher codimension. An
$n$-dimensional integrable submanifold of
$T_1S^{n+1}$ is called a Legendre
submanifold. 

Pinkall (1985a) investigated
the Legendre submanifold of $T_1S^{n+1}$
that he called Lie geometric hypersurfaces
of the sphere $S^{n+1}$. Each oriented
hypersurface $M$ of $S^{n+1}$ gives rise to
a Legendre submanifold $L_M$ of $T_1S^{n+1}$
by associating to $M$ the set $L_M$ of
oriented unit normal vectors along $M$.

The image of a Legendre submanifold is
called a wavefront. Pinkall  showed how the
basic theory of hypersurfaces can be
extended to wavefronts. A contact
transformation is a diffeomorphism $F$ of
$T_1S^{n+1}$ which satisfies the property
$dF(\hbox{ker}\,\omega)=\hbox{ker}\,\omega$.
In particular, if $F$ leaves the
class of Legendre submanifolds that come
from lifting an oriented hypersurface of
$S^{n+1}$ invariant, then $F$ is called a
Lie sphere transformation. The class of
transformations of $S^{n+1}$ that map
spheres to spheres are called M\"obius
transformations; they are exactly the
conformal automorphisms of $S^{n+1}$,
according to a result of J. Liouville. 

The classification of Lie sphere
transformations is done by mapping the space
of oriented spheres of $S^{n+1}$ onto the
Lie quadric, that is, the quadric of type
$(n+2,2)$ in $RP^{n+3}$. Then the Lie sphere
transformations correspond to the
projective transformations of $RP^{n+3}$
which leave the Lie quadric invariant.

Pinkall pointed out that, for each Legendre
submanifold $L$ of $T_1S^{n+1}$, the concept 
of a principal direction can be defined at a
point $p\in L$ as a direction in which $L$
has higher order contact at $p$ with the
lift of a hypersphere, called an osculating
sphere, to $T_1S^{n+1}$. Define the
principal radii of $L$ at $p$ to be the
radii of the corresponding osculating
spheres. The tangent space $T_pL$ then
decomposes into $E_1\oplus\cdots\oplus
E_g$, where each $E_i$ is a maximal
subspace which consists of principal
directions. The multiplicity of a principal
radius is nothing but the dimension of the
corresponding space $E_i$. 

Let $L$ be a Legendre submanifold of
$T_1S^{n+1}$ and let $S$ be a submanifold of
$L$ such that $T_pS$ is one of the spaces
$E_1,\ldots,E_g$ in the decomposition. $S$
is called a curvature surface of $L$
according to Pinkall (1985a). When $L$ is a
Legendre submanifold of $T_1S^{n+1}$ such
that a continuous principal radius function
is constant along its corresponding 
curvature surface, then $L$ is called a
Dupin hypersurface  [Pinkall 1985a]. 

In [Pinkall 1985a] a proper Dupin
hypersurface is defined as a Dupin
hypersurface $L$ which satisfied the
property that the multiplicities of the
principal radii at each point $p\in L$ are
independent of the point $p$. Pinkall showed
that the class of proper Dupin hypersurfaces
is invariant under the Lie sphere
transformations. 

As a generalization of the classical cyclides
in $E^3$, Dupin hypersurfaces in $T_1S^{n+1}$
with two distinct principal radii at each
point are called cyclides of Dupin. Pinkall
classified the cyclides of Dupin in
$T_1S^{n+1}$ and proved that they are Lie
equivalent to the Lagrangian submanifold of
$T_1S^{n+1}$ obtained by the standard product
embedding of $S^k(1/\sqrt{2})\times 
S^{n-k}(1/\sqrt{2})$, where $k$ and $n-k$ are
the multiplicities of the principal
curvatures of the cyclide of Dupin. 

Pinkall (1985b) proved that a
Dupin hypersurface with three distinct
principal curvatures in  $T_1S^4$ is either
reducible or Lie equivalent to a piece of
Cartan's isoparametric hypersurface in
$S^4$ with three distinct principal
curvatures. T. E. Cecil and G. R. Jensen 
(1998) extended Pinkall's result and
showed that if a proper Dupin hypersurface
in $T_1S^{n+1}$ contains no reducible open
subset, then it is Lie equivalent to a
piece of an isoparametric hypersurface with
three principal curvatures.

\subsection{Hypersurfaces with constant
mean curvature}

Surfaces of constant mean curvature occur
naturally in physics. S. D. Poisson
(1781--1840) showed in 1828 that if a surface
in $E^3$ is the interface between two media
in equilibrium, then the mean curvature $H$
of the surface is constant and equal to
$H=k(p_1-p_2)$, where
$p_1$ and $p_2$ denote the pressures in the
media and $\lambda=1/k$ is called the
coefficient of surface tension.
Furthermore, a hypersurface with constant mean
curvature is a solution to a variational
problem; namely, with respect to any
volume-preserving normal variation of a
domain $D$ in a Euclidean space, the mean
curvature of $M=\partial D$, the boundary of
$D$, is constant if and only if the area of
$M$ is critical, that is, it satisfies
$A'(0)=0$. 
\vskip.1in

\noindent{\bf 10.5.1. Hopf's problem
and Wente's tori}

H. Hopf (1894--1971) proved in 1951 that any
immersion of a surface, topologically a
sphere, with  constant mean curvature in  $
E^3$ must be a round sphere. A. D.
Alexandrov (1912-- ) showed in 1958 that the
only compact embedded surfaces in $ E^3$ of
constant mean curvature are  round spheres.
That left open the possibility of immersed
constant mean curvature surfaces of higher
genus.

Surprisingly, H. C. Wente constructed in
1984  examples of tori of constant mean
curvature  in $ E^3$. Wente's examples
solved the long-standing problem of Hopf: Is
a compact constant mean curvature immersed
surface in $ E^3$ necessarily a round
sphere?

Wente's work inspired a string of further
research on compact tori of constant mean
curvature [Abresch 1987, Spruck 1986,
Pinkall-Sterling 1989]. For instance, 
U. Abresch classified all constant mean
curvature tori having one family of planar
curvature lines. 

 N. Kapouleas (1991) showed that there
also exist  compact  constant mean
curvature  surfaces of  every genus $\geq
3$ in $E^3$. Also Kapouleas (1991)
provided a general construction method for
complete surfaces of constant mean
curvature in $ E^3$.

W. Y. Hsiang (1982b,1982c) constructed 
infinitely many noncongruent immersions
of topological $n$-sphere
$S^n$ into $E^{n+1}$ with constant mean
curvature for each $n\geq 3$.
\vskip.1in

\noindent{\bf 10.5.2. Delaunay's
surfaces and generalizations}

Delaunay's surfaces in $ E^3$ introduced
by C. E. Delaunay (1916--1872) are surfaces
of revolution of constant mean curvature. N.
J. Korevaar, R. Kusner and B.  Solomon (1989)
proved that, besides the circular cylinder,
Delaunay surfaces are the only
doubly-connected surfaces properly embedded
in $ E^3$ with nonzero constant mean
curvature.

W. Y. Hsiang and W. C. Yu (1981) and W.
Y. Hsiang (1982a,1982c) constructed 
one-parameter family of hypersurfaces of
revolution, symmetric under $O(n - 1)$,
having constant mean curvature in the 
$n$-sphere $S^n(1)$, in the Euclidean
$n$-space $E^n$, or in the hyperbolic
$n$-space $H^n(-1)$.
\vskip.1in

\noindent{\bf 10.5.3. Hypersurfaces
with $K\leq 0$ or with $K\geq 0$}

T. Klotz and R. Osserman (1966) proved that a
complete surface of nonzero constant mean
curvature is a circular cylinder if its Gauss
curvature
$K$ is $\leq 0$.  B. Smyth and K. Nomizu
(1969) showed that the only compact
hypersurfaces of constant mean curvature in
$ E^{n+1}$ with non-negative sectional
curvatures are  round $n$-spheres.

B. Smyth and K. Nomizu (1969) also proved
that  totally umbilical hypersurfaces and a
standard product embedding of the product
of two spheres are the only compact
hypersurfaces of constant mean curvature
in $S^{n+1}$ with nonnegative sectional
curvature. 

 Let $M$ be a compact hypersurface 
of constant mean curvature in $S^{n+1}$.
Denote by $B=\sqrt{S}$ the norm of the
second fundamental form of $M$ in $S^{n+1}$.
Z. H. Hou (1997) proved the following: 

(1) If $B<2\sqrt{n-1}$, $M$ is a small
hypersphere $S^{n}(r)$ of radius
$r=\sqrt{n/(n+B)}$.

(2) If $B=2\sqrt{n-1}$, $M$ is either
$S^n(r_0)$ or $S^1(r)\times S^{n-1}(s)$, 
where $$r_0^2={n\over{n+2\sqrt{n-1}}},\;\;
r^2={1\over{1+\sqrt{n-1}}},\;\;
s^2={\sqrt{n-1}\over{1+\sqrt{n-1}}}.$$
\vskip.1in

\noindent{\bf 10.5.4. Hypersurfaces of constant
mean curvature in hyperbolic
spaces}

In the hyperbolic 3-space $H^3$ of sectional
curvature $-1$, the behavior of a surface of
constant mean curvature $H$ depends on the
value of $H$. If $|H|<1$, the area of the
surface grows exponentially. In the case
where $|H|>1$, the surface can be compact,
like geodesic spheres. In the border case
$|H|=1$, there exist examples such that the
area grows polynomially, and it is known that
some properties similar to those of minimal
surfaces in Euclidean space hold. 

N. Korevaar, R. Kusner, W. H. Meeks and B.
Solomon (1992) studied constant mean
curvature surfaces  in  $H^3$ and proved the
following: 

Let $M$ be a complete
properly embedded surface  in $H^3$ with
constant mean curvature greater than that of
a horosphere. Then 

(1) $M$ is not homeomorphic to a closed
surface punctured in one point. 

(2) If $M$ is homeomorphic to a closed
surface punctured in two points, then
$M$ is Delaunay, that is, $M$ is a
constant mean curvature surface of
revolution.

(3) If $M$ is
 homeomorphic to a closed surface punctured
in three points, then
$M$ remains at a bounded distance from 
a geodesic plane of reflective symmetry and
each half of
$M$ determined by the geodesic plane is
 a graph over this plane with respect to the
distance function to the plane. 

Moreover, they showed that the annular ends
of $M$ must exponentially converge to
Delaunay surfaces, which are constant mean
curvature surfaces of revolution in $H^3$.

For hypersurfaces of constant mean 
curvature in  $H^{n+1}$, M. do Carmo
and H. B. Lawson (1983) proved  the following: 

Let $M$ be a complete  hypersurface properly
embedded in
$H^{n+1}$ with constant mean curvature, and let
$\partial_\infty M\subset S^n (\infty)$ be its
asymptotic boundary. Then

 (a) if $M$ is
compact, it is a sphere;
if $\partial_\infty M$ consists of
exactly one point, then
$M$ is a horosphere; and

 (b) if $\partial_\infty M$ is a sphere and
$M$ separates poles, then $M$ is a
hypersphere.

As a consequence,  if $M$ is a
hypersurface of constant mean curvature in
$H^{n+1}$, admitting a one-to-one projection onto a
geodesic hyperplane, then $M$ is a hyperplane. This
result can be restated as follows: A nonparametric
entire hypersurface, that is, the graph of a
function $f$ defined in some
$H^n$, of constant mean curvature is a
hypersphere, which is a close analogue of the
 Bernstein theorem. 

\vskip.1in

\noindent{\bf 10.5.5. Weierstrass type 
representation for surfaces with
constant mean curvature}

K. Kenmotsu (1979) established an integral
representation formula for arbitrary surfaces
in $E^3$ with nonvanishing mean
curvature $H$ which describes the
surface as a branched conformal
immersion in terms of its mean curvature and
its Gauss map. Specifically, let
$\psi$ denote the complex-valued function on
the surface obtained by composing the Gauss
map with stereographic projection. Then,
letting $f= -\psi_z/H(1+\psi^2),$ one has
$$x_1=\hbox{Re} {\int} (1 - \psi^2)f
dz,\;\; x_2=\hbox{Re} {\int}i(1+ \psi^2)f
dz,\;\; x_3=2\hbox{Re} {\int}\, \psi f
dz.$$ Thus, given a real function $H$ and a
complex function $\psi$ on a simply-connected
domain D, a necessary and sufficient
condition that there exists a map $x: D\to
E^3$ defining a branched surface in
isothermal coordinates having $H$ as mean
curvature and $\psi$ as the Gauss map is that
the differential relations above between $H$
and $\psi$ be satisfied. As a special case, if
$H$ is a nonzero constant, then the resultant
equation for $\psi$ is the precise condition
for $\psi$ to define a harmonic map into the
unit sphere.  By virtue of this representation
formula, if a  harmonic map $\psi$ from a
Riemann surface $\Sigma$ into $S^2$ is given,
then one can construct a branched immersion of
a constant mean curvature surface whose Gauss
map is $\psi$.

J. Dorfmeister, F. Pedit and H. Wu (1997) showed
that every constant mean
curvature immersion $\Phi:D\to E^3$, $D$ the
whole complex plane or the open unit disk in
{\bf C}, can be produced from a meromorphic
matrix valued 1-form
$$\xi=\lambda^{-1}
\begin{pmatrix} 0& f(z)\\g(z)& 0\end{pmatrix} dz,$$
$\lambda\in S^1$, the so called meromorphic
potential. Here $f$ and $g$ are meromorphic
functions of $z\in D$ and $f(z)g(z)=E(z)$, where
$E(z)dz^2$ is, up to a constant factor, the
Hopf differential of the surface. J. Dorfmeister,
and G. Haak (1997)  gave an explicit
characterization of the zero and pole order of
the meromorphic functions for the
branched constant mean curvature surface to
be a smooth immersion.

R. Bryant (1987b) showed that there also
exists a  Weierstrass type  representation 
for surfaces of constant mean curvature
$H=c>0$ in $H^3(-c^2)$.  In particular, any 
constant mean curvature one surface in 
$H^3(-1)$ can be constructed from an ${\frak
{sl}}(2,{\bold C})$-valued  holomorphic
1-form satisfying some  conditions (or
equivalently a pair of a meromorphic function
and a holomorphic  1-form) on a Riemann
surface. 

Bryant's
representation can be stated as follows.

We identify each point $(t,x_2,x_3,x_4)$ of
$L^4$ with a $2\times 2$ Hermitian matrix
$$\begin{pmatrix} t+x_4 & x_2+ix_3\\ x_2-ix_3& t-x_4\end{pmatrix} \in\;\hbox{Herm}\,(2).\leqno(10.5)$$
Then $H^3(-1)$ is identified with
$$\aligned H^3(-1)&=\{X\in \hbox{Herm}\,
(2):\det (X)=1,\;\hbox{trace}\,(X)>0\}\\
&=\{a\cdot a^*:a\in SL(2,
C)\},\endaligned \leqno(10.6)$$
where $a^*=\bar a^T$. Under this
identification, each element $a$ of the
group $PSL(2, C):= SL(2, C)/\{\pm
1\}$ acts isometrically on $H^3(-1)\ni
X\mapsto a\cdot X\cdot a^*$.

Bryant proved the following result.

 Let $M$ be a
simply-connected Riemann surface and $z_0\in
M$ a fixed point. Take a meromorphic
function $\psi$ and a holomorphic 1-form
$\omega$ on $M$ such that
$ds^2:=(1+|\psi|^2)^2\omega\cdot\bar\omega$
is positive definite on $M$. Then 
there exists a unique holomorphic immersion
$F:M\to PSL(2, C)$ that satisfies

(1) $F(z_0)=\pm\,$id.;

(2) $F^{-1}\cdot dF=\begin{pmatrix} \psi &
-\psi^2\\ 1 & -\psi\end{pmatrix}\omega$;

(3) $f=F\cdot F^*:M\to H^3(-1)$ is a conformal immersion with constant man
curvature 1  whose first fundamental form is $ds^2$.

Conversely, any conformal constant mean
curvature 1 immersions in
$H^3(-1)$ are obtained as above.

 M. Umehara and
K. Yamada (1996) showed that Bryant's
representation formula for surfaces of
constant mean curvature $c$ in $H^3(-c^2)$
can be deformed to the Weierstrass
representation formula as $c$ tends to 0. 

R. Aiyama and K. Akutagawa (1997a) gave 
representation formulas for surfaces of
constant mean curvature $H$ in the
hyperbolic 3-space $H^3(-c^2)$ with $H>c>0$.
R. Aiyama and K. Akutagawa (1997b) also gave 
representation formulas for surfaces of
constant mean curvature in the 3-sphere
$S^3(c^2)$. Their formulas show that every
such surface in $H^3(-c^2)$ or in $S^3(c^2)$
can be represented locally by a harmonic map
to the unit 2-sphere. 

Further results on surfaces of  constant mean
curvature  in $H^3$ were obtained by M.
Umehara and K. Yamada
(1992,1993,1996,1997a, 1997b).
\vskip.1in

\noindent{\bf 10.5.6. Stability of surfaces with
constant mean curvature}

Since a compact constant mean curvature
surface in $E^3$ is a critical point of the
area functional with respect to
volume-preserving normal variations, one can
define the stability of such surfaces:   A
compact  constant mean curvature surface in
$E^3$ is called stable if $A''(0)> 0$ with
respect to the  class of volume-preserving
normal variations. 

M. do Carmo and A. Da Silveira (1990) proved
that the index of $\Delta-2K$ is finite if
and only if the total curvature is finite
for a complete surface of constant mean
curvature one in the hyperbolic 3-space
$H^3$.

J. A. Barbosa and M. do Carmo (1984) proved
that the spheres are the only compact stable
hypersurfaces of constant mean curvature in
$E^{n+1}$.  This result was generalized
for closed constant mean
hypersurfaces in $S^{n+1}$ and $H^{n+1}$ by
Barbosa, do Carmo and J. Eschenburg in 1988.
For surfaces this result was
extended by H. Mori (1983), B. Palmer (1986)
and F. J. L\'opez and A. Ros (1989) to
complete surfaces, where the stability
assumption applied to every compact
subdomain and the surface is assumed to have
nonzero constant mean curvature. 

Also A. Da Silveira (1987) studied complete
noncompact surfaces which are immersed as
stable constant mean curvature surfaces in
$ E^3$ or in the hyperbolic space 
$H^3$. In the case of $E^3$ he
proved that the immersion is a plane.  For
$H^3$ he showed, under the condition of
nonnegative mean curvature $H$, that for
$H\geq 1$ only horospheres can occur, while
for $H<1$ there exists a one-parameter family
of stable nonumbilic embeddings. His theorems
are generalizations and improvements of
previous results by Barbosa, do Carmo and
Eschenburg, Fischer-Colbrie, Peng, and
Schoen. His notion of stability is slightly
weaker than the previous one for minimal
immersions. Examples show that the two
definitions in general do not agree. For
three-dimensional simply-connected complete
Riemannian manifolds with positive constant
sectional curvature he proved that there
exist no complete and noncompact stable
immersions with constant mean curvature.

For hypersurfaces in $ E^{n+1}$ H. P.
Luo (1996) proved that if a complete
noncompact stable hypersurface has
nonnegative Ricci curvature, it is minimal.

\subsection{Hypersurfaces with
constant higher order mean
curvature}

The $r$-th mean curvature $H_r$ of a
hypersurface $M$ is defined as the 
elementary symmetric polynomial of degree $r$
in the principal curvatures
$\kappa_1,\ldots,\kappa_n$ of $M$, that is, 
$$H_r=\sum_{i_1<\cdots<i_r}
\kappa_{i_1}\cdots \kappa_{i_r}.\leqno (10.7)$$

For a hypersurface in $ E^{n+1}$,
$H_1,H_2$ and $H_n$ are the mean curvature,
the scalar curvature, and the Gauss-Kronecker
curvature, respectively (up to suitable
constants). 

For a compact oriented hypersurface $f:M\to
 E^{n+1}$, the $r$-th mean curvatures are
related by the following formulas of H.
Minkowski (1903):
$${n\choose {r-1}}\int_M H_{r-1}dV=-
{n\choose {r}}\int_M
\left<f,\xi\right>H_rdV,\quad 1\leq r\leq
n,\leqno (10.8)$$ where $\xi$ is a unit normal
vector field of $M$ in $ E^{n+1}$.

A. Ros (1987) and, 
independently, N. J. Korevaar (1988) proved
that, for any $r, 1\leq r\leq n$, the round
sphere is the only compact hypersurface with
constant $r$-th mean curvature $H_r$ embedded
in $ E^{n+1}$. 

For a compact hypersurface $M$ embedded in
hyperbolic space, N. J. Korevaar (1988) and
Montiel and Ros (1991) proved that if any of
the higher order mean curvatures is constant,
then $M$ must be a geodesic hypersphere. 
The same is true if $M$
lies in a hemisphere of $S^{n+1}$. 

The above results
are not true in general, since all the
isoparametric hypersurfaces of $S^{n+1}$ have
all mean curvatures constant. 

For a compact immersed hypersurface $f:M\to
R^{n+1}(c)$ in a complete simply-connected
real space form
$R^{n+1}(c)$, let $F_r(H_1,H_2,\ldots,H_r)$
be the function defined inductively by
$$F_0=1,\quad F_1=H_1,\quad
F_r=H_r+{{(n-r+1)c}\over {r-1}}F_{r-2},\quad
2\leq r\leq n-1.\leqno (10.9)$$
A variation of $f$ is a differentiable map
$X:I\times M\to R^{n+1}(c)$ such that 
$X_0=f$ and, for each $t\in I,\,
X_t(x)=X(t,x),\, x\in M,$ is an immersion.

The balance of volume is defined to be the
function $V:I\to
\hbox{\bf R}$ given by
$V(t)=\int_{[0,t]\times M} F^*(d\bar V),$
where $d\bar V$ denotes the volume element
of $R^{n+1}(c)$. In Euclidean case, it
measures the balance of the volume of the
enclosed domain from the time 0 to time
$t$. So, in this case,
$V\equiv 0$ means that the volume of the
domain bounded by the hypersurface is kept
constant while the time changes. 

A variation of $f$ is said to be
volume-preserving if $V(t)\equiv 0$.

Put  $$A_r=\int_M
F_r(H_1,H_2,\ldots,H_r)dV.\leqno (10.10)$$

In the class of volume-preserving variations
of $f:M\to R^{n+1}(c)$, the first
variational formula of $f$ is given by
$$A'_r(\phi)=\int_M\{-(r+1)H_{r+1}+\kappa\}\phi
dV,
\leqno (10.11)$$
where $\kappa$ stands for a constant and
$\phi$ is the normal projection of the
variation vector field $\xi$. Thus, immersions
with constant
$(r+1)$-th mean curvature  arise as critical
points for the variational problem of
minimizing $A_r$, keeping the balance of
volume zero (cf. [Barbosa-Colares 1997,
Reilly 1973]).

An immersion $f:M\to R^{n+1}(c)$ with
constant $(r+1)$-th mean curvature
is said to be $r$-stable if its second
variation $A''_r(\phi)$ is $> 0$, for any
compact support function $\phi:M\to
\hbox{\bf R}$ that satisfies $\int_M \phi
dV=0$.

H. Alencar, M. do Carmo and H. Rosenberg
(1993) proved that hyperspheres are the only
$r$-stable immersed compact orientable
hypersurfaces in Euclidean space.

J. L. Barbosa and A. G. Colares (1997)
showed that geodesic hyperspheres are the
only $r$-stable immersed  compact orientable
hypersurfaces in an open hemisphere of
$S^{n+1}$ or in the hyperbolic space
$H^{n+1}$. When $r=1$, this is due to
[Alencar-do Carmo-Colares 1993].

\subsection{Harmonic spaces and
Lichnerowicz conjecture}

A Riemannian manifold $M$ is called a
harmonic space if all sufficiently small
geodesic hyperspheres have constant mean
curvature.  

Let $(x_1,\ldots,x_n)$ denote a system of
Cartesian coordinates of Euclidean
$n$-space $ E^n$ centered at a point 0, and
 $\Delta$ denote the Laplacian of
$ E^n$. Then it is well-known that the
Laplace equation $\Delta \phi=0$ admits a nice
solution given by
$$\phi(x)=\begin{cases} r^{2-n}\quad\hbox{if}
\quad n>2,\\ \ln r\quad \hbox{if} \quad n=2,\end{cases}$$
where $r$ is the distance of the point $x$ from
the origin 0. This implies that the Laplace
equation has a solution which is constant on each
hypersphere centered at 0.  Clearly, this is
true for any arbitrary center in $ E^n$.

In his 1930 doctoral thesis at Oxford
University, H. S. Ruse (1905--1974) made an
attempt to solve Laplace's equation on a
general Riemannian manifold and to find a
solution which depends only on the geodesic
distance. He realized later that it was
implicitly assumed in his thesis that such
a solution always exists and this is not
the case. Stimulated by this incident, E.
C. Copson (1901--1980) and Ruse started in
1939 the study of the class of Riemannian
manifolds which admit such a solution. This
is the beginning of the study of harmonic
manifolds.

In 1944 A. Lichnerowicz (1915-- ) showed
that a harmonic manifold of dimension $\leq
4$ is either a flat space or a rank one
locally symmetric space. From this one
conjectures  that the same conclusion holds
true without the dimension  hypotheses;
which is known as  Lichnerowicz's conjecture.

Lichnerowicz's conjecture has been
proved by Z. I. Szab\"o (1990) for compact
harmonic manifolds with finite
fundamental groups. 

A Riemannian manifold of negative curvature
 is said to be asymptotically harmonic if the
mean curvatures of the geodesic horospheres
are constant. P. Foulon and F. Labourie (1992)
proved that if $M$ is a compact ($C^\infty$-) 
 negatively curved asymptotically harmonic
manifold, then the geodesic flow of $M$ is
$C^\infty$ conjugate to that of a rank one
locally symmetric space.
On the other hand, G. Besson, G. Courtois
and S. Gallot (1995) proved that a Riemannian
manifold whose geodesic flow is
$C\sp 1$ conjugate to that of a compact locally
symmetric manifold $N$ is isometric to $N$.
Thus, a compact negatively curved
asymptotically  harmonic Riemannian manifold
 is locally symmetric; this in
particular proves 
Lichnerowicz's conjecture for the compact
negatively curved case.

On the other hand, E. Damek
and F. Ricci showed in 1992 that 
there exist noncompact counterexamples to
the conjecture; namely, there exists a
class of harmonic homogeneous simply
connected manifolds of negative curvature
which are not symmetric.

Damek and Ricci's examples are given as
follows: Let $\frak n$ be a two-step nilpotent
Lie algebra
 with inner product $\left<\,\;,\;\right>$,
such that if $\frak z$ is the center of
$\frak n$ and $\frak o=\frak  z^\perp$, the map 
$J_Z: \frak  o\to\frak  o$ defined by
$\left< J_ZX,Y\right>=\left<[X,Y],Z\right>$
satisfies $J^2_Z=-\vert Z\vert^2$ for
all $X,\,Y\in\frak  o$ and $Z\in\frak  z$. The
connected and simply connected Lie group $N$
generated by this algebra $\frak  n$ is
classically referred to as the Heisenberg group.
Let $\frak  n$ be solvably extended to $\frak 
s=\frak  o\oplus\frak  z\oplus \hbox{\bf
R}T$ by adding the rule
$[T,X+Y]=X/2+Z$, and denote by $S=NA$ 
$(A=\exp_S(\hbox{\bf
R}T))$, the corresponding
connected and simply connected Lie group. By the
use of the admissible invariant metric, $S$ is
then made into a Riemannian manifold. Effecting on
$S$ a suitable Cayley transform one can introduce
the normal coordinates $(r,w)$ on the ball
$$B=\{(X,Y,t)\in\frak  o\oplus\frak  z\oplus
\hbox{\bf
R}T,\ r^2=\vert X\vert^2+\vert
Y\vert^2+\vert t\vert^2<1\}$$ around the
identity element
$\frak  e=(0,0,1)$, with respect to which the
volume element on $B$ is computed as
$$2^{m+k}\left({
\cosh}\left({\rho\over 2}\right)\right)^
k\left({\sinh}\left({\rho\over
2}\right)\right)^{m+k}d\rho\, d\sigma(w),$$
where
$$m=\dim\frak  z,\quad k=\dim\frak 
o,\quad\rho=\log\left({{1+r}\over{1-r}}
\right)$$
and
$d\sigma(w)$ denotes the surface element of the 
sphere $S^{m+k}$. Thus $S$ turns out to be a
harmonic space. Note that
$m$ $(=\dim\frak  z)$ is quite arbitrary. 

Consider the symmetric space $M=G/K$ of
 noncompact type and let $G=NAK$ be the
Iwasawa decomposition. Then the map $s\to
sK$ is an isometry of $S=NA$ onto $G/K$, and
it is known that if $M$ is a symmetric space
of rank one, the dimension $k$ of the center
$\frak  z$ of $N$ equals $1,3$ or 7 only.
This leads to infinitely many harmonic
spaces $S$ that are not rank-one symmetric.

\vfill\eject

\section{Totally geodesic submanifolds}

The notion of totally geodesic submanifolds
was introduced in 1901 by J. Hada\-mard
(1865--1963). Hadamard defined (totally)
geodesic submanifolds of a Riemannian
manifold as submanifolds such that each
geodesic of them is a geodesic of the
ambient space. This condition is equivalent
to the vanishing on the second fundamental
form on the submanifolds. 1-dimensional
totally geodesic submanifolds are nothing
but geodesics. Totally geodesic submanifolds
are the simplest and the most fundamental
submanifolds of Riemannian manifolds. 

It is easy to show that every connected
component of the fixed point set of an
isometry on a Riemannian manifold is a
totally geodesic submanifold.

Totally geodesic submanifolds of
a Euclidean
space  are affine subspaces  and totally
geodesic  submanifolds of a Riemannian sphere are
the greatest spheres.

It is much more difficult to classify  totally
geodesic submanifolds of a  Riemannian
manifold in general.

\subsection{Cartan's theorem}

Let $M$ be a Riemannian $n$-manifold
with $n\geq 3$. For a vector $v$ in the tangent
space $T_pM$ at $p\in M$, denote by $\gamma_v$ the
geodesic through $p$ whose tangent vector at $p$
is $v$. Denote by $R_v(t)$ the $(1,3)$-tensor on
$T_pM$ obtained by the parallel translation of
the curvature tensor at $\gamma_v(t)$ along the
geodesic $\gamma_v$. Also define a (1,2)-tensor
$r_v(t)$ on $T_pM$ by
$$r_v(t)(x,y)=R_v(t)(v,x)y,\quad x,y\in T_pM.$$

 The following result of \'E. Cartan
provides  necessary and sufficient
conditions for the existence of totally
geodesic submanifolds in Riemannian
manifolds in general.

 Let $V$ be a subspace of the tangent space
$T_pM$ of a Riemannian manifold
$M$ at a point $p$. Then the following three
conditions are equivalent.

(1) There is a totally geodesic
submanifold of $M$ through $p$ whose tangent
space at $p$ is $V$.

(2) There is a positive number $\epsilon$
such that for any unit vector $v\in V$ and
any $t\in (-\epsilon,\epsilon)$,
$R_v(t)(x,y)z\in V$ for any
$x,y,z\in V$.

(3) There is a positive number $\epsilon$
such that for any unit vector $v\in V$ and
any $t\in (-\epsilon,\epsilon)$,
$r_v(t)(x,y)\in V$ for any $x,y\in V$.

\subsection{Totally geodesic submanifolds of symmetric spaces}

The class of Riemannian manifolds with parallel
Riemannian curvature tensor, that is, $\nabla R=0$,
was first introduced independently by P. A.
Shirokov (1895--1944) in 1925 and by H. Levy in
1926. This class is known today as the class of
locally symmetric Riemannian spaces. \`E. Cartan
noticed in 1926 that irreducible spaces of this
type are separated into ten large classes each
of which depends on one or two arbitrary
integers, and in addition there exist twelve
special classes corresponding to the
exceptional simple groups. Based on this,
Cartan created his theory of symmetric
Riemannian spaces in his famous  papers
``Sur une classe remarquable d'espaces de
Riemann'' [Cartan 1926/7].

An isometry $s$ of a Riemannian manifold $M$ is
called involutive if its iterate $s^2=s\circ s$
is the identity map. A Riemannian manifold $M$
is called a symmetric space if, for each point
$p\in M$, there exists an involutive 
isometry $s_p$ of $M$ such that $p$ is an
isolated fixed point of $s_p$. The $s_p$ is
called the (point) symmetry of $M$ at the
point $p$. 

Denote by $G_M$, or simply by
$G$, the closure of the group of isometries
generated by
$\{s_p:p\in M\}$ in the compact-open topology.
Then $G$ is a Lie group which acts transitively on
the symmetric space; hence the typical isotropy
subgroup $H$, say at $o$, is compact and $M=G/H$.

Every complete totally geodesic submanifold of a
symmetric space is a symmetric space. For a
symmetric space $M$, the dimension of a maximal
flat totally geodesic submanifold of $M$ is a
well-defined natural number which is called the
rank of $M$, denoted by $rk(M)$. 

It follows from
the equation of Gauss that $rk(B)\leq rk(M)$
for each totally geodesic submanifold $B$ of a
symmetric space $M$.
\vskip.1in

\noindent{\bf 11.2.1. Canonical decomposition and
Cartan's criterion}

If $M=G/H$ is a symmetric space and $o$ is a point
in $M$, then the map $$\sigma:G\to G$$ defined
by $\sigma(g)=s_o g s_o$ is an involutive
automorphism of $G$. Let $\frak g$ and $\frak h$
be the Lie algebras of $G$ and $H$, respectively.
Then $\sigma$ gives rise to an involutive
automorphism of $\frak g$, also denoted by
$\sigma$.  $\frak h$ is the
eigenspace of $\sigma$ with eigenvalue 1. 

Let $\frak m$ denote the eigenspace of $\sigma$
on $\frak g$ with eigenvalue $-1$. One has the
 decomposition:
$${\frak g}={\frak h}+{\frak m},$$ which is
called the Cartan decomposition or the
canonical decomposition of
${\frak g}$ with respect to $\sigma$. 

The subspace $\frak m$ can
be identified with the tangent space of the
symmetric space $M$ at $o$ in a natural way.
A linear subspace $\Cal L$
of $\frak m$ is called a Lie triple system if it
satisfies
$[\,[{\Cal L},{\Cal L}],{\Cal L}]\subset {\Cal
L}$.

The following result of \'E. Cartan provides
a simple relationship between totally
geodesic submanifolds and Lie triple systems
of a symmetric space: Let $M$
be a symmetric space. Then a subspace $\Cal L$ of
$\frak m$ forms a Lie triple system if and only if
$\Cal L$ is the tangent space of a totally
geodesic submanifold of $M$ through $o$.
\vskip.1in

\noindent{\bf 11.2.2. Totally geodesic
submanifolds of  rank one symmetric
spaces}

Applying Cartan's criterion, J. A. Wolf
 completely classified in 1963 totally geodesic
submanifolds in rank one symmetric
spaces and obtained the following:

(1) The
maximal totally geodesic submanifolds of
the real projective $m$-space $RP^{m}$ are
$RP^{m-1}$; 

(2) The
maximal totally geodesic submanifolds of the
complex projective $m$-space $CP^m$ are $RP^m$
and $CP^{m-1}$;

(3) The maximal totally geodesic submanifolds
of the quaternionic projective $m$-space
$HP^m$ are $HP^{m-1}$ and $ CP^m$; and 

(4) The maximal totally geodesic submanifolds
of the Cayley plane $\Cal OP^2$ are $ HP^2$ and
$\Cal OP^1=S^8$. 
\vskip.1in

\noindent{\bf 11.2.3. Totally geodesic
submanifolds of complex quadric}

Applying Cartan's  criterion,  Chen and H. S.
Lue (1975a) classified totally geodesic
surfaces in the complex quadric:
$Q_m=SO(m+2)/SO(2)\times SO(m),$  $m>1$. The
complete classification of totally geodesic
submanifolds of $Q_m$ was obtained by   Chen
and T. Nagano in 1977. More precisely, they
proved the following.

(1) If $B$ is a maximal totally geodesic
submanifold of $Q_m$, $B$ is one of the following
three spaces: 

(1-a) $Q_{m-1}$, embedded as a K\"ahler
submanifold; 

(1-b) a local Riemannian product,
$\left(S^p\times
S^q\right)/\{\pm\,\hbox{id}\,\}$, of two
spheres
$S^p$ and $S^q$, $p+q=m$, of the same radius,
embedded as a Lagrangian submanifold; and 

(1-c) the complex projective
space $CP^n$ with $2n=m$, embedded as a K\"ahler
submanifold.

(2) If $B$ is a non-maximal totally geodesic
submanifold of $Q_m$, $M$ is either contained in
$Q_{m-1}$ in an appropriate position in $Q_m$, or
the real projective space $RP^n$ with $2n=m$.

Chen and Nagano (1977) also proved the
following:

(3) Each homology group $H_k(Q_m;\hbox{\bf
Z}),\, k<2m$, is spanned by the classes of
totally geodesic submanifolds of $Q_m$;

(4) The cohomology ring $H^*(Q_m;\hbox{\bf
Z})$ is generated by the Poincar\'e duals of
totally geodesic submanifolds of $Q_m$; and

(4) There is a maximal totally geodesic
submanifold $M$ of $Q_m$ such that the
differentiable manifold $Q_m$ is the union of
the normal bundles to $M$ and to its focal
manifold  with the nonzero vectors
identified in some way.
\vskip.1in

\noindent{\bf 11.2.4. Totally geodesic
submanifolds of compact Lie groups}

Totally geodesic submanifolds of  compact Lie
groups equipped with  biinvariant metrics have
been  determined in [Chen-Nagano 1978].

 Let $M$ be a compact Lie group with a
biinvariant metric. Then the local isomorphism
classes of totally geodesic submanifolds of
$M$ are those of symmetric space $B=G_B/H_B$
such that $G_B$ are subgroup of $G_M=M\times
M$.

\vskip.1in

\noindent{\bf 11.2.5. $(M_+,M_-)$-method}

In general, it is quite difficult to classify
totally geodesic submanifolds of a given
symmetric space with rank $\geq 2$ by
classifying the Lie triple systems via
Cartan's criterion. For this reason a  new
approach to compact symmetric spaces was
introduced by  Chen and  Nagano
[Chen-Nagano 1978, Chen 1987a]. Using
their method, totally geodesic submanifolds
in compact symmetric spaces were
systematically investigated. 

The method of
Chen and Nagano works as follows: A pair of
points $\{o,p\}$ in a compact symmetric
space $M$ is called an antipodal pair if
there exists a smooth closed geodesic
$\gamma$ in $M$ such that $p$ is the
midpoint of $\gamma$ from
$o$. For each pair
$\{o,p\}$ of antipodal points in a compact
symmetric space $M=G/H$, they
introduced a pair of orthogonal totally
geodesic submanifolds
 $M_+^o(p),M_-^o(p)$ through $p$ such
that $$\dim M_+^o(p)+\dim M_-^o(p)=\dim M,\quad
rk(M_-^o(p))=rk(M), \quad M_+^o(p)=H(p).$$
 
 The totally geodesic submanifolds 
$\,M_+$'s and $\,M_-$'s are called  polars
and meridians of $M$, respectively. 

A compact symmetric space $M$ is globally
determined by its polars and meridians. In
fact, two compact symmetric spaces $M$ and $N$
are isometric if and only if some pair
$(M_+(p),M_-(p))$ of $M$ is isometric to some
pair $(N_+(q),N_-(q))$ of $N$ pairwise
[Chen-Nagano 1978, Nagano 1992]. 

If $B$ is a
complete totally geodesic submanifold of a
compact symmetric space $M$, then, for any pair
$(B_+(p),B_-(p))$ of $B$, there is a pair
$(M_+(q),M_-(q))$ of $M$ such that $B_+(p)$
and
$B_-(p)$ are totally geodesic submanifolds of
$M_+(q)$ and $M_-(q)$, respectively. Since
the same argument applies to the totally
geodesic submanifold $B_+(p)\subset M_+(q)$
and to the totally geodesic submanifold
$B_-(p)\subset M_-(q)$, one obtains strings
of conditions.  

Also, given a pair of antipodal points
$\{o,p\}$ in a compact symmetric space
$M$, one obtains an ordered pair
$(M^o_+(p),M^o_-(p))$ as above.  Two pairs 
$(M^o_+(p),M^o_-(p))$ and 
$(M^{o'}_+(p'),M^{o'}_-(p'))$ are
called equivalent if there is an isometry on
$M$ which carries one to the other. Let
$P(M)$ denote the corresponding moduli space.
Then $P(M)$ is a finite set which is a global
Riemannian invariant of $M$. 

In general, one has $\# P(M)\leq 2^{rk(M)}-1$.
A compact irreducible symmetric space
satisfying $\# P(M)= 2^{rk(M)}-1$ if and only
if it is a rank one symmetric space.

Every isometric totally geodesic embedding
$f:B\to M$ of a compact symmetric space into
another induced a pairwise totally geodesic
immersion
$P(f):P(B)\to P(M)$. In particular, if $B$
and $M$ have the same rank, then $P(f)$ is
surjective, hence, $\#P(B)\geq \#P(M)$;
this provides us a useful obstruction to
totally geodesic embeddings as well.

In particular, Chen and Nagano's results imply
the following:

(1) Any compact symmetric space $M$
of dimension $\geq 2$ admits a totally
geodesic submanifold $B$ satisfying $\frac12
\dim M\leq \dim B<\dim M$;

(2) Spheres and hyperbolic spaces are the only
simply-connected irreducible symmetric spaces
admitting a totally geodesic hypersurface.

The result (2) was extended by K. Tojo (1997a)
to the following: Let $G$ be a compact simple
Lie group and $K$ a closed subgroup of $G$. If
the normal homogeneous space $M=G/K$ contains
a totally geodesic hypersurface, then $M$ is
a space with constant sectional curvature.

Further information on polars and meridians
and on their applications to both geometry and
topology can be found in [Chen 1987;
Chen-Nagano 1978; Nagano 1988; Nagano
1992; Nagano-Sumi 1989; Burns 1992,1993; 
Peterson 1987;  Burns-Clancy 1994].

There is a duality between totally geodesic
submanifolds of symmetric spaces of compact type
and  of their non-compact duals.

For the investigation of totally geodesic
submanifolds in some symmetric spaces of
non-compact type  with rank $\geq 2$, see
also [Berger 1957].
\vskip.1in

\noindent{\bf 11.2.6. The 2-number $\#_2M$}

The notion of 2-number was introduced by 
 Chen and  Nagano in
1982. The notion of
2-number can also be applied to determine
totally geodesic embeddings in symmetric
spaces. 

For a compact symmetric space $M$, the
2-number, denoted by $\#_2M$, is defined as
the maximal possible cardinality $\#_2A_2$
of a subset $A_2$ of $M$ such that the point
symmetry $s_x$ fixes every point of $A_2$
for every $x\in A_2$. 

The 2-number $\#_2M$
is finite. The definition is equivalent to
saying that $\#_2M$ is a maximal possible
cardinality $\#A_2$ of a subset $A_2$ of $M$
such that, for every pair of points $x$ and
$y$ of $A_2$, there exists a closed geodesic
of $M$ on which $x$ and $y$ are antipodal to
each other. Thus, the invariant can also be
defined on any Riemannian manifold.

The geometric invariant $\#_2M$ is an
obstruction to the existence of a totally
geodesic embedding $f:B\to M$, since the
existence of $f$ clearly implies the
inequality $\#_2B\leq \#_2M$. For example,
while the complex Grassmann manifold
$G^{C}(2,2)$ of the 2-dimen\-sional
complex subspaces of the complex vector space
$ C^4$ is obviously embedded into $G^{
C}(3,3)$ as a totally geodesic submanifold,
the ``bottom space'' $G^{C}(2,2)^*$
 obtained by identifying every
member of $G^{C}(2,2)$ with its
orthogonal complement in $ C^4$,
however, cannot be totally geodesically
embedded into $G^{C}(3,3)^*$, simply
because $\#_2G^{
C}(2,2)^*=15>12=\#_2G^{C}(3,3)^*$. 

The 2-number is not an obstruction to a
topological embedding; for instance, the
real projective space $ RP^n$ can
be topologically embedded  in a sufficiently
high dimensional sphere, but the 2-number
$\#_2 RP^n=n+1 >2$ simply prohibits a
totally geodesic embedding of  $RP^n$ into
any sphere whose 2-number is 2, regardless of
dimension.

The invariant $\#_2M$ has certain bearings
on the topology of $M$ in other aspects; for
instance, Chen and Nagano proved that
$\#_2M$ is equal to the Euler number $\Cal
X(M)$ of $M$, if $M$ is a semisimple
Hermitian symmetric space. And in general
they proved that the inequality
$\#_2M\geq \Cal X(M)$ holds for any compact
symmetric space $M$ (cf. [Chen-Nagano 1988]
for details). 

M. Takeuchi (1989)
proved that $\#_2M=\dim H(M;\hbox{\bf Z}_2)$
for any symmetric $R$-space.  This formula
is actually correct for every compact
symmetric space which have been checked.

For a group manifold $G$, Chen and Nagano
showed that $\#_2G=2^{r_2}$, where
$r_2$ is the 2-rank of $G$, which by
definition is the maximal possible rank of
the elementary 2-subgroup $\hbox{\bf
Z}_2\times\cdots\times \hbox{\bf Z}_2$ of
$G$. An immediate application of this fact
is that the algebraic notion of 2-rank of a
compact Lie group $G$, studied by Borel and
Serre (1953), can be determined by
computing the 2-number
$\#_2 G$ of the group manifold $G$ via the
theory of submanifolds.

For the determination of $\#_2M$ of compact
symmetric spaces and of group manifolds, and
 the relationship between $\#_2M$ and
$(M_+,M_-)$'s, see [Chen-Nagano 1988].

The notion of 2-number was extended in 1993 to
$k$-number for compact $k$-symmetric spaces
by C. U. S\'anchez. She also obtained a
formula for $k$-number similar to Takeuchi's
for flag manifolds. 

\subsection{Stability of totally geodesic
submanifolds}

A minimal submanifold $N$ of a Riemannian
manifold $M$ is called stable if its second
variation for the volume functional of $M$ is
positive for every variation of $N$ in $M$.
It is an interesting and important problem to
find all stable minimal submanifolds in each
symmetric space, in particular, to determine all
stable totally geodesic submanifolds. 
\vskip.1in

\noindent{\bf 11.3.1. Stability of submanifolds in
compact rank one symmetric spaces}

Stability of compact totally geodesic 
submanifolds in compact rank one symmetric spaces
have been completely determined in [Simons
1968, Lawson-Simons 1973, Ohnita 1986a]:

(1) Compact totally geodesic submanifolds of
$S^m$ are unstable.

(2) Compact totally geodesic submanifolds of
$RP^m$ are stable.

(3) A compact totally geodesic
submanifold  of $CP^m$ is stable if and only if
it is a  complex projective subspace.

(4) A compact totally geodesic submanifold 
of $HP^m$ is stable if and only if it is a 
quaternionic projective subspace.

(5) A compact totally geodesic
submanifold  of the Cayley plane $\Cal OP^2$ is
stable if and only if
it is a  Cayley projective line $\Cal OP^1=S^8$
of $\Cal OP^2$.
\vskip.1in

\noindent{\bf 11.3.2. An algorithm for
determining the stability of
totally geodesic submanifolds in
symmetric spaces}

There is a general algorithm, 
discovered by B. Y. Chen, P. F. Leung and T.
Nagano in 1980, for determining
the stability of totally geodesic
submanifolds in compact symmetric spaces.  

Let $N$ be a  totally geodesic submanifold
of a compact symmetric  space
$M$. There is a finitely covering group
$G_N$ of the connected isometry group $G_{N}^o$ of
$N$ such that $G_N$ is a subgroup of the
connected isometry  group $G_M$ of $M$ which
leaves $N$ invariant, provided that 
$G_{N}^o$ is semisimple. Let $\Cal P$ denote
the orthogonal complement of the Lie algebra
${\frak g}_{_{N}}$ in the Lie algebra ${\frak
g}_{_{M}}$ with respect to the biinvariant
inner product on ${\frak g}_{_{M}}$ which is
compatible with the metric of $M$. Every
member of ${\frak g}_{_{M}}$  is thought of
as a Killing vector field because of the
action of $\, G_M$ on $\, M$. 

Let $\hat P$ denote the space of the vector
fields corresponding to the member of $\Cal
P$ restricted to the submanifold $N$. Then,
to every member of $\, \Cal P\,$ there
corresponds a unique (but not canonical)
vector field $v\in {\hat P}$ which is a
normal vector field and hence
${\hat P}$ is a $G_N$-invariant subspace of the
space $\Gamma (T^{\perp}N)$ of the sections
of the normal bundle to $N$. Moreover, ${\hat P}$ is
homomorphic with $\Cal P$ as a $G_N$-module.
The group $G_N$ acts on 
$\Gamma(T^{\perp}N)$ and
hence on the differential operators: 
$\Gamma(T^{\perp}N)
\to\Gamma(T^{\perp}N)$. $G_N$ leaves $L$ fixed,
since $L$ is defined with $N$ and the metric of $M$
only. Therefore, each eigenspace of $L$ is left
invariant by $G_N$. 

Let $V$ be one of its
$G_N$-invariant irreducible subspaces. One has  a
representation $\rho : G_{N} \rightarrow GL(V).$ 
Denote by $c(V)$ or $c(\rho)$ the eigenvalue of the
corresponding  Casimir operator. 
To define $c(V)$ one fixes an orthonormal basis
$(e_{\lambda})$ for
${\frak g}_{_{N}}$ and consider the linear
endomorphism $C$ or $C_V$ of $V$ defined by 
$$C=-\sum\,\rho(e_{\lambda})^{2}.\leqno(11.1)$$
Then $C$ is $c(V)I_V$, where $I_V$ is
the identity map on $V$. The 
Casimir operator is given by $C_{V}=-\sum
[e_{\lambda} ,[e_{\lambda},V]\,]$ for every member
$v$ of $V$ (after extending to a neighborhood of
$N$). 

A compact totally geodesic submanifold $N $
$ (=G_{N}/K_{N})$ of a compact symmetric
space $M
\,(=G_{M}/K_{M})$ is stable as a minimal
submanifold if and only if one has $c(V)\geq c(P')$
for the eigenvalue of the Casimir operator of every
simple $G_N$-module $V$ which shares as a
$K_N$-module some simple $K_N$-submodule of the
$K_N$-module $T_{o}^{\perp}N$ in common with some
simple $G_N$-submodule $P'$ of $\hat P$.

Roughly speaking, the algorithm says that
$N$ is stable if and only if
$c(V)\geq c({\Cal P})$ for every $G_N$-invariant
irreducible space $V$ (cf. [Chen 1990]).

Applying their algorithm, Chen, Leung
and Nagano obtained in 1980 the following
results:

(1) A compact subgroup $N$ of a compact
Lie group $M$  with  a biinvariant metric is
stable if $N$ has the same rank as $M$ and
$M$ has nontrivial center.

(2) Every meridian $M_-$ of a compact group
manifold $M$ is stable if $M$ has
nontrivial center.

(3)  Let $G^R(p,q)=SO(p+q)/SO(p)\times SO(q)$
be a real Grassmann manifold isometrically
immersed in a complex Grassmann manifold
$G^C(p,q)$ as a totally real totally geodesic
submanifold in a natural way. Then
$G^R(p,q)$ is unstable in
$G^C(p,q)$.

When $p=1$, statement (3) reduces  to
a result of Lawson and Simons (1973).

Applying the algorithm, M. Takeuchi
(1984)  completely determined  the stability
of totally geodesic Lagrangian submanifolds
of compact Hermitian symmetric spaces. He
proved that if $M$ is a compact Hermitian
symmetric space and $B$ a compact Lagrangian
totally geodesic submanifold of $M$, then $B$
is stable if and only if $B$ is
simply-connected.

K. Mashimo and H. Tasaki (1990b) applied the
same algorithm to determine  the stability of
maximal tori of compact Lie groups and obtained
the following:

(1) Let $G$ be a connected closed
subgroup of maximal rank in a compact Lie
group $U$ equipped with biinvariant metric.
If a maximal torus of $U$ is stable, then
$G$ is also stable.

(2) Let $U$ be a compact connected simple Lie
group and $T$ be a maximal torus. Then $T$ is
unstable if and only if $U$ is isomorphic to
$$SU(r+1),\;\; Spin(5),\;\;
Spin(7),\;\;Sp(r)\quad
\hbox{or}\quad G_2.$$

Further results concerning the stability of
certain subgroups of compact Lie groups equipped
with  biinvariant metrics can also be found in
[Fomenko 1972, Thi 1977, Brothers 1986,
Mashimo-Tasaki 1990a].
 
The stabilities of all the $M_+$'s (polars) and the
$M_-$'s (meridians) of a compact irreducible
symmetric space
$M$ were determined by M. S. Tanaka (1995). In
particular, she proved that all
polars and meridians of a compact Hermitian
symmetric space are stable.

Let $G$ be a compact connected Lie group, 
$\sigma$ an automorphism of $G$ and $K=\{k\in
G\colon
\sigma(k)=k\}$. A mapping $\Sigma\colon G\to G$,
 $g\mapsto g\sigma g^{-1}$, induces the Cartan
embedding of $G/K$ into $G$ in a natural way.  
If $M$ is a compact simple Lie
group $G$, then the $G_+$'s are images of Cartan
 embeddings and the $G_-$'s are the sets of
fixed points of involutive automorphisms.

 K. Mashimo (1992)
proved that, if
$G$ is simple and $\sigma$ is involutive, the image $\Sigma(G/K)$ is unstable only if either
$G/K$ is a Hermitian symmetric space or the pair $(G,K)$ is one of the
four cases: 
$$ (SU(n),SO(n))\;\; (n\geq
3),\quad (SU(4m+2)/\{\pm\,I\}, SO(4m+2)/
\{\pm\,I\})\;\;
(m\geq 1),$$ $$ (Spin(n),(Spin(n-3)\times 
Spin(3))/{\hbox{\bf Z}}_2)\;\;(n\geq
7),\quad (G_2,SO(4)).$$ 

\vskip.1in

\noindent{\bf 11.3.3. Ohnita's formulas}

Y. Ohnita (1987) improved the above
algorithm to include the formulas for the
index, the nullity and the Killing nullity
of a compact totally geodesic submanifold in
a compact symmetric space.

Let \ $f:N\to M$ \ be a compact totally geodesic
submanifold of a compact Riemannian symmetric
space. Then  $f:N\to M$ is expressed as follows:
There are compact symmetric pairs
$(G,K)$ and $(U,L)$ with $N=G/K,\, M=U/L$ so
that  $f:N\to M$ is
given by $uK\mapsto \rho(u)L$, where  $\rho:G\to
U$ is an analytic homomorphism
with $\rho(K)\subset L$ and the injective
differential $\rho:\frak g\to\frak u$ which
satisfies $\rho(\frak m)\subset \frak p$. Here
$\frak u=\frak l +\frak p$ and $\frak g=\frak
k+\frak m$ are the canonical decompositions of
the Lie algebras  $u$ and $g$, respectively.

Let $\frak m^\perp$
denote the orthogonal complement of $\rho(\frak
m)$ with $\frak p$ relative to the
ad$(U)$-invariant inner product $(\;,\;)$ on
$\frak u$ such that $(\;,\;)$ induces the
metric of $M$. Let $\frak k^\perp$ be the
orthogonal complement of $\rho(\frak k)$ in
$\frak l$. Put $\frak g^\perp=\frak k^\perp
+\frak m^\perp$. Then $\frak g^\perp$ is the
orthogonal complement of $\rho(\frak g)$ in
$\frak u$ relative to $(\;,\;)$, and $\frak
g^\perp$ is ad$_\rho (G)$-invariant. 

Let $\theta$
be the involutive automorphism of the symmetric
pair $(U,L)$. Choose an orthogonal decomposition\
$\frak g^\perp=\frak g_1^\perp\oplus\cdots
\oplus\frak g_t^\perp$\ such that each $\frak
g^\perp_i$ is an irreducible ad$_\rho
(G)$-invariant subspace with $\theta(\frak
g^\perp_i)=\frak g^\perp_i$. Then, by Schur's
lemma, the Casimir operator $C$ of the
representation of $G$ on each $\frak g^\perp_i$ is
$a_i I$ for some $a_i\in\hbox{\bf C}$. 

Put $\frak m_i^\perp =\frak m\cap \frak
g^\perp_i$ and let $D(G)$ denote the set of all
equivalent classes of finite dimensional
irreducible complex representations of $G$. For
each $\lambda\in D(G)$,
$(\rho_\lambda,V_\lambda)$ is a fixed
representation of $\lambda$. 

For each $\lambda\in
D(G)$, we assign a map $A_\lambda$ from
$V_\lambda\otimes \hbox{Hom}_K(V_\lambda,W)$ to
$C^\infty(G,W)_K$ by the rule
$A_\lambda(v\otimes
L)(u)=L(\rho_\lambda(u^{-1})v)$. Here
Hom$_K(V_\lambda,W)$ denotes the space of
all linear maps $L$ of $V_\lambda$ into $W$
so that $\sigma(k)\cdot L=L\cdot
\rho_\lambda(k)$ for all
$k\in K$. 

Y. Ohnita's formulas for the index $i(f)$,
the nullity $n(f)$, and the Killing nullity
$n_k(f)$ are given respectively by

 $i(f)=\sum_{i=1}^t\sum_{\lambda\in D(G),
a_\lambda <a_i} m(\lambda)d_\lambda$,

 $n(f)=\sum_{i=1}^t\sum_{\lambda\in D(G),
a_\lambda =a_i} m(\lambda)d_\lambda$,

 $n_k(f)=\sum_{i=1,\frak
m_i^\perp\ne \{0\}}^t\dim \frak g_i^\perp$,

\noindent where $m(\lambda)=\dim\,
\hbox{Hom}_K(V_\lambda,(\frak m_i^\perp)^{
C})$ and $d_\lambda$ denotes the dimension
of the representation $\lambda$.

 By applying his formulas,
Ohnita  determined the  indices, the
nullities and the Killing nullities for all
totally geodesic submanifolds in compact
rank one  symmetric spaces and Helgason spheres
in all compact irreducible symmetric spaces.

Q. Zhao applied Ohnita's formulas in
[Zhao 1996] to determine the  indices, the
nullities and the Killing nullities of maximal
totally geodesic submanifolds of the complex
quadric
$Q_m$.

\subsection{Helgason's spheres}

Let $M$ be a compact irreducible symmetric space
and let $\kappa$ denote the maximum of the
sectional curvatures of $M$. By a theorem of
\'E. Cartan, the same dimensional maximal
totally geodesic flat submanifolds of $M$,
that is, the maximal tori, are all conjugate
under the largest connected Lie group
$I_0(M)$ of isometries.

S. Helgason (1927-- ) proved in 1966 an
analogous result for the maximum curvature
$\kappa$ as follows:

(1) A compact irreducible symmetric space
$M$ contains totally geodesic submanifolds of
maximum constant curvature $\kappa$.

(2) Any two such totally geodesic
submanifolds of the same dimension are
conjugate under
$I_0(M)$.

(3) The maximal dimension of such
submanifolds is $1+m(\bar\delta)$, where
$m(\bar\delta)$ is the multiplicity of the
highest restricted root. Also
$\kappa=||\bar\delta||^2$, where $||\;\; ||$
denotes length.

(4) If $M$ is a simply-connected
compact irreducible symmetric space, then the
closed geodesics in $M$ of minimal length are
permuted transitively by $I_0(M)$.

A maximal dimensional totally geodesic sphere
with maximal possible sectional curvature
$\kappa$ in a compact irreducible symmetric
space is known as a Helgason's sphere.

The stability of Helgason's spheres was
determined by  Y. Ohnita in 1987. He proved
that every Helgason's sphere in a compact
irreducible symmetric space is stable as
totally geodesic submanifolds. 

The Helgason sphere in a compact simple Lie
group is the compact simple 3-dimensional Lie
subgroup associated with the highest root.
In 1977 C. T. Dao showed that
Helgason's spheres in some compact classical
Lie groups are homologically volume
minimizing. H. Tasaki proved in 1985 that the
Helgason spheres in any compact simple Lie
groups are homologically volume minimizing
using the canonical 3-form
$\left<[X,Y],Z\right>$ divided by the length
of the highest root as a calibration and by
root systems. 

Le Khong Van showed in 1993 the volume
minimizing property of Helgason spheres in
all compact simply-connected irreducible
symmetric spaces. In 1995 H. Tasaki obtained
 estimates of the volume of Helgason's
spheres and of cut loci using the
generalized Poincar\'e formula proved by R.
Howard (1993) and established the inequality
of volume minimizing of the Helgason
spheres. He then proved that under certain
suitable conditions a
$k$-dimensional Helgason sphere of a compact
symmetric space is volume minimizing in the
class of submanifolds of dimension $k$ whose
inclusion maps are not null homotopic.
Tasaki's conditions hold automatically for
compact symmetric spaces of rank one, compact
Hermitian symmetric spaces, and quaternionic
Grassmann manifolds. 

The maximal dimension of totally geodesic
submanifolds of positive constant
sectional curvature  in an arbitrary compact
irreducible symmetric space have been 
completely determined in [Chen-Nagano 1978,
Nagano-Sumi 1991] by applying the
$(M_+,M_-)$-method. Such dimensions play an
important role in the study of totally
umbilical submanifolds, in particular, in
the study of extrinsic spheres, in locally
symmetric spaces. 

\subsection{Frankel's theorem}

T. J. Frankel (1961) proved that if $M$ is a
complete Riemannian manifold of positive
sectional curvature and $V$ and $W$ are
two compact totally geodesic submanifolds
with $\dim V+\dim W\geq \dim M$, then $V$ and
$W$ have a nonempty intersection. Applying
this result, Frankel showed in 1966 that if
$M$ is a complete Riemannian manifold of
strictly positive sectional curvature and if
$V$ is a compact totally geodesic submanifold
of $M$ with $\dim V\geq\frac12\dim M$, then
the homomorphism of fundamental groups:
$\pi_1(V)\to\pi_1(M)$ is surjective.

In 1961 Frankel also showed that if $M$ is a
complete K\"ahler manifold of positive
sectional curvature, then any two compact
K\"ahler submanifolds must intersect if their
dimension sum is at least that of $M$. S. I.
Goldberg and S. Kobayashi (1967) extended this
Frankel's result to complete K\"ahler manifolds
of positive bisectional curvature.

Frankel's theorems were extended in 1996 by K.
Kenmotsu and C. Xia to complete
Riemannian manifolds which has positive
$k$-Ricci curvature or to K\"ahler
manifold with partially positive bisectional
curvature. For instance, they proved that if
$M$ is a complete Riemannian manifold with
nonnegative
$k$-Ricci curvature and $V$ and $W$ are
complete totally geodesic submanifolds
satisfying (a) $V$ and $W$ are immersed as 
closed subsets, (b) one of $V$ and $W$ is
compact, (3)
$M$ has positive $k$-Ricci curvature either
at all points of $V$ or at all points of $W$,
and (4) $\dim V+\dim W\geq\dim M+k-1$, then
$V$ and $W$ must intersect.

In 1966 Frankel also proved that two
compact minimal hypersurfaces of a compact
Riemannian manifold with positive Ricci
curvature must intersect.

\vfill\eject

\section{Totally umbilical submanifolds}

A submanifold $N$ of a Riemannian manifold $M$ is
called totally umbilical if its second
fundamental form $h$ is proportional to its first
fundamental form $g$, that is,
$h(X,Y)=g(X,Y) H$ for vectors $X,Y$
tangent to $N$, where $ H$ is the
mean curvature vector. Total
umbilicity is a conformal invariant in the
sense that if $N$ is a totally umbilical
submanifold of a Riemannian manifold $M$, then
$N$ is also a totally umbilical submanifold of
$M$ endowed with another Riemannian metric
which is conformal equivalent to the original
Riemannian metric on $M$. Each connected
component of the fixed point set of a conformal
transformation on a Riemannian manifold is a
totally umbilical submanifold. From Riemannian
geometric point of views, totally umbilical
submanifolds are the simplest submanifolds
next to totally geodesic ones. 

 A totally umbilical
submanifold $N$ in a Riemannian manifold $M$ is
called an extrinsic sphere if its mean curvature
vector field is a nonzero  parallel
 normal vector field. 1-dimensional extrinsic spheres
in Riemannian manifolds are called circles.

Hyperspheres in  Euclidean space are the most
well-known examples of totally umbilical
submanifolds and also of extrinsic spheres.

Since every curve in a Riemannian manifold is
totally umbilical, we shall only consider
totally umbilical submanifolds of dimension
$\geq 2$.

Totally umbilical surfaces in $E^3$ are
open parts of planes and round spheres.
This result was first proved by J.
Meusnier in 1785 who showed that open
parts of planes and spheres are the
only surfaces in $E^3$ satisfying the
property that the curvature of the plane
sections through each point of the
surface are equal.

In 1954 J. A. Schouten proved that every
totally umbilical submanifold of dimension
$\geq 4$ in a conformally flat space is
conformally flat.

An $n$-dimensional submanifold $M$ of $E^m$
satisfies $S\geq n H^2$, with the equality
holding identically if and only if $M$ is
totally umbilical, where $S$ and $H^2$ denote
the squared norm of the second fundamental
form and the squared mean curvature
function.

\subsection{Totally umbilical submanifolds of real space forms}

Totally umbilical submanifolds in real space
forms have been completely classified. Totally
umbilical submanifolds of dimension $\geq 2$
in real space forms are either totally
geodesic submanifolds or extrinsic spheres.  An
$n$-dimensional  non-totally geodesic, totally
umbilical submanifold of a Euclidean
$m$-space  $ E^m$ is an ordinary
hypersphere which is contained in an affine
$(n+1)$-subspace of $ E^m$. 

An $n$-dimensional totally umbilical
submanifold of a Riemannian $m$-sphere
$S^m$ or of a hyperbolic $m$-space $H^m$ is
contained in an $(n+1)$-dimensional totally
geodesic submanifold as a totally
umbilical hypersurface. 

An $n$-dimensional non-totally geodesic
totally umbilical submanifold of a real
projective $m$-space $RP^m$ is contained in
a totally geodesic $RP^{n+1}$ of $RP^m$. Such a
submanifold is obtained from a totally
umbilical submanifold of a Riemannian
$m$-sphere via the two-fold Riemannian
covering map $\pi:S^m\to RP^m$.

\subsection{Totally umbilical submanifolds of complex space forms}

Totally umbilical submanifolds in other rank
one symmetric spaces are also  known.
Totally umbilical submanifolds in
$CP^m$ and in its non-compact dual  are
classified in [Chen-Ogiue 1974c].

Let $N$ be an $n$-dimensional, $(n\geq 2)$,
totally umbilical submanifold of a real
$2m$-dimensional K\"ahler manifold $M$ of
constant holomorphic sectional curvature $4c$,
$c\ne0$. Then $N$ is one of the following
submanifolds:

(1) a complex space form isometrically immersed
in $M$  as a totally geodesic complex
submanifold;

(2) a real space form isometrically immersed in
$M$ as a totally real, totally geodesic
submanifold; 

(3) a real space form isometrically immersed in
an $(n+1)$-dimensional totally real totally
geodesic submanifold of $M$ 
 as an extrinsic hypersphere.

\subsection{Totally umbilical submanifolds of quaternionic space forms}

Totally umbilical submanifolds in a
quaternionic projective space
$HP^m$ and in its non-compact dual have been
classified in [Chen 1978].

Let $N$ be an $n$-dimensional, $(n\geq 4)$,
totally umbilical submanifold of a real
$4m$-dimensional quaternionic space form $M$ of
constant quaternionic sectional curvature $4c$,
$c\ne0$. Then $N$ is one of the following
submanifolds:

(1) a quaternionic space form
isometrically immersed in 
$M$  as a totally geodesic quaternionic
submanifold;

(2) a complex space form isometrically immersed
in  $M$  as a totally geodesic, totally complex
submanifold;

(3) a real space form isometrically immersed in
$M$ as a totally real, totally geodesic
submanifold; 

(4) a real space form isometrically immersed
 in an $(n+1)$-dimensional totally real totally
geodesic submanifold of $M$
 as an extrinsic hypersphere.

\subsection{Totally umbilical submanifolds of the Cayley plane}

Totally umbilical submanifolds of the Cayley
plane $\Cal OP^2$ and of its noncompact dual
have also been classified [Chen 1977b,
Nikolaevskij 1994]:

(1) The maximum dimension of totally
umbilical submanifolds in the Cayley plane
$\Cal OP^2$ is 8.

(2) A maximal totally umbilical submanifold of
$\Cal OP^2$ is one of the following:

(2-a) a totally geodesic quaternionic
projective plane $HP^2$;

(2-b) a totally geodesic Cayley line $\Cal
OP^1=S^8$;

(2-b) an extrinsic hypersurface of a
totally geodesic Cayley line
$\Cal OP^1=S^8$; or

(2-c) a non-totally geodesic totally umbilical
submanifold of a totally geodesic quaternionic
projective plane $HP^2$.

The corresponding result also holds for
totally umbilical submanifolds in the
non-compact dual of the Cayley plane.

\subsection{Totally umbilical submanifolds in complex quadric}

Yu. A. Nikolaevskij (1991) classified totally
umbilical submanifolds in complex quadric
$Q_m=SO(m+2)/SO(m)\times SO(2)$.

 An $n$-dimensional $(n\geq 3)$ totally
umbilical submanifold $N$ of the complex
quadric $Q_m$ is one of the following:

(1) a totally geodesic submanifold;

(2) an extrinsic sphere;

(3) a totally umbilical submanifold with
nonzero and nonparallel mean curvature
vector.

The last case  can be
described in the following two ways: 

(3-a) $N$ is an umbilical hypersurface of
non-constant mean curvature lying in the totally
geodesic $S\sp{p}\times S\sp{1}\in Q_m$; or 

(3-b) $N$ is a diagonal of the product of two
small spheres lying in the totally geodesic
$S\sp{l+1}\times S\sp{l+1}\in Q_m$; moreover, 
the mean and sectional curvatures of $N$ are
both constant.

\subsection{Totally umbilical submanifolds of locally symmetric spaces}

T. Miyazawa and G. Chuman (1972)
studied totally umbilical submanifolds
of locally symmetric spaces and obtained the
following. 

(1) A totally
umbilical submanifold of a locally symmetric
space is locally symmetric if and only if
its mean curvature is constant.

(2) Let $N$ be a totally umbilical
submanifold of dimension $\geq 4$ in a locally
symmetric space. If the mean curvature is
nowhere zero, then $M$ is conformally flat.

A Riemannian manifold is called  reducible if it
is locally the Riemannian product of two
Riemannian manifolds of positive dimensions. 

The following results of Chen (1981a)
determines  reducible totally umbilical
submanifolds of locally symmetric spaces:

(1) Every reducible totally umbilical
submanifold of a locally symmetric space has
constant mean curvature.

(2) A reducible totally umbilical submanifold
of a locally symmetric space $M$ is one
of the following locally symmetric spaces:

(2.1) A totally geodesic submanifold,

(2.2) a locally Riemannian product of a curve
and a Riemannian space form of
constant curvature,

(2.3) a locally Riemannian product of
two Riemannian space forms of
constant sectional curvatures
$c$ and $-c$, $c\ne0$, respectively.
\smallskip

Irreducible totally
umbilical submanifolds in locally symmetric spaces
do not have constant mean curvature in general.
For irreducible totally umbilical
submanifolds with constant mean curvature, 
we have the following [Chen 1980a]: If $N$ is
an $n$-dimensional $(n\geq 2)$ irreducible
totally umbilical submanifold with constant
mean curvature in a locally symmetric space
$M$, then $N$ is either a totally geodesic
submanifold or a real space form. Furthermore,
$N$ is either a  totally geodesic submanifold
or an extrinsic sphere, unless
$\dim N<{1\over 2}\dim M$.

Chen and P. Verheyen (1983) studied totally
umbilical submanifolds in locally Hermitian
symmetric spaces and obtained the following:
Let $N$ be an $n$-dimensional $(n\geq 4)$ totally
umbilical submanifold of a locally
Hermitian symmetric space
$M$. If $n>\dim_C M$, then either

(1) $N$ is a totally geodesic submanifold, or

(2) $rk(M)>n$ and $N$ is an extrinsic sphere
in a maximal flat totally geodesic
submanifold of $M$.

 This result implies, in
particular, that if $N$ is an $n$-dimensional
$(n\geq 4)$ totally umbilical submanifold in a
locally Hermitian symmetric space
$M$ of compact or non-compact type and if 
$n>\dim_{ C} M$, then $N$ is totally
geodesic.

For totally umbilical submanifolds in symmetric
spaces, Yu. A. Nikolaevskij (1994) proved
that an $n$-dimensional $(n\geq 3)$ totally
umbilical submanifold in a globally
symmetric space $M$ is either totally
geodesic or totally umbilical and complete
in a totally geodesic submanifold of $M$
which is isometric to the product of some
Riemannian space forms. In particular, he
proved that if $N$ is an $n$-dimensional
totally umbilical submanifold in an
irreducible symmetric space of compact
type, then it is either totally geodesic or
totally umbilical in the totally geodesic
product $\bar M$ of flat torus and spheres. 

\subsection{Extrinsic spheres in locally symmetric spaces}

Extrinsic spheres in Riemannian manifolds can be
characterized as follows:
Let $N$ be an $n$-dimensional $(n\geq
2)$ submanifold of a Riemannian
manifold $M$. If, for some $r>0$, every
circle of radius $r$ in $N$ is a circle in $M$,
then $N$ is an extrinsic sphere in $M$.
Conversely, if $N$ is an extrinsic sphere in $M$,
then every circle in $N$ is a circle in $M$ 
[Nomizu-Yano 1974].

 Extrinsic spheres
in locally symmetric spaces were completely
classified by Chen (1977a) as follows:

Let $N$ be an $n$-dimensional $(n\geq 2)$
extrinsic sphere of a locally symmetric
space $M$. Then $N$ is an extrinsic
hypersphere of an $(n+1)$-dimensional
totally geodesic  submanifold
of $M$ with constant sectional curvature.

Conversely, every extrinsic sphere of dimension
$\geq 2$ in a locally symmetric space is obtained
in such way.
\smallskip

The above result implies, in particular, that

(1) every extrinsic sphere in a locally
symmetric space is a real space form, and

(2) real space forms are the only
locally symmetric spaces of dimension $\geq
3$ which admit an  extrinsic hypersphere.

An extrinsic sphere in a Riemannian manifold
is not necessary a Riemannian sphere in general.
In contrast, for extrinsic spheres in  K\"ahler
manifolds, we have
the following result of [Chen 1976a].
 
Every complete simply-connected
even-dimensional extrinsic sphere  in a
K\"ahler manifold is isometric to a
Riemannian sphere if it has flat normal
connection. 

This result is false if the
extrinsic sphere is odd-dimensional. 
In fact, there exist many complete 
odd-dimensional simply-connected extrinsic
spheres with flat normal connection  in 
K\"ahler manifolds which are not Riemannian
spheres, even not homotopy spheres 
[Chen 1981a].

In 1984 S. Yamaguchi, H. Nemoto and N.
Kawabata showed that if a complete,
connected and simply-connected extrinsic
sphere  in a K\"ahler manifold is not
isometric to an ordinary sphere, then it is
homothetic to either a Sasakian manifold or a
totally real submanifold. 

Extrinsic spheres in locally conformally
K\"ahler manifolds were treated in
[Drago\-mir-Ornea 1998].

\subsection{Totally umbilical hypersurfaces}

 O. Kowalski (1972) proved
that every totally umbilical hypersurface of
an Einstein manifold of dimension $\geq 3$
is either a totally geodesic hypersurface or
an extrinsic hypersphere.

Not every Riemannian manifold admits a
totally umbilical hypersurface. The
following result
 of Chen (1981a) determined all locally
symmetric spaces which admit a non-totally
geodesic totally umbilical hypersurface:

A locally symmetric space  admits a non-totally
geodesic, totally umbilical hypersurface if and
only if, locally, it is one of the following
spaces:

(1) a real space form;

(2) a  Riemannian product of a line
and a Riemannian space form;

(3) a  Riemannian product of two
real space forms of constant curvatures $c$
and $-c$, respectively.

Thus, if an irreducible locally
symmetric space admits a totally umbilical
hypersurface, then it is a real space form.

It follows from the above result that 
non-totally geodesic, totally umbilical
hypersurfaces of a locally symmetric space are
given locally by the fixed point sets of some
conformal mappings. More precisely,  if $N$ is
a non-totally geodesic totally umbilical
hypersurface of a locally symmetric space $M$,
then, for each point
$x\in N$, there is a neighborhood $U$ of $x$ in
$M$ and a conformal mapping $\phi$ of $U$ into $M$
such that $U\cap N$ lies in the fixed point set of
$\phi$. 

K. Tojo (1997b) studied Riemannian
homogeneous spaces which admit extrinsic
hyperspheres and proved that if a natural
reductive homogeneous Riemannian manifold
admits an extrinsic hypersphere, then it
must be a real space form. 

K. Tsukada (1996) proved that if $M=G/H$ is a
Riemannian homogeneous space such that the
identity component of $H$ acts irreducibly
on the tangent space and if $\dim M\geq 3$,
then $M$ admits no totally umbilical
hypersurfaces unless $M$ is a real space
form.

\vfill\eject

\section{Conformally flat submanifolds}

A Riemannian $n$-manifold $M$ is called
conformally flat if, at each point $x\in M$,
there is a neighborhood of $x$ in $M$ which
is conformal to the Euclidean $n$-space.  An
immersed submanifold $f:M\to E^m$ is called a
conformally flat submanifold if the
submanifold is conformally flat with
respect to the induced metric.

Since every Riemannian 2-manifold is
conformally flat, due to the existence of
local isothermal coordinate system, we only
consider conformally flat manifolds of
dimension greater than or equal to 3.

According to a
well-known result of H. Weyl (1918), a
Riemannian manifold of dimension $n\geq 4$ is
conformally flat if and only if its Weyl
conformal curvature tensor $W$ vanishes
identically. The Weyl conformal curvature
tensor $W$ vanishes identically for $n=3$. 

Define a tensor $L$ of type $(0,2)$ on a
Riemannian $n$-manifold by
$$L(X,Y)=-\left({1\over
{n-2}}\right)Ric(X,Y)+\left({\rho\over{2
(n-1)(n-2)}}\right)g(X,Y),\leqno (13.1)$$ where
$Ric$ is the Ricci tensor and
$\rho=\,\hbox{trace}\, Ric$. Let $D$ be the
 (0,3)-tensor defined by
$$D(X,Y,Z)=(\nabla_X L)(Y,Z)-(\nabla_Y
L)(X,Z). \leqno(13.2)$$

H. Weyl (1918) proved that the tensor $D$
vanishes identically for a conformally flat
manifold of dimension
$n\geq 4$; and a Riemannian
3-manifold is conformally flat if and only
if $D$ vanishes identically..

N. H. Kuiper (1949) proved that a compact
simply-connected conformally flat manifold of
dimension $\geq 2$ is conformally equivalent
to $S^n$.

\subsection{Conformally flat
hypersurfaces} .

\noindent{\bf 13.1.1. Quasi-umbilicity of
conformally flat hypersurfaces} 

The study of nonflat conformally flat
hypersurfaces of dimension $n\geq 4$ was
initiated by \'E. Cartan around  1918. He
proved that a hypersurface of dimension $\geq
4$ in Euclidean space is conformally flat if
and only if it is quasi-umbilical, that is,
it has a principal curvature with
multiplicity $\geq n-1$. 

The Codazzi equation implies that every
quasi-umbilical hypersurface of dimension 3
in a conformally flat $4$-manifold is
conformally flat. The converse of this is
not true in general. In fact, G. M. Lancaster
(1973) showed that there exist conformally flat
hypersurfaces in $ E^4$ which have three
distinct principal curvatures;  hence
there exist conformally flat hypersurfaces
in $ E^4$ which are not quasi-umbilical.

\noindent{\bf 13.1.2. Canal
hypersurfaces} 

A hypersurface of  Euclidean space is
called a canal hypersurface if it is the
envelope of one-parameter family of hyperspheres.
\'E. Cartan (1918) proved that a canal
hypersurface of dimension $n\geq 4$ in
Euclidean space is conformally flat. 

\noindent{\bf 13.1.3. Conformally flat
hypersurfaces as loci of spheres} 

A conformally flat hypersurface of dimension
$\geq 4$ in a real space form $R^{n+1}(c)$
is a locus of $(n-1)$-spheres, in the sense
that it is given by smooth gluing of some
$n$-dimensional submanifolds of $M$ (possibly
with boundary) such that each of the
submanifolds is foliated by totally umbilical
$(n-1)$-submanifolds of $R^{n+1}(c)$ (cf.
[Chen 1973b]). 

D. E. Blair (1975) proved that the
generalized catenoid and the hyperplanes are
the only conformally flat minimal
hypersurfaces in $ E^{n+1}$ with
$n\geq 4$. 

\noindent{\bf 13.1.4. Intrinsic properties
of conformally flat hypersurfaces} 

 An intrinsic characterization of 
conformally flat manifolds admitting
isometric immersions in real space forms as 
hypersurfaces was given by Chen and Yano.

A conformally flat manifold $M$ of dimension
$n\geq 4$ is called special if there exist
three functions $\alpha,\beta$ and $\gamma$
on $M$ such that the tensor $L$ defined by
(13.1) takes the form:
$$L=-{1\over
2}\left(k+\alpha^2\right)g-\alpha\beta\,\omega\otimes
\omega,\leqno(13.3)$$
for some constant $k$, where $\omega$ is a
unit 1-form satisfying $d\alpha=\gamma\omega$
on the open subset $U=\{x\in M:\beta(x)\ne
0\}$, where $L=-{1\over
2}\left(k+\alpha^2\right)g$ on $\{x\in
M:\beta(x)=0\}$.

For a special conformally flat space $M$, we
define a real number $i_M$ by
$$\aligned i_M=&\sup\{k\in \hbox{\bf R}:
L=-{{k+\alpha^2}\over
2}g-\alpha\beta\,\omega\otimes
\omega\;\;\\ &\hbox{for some functions }
\alpha,\;\beta \hbox{ on } M\},\endaligned$$
which is called the index of the special
conformally flat manifold. 

Chen and  Yano proved the following (cf. [Chen
1973b]):

(1) Every conformally flat hypersurface of dimension $n\geq 4$ in a real space form is special.

(2) Conversely, every simply-connected special conformally flat manifold of dimension $n\geq
4$ with index $i_M$ can be isometrically
immersed in every real space form of
curvature $k<i_M$ as a hypersurface and it
cannot be isometrically immersed in any real
space form of curvature $k>i_M$ as a
hypersurface.
\smallskip

In 1982, M. do Carmo and M. Dajczer 
showed that conformally flat manifolds are
the only Riemannian $n$-manifolds, $n\geq
4$, which can be isometrically immersed  as a
hypersurface in two real
space forms of different curvatures.

U. Pinkall (1988) proved that every
compact conformally flat hypersurface in $
E^{n+1},\,n\geq 4$ is conformally equivalent to
a classical Schottky manifold. Pinkall's
result improves a result
of [do Carmo-Dajczer-Mercuri 1985].

It is not known whether every classical
Schottky manifold admits a conformal
immersion into $ E^{n+1}$.

\noindent{\bf 13.1.5. Taut conformally flat
hypersurfaces} 

 Let $f:M\to E^m$ be an immersion and $p\in 
E^m$. Denote the function $$x\in M\to
|f(x)-p|^2$$ by $L_p$. Suppose $\phi:M\to
\hbox{\bf R}$ is a Morse function on 
manifold $M$. If for all real $r$,
$M_r=\phi^{-1}(-\infty,r]$ is compact, then
the Morse inequality $\mu_k\geq
\beta_k$ holds, where $\mu_k$ is the number of
critical points of index $k$ which $\phi$ has
on $M_r$, and $\beta_k$ is the $k$-th Betti
number of $M_r$ over any field $F$. 

The function $\phi$ is called a $T$-function if
there is a field $F$ such that the Morse
inequality is an equality for all $r$ and $k$.
An immersion $f:M\to E^m$ is said to be
taut if every function of the form $L_p,\,
p\in E^m$, is a $T$-function. 

T. E. Cecil and P. Ryan (1980) proved that a
complete conformally flat hypersurface of $
E^{n+1},\, n\geq 4$, is taut if and only if it
is one of the following:

(a)  a hyperplane or a round sphere;

(b)  a cylinder over a circle or round
$(n-1)$-sphere; \

(c)  a ring cyclide
(diffeomorphic to $S^1\times S^{n-1}$);

(d) a parabolic cyclide (diffeomorphic
to $S^1\times S^{n-1}$ with a point
removed). 

The proofs of all of the above results rely on
Cartan's condition of quasiumbilicity on 
conformally flat hypersurfaces.

\subsection{Conformally flat 
submanifolds}. 
\vskip.1in

\noindent{\bf 13.2.1. Totally quasi-umbilical
submanifolds} 

An $n$-dimen\-sional submanifold $M$ of a
Riemannian $(n+p)$-manifold $N$ is called
quasiumbilical (respectively, umbilical) with
respect to a  normal vector field $\xi$ if
the shape operator $A_\xi$ has an eigenvalue
with multiplicity $\geq n-1$ (respectively,
multiplicity $n$). In this case, $\xi$ is
called a quasiumbilical  (respectively,
umbilical) normal section of $M$. 

An $n$-dimensional submanifold $M$ of  a
Riemannian $(n+p)$-manifold
is called totally quasiumbilical if there
exist $p$ mutually orthogonal quasiumbilical
normal sections on $M$. 

A result of Chen and Yano (1972) states that a
totally quasiumbilical submanifold of
dimension $\geq 4$ in a conformally flat
manifold is conformally flat.

The property of being totally
quasiumbilicity is a conformal invariant, that
is, the property remains  under every
conformal change of the metric of the ambient
space. 

For conformally flat submanifolds of higher
codimension, Chen and L. Verstraelen (1977)
proved that an $n$-dimensional ($n\geq 4$)
conformally flat submanifold $M$ with flat
normal connection in a conformally flat
$(n+p)$-manifold is
totally quasiumbilical if  $p<n-2$. 

\"U. Lumiste and M. V\"aljas (1989) showed
that Chen-Verstraelen's result is sharp in
the sense that there exists a conformally
flat  submanifold $M$ of dimension $n\geq 4$
with flat normal connection in a conformally
flat $(2n-2)$-manifold which is not
totally quasiumbilical. 

 J. D. Moore and J. M. Morvan (1978) showed
that an $n$-dimensional conformally flat
submanifold $M$ in $ E^{n+p}$ is
totally quasiumbilical if  $n>7$ and $p<4$.

 A relationship between the notion
of quasiumbilicity and focal points of
conformally flat submanifolds was given in
[Morvan-Zafindratafa 1986].
\vskip.1in

\noindent{\bf  13.2.2. Submanifolds admitting an
umbilical normal section} 

A normal vector field $\xi$ of a submanifold
is said to be parallel if $D_X\xi=0$ for each
tangent vector $X$ of $M$, where $D$ denotes
the normal connection. The length of a
parallel normal vector field is constant.
A normal vector field
$\xi$ of $M$ is called nonparallel if, for
each $x\in M$, there is a tangent vector
$X\in T_xM$ such that $D_X\xi\ne 0$.

A submanifold $f:M\to E^{n+p}$ of a
Euclidean space admits an umbilical 
parallel unit normal section if and only if
$f$ is spherical, that is, $f(M)$ is
contained in a hypersphere of $ E^{n+p}$.
On the other hand, if a submanifold $M$ of
codimension 2 in a Euclidean space admits an
umbilical  nonparallel unit normal section
$\xi$, then $M$ must be quasiumbilical with
respect to each unit normal vector field
perpendicular to
$\xi$. Thus, $M$ is a totally quasiumbilical
submanifold. In particular, if $\dim M\geq
4$, then the submanifold is conformally flat
[Chen-Yano 1973].
\vskip.1in

\noindent{\bf  13.2.3. Normal curvature
tensor as a conformal invariant} 

Another general conformal property of
submanifolds is that the normal curvature
tensor of a submanifold of a Riemannian
manifold is a conformal invariant [Chen 1974].

If $M$ is a submanifold of a Riemannian
manifold $N$ endowed with Riemannian
metric
$\tilde g$ and if
$\tilde g^*=e^{2\rho}\tilde g$ is a conformal
change of the metric $\tilde g$, then the
normal curvature tensor $R^D$ of $M$ in
$(N,\tilde g)$ and  the normal
curvature tensor
$R{^*}^{D}$ of $M$ in $(N,\tilde g^*)$
satisfy $R^*{}^D(X,Y)=R^D(X,Y)$ for any
vectors $X,Y$ tangent to $M$. 
\vskip.1in

\noindent{\bf  13.2.4. $CW$-decomposition
of a conformally flat submanifold} 

Applying Morse theory, J. D. Moore (1977)
proved that  a compact $n$-dimensional
conformally flat submanifold of $ E^{n+p}$
possesses a $CW$-decomposition with no cells of
dimension $k$ with $p<k<n-p$.
\vskip.1in

\noindent{\bf  13.2.5. A sufficient
condition for a locus of spheres
to be conformally flat} 

Chen (1973b) proved that a
locus  $M^n$  of
$(n-1)$-spheres ($n\geq 4$) in an
$(n+p)$-dimensional real space form
$R^{n+p}(c)$ with arbitrary  $p$ is
conformally  flat if the unit normal vector
field of each leaf in
$M^n$ is a parallel vector field in  the normal
bundle of the leaf in $R^{n+p}(c)$.

\vskip.1in

\noindent{\bf  13.2.6. Conformally flat
manifolds and hypersurfaces of light
cone} 

A simply-connected
Riemannian manifold of dimension $n\geq 3$ is
conformally flat if and only if it can be
isometrically immersed as a hypersurface of
the light cone $$V^{n+1}=\{X\in
L^{n+2}:\left<X,X\right>=0, X\ne 0\},$$
where
$\left<\;,\right>$ is the semi-definite
metric on $V^{n+1}$ induced from the standard 
Lorentzian metric on the flat
$(n+2)$-dimensional Lorentzian space
$L^{n+2}$ [Brinkmann 1923, Asperti-Dajczer
1989]. 
\vskip.1in

\noindent{\bf  13.2.7. Conformally flat
submanifolds with constant index
of conformal nullity}

When $f:M^n\to  E^{n+p},\, p\leq n-3$, is
a conformally flat submanifold, at each point
$x\in M$, there is an umbilical subspace $\Cal
U(x)
\subset T_xM$ with $\dim {\Cal U}(x)\geq n-p$.
Hence, there is a unit vector $\eta\in
T^\perp_xM$ and a nonnegative number $\lambda$
so that the second fundamental form satisfies
$h(Z,X)=\lambda\left<Z,X\right>\eta$, for
each
$Z\in\Cal U(x)$ and each $X\in T_xM$. 

The umbilical distribution $\Cal U$ is
integrable on any open subset where the
dimension of
$\Cal U(x)$ is constant, which is
denoted by $\nu_f^c(x)$.
$\nu_f^c(x)$ is called the index of conformal
nullity. 

The leaves of the umbilical
distribution are extrinsic spheres in $M$,
hence they are totally umbilical submanifolds
with parallel mean curvature vector.

We  say that an isometric immersion
$F:N^{n+1}\to\tilde N^{n+p}$ extends an
isometric immersion $f:M^n\to\tilde N^{n+p}$
when there exists an isometric embedding of
$M^n$ into $N^{n+1}$ such that $F|_M=f$.

M. Dajczer and L. A. Florit (1996) proved the
following: Let $f:M^n\to  E^{n+p},\, n\geq
5,\, p\leq n-3$, be a simply-connected
conformally flat submanifold without flat
points. If $f$ has constant index of conformal
nullity $\nu^c_f$, say $\ell$, then there exist
an extension $F:N^{n+1}\to E^{n+p}$ of $f$
and an isometric immersion $G:N^{n+1}\to
L^{n+2}$ so that $M^n=G(N^{n+1})\cap V^{n+1}$.
Moreover, $F$ and $G$ carry a common
$(\ell+1)$-dimensional relative nullity
foliation.

A conformally flat submanifold $f:M^n\to
E^{n+2},\,n\geq 5$, is called generic when
its umbilical direction $\eta\in T^\perp M$
possesses everywhere a nonzero principal
curvature $\lambda$ of multiplicity $n-2$.  An
immersion
$f:M^n\to E^{n+2}$ is called a composition
if there exist an open subset $U\subset 
E^{n+1}$ and isometric immersions $\tilde
f:M^n\to U$ and $H:U\to E^{n+2}$ such that
$f=H\circ\tilde f$. 

Dajczer and Florit
(1996)  proved that any conformally flat
submanifold $f:M^n\to
E^{n+2},\,n\geq 5$, without flat points is
locally along an open dense subset either
generic or a composition.

\vskip.1in

\noindent{\bf  13.2.8. A non-immersion
theorem} 

H. Rademacher (1988) proved that if there exists a constant $c$ such that the Ricci
curvature of a compact conformally flat
$n$-manifold $M$, $n\geq 4$, satisfies $$-c\sp 2\leq{\rm Ric}(X)\leq-\left(
{{n-1}\over{2n-3}}\right)
c^2$$  for any unit vector $X$ tangent to $M$, then $M$ cannot
be conformally immersed in $S^{2n-2}$.

\vfill\eject

\section{Submanifolds with parallel mean
curvature vector}

A submanifold of a Riemannian manifold is said to
have parallel mean curvature vector  if the mean
curvature vector field $H$ is parallel as a
section of the normal bundle.
Trivially, every minimal submanifold of a
Riemannian manifold  has parallel mean
curvature vector.  A minimal
submanifold of a hypersphere in Euclidean space
has nonzero parallel mean curvature vector
when it  is considered as a submanifold of
the ambient Euclidean space. Furthermore, a
hypersurface of a Riemannian manifold  has
parallel mean curvature vector if and only
if it has constant mean curvature. 
By introducing the notion of twisted
product, Chen (1981a) showed that every
Riemannian manifold can be embedded in some
twisted product Riemannian manifold as a
submanifold with nonzero parallel mean
curvature vector. H. Reckziegel (1974)
proved if $M$ is a compact submanifold of a
manifold $N$ such that $TN|_M$ has a
metric  $\gamma$ and $\eta$ is a nonzero
normal vector field of constant length in
$(TM)^\perp$, then $\gamma$ can be
extended to a Riemannian metric on $N$
such that $M$ is an extrinsic sphere with
parallel mean curvature vector $\eta$.

It was S. S. Chern who first suggested in
the mid 1960s that the notion of parallel
mean curvature vector as the natural
extension of constant mean curvature for
hypersurfaces.

\subsection{ Gauss map and mean
curvature vector} 

For an isometric immersion $f: M\to 
E^m$  of  an oriented $n$-dimensional
Riemannian manifold into Euclidean m-space,
the Gauss map
$$G : M\to G^R(m-n,m)$$ of $f$ is a smooth
map which carries a point $x\in M$ into the
oriented $(m-n)$-plane in $ E^m$, which is
obtained from the parallel translation of
the normal space of $M$ at $x$ in $ E^m$,  
where $G^R(m-n,m)$ denotes the  Grassmannian
manifold consisting of oriented 
$(m-n)$-planes in $E^m$.

E. A. Ruh and J. Vilms (1970)
characterized submanifolds of Euclidean space with
parallel mean curvature vector  as follows: A
submanifold $M$ of a Euclidean $m$-space
$ E^m$ has parallel mean curvature vector
if and only if its Gauss map $G$ is
harmonic.

\subsection{Riemann sphere
with parallel mean curvature
vector} 

One-dimensional submanifolds with parallel
mean curvature vector are nothing but
geodesics and circles. 

D. Ferus (1971b) and E. A. Ruh (1971)
determined completely closed surfaces of
genus zero with parallel mean curvature
vector in Euclidean space:

 Let $M$ be a closed oriented 
surface of genus zero in $E^m$. If $M$
has parallel mean curvature vector,
then $M$ is contained in a hypersphere
of $ E^m$ as a minimal surface.

\subsection{Surfaces with parallel mean curvature
vector} 

Surfaces in $E^4$ with parallel mean
curvature vector were classified by D.
Hoffman in his doctoral thesis (Stanford
University 1971). 

The complete classification of 
surfaces in Euclidean
$n$-space, $n\geq 4$, with parallel mean
curvature vector was obtained by Chen
(1972) and, independently by Yau (1974).

 A surface $M$ of  Euclidean $m$-space
$ E^m$ has parallel mean curvature vector
if and only if  it is one of the following
surfaces:

(1)a minimal surface of $ E^m$; 

(2)  a minimal surface of a hypersphere of $ E^m$;

(3) a surface of $ E^3$ with constant mean curvature;

(4) a surface of constant mean curvature
lying in a hypersphere of an affine
4-subspace of $ E^m$.

Similar results hold for surfaces with parallel
mean curvature vector in spheres and
in real hyperbolic spaces as well (cf.
[Chen 1973b]).

For surfaces of constant Gaussian
curvature, Chen and G. D. Ludden (1972) and
D. Hoffman (1973) proved that minimal
surfaces of a small hypersphere, open
pieces of the product of two plane
circles, and open pieces of a circular
cylinder are the only non-minimal
surfaces in Euclidean space with
parallel mean curvature vector and with
constant Gaussian curvature.

For  a compact surface $M$ with positive
constant Gaussian curvature, 
 K. Enomoto (1985) proved that if 
$f:M\to E^m$ is an isometric embedding
with constant  mean curvature and flat
normal  connection, then $f(M)$ is a round
sphere in an affine 3-subspace of $ E^m$.
 Enomoto also proved that if $f:M\to
E^{n+2}$ is an isometric embedding of a compact
Riemannian $n$-manifold, $ n\geq 4$, with
positive constant sectional curvature and
with constant mean curvature, then $f(M)$ is
a round $n$-sphere in a hyperplane of $
E^{n+2}$.

For codimension two using the method of
equivariant differential geometry, W. T.
Hsiang, W. Y. Hsiang and I. Sterling
(1985)  proved the following: 

(a) There exist infinitely many
codimension two embeddings of distinct
knot types of $S^{4k+1}$ into
$S^{4k+3}(1)$ with parallel mean
curvature vector of arbitrarily small
constant length. 

(b) There exist infinitely many
codimension two embeddings of distinct
knot types of the Kervaire exotic sphere 
$\Sigma^{4k+1}_0$ into $S^{4k+3}(1)$
with parallel mean curvature vector having
length of arbitrarily small constant
value. 

(c) There exist infinitely many constant
mean curvature embeddings of
$(4k-1)$-dimensional generalized lens spaces into
$S^{4k+1}(1)$. 

It remains as an open problem to 
completely classify submanifolds
of dimension $\geq 3$ with parallel mean
curvature vector in real space forms.

\subsection{Surfaces with parallel
normalized mean curvature vector} 

Chen (1980b) defined a submanifold  in a
Riemannian manifold to have parallel
normalized mean curvature vector field if
there exists a unit parallel vector field
$\xi$ which is parallel to the mean curvature
vector field $H$, that is, $H=\alpha \xi$ for
some unit parallel normal vector field $\xi$.

Submanifolds with nonzero parallel mean
curvature vector field also have parallel
normalized mean curvature vector field.  The
condition to have parallel
normalized mean curvature vector field is much
weaker than the condition to have parallel
mean curvature vector field. For instance,
every hypersurface in a Riemannian manifold
always has parallel normalized mean curvature
vector field.

For surfaces with parallel normalized
mean curvature vector field, we have the
following results from [Chen 1980b]: 

(1) Let $M$ be a Riemann sphere immersed in a
Euclidean $m$-space $ E^m$.  Then $M$
has  parallel normalized mean curvature
vector field if and only if either

(1-a) $M$ is immersed in a  hypersphere of
$ E^m$ as a minimal surface, or 

(1-b) $M$ is immersed in
an affine 3-subspace of $E^m$.

(2) A surface $M$ of class $C^\omega$ in a
Euclidean $m$-space $ E^m$ has parallel
normalized mean curvature vector field if and
only if $M$ is one of the following surfaces:

(2-a) a minimal surface of a hyperplane of
$ E^m$,

(2-b) a  surface in an affine 4-subspace of $
E^m$ with  parallel normalized mean curvature
vector.

\smallskip
Every surface in a Euclidean 3-space has
parallel normalized mean curvature
vector field. Moreover, there exist abundant
examples of surfaces which lie fully in a
Euclidean 4-space with parallel normalized
mean curvature vector field, but  not with
parallel mean curvature vector field.

\subsection{ Submanifolds
satisfying additional conditions} 

It is a classical theorem of Liebmann that the
only closed convex surfaces in Euclidean
3-space having constant mean curvature are 
round spheres. B. Smyth (1973) extended
Liebmann's result to the following:
 Let $M$ be a compact $n$-dimensional
submanifold with nonnegative sectional
curvature in Euclidean $m$-space. If  $M$
has parallel mean curvature vector, then
$M$ is a product submanifold
$M_1\times\cdots\times M_k$, where
each $M_i$ is a minimal submanifold
in a hypersphere of an affine
subspace of $ E^m$.

Further, K. Yano (1912--1993) and S.
Ishihara (1922-- ) in 1971 and J. Erbacher
in 1972 extended Liebmann's result to the
following: Let
$M$ be an $n$-dimensional submanifold in
Euclidean $m$-space with nonnegative sectional
curvature. Suppose that the
 mean curvature vector  is parallel in the
normal bundle and the normal connection is flat.
If $M$ is either compact or has constant
scalar curvature, then  
$M$ is the standard product immersion of the
product 
$S^{n_1}(r_1)\times\cdots\times S^{n_k}(r_k)$ of
some spheres.

Recently, Y. Zheng (1997) proved the following:
Let $M$ be a compact
orientable submanifold with constant scalar
curvature and with nonnegative sectional
curvature immersed in  a real space form of
constant sectional curvature $c$. Suppose
that $M$ has flat normal connection. If the
normalized scalar curvature of $M$ is
greater than $c$, then $M$ is either
totally umbilical or locally the Riemannian
product of several totally umbilical
constantly curved submanifolds. 

For complete submanifolds $M^n$ of dimension
$\geq 3$ in Euclidean space,   Y. B. Shen
(1985) proved the following: Let $M^n$
$(n\ge 3)$ be a  complete submanifold in the
Euclidean space $ E^{m}$ with parallel
 mean curvature vector. If the squared mean
curvature $H^2$ and the squared length $S$
of the second fundamental form of
$M$ satisfies
$$(n-1)S\leq n^2H^2,$$ then $M^n$
is an $n$-plane, an $n$-sphere $S^n$, or a
circular cylinder $S^{n-1}\times E^1$.
This extended some results of
[Chen-Okumura 1973, Okumura 1973].

G. Chen and X. Zou (1995) studied compact
submanifolds of spheres with parallel mean
curvature vector and proved the following:

 Let $M$ be an $n$-dimensional compact
submanifold with nonzero parallel mean
curvature vector in the unit $(n+p)$-sphere.
Then 

(1) $M$ is totally geodesic, if one of the
following two conditions hold:
$$S\leq \min\left\{\frac23
n,\frac{2n}{1+\sqrt{\frac n 2}}\right\},\quad
p\geq 2\quad
\hbox{and}\quad n\ne 8;$$
$$S\leq \min\left\{\frac n{2-\frac 1{p-1}}
,\frac{2n}{1+\sqrt{\frac n 2}}\right\},\quad
p\geq 1\quad
\hbox{and}\quad (n,p)\ne (8,3);$$

(2) $M$ is totally umbilical, if $2\leq n\leq
7,\; p\geq 2,$ and $S\leq \frac23 n.$

\subsection{ Homogeneous
submanifolds  with parallel mean 
curvature vector}

C. Olmos (1994,1995) studied homogeneous
submanifolds of Euclidean space and proved
the following.

(a) If $M$ is a compact homogeneous
submanifold of a Euclidean space with
 parallel mean curvature vector which is not
minimal in a sphere, then $M$ is an orbit of
the isotropy representation of a simple symmetric
space; 

 (b) A homogeneous irreducible
submanifold of Euclidean space with parallel mean 
curvature vector is either minimal, or minimal
in a sphere, or an orbit of the isotropy 
representation of a simple
Riemannian symmetric space.

\vfill\eject

\section{ K\"ahler
submanifolds of K\"ahler manifolds}

 According to the behavior of the
tangent bundle of a submanifold with respect to the
action of the almost complex structure $J$ of 
the ambient
manifold, there are several typical classes of
submanifolds, namely,  K\"ahler submanifolds, 
totally real submanifolds, $CR$-submanifolds and
slant submanifolds.

In this section, the dimensions of  complex
manifolds always mean the complex dimensions,
unless mentioned otherwise.

The theory of submanifolds of a  K\"ahler
manifold began as a separate area of study
in the last century with the investigation
of algebraic curves and algebraic surfaces
in classical algebraic geometry. The study
of complex submanifolds of  K\"ahler
manifolds from differential geometrical
points of view (that is, with emphasis on
the Riemannian metric) was initiated by E.
Calabi in the early of 1950's. 

\subsection{Basic properties of
K\"ahler submanifolds} 

A submanifold of a complex manifold is called a
complex submanifold  if each of its tangent 
spaces is invariant under the almost complex
structure of the ambient manifold. 
A complex submanifold of a
K\"ahler manifold is itself a K\"ahler manifold with
respect to its induced metric. By a  K\"ahler
submanifold we mean a complex submanifold with
the induced K\"ahler structure. 

It was proved by Calabi (1953) that K\"ahler
submanifolds of K\"ahler manifolds always
have rigidity. Thus, for any two full
K\"ahler immersions $f$ and $f'$ of the same
K\"ahler manifold $M$ into $CP^m$ and into
$CP^N$, respectively, we have $m=N$ and,
moreover, there exists a unique holomorphic
isometry $\Psi$ of $CP^m$ onto itself such
that $\Psi\circ f=f'$.

The second fundamental form of a K\"ahler
submanifold $M$ of a K\"ahler manifold with
the almost complex structure $J$ satisfies
$$h(JX,Y)=h(X,JY)=Jh(X,Y),$$  for $X,Y$
tangent to $M$. From this it follows that
K\"ahler submanifolds of K\"ahler manifolds
are always minimal. 

Compact K\"ahler submanifolds of K\"ahler
manifolds are also stable and  have
the property of being absolutely volume
minimizing inside the homology class.
Moreover, a compact K\"ahler submanifold $M$
of a K\"ahler manifold $\tilde M$ can never
be homologous to zero, that is, there
exists no submanifold $M'$ of
$\tilde M$ such that $M$ is the boundary of $M'$
[Wirtinger 1936]. 

A K\"ahler submanifold of a K\"ahler manifold is said
to be of degree $k$ if the pure part of the
$(k-1)$-st covariant derivative of the second
fundamental form is identically zero but the pure
part of the $(k-2)$-nd one is not identically zero.
In particular, degree 1 is nothing but
totally geodesic and  degree 2 is equivalent
to parallel second fundamental form but not
totally geodesic.

Let $CP^m(c)$ denote the complex projective
$m$-space equipped with the Fubini-Study metric of
constant holomorphic sectional curvature $c$
and  $M$ be an $n$-dimen\-sional K\"ahler
submanifold of $CP^m(c)$.  

Denote by $g,K, H, Ric,\rho$ and $h$ the
metric tensor, the sectional curvature, the
holomorphic sectional curvature, the Ricci
tensor, the scalar curvature and the second
fundamental form of $M$, where
$\rho=\sum_{i\ne j} K(e_i\wedge e_j)$
and $\{e_1,\ldots,e_n\}$ is an orthonormal
frame on $M$.

It follows from the equation of Gauss that an
$n$-dimensional K\"ahler submanifold of
a K\"ahler $m$-manifold $\tilde M^m(c)$ of
constant holomorphic sectional curvature $c$
satisfies the following curvature properties
in general:

(1) $H\leq c$, with equality holding
identically if and only if $M$ is a totally
geodesic K\"ahler submanifold.

(2) $Ric\leq {c\over 2}(n+1)g$, with
equality holding identically if and only if $M$
is  a totally geodesic K\"ahler submanifold.

(3) $\rho\leq n(n+1)c$, with equality holding
identically if and only if $M$ is  a totally
geodesic K\"ahler submanifold.

Let $z_0,z_1,\ldots,z_{n+1}$ be a
homogeneous coordinate system of
$CP^{n+1}(c)$ and $Q_n$ be the complex
quadric hypersurface of $CP^{n+1}(c)$
defined by
$$Q_n=\{(z_0,z_1,\ldots,z_{n+1})\in
CP^{n+1}(c):
\sum z_i^2=0\}.$$
Then $Q_n$ is complex analytically
isometric to the compact Hermitian symmetric
space
$SO(n+2)/SO(n)\times SO(2)$. 

With respect to the induced K\"ahler metric $g$,
$Q_n$ satisfies the following:

(1) $0\leq K\leq c$ for $n\geq 2$ and
$K={c\over 2}$ for $n=1$.

(2) ${c\over 2}\leq H\leq c$ for $n\geq 2$
and $H={c\over 2}$ for $n=1$.

(3) $Ric ={n\over 2}cg$.

(4) $\rho=n^2c$.
\medskip

P. B. Kronheimer and T. S. Mrowka
(1994) proved the Thom conjecture concerning
the genus of embedded surfaces in $CP^2$,
namely, they proved that if
$C$ is a smooth holomorphic curve in
$CP^2$, and $C'$ is a smoothly embedded
oriented 2-manifold representing the same
homology class as $C$, then the genus of
$C'$ satisfies $g(C')\geq g(C)$. The proof
of this result uses Seiberg-Witten's
invariant on 4-manifolds.

\subsection{Complex space forms and Chern
classes} 

A K\"ahler manifold is called a complex space
form if it has constant holomorphic
sectional  curvature. The universal covering
of a complete complex space form $\tilde
M^n(c)$ is the complex projective $n$-space
$CP^n(c)$, the complex Euclidean $n$-space
$ C^n$, or the complex hyperbolic space
$CH^n(c)$, according to $c>0,\,c=0$, or $c<0$.

 Complex space forms can be
characterized in terms of the first and the
second Chern classes. In fact, B. Y. Chen
and K. Ogiue (1975) proved the following sharp
inequality between the first and the second
Chern numbers for Einstein-K\"ahler
manifolds. They also applied their
inequality to characterize complex space
forms:

 Let $M$ be an $n$-dimensional
compact Einstein-K\"ahler manifold. Then the
first and the second Chern classes of $M$
satisfy
$$\varepsilon^{n}\int_M 2(n+1)c_1^{n-2}c_2\geq
\varepsilon^{n}
\int_M nc_1^n\quad (\varepsilon \;\hbox{the sign
of } \rho),\leqno(15.1)$$ with the equality
holding if and  only if $M$ is either a
complex space form or a Ricci-flat K\"ahler
manifold. 
\smallskip
Chen-Ogiue's inequality was extended by M.
L\"ubke in 1982 to Einstein-Hermitian vector
bundle over compact K\"ahler manifolds in the
sense of S. Kobayashi (1932-- )
(cf. [Kobayashi 1987]).

In this respect, we mention that T. Aubin
(1976) proved that if $M$ is a compact
 K\"ahler manifold with $c_1<0$ (that is,
$c_1$ is represented by a negative definite
real
$(1,1)$-form), then there exists a unique
Einstein-K\"ahler metric on $M$ whose K\"ahler
form is cohomologous to the K\"ahler form of
the initial given metric.
Consequently, by combining these two results,
it follows that every compact K\"ahler
manifold with $c_1<0$ satisfies inequality
(15.1), with the equality holding if and only
if $M$ is covered by the complex hyperbolic
$n$-space. 

S. T. Yau (1977) proved that if $M$ is a
compact K\"ahler manifold with $c_1=0$, then
it admits a Ricci-flat K\"ahler metric. 

\subsection{ K\"ahler immersions 
of complex space forms in complex space forms} 

M. Umehara (1987a) studied K\"ahler immersions
between complex space forms and obtained the
following.

(1) A K\"ahler submanifold 
of a complex Euclidean space cannot be a
K\"ahler submanifold of any complex hyperbolic
space; 

(2) A K\"ahler submanifold of a complex
Euclidean space  cannot be a K\"ahler
submanifold of any  complex projective
space, and 

(3) A K\"ahler submanifold of a complex
hyperbolic space cannot be a K\"ahler
submanifold of any complex projective space.

K\"ahler immersions of complex space forms
in complex space forms are completely classified by
H. Nakagawa and K. Ogiue in 1976 as follows.

Let $M^n(c)$ be an $n$-dimensional complex space
form isometrically immersed  in an
$m$-dimensional complex space form $\tilde
M^m(\bar c)$ as a K\"ahler submanifold such
that the immersion is full. Then

(1) if $\bar c\leq 0$, then the
immersion is totally geodesic;

(2) if $\bar c>0$, then $\bar c=\mu c$ and
$m={{n+\mu}\choose \mu}-1$ for some positive
integer $\mu$. Moreover, either  the immersion
is totally geodesic or locally the immersion
is given by one of the Veronese embeddings.

This result is due to Calabi (1953)
when both $M^n(c)$ and $\tilde M^m(\bar c)$ are
complete simply-connected complex space forms.

An immersion $f\:M\to\tilde M$ between
Riemannian manifolds is called proper
$d$-planar geodesic if every geodesic in $M$
is mapped into a $d$-dimensional totally
geodesic subspace of $\tilde M$, but not into
any $d-1$-dimensional totally
geodesic subspace of $\tilde M$.

J. S. Pak and K. Sakamoto (1986,1988) proved
that and $f:M^n\to CP^m$ a proper $d$-planar
geodesic K\"ahler immersion from a K\"ahler
manifold into $CP^m$, $d$ odd or $d\in
\{2,4\}$, then $f$ is equivalent to the $d$-th
Veronese embedding of $CP^n$ into $CP^m$. 

\subsection{Einstein-K\"ahler
submanifolds and K\"ahler submanifolds 
$\tilde M$ satisfying $Ric(X,Y)={\widetilde
Ric}(X,Y)$}

For complex hypersurfaces of complex space
forms we have  the
following:

(1) Let $M$ be a K\"ahler hypersurface of
an $(n+1)$-dimensional complex space form
$\tilde M^{n+1}(c)$. If $n\geq 2$ and  $M$ is
Einstein, then either $M$ is totally geodesic or
$Ric={n\over 2}cg$. 

The latter case occurs only
when $c>0$.  Moreover, the immersion is
rigid [Smyth 1968, Chern 1967].

(2)   Let $M$ be a compact K\"ahler
hypersurface embedded in $CP^{n+1}$. If
$M$ has constant scalar curvature, then $M$ is
either totally geodesic in $CP^{n+1}$ or 
holomorphically isometric to $Q_n$ in
$CP^{n+1}$. [Kobayashi 1967a].

(3) Let $M$ be a  K\"ahler hypersurface of
an $(n+1)$-dimensional complex space form
$\tilde M^{n+1}(c)$. If the Ricci tensor of $M$
is parallel, then $M$ is an Einstein space
 [Takahashi 1967]. 

J. Hano (1975) proved that, besides linear
subspaces, $Q_n$ is the only
Einstein-K\"ahler submanifold of a complex
projective space which is a complete
intersection.

M. Umehara (1987b) proved that  every
Einstein-K\"ahler submanifold of $C^m$ or
$CH^m$ is totally geodesic.

B. Smyth (1968) proved that the normal
connection of a K\"ahler hypersurface $M^n$ in a
K\"ahler manifold $\tilde M^{n+1}$ is flat if and
only if the Ricci tensors of $M^n$ and $\tilde
M^{n+1}$ satisfy
$Ric(X,Y)={\widetilde Ric}(X,Y)$ for $X,Y$
tangent to $M^n$. 

For general K\"ahler submanifolds  
Chen and Lue (1975b) proved the
following.

Let $M$ be a compact K\"ahler submanifold
of a compact K\"ahler manifold $\tilde M$. Then

(1) if the normal connection is flat, the
Ricci tensors of $M$ and $\tilde M$
satisfy
$Ric(X,Y)={\widetilde Ric}(X,Y)$ for $X,Y$
tangent to $M$;

(2) if $Ric(X,Y)={\widetilde Ric}(X,Y)$ for any
$X,Y$ tangent to $M$, then the first Chern
class of the normal bundle  is trivial,
that is, $c_1(T^\perp M)=0$;

(3) if $\tilde M$ is flat, then the first
Chern class of the normal bundle is trivial
if and only if the normal connection is flat.

\subsection{ Ogiue's conjectures and curvature pinching}

An $n$-dimensional complex projective space
of constant holomorphic sectional
curvature $c$ can be analytically
isometrically embedded into an
$\left({{n+\mu} \choose
\mu}-1\right)$-dimen\-sional complex
projective space of constant holomorphic
sectional curvature $\mu c$. Such an
embedding is given by all homogeneous
monomials of degree $\mu$ in homogeneous
coordinates, which is called  the $\mu$-th
Veronese embedding of
$CP^n(c)$. The degree of the $\mu$-th Veronese
embedding is $\mu$.

The Veronese embeddings were
characterized  by A. Ros (1986) in terms of
holomorphic sectional curvature in the
following theorem:
 If a
compact $n$-dimensional K\"ahler submanifold $M$
immersed in $CP^m(c)$ satisfies 
$${c\over{\mu+1}}< H\leq
{c\over\mu},$$ then $M=CP^n({c\over\mu})$
and the immersion is given by the $\mu$-th
Veronese embedding.

K\"ahler submanifolds of degree $\leq 2$ are
characterized by Ros (1985b) as
follows: If a compact K\"ahler submanifold $M$
immersed in $CP^m(c)$ satisfies $ H\geq
1/2$, then degree $\leq 2$; and, moreover, the
K\"ahler submanifold is congruent to one of the
following six K\"ahler manifolds: $$
CP^n(c),\;\; CP^{n}\left({c\over 2}\right),\;
Q_n=SO(n+2)/SO(n)\times SO(2),$$ $$
SU(r+2)/S(U(r)\times U(2)),\; r\geq 3,\;\;
SO(10)/U(5),$$ $$ E_6/Spin(10)\times
SO(2).$$

A. Ros and L. Verstraelen (1984) and  Liao
(1988) characterized the second Veronese
embedding in terms of
sectional curvature in the following theorem.

If a compact $n$-dimensional
$(n\geq 2)$ K\"ahler submanifold
$M$ immersed in $CP^m(c)$ satisfies $K\geq
1/8$, then either $M$ is totally geodesic or
$M=CP^n({c\over 2})$ and the immersion is the
second Veronese embedding.

A. Ros (1986) characterized all Veronese
embeddings by curvature pinching as follows:
If a compact
$n$-dimen\-sional,
$n\leq 2$, K\"ahler submanifold
$M$ immersed in $CP^m(c)$ satisfies $${c\over
{4(\mu+1)}}\leq K\leq {c\over\mu},$$ then either
$M=CP({c\over\mu})$ and the immersion is given by
the $\mu$-th Veronese embedding
or $M=CP^n({c\over{\mu+1}})$ and the immersion is
given by the $(\mu+1)$-st Veronese embedding.

J. H. Cheng (1981) and R. J. Liao (1988)
characterized complex quadric hypersurface  in
terms of its scalar curvature:

If a compact $n$-dimensional K\"ahler
submanifold $M$ immersed in $CP^m(c)$ satisfies
$\rho\geq n^2$, then either $M$ is totally
geodesic or $M=Q_n$.

The above results gave affirmative
answers to some of Ogiue's conjectures.

Totally geodesic K\"ahler submanifolds of
$CP^m(c)$ are  characterized in terms of
Ricci curvature by K. Ogiue (1972a): If a
compact $n$-dimensional K\"ahler
submanifold
$M$ immersed in $CP^m(c)$ satisfies $Ric >{n\over
2}c$, then  $M$ is totally geodesic.

 Chen and  Ogiue (1974a) proved
that if an $n$-dimensional
 K\"ahler submanifold $M$ of $CP^m(c)$ satisfies
$Ric={n\over 2}cg$, then $M$ is an open piece
of  $Q_n$ which is embedded in some totally
geodesic $CP^{n+1}(c)$ in $CP^m(c)$.

For compact K\"ahler submanifolds with
positive sectional curvature, Y. B. Shen
(1995) proved the following.

  Let $M^n$ be an $n$-dimensional, $n\geq 2,$
compact K\"ahler submanifold immersed in $C
P^{n+p}$ with $p<n$. Then $M^n$ has nonnegative
sectional curvature  if and only if
$M^n$ is one of the following:

(1) a totally geodesic K\"ahler submanifold in
$CP^{n+p}$; or 

(2) an embedded submanifold
congruent to the standard full embedding of

(2-a)  the complex quadric
$Q_n$ as a hypersurface; or

(2-b) of $CP^{n-1}\times CP^1$ with
codimension $n-1$; or

(2-c)  of $U(5)/U(3)\times U(2)$ of
codimension 3; or

(2-d)  of $ U(6)/
U(4)\times  U(2)$ with codimension 6; or

(2-e)  of ${SO}(10)/ U(5)$ with codimension
5, or

(2-f)  of $E_6/{ Spin}(10)\times T$ with
codimension 10.

F. Zheng (1996) studied K\"ahler submanifolds
of complex Euclidean spaces and proved
that if $M$ is an
$n$-dimensional K\"ahler submanifold 
with   nonpositive sectional curvature in
$C^{n+r}$ with $r\leq n$, then 

(a) if $r<n$,  $M$ is ruled, that is, 
there exists a holomorphic
foliation ${\Cal F}$ on an open dense 
subset $U$ of $M$ with each leaf of ${\Cal F}$
totally geodesic in $ C^{n+r}$;

(b) if $r=n$, then either $M$ is ruled,
or $M$ is locally holomorphically
isometric to the product of $n$ complex
plane curves. 

The dimension bound is sharp, as there exist 
examples of negatively curved submanifolds
$M^n$ in $C^{2n+1}$ which are not ruled or
product manifolds. 

\subsection{Segre embedding}

Let $(z_0^i,\ldots,z_{N_i}^i)$ $(1\leq i\leq
s)$ denote the homogeneous coordinates of
$CP^{N_i}$. Define a map:
$$S_{N_1\cdots N_s}:CP^{N_1}\times\cdots\times
CP^{N_s}\to CP^N,\quad N=\prod_{i=1}^s (N_i+1)-1,$$
which maps a point
$((z_0^1,\ldots,z_{N_1}^1),\ldots,
(z_0^s,\ldots,z_{N_s}^s))$ of the product
K\"ahler manifold
$CP^{N_1}\times\cdots\times CP^{N_s}$ to
the point $(z^1_{i_1}\cdots
z^s_{i_j})_{1\leq i_1\leq
N_1,\ldots,1\leq i_s\leq N_s}$ in
$CP^N$. The map
$S_{N_1\cdots N_s}$ is a K\"ahler embedding
which is called the Segre embedding [Segre
1891].

Concerning Segre embedding, B. Y. Chen
(1981b) and B. Y. Chen and W. E. Kuan
(1985) proved the following:

 Let $M_1,\ldots,M_s$ be K\"ahler
manifolds of dimensions $N_1,\ldots,$ $N_s$,
respectively. Then every K\"ahler immersion
$$f:M_1\times\cdots\times M_s\to
CP^N,\quad N=\prod_{i=1}^s (N_i+1)-1,$$ of
$M_1\times\cdots\times M_s$ into $CP^N$   is
locally  the Segre embedding, that is,
$M_1,\ldots,M_s$ are open portions of
$CP^{N_1},\ldots, CP^{N_s}$, respectively, and
moreover, the K\"ahler immersion $f$ is congruent to
the Segre embedding.

This theorem was proved in [Nakagawa-Takagi
1976] under two additional assumptions;
namely, $s=2$ and the K\"ahler immersion $f$
has parallel second fundamental form.

\subsection{Parallel K\"ahler
submanifolds}

Ogiue (1972b) studied complex space forms in
complex space forms with parallel second
fundamental form and proved the following:
 Let $M^n(c)$ be a complex space form
analytically isometrically immersed in  another
complex space form $M^m(\bar c)$. If the second
fundamental form of the immersion is parallel, then
either the immersion is totally geodesic or $\bar
c>0$ and the immersion is given by the second
Veronese embedding.

Complete K\"ahler submanifolds in a complex
projective space with parallel second fundamental
form were  completely classified by
H. Nakagawa and R. Tagaki in 1976. Their
result states as follows.

 Let $M$ be a complete K\"ahler submanifold
embedded in $CP^m(c)$.  If $M$ is irreducible
and has parallel second fundamental form,
then $M$ is congruent to one of the following
six kinds of K\"ahler submanifolds:
$$CP^n(c),\;\; CP^{n}\left({c\over
2}\right),\; Q_n=SO(n+2)/SO(n)\times SO(2),$$ $$
SU(r+2)/S(U(r)\times U(2)),\; r\geq
3,\;\; SO(10)/U(5),\\ E_6/Spin(10)\times
SO(2).$$

If $M$ is reducible and has parallel second
fundamental form, then $M$ is congruent to
$CP^{n_1}\times CP^{n_2}$ with $n=n_1+n_2$ and the
embedding is given by the Segre embedding.

K. Tsukada (1985a) studied parallel
K\"ahler submanifolds of Hermitian symmetric
spaces and obtained the following: Let
$f: M\to\tilde M$ be a  K\"ahler immersion
of a  complete K\"ahler manifold $M$ into a
simply-connected Hermitian symmetric space
$\tilde M$. If $f$ has  parallel second
fundamental form, then $M$ is the direct
product of a complex Euclidean space and
semisimple Hermitian symmetric spaces.
Moreover,
$f=f_ 2\circ f_ 1$, where $f\sb 1$ is a
direct product of identity maps and (not totally
geodesic) parallel K\"ahler embeddings into
complex projective spaces, and
$f_ 2$ is a totally geodesic K\"ahler
embedding.

\subsection{Symmetric and
homogeneous K\"ahler submanifolds}

Suppose $f_i:M_i\to CP^{N_i},\;
i=1,\ldots,s,$ are full K\"ahler embeddings
of irreducible Hermitian symmetric spaces of
compact type and
$N=\prod_{i=1}^s (N_i+1)-1$. Then the composition
$$S_{N_1\cdots N_s}\circ
(f_1\times\cdots\times
f_s):M_1\times\cdots\times M_s\to
CP^N$$ is a full K\"ahler embedding, which is
called the tensor product of $f_1,\ldots,f_s$.

 H. Nakagawa and R. Tagaki
(1976) and R. Tagaki and M. Takeuchi (1977)
had obtained a close relation between the
degree and the rank of a symmetric K\"ahler
submanifold in complex projective space;
namely, they proved the following.

 Let $f_i:M_i\to CP^{N_i},\; i=1,\ldots,s,$ are
$p_i$-th full K\"ahler embeddings of irreducible
Hermitian symmetric spaces of compact type. Then
the degree of the tensor product of
$f_1,\ldots,f_s$ is given by $\sum_{i=1}^s r_ip_i$,
where $r_i=rk(M_i)$.

M. Takeuchi (1978) studied K\"ahler
immersions of homogeneous K\"ahler
manifolds and proved the following:

Let $f:M\to CP^N$ be a K\"ahler immersion of a
globally homogeneous K\"ahler manifold $M$.
Then

(1) $M$ is compact and simply-connected;

(2) $f$ is an embedding; and 

(3) $M$ is the orbit in $CP^N$ of the highest
weight in an irreducible unitary
representation of a compact semisimple Lie
group.

A different characterization of homogeneous
K\"ahler submanifolds in $CP^N$ was given by
S. Console and A. Fino (1996).

H. Nakagawa and R. Tagaki (1976) proved
that there do not exist K\"ahler immersions
from locally symmetric Hermitian manifolds
into complex Euclidean spaces and complex
hyperbolic spaces except the totally
geodesic ones. 

\subsection{Relative nullity
of K\"ahler submanifolds and reduction
theorem}

For a submanifold $M$ in a Riemannian manifold,
the subspace $$N_p=\{X\in T_pM:h(X,Y)=0,
\hbox{ for all } Y\in T_pM\},\quad p\in M$$
is called the relative nullity space of $M$
at $p$. The dimension $\mu(p)$  of $N_p$ is
called the relative nullity of $M$ at $p$.
The subset $U$ of $M$ where $\mu(p)$ assumes
the minimum, say $\mu$, is open in $M$, and
$\mu$ is called the index of relative
nullity.

K. Abe (1973) studied the index of relative
nullity of K\"ahler submanifolds and obtained the
following:

(1) Let $M$ be an $n$-dimensional K\"ahler
submanifold of a complex projective
$m$-space $CP^m$. If $M$ is complete, then
the index of relative  nullity is either $0$
or $2n$. In particular, if $\mu>0$, then
$M=CP^n$ and it is embedded as a totally
geodesic submanifold.

(2) Let $M$ be an $n$-dimensional complete
K\"ahler submanifold of the complex Euclidean
$m$-space $ C^m$. If the relative nullity
is greater than or equal to $n-1$, then $M$
is $(n-1)$-cylindrical. 

Let $M$ be an $n$-dimensional K\"ahler
submanifold of an $(n+p)$-dimensional complex
space form $\tilde M^{n+p}(c)$. A subbundle
$E$ of the normal bundle $T^\perp M$ is called
holomorphic if $E$ is invariant under the
action of the almost complex structure $J$ of
$\tilde M^{n+p}(c)$. 
For a holomorphic subbundle $E$ of $T^\perp
M$, let
$$\nu_E (x)=\dim_{ C}\{X\in T_xM: A_\xi X=0
\hbox{ for all } \xi\in E_x\}.$$
Put $\nu_E=$ Min$_{x\in M} \nu_E (x)$, which
is called the index of relative nullity with
respect to $E$. 

 Chen and  Ogiue (1973a) proved the
following: 
Let $M$ be an $n$-dimen\-sional
K\"ahler submanifold of a complex space form
$\tilde M^{n+p}(c)$. If there exists an
$r$-dimensional parallel normal subbundle $E$
of the normal bundle such that $\nu_E(x)\equiv
0$, then $M$ is contained in an
$(n+r)$-dimensional totally geodesic
submanifold of $\tilde M^{n+p}(c)$.

This result implies the following result of
Chen and Ogiue (1973a) and T. E. Cecil
(1974): If $M$ is  an $n$-dimensional K\"ahler
submanifold of a complete complex space form
$\tilde M^{n+p}(c)$ such that the first normal
space, Im $h$, defines an $r$-dimensional
parallel subbundle of the normal bundle, then
$M$ is contained in an $(n+r)$-dimensional
totally geodesic submanifold of $\tilde
M^{n+p}$.

 \vfill\eject

\section{Totally real and Lagrangian submanifolds of K\"ahler manifolds}
 
The study of totally real submanifolds  of a K\"ahler manifold from
differential geometric points of views was initiated in the early 1970's. A totally real
submanifold $M$ of  a K\"ahler manifold  $\tilde M$
is a submanifold such that the almost complex structure
$J$ of the ambient manifold $\tilde M$ carries each
tangent space of $M$ into the corresponding normal space of $M$,
that is, $J(T_{p}M)\subset T_{p}^{\perp}M$ for any
point $p\in M$. In other words, $M$  is a totally
real submanifold  if and only if, for any nonzero vector
$X$ tangent to $M$ at any point $p \in M$, the
angle between $JX$ and the tangent plane $T_{p}M$
is equal to ${\pi\over 2}$, identically. 
A totally real submanifold $M$ of a K\"ahler
manifold $\tilde M$ is called Lagrangian if
$\dim_R M=\dim_{C} \tilde M$. 

1-dimensional submanifolds, that
is, real curves, in a K\"ahler manifold are
always totally real. For this reason, we only 
consider totally real submanifolds of
dimension $\geq 2$. 

A submanifold $M$ of dimension $\geq 2$ in a
non-flat complex space form $\tilde M$ is
curvature invariant, that is, the Riemann
curvature tensor $\tilde R$ of $\tilde M$
satisfies $\tilde R(X,Y)TM\subset TM$ for
$X,Y$ tangent to $M$, if and only if $M$ is
either a K\"ahler submanifold or a totally
real submanifold [Chen-Ogiue 1974b].

 For a Lagrangian submanifold $M$ of a K\"ahler
manifold $(\tilde M,g,J)$, the tangent bundle
$TM$ and the normal bundle $T^\perp M$ are
isomorphic via the almost complex structure
$J$ of the ambient manifold. In particular,
this implies that the Lagrangian submanifold
has flat normal connection if and only if
the submanifold is a flat Riemannian
manifold. 

Let $h$ denote the second fundamental form of
the  Lagrangian submanifold in $\tilde M$ and
let $\alpha=Jh$. Another important property of
Lagrangian submanifolds is that
$g(\alpha(X,Y),JZ)$ is totally symmetric,
that is, we have
$$g(\alpha(X,Y),JZ)=g(\alpha(Y,Z),JX)=
g(\alpha(Z,X),JY)\leqno (16.1)$$
for any vectors $X,Y,Z$ tangent to $M$.

A result of M. L. Gromov (1985) implies
that every compact embedded Lagrangian
submanifold of $C^n$ is not simply-connected
(see [Sikorav 1986] for a complete proof of
this fact). This result is not true when the
compact Lagrangian submanifolds were
immersed but not embedded.

\subsection{ Basic properties of
Lagrangian submanifolds} 

A general K\"ahler manifold may not have any
minimal Lagrangian submanifold. 
Also, the only minimal Lagrangian immersion
of a topological 2-sphere into $CP^2$ is the
totally geodesic one. In contrast, minimal
Lagrangian submanifolds in an Einstein-K\"ahler
manifold exist in abundance, at least locally
(cf. [Bryant 1987]).

For surfaces in $E^4$  Chen and J. M.
Morvan (1987) proved that  an orientable
minimal surface $M$ in $ E^4$ is Lagrangian
with respect to an orthogonal almost complex
structure on $ E^4$ if and only if it is
holomorphic with respect to some orthogonal
almost complex structure on $E^4$.

A simply-connected Riemannian 2-manifold
$(M,g)$ with Gaussian curvature $K$
less than a constant $\,c\,$ admits a
Lagrangian isometric minimal immersion into a
complete simply-connected complex space form
$\tilde M^2(4c)$ if and only if it satisfies
the following differential equation [Chen
1997c; Chen-Dillen-Verstraelen-Vrancken
1995b]:
$$\Delta\ln (c-K)=6K,\leqno (16.2)$$
where $\Delta$ is the Laplacian on $M$
associated with the metric $g$.

The intrinsic and extrinsic structures of
Lagrangian minimal surfaces in
complete simply-connected complex space forms
were determined in [Chen 1997c] as follows:

 Let $f:M\to \tilde M^2(4c)$ be a minimal
Lagrangian surface without totally geodesic
points. Then,  with respect to a suitable
coordinate system $\{x,y\}$, we have

(1) the metric tensor of
$M$ takes the form of $$g=E\left(dx^2 +
dy^2\right)\leqno (16.3)$$
for some positive function $E$ satisfying
$$ \Delta_0(\ln E) =4E^{-2}-2cE,\leqno(16.4)$$
where $\Delta_0={{\partial^2
}\over{\partial x^2}} +{{\partial^2
}\over{\partial y^2}}$, and

(2) the second fundamental form of $L$ is
given by $$ h\left({\partial\over{\partial x}},{\partial\over{\partial x}}\right)=-{1\over E}J\left(
{\partial\over{\partial x}}\right),\quad h\left({\partial\over{\partial x}},{\partial\over{\partial
y}}\right)={1\over E}J\left( {\partial\over{\partial y}}\right),
$$ $$h\left({\partial\over{\partial y}},{\partial\over{\partial y}}\right)={1\over E}J
\left({\partial\over{\partial x}}\right).\leqno(16.5)$$

Conversely, if $E$ is a positive
function defined on a simply-connected domain
$U$ of $ E^2$ satisfying (16.4) and  
$g=E (dx^2 + dy^2)$ is the
metric tensor on $U$, then, up to rigid motions
of $\tilde M^2(4c)$, there is a unique minimal
Lagrangian isometric immersion of
$(U,g)$ into a complete
simply-connected complex space form $\tilde
M^2(4c)$ whose second fundamental form is
given by (16.5).

R. Harvey and H. B. Lawson (1982) studied
the so-called special Lagrangian
submanifolds in $C^n$, which are calibrated
by the $n$-form $Re(dz_1\wedge\cdots\wedge
dz_n)$. 

Being calibrated implies volume
minimizing in the same homology class. So, in
particular, the special Lagrangian
submanifolds are oriented minimal Lagrangian
submanifolds. In fact, they proved that a
special Lagrangian submanifold $M$ (with
boundary $\partial M$) in $ C^n$ is volume
minimizing in the class of all submanifolds
$N$ of $ C^n$ satisfying $[M]=[N]\in
H_n^c( C^n;\hbox{\bf R})$ with $\partial
M=\partial N$. Harvey and Lawson (1982) 
constructed many examples of special
Lagrangian submanifolds in $C^n$.

Using the idea of calibrations, one can show
that every Lagrangian minimal submanifold in
an Einstein-K\"ahler manifold $\tilde M$ with
$c_1(\tilde M)=0$ is volume minimizing. It
is false for the case $c_1=\lambda\omega$
with $\lambda>0$, where $\omega$ is the
canonical symplectic form on $\tilde M$. It
is unknown for the case
$c_1=\lambda\omega$ with $\lambda<0$ (cf.
[Bryant 1987c]).  J. G. Wolfson (1989)
showed that a compact minimal surface
without complex tangent points in an
Einstein-K\"ahler surface with $c_1<0$ is
Lagrangian.

Y. I. Lee (1994) studied embedded surfaces
which represent a second homology class
of an Einstein-K\"ahler surface. She obtained
the following:  

Let $N$ be  an Einstein-K\"ahler surface with
$c_1(N)<0$. Suppose $[A]\in H_2(N,\hbox{\bf
Z})$ and there exists an embedded  surface
without complex tangent points of genus $r$
which represents $[A]$. Then every connected
embedded minimal surface in $[A]$ has genus at
least $r$. Moreover, the equality occurs if
and only if the embedded minimal surface is
Lagrangian.

Notice that by using an adjunction inequality
for positive classes obtained by P. B.
Kronheimer and T. S. Mrowka (1994), the
minimality condition in her result can be
dropped; namely, under the same hypothesis,
one can conclude that every connected
embedded surface in [A] has genus at least
$r$.  Moreover, equality occurs if and
only if the embedded surface is Lagrangian.

Recently, Y. I. Lee also obtained the following
result:

Let $(N,g)$ be an Einstein-K\"ahler surface
with $c_1<0$. If an integral homology class
 $[A]\in H_2(N,\hbox{\bf Z})$ can be represented
by a union of Lagrangian branched minimal
surfaces with respect to $g$, then, for any other
Einstein-K\"ahler metric $g'$ on $N$ which
can be connected to $g$ via a family of
Einstein-K\"ahler metrics on $N$, $[A]$ can also be
represented by a union of branched minimal
surfaces with respect to $g'$. 

\subsection{ A vanishing theorem and its applications}

For compact Lagrangian submanifolds in
Einstein-K\"ahler manifolds, there is the
following vanishing theorem  [Chen 1998a]: 

Let $M$ be a compact manifold with 
finite fundamental group $\pi_1(M)$ or
vanishing first Betti number
$\beta_1(M)$. Then every Lagrangian
immersion from $M$ into any
Einstein-K\"ahler manifold must
have some minimal points.

This vanishing theorem has the following
interesting geometrical consequences:

(1) There do not exist Lagrangian isometric
immersions from a compact Riemannian
$n$-manifold with positive Ricci curvature
into any flat K\"ahler $n$-manifold or into any
complex hyperbolic $n$-space;

(2) Every Lagrangian isometric immersion of
constant mean curvature from a compact
Riemannian manifold with positive Ricci
curvature into any Einstein-K\"ahler
manifold is a minimal immersion; and

(3) Every Lagrangian isometric immersion of
constant mean curvature from a spherical
space form into a complex projective
$n$-space $CP^n$ is a totally geodesic
immersion.

This vanishing theorem is sharp in the
following sense: 

(a) The conditions on $\beta_1$ and
$\pi_1$ given in the vanishing theorem cannot
be removed, since the standard Lagrangian
embedding of $T^n=S^1\times\cdots\times S^1$
into $C^1\times\cdots\times C^1=C^n$ is a
Lagrangian embedding with nonzero constant
mean curvature; and 

(b) ``Lagrangian immersion'' in the theorem
cannot be replaced by the weaker condition
``totally real immersion'', since $S^n$ has
both trivial first Betti number and trivial
fundamental group; and the standard totally
real embedding of $S^n$ in $E^{n+1} \subset
C^{n+1}$ is a totally real submanifold with
nonzero constant mean curvature.

\subsection{The Hopf lift of Lagrangian
submanifolds of nonflat complex space forms
}

There is a general method for constructing 
Lagrangian submanifolds both in complex
projective spaces and in complex hyperbolic
spaces.

Let $S^{2n+1}(c)$ be the  hypersphere
of $ C^{n+1}$ with constant sectional
curvature $c$ centered at the origin. We
consider the Hopf fibration
$$\pi \colon\;S^{2n+1}(c)\to CP^n(4c).\leqno
(16.6)$$ Then $\pi$ is a Riemannian submersion,
meaning that $\pi_*$, restricted to the
horizontal space, is an isometry. Note that
given $z\in S^{2n+1}(c)$, the horizontal space
at $z$ is the orthogonal complement of $i z$
with respect to the metric induced on
$S^{2n+1}(c)$ from the usual Hermitian
Euclidean metric on $ C^{n+1}$.  Moreover,
given a horizontal vector $x$, then $i x$ is
again horizontal (and tangent to the sphere)
and $\pi_*(ix)=J (\pi_*(x))$, where $J$ is the
almost complex structure on $CP^n(4c)$. 

The main result of H. Reckziegel (1985) is
the following: Let $g\colon\; M\to CP^n(4c)$
be a Lagrangian isometric immersion. Then
there exists an isometric covering map
$\tau\colon\; \widehat M \to M$, and a horizontal
isometric immersion $f\colon\; \widehat  M\to
S^{2n+1}(c)$ such that $g(\tau)=\pi(f)$. 
Hence every Lagrangian immersion can be
lifted locally (or globally if we assume the
manifold is simply connected) to a horizontal
immersion of the same Riemannian manifold. 

Conversely, let $f\colon\; \widehat M\to
S^{2n+1}(c)$ be a horizontal isometric immersion.  Then
$g=\pi(f)\colon\;  M\to CP^n(4c)$ is again an
isometric immersion, which is Lagrangian. 
Under this correspondence, the second
fundamental forms $h^f$ and $h^g$ of $f$ and
$g$ satisfy $\pi_*h^f=h^g$.  Moreover, $h^f$
is horizontal with respect to $\pi$. 

In the complex hyperbolic case, we consider  the
complex number $(n+1)$-space $ C^{n+1}_1$
endowed with the pseudo-Euclidean metric
$g_0$ given by $$g_0=-dz_1d\bar z_1
+\sum_{j=2}^{n+1}dz_jd\bar z_j \leqno(16.7)$$ 

Put $$H^{2n+1}_1(c)=\left\{z=(z_1,z_2,
\ldots,z_{n+1})\colon\;
g_0(z,z)=\tfrac1{c}<0\right\}.\leqno(16.8)$$
$H^{2n+1}_1(c)$ is known as the
anti-de$\,$Sitter space-time.

Let $$T'_z=\{ u\in  C^{n+1}_1\colon\;
\left<u,z\right>=0\}, \quad
H_1^1=\{\lambda\in 
\hbox{\bf C}\colon\;
\lambda\bar\lambda=1\},$$ where
$\left<\;,\;\right>$ denotes the Hermitian
inner product on
$ C^{n+1}_1$ whose real part is $g_0$.
Then we have an $H^1_1$-action on
$H_1^{2n+1}(c)$,  $z\mapsto \lambda z$ and at each point $z\in
H^{2n+1}_1(c)$, the vector $iz$ is tangent to the flow of the
action. Since the metric $g_0$ is Hermitian, we have
$g_0(iz,iz)={1\over c}$. Note that the orbit is given by
$x_t=(\cos t + i\sin t)z$ and ${{dx_t}\over{dt}}=iz_t$. Thus the
orbit lies in the negative definite plane spanned by $z$ and
$iz$. The quotient space $H^{2n+1}_1/\sim$, under the
identification induced from the action, is the complex
hyperbolic space $CH^n(4c)$ with constant
holomorphic sectional curvature $4c$,
 with the complex structure $J$ induced from the canonical complex
structure $J$ on $ C^{n+1}_1$ via the
following totally geodesic fibration:
$$\pi\colon\; H^{2n+1}_1(c)\rightarrow 
CH^n(4c).\leqno(16.9)$$

Just as in the case of complex projective
spaces, let
$g\colon\; M\to CH^n(4c)$ be a Lagrangian
isometric immersion. Then there exists an
isometric covering map
$\tau\colon\; \widehat M \to M$, and a horizontal
isometric immersion
$f\colon\; \widehat M\to H_1^{2n+1}(c)$ such that
$g(\tau)=\pi(f)$.  Hence every totally real immersion can be
lifted locally (or globally if we assume the manifold is simply
connected) to a horizontal immersion. 

Conversely, let
$f:\widehat M\to H_1^{2n+1}(c)$ be a horizontal isometric
immersion.  Then $g=\pi(f)\colon\;  M\to
 CH^n(4c)$ is again
an isometric immersion, which is Lagrangian.

  Similarly, under
this correspondence, the second fundamental forms $h^f$ and
$h^g$ of $f$ and $g$ satisfy $\pi_*h^f=h^g$.  Moreover, $h^f$ is
horizontal with respect to  $\pi$.

This construction method have been used by
various authors. For instance, M.  Dajczer and
R. Tojeiro (1995b)  applied  this technique
to show that the projection of the Hopf
fibration provides a one-to-one
correspondence between the set of  totally
real flat submanifolds in complex
projective space
 and the set of symmetric flat submanifolds
in Euclidean sphere with the same codimension,
where a flat submanifold $F:M^{n+1}_0\to
S^{2n+1}$ is called symmetric if its
corresponding principal coordinates generate a
solution of the symmetric generalized wave
equation. 

\subsection{Totally real minimal
submanifolds of complex space
forms}

The equation of Gauss
implies that an
$n$-dimensional totally real minimal submanifold
of a complex space form 
$\tilde M^n(4c)$ satisfies the following
properties:

(1) $Ric\leq (n-1)cg$, with equality holding
if and only if it is totally geodesic.

(2) $\rho\leq n(n-1)c$, with equality
holding if and only if it is totally geodesic.

B. Y. Chen and C. S. Houh (1979) proved
that if the sectional curvature of  an
$n$-dimensional compact totally real
submanifold $M$ of $CP^m(4c)$ satisfies
$K\geq (n-1)c/(2n-1)$, then either $M$ is
totally geodesic or $n=2,\,m\geq 3,$
and $M$ is of constant Gaussian
curvature $c/3$. In both cases, the
immersion is rigid.

 K. Kenmotsu (1985) gave a
canonical description of all totally real
isometric minimal immersions of
$ E^2$ into $CP^m$ as follows.

Let $\Omega$ be a simply-connected domain in the
Euclidean plane $ E^2$ with metric
$g=2|dz|^2$ where $z=x+\sqrt{-1}y$ is the
standard complex coordinate on $E^2$, and
let $\psi:\Omega\to CP^m$ be a totally real
isometric minimal immersion. Then, up to a
holomorphic isometry of $CP^m$, $\psi$ is given by
$$\psi(z)=[r_0e^{\mu_{0}z-\bar \mu_0
\bar z},\ldots, r_me^{\mu_{m}z-\bar\mu_m
\bar z}],\leqno (16.10)$$ where $r_0,\ldots,r_m$
are non-negative real numbers and
$\mu_0,\ldots,\mu_m$ are complex numbers of unit
modulus satisfying the following three
conditions:
$$\sum_{j=0}^m r^2_j=1,\quad \sum_{j=0}^m
r^2_j\mu_j=0,\quad \sum_{j=0}^m
r^2_j\mu_j^2=0.\leqno (16.11)$$

Furthermore, $\psi$ is linearly full, that is,
$\psi(\Omega)$ does not lie in any linear
hyperplane of $CP^m$, if and only if $r_0,\ldots,
r_m$ are strictly positive and
$\mu_0,\ldots,\mu_m$ are distinct, in which case
the complex numbers $r_0\mu_0,\ldots,r_m\mu_m$
are uniquely determined up to permutations.

Conversely, any map $\psi$ defined by
$(16.10)$, where
$r_0,\ldots,r_m,\mu_0,\ldots,\mu_m$ satisfy
$(16.11)$, is a totally real isometric minimal
immersion of $\Omega$ into $CP^m$.

\subsection{ Lagrangian real space form
in  complex space form} 

The simplest examples of Lagrangian (or more
generally totally real) submanifolds of complex
space forms are  totally geodesic Lagrangian
submanifolds. A totally geodesic
Lagrangian submanifold
$M$ of a complex space form $\tilde M^n(4c)$ of
constant holomorphic sectional curvature $4c$ is a
real space form of constant curvature $c$. 

The real projective
$n$-space \ $RP^n(1)\;$ (respectively, the
real hyperbolic $n$-space $H^n(-1)$) can be
isometrically embedded in
$CP^n(4)$ (respectively, in 
complex hyperbolic space $CH^n(-4)$) as a
Lagrangian totally geodesic submanifold. 

Non-totally
geodesic  Lagrangian  isometric immersions
from  real space forms of constant curvature
$c$ into a complex space form 
$\tilde M^n(4c)$ were determined by 
Chen, Dillen, Verstraelen, and Vrancken
(1998). Associated with each twisted product
decomposition of a real space form, they
 introduced a canonical 1-form, called
the twistor form of the twisted product
decomposition. Their result says that if the
twistor form of a twisted product
decomposition of a simply-connected real
space form of constant  curvature $c$ is
twisted closed, then, up to motions, it
admits a unique ``adapted'' Lagrangian
isometric immersion  into the  complex space
form $\tilde M^n(4c)$. 

Conversely, if $\; L:\;M^n(c)\to
\tilde M^n(4c)\;$ is a non-totally geodesic 
Lagrangian  isometric immersion of a
real space form $M^n(c)$ of constant
sectional curvature $c$ into a complex space
form $\tilde M^n(4c)$, then
$M^n(c)$ admits an appropriate twisted product
decomposition with twisted closed twistor form
and, moreover, the Lagrangian  immersion $L$ is
given by the corresponding adapted Lagrangian
isometric immersion of the twisted product. 

Chen and Ogiue (1974b) proved that a
Lagrangian minimal submanifold of constant
sectional curvature $c$ in a complex space
form $\tilde M^{n}(4\tilde c)$ is either
totally geodesic or $c\leq 0$. 
N. Ejiri (1982) proved that the only
Lagrangian  minimal submanifolds of constant
sectional curvature $c\leq 0$ in a complex
space form are the flat ones. Ejiri's result
extends the corresponding result of Chen and
Ogiue (1974b) for $n=2$ to $n\geq 2$.

A  submanifold $M$ of a  Riemannian manifold
is called a Chen submanifold if $$\sum_{i,j}
\left<{}\right. h(e_i,e_j),
H\left. {}\right>h(e_i,e_j)\leqno (16.12)$$ is
parallel to the mean curvature vector
$ H$, where $\{e_i\}$ is an
orthonormal frame of the submanifold $M$
(for general properties of Chen
submanifolds, cf.
[Gheysens-Verheyen-Verstraelen 1981, Rouxel
1994]).

 M. Kotani (1986) studied Lagrangian Chen
submanifolds of constant curvature in complex
space forms and obtained the following:

If $M$ is a Lagrangian Chen submanifold with
constant sectional curvature  in a complex
space form $\tilde M^n(4\tilde c)$ with
$c<\tilde c$, then either $M$ is minimal, or
locally, $M=I\times \tilde L^{n-1}$ with
metric $g=dt^2+f(t)\tilde g$, where $I$ is an
open interval, $(\tilde L^{n-1},\tilde g)$ is
the following submanifold in $\tilde
M^n(4\tilde c)$:
$$\aligned &\tilde L^{n-1}\;\subset 
S^{2n-1}\subset \tilde M^{n}(4\tilde c)\\
\pi &\downarrow\hskip.48in \downarrow\pi
 \\&
L^{n-1}\;\subset CP^{n-1},
\endaligned\leqno (16.13)$$
$S^{2n-1}$ is a geodesic hypersphere in
$\tilde M^n(4\tilde c)$, and $\tilde L$ is the
horizontal lift of a Lagrangian minimal flat
torus $L^{n-1}$ in $CP^{n-1}$. 

A flat torus $T^n$ can be
isometrically immersed in $CP^n$ as Lagrangian
minimal submanifold with parallel second
fundamental form, hence with parallel
nonzero mean curvature vector. 

M. Dajczer and R. Tojeiro (1995b) proved that a
complete flat Lagrangian submanifold in $CP^n$
with constant mean curvature is a flat torus
$T^n$ (with parallel second fundamental
form).

\subsection{Inequalities for Lagrangian
submanifolds} 

For any $n$-dimensional Lagrangian submanifold
$M$ in a complex space form $\tilde M^n(4c)$
and for any $k$-tuple $(n_1,\ldots,n_k)\in
\Cal S(n)$, the $\delta$-invariant
$\delta(n_1,\ldots,n_k)$ must satisfies
following inequality [Chen 1996f]:
$$\delta(n_1,\ldots,n_k)\leq
b(n_1,\ldots, n_k)H^2+a(n_1,\ldots,n_k)c,\leqno
(16.14)$$ where
$\delta(n_1,\ldots,n_k),\,a(n_1,\ldots,n_k)$
and
$b(n_1,\ldots,n_k)$ are defined by
(3.14), (3.15) and (3.16), respectively.

There exist abundant examples of Lagrangian
submanifolds in complex space forms
which satisfy the equality case of (16.14). 

Inequality (16.14) implies that
for any Lagrangian submanifold $M$ in a
complex space form
$\tilde M^n(4c)$, one has

$$\delta(2) \leq  {{n^2(n-2)}\over{2(n-1)}}
H^2 +\tfrac12(n+1)(n-2) c.\leqno (16.15)$$

Chen, Dillen, Verstraelen and Vrancken (1994)
proved that if a Lagrangian submanifold $M$ of
$\tilde M^n(4c)$ satisfies the equality case
of (16.17) identically, then $M$ must be
minimal in $\tilde M^n(4c)$. 

There exist many Lagrangian minimal
submanifolds of $\tilde M^n(4c), \, c\in
\{-1,0,1\},$ which satisfy the equality case
of (16.15) identically. 
However, Chen, Dillen, Verstraelen and Vrancken
(1996) proved that  if $M$ has constant scalar
curvature, then the equality occurs when
and only when either $M$ is totally geodesic
in $\tilde M^n(4c)$ or $n=3,\, c=1$ and the
immersion is locally congruent to a special
non-standard immersion $\psi:S^3\to CP^3$ of a
topological 3-sphere into $CP^3$ which is
called an exotic immersion of $S^3$. The
classification of Lagrangian submanifolds of
$CP^3$ satisfying the equality case of
(16.15) was given in
[Bolton-Scharlach-Vrancken-Woodward 1998]

Lagrangian submanifolds of the complex
hyperbolic $n$-space $CH^n$, $n\geq 3$,
satisfying the equality case of (16.15) were
classified in [Chen-Vrancken 1997c].

We remark that inequality (16.14) holds for
an arbitrary $n$-dimensional submanifold in
$CH^n(4c)$ for $c<0$ as well.

\subsection{ Riemannian and
topological obstructions to Lagrangian
immersions}

M. L. Gromov (1971) proved that a compact
$n$-manifold $M$ admits a Lagrangian
immersion into $C^n$ if and only if the
complexification of the tangent bundle of
$M$, $TM\otimes \hbox{\bf C}$, is trivial. 
Since the tangent bundle of a 3-manifold is
always trivial, Gromov's result implies that
there does not exist topological obstruction
to Lagrangian immersions for compact
3-manifolds.

In contrast, by applying inequality (16.14)
and the vanishing theorem mentioned in
\S16.2, one obtains the following sharp
obstructions to isometric Lagrangian
immersions of compact Riemannian
manifolds into complex space forms [Chen
1996f]:

 Let $M$ be a compact Riemannian manifold with
finite fundamental group $\pi_1(M)$ or
$b_1(M)=0$. If there exists a $k$-tuple
$(n_1,\ldots,n_k)\in \Cal S(n)$ such that
$$\delta(n_1,\ldots,n_k)>
{1\over2}\Big({{n(n-1)}}-\sum_{j=1}^k
{{n_j(n_j-1)}}\Big)c,\leqno(16.16)$$ then $M$
admits no Lagrangian isometric immersion
into a complex space form of constant
holomorphic sectional curvature $4c$.
\smallskip

An immediate important consequence of the
above result is the first necessary
intrinsic condition for compact Lagrangian
submanifolds in $C^n$; namely, the Ricci
curvature of every compact Lagrangian
submanifold $M$ in
$C^n$ must satisfies $\,\inf_u Ric(u)\leq
0$, where $u$ runs over all unit tangent
vectors of $M$ [Chen 1997d]. For Lagrangian
surfaces, this means that the Gaussian
curvature of every compact Lagrangian
surface $M$ in
$C^2$ must be nonpositive at some points
on $M$. Another immediate consequence is
that every compact irreducible symmetric
space cannot be isometrically immersed in a
complex Euclidean space as a Lagrangian
submanifold.
\smallskip

Let $f: E^{n+1}\to C^n$ be the map
defined
by $$f(x_0,\ldots,x_n)={1\over{1+x_0^2}}
(x_1,\ldots, x_n,x_0x_1,\ldots, x_0x_n).
\leqno(16.17)$$ Then $f$ induces an immersion
$w:S^{n}\to  C^n$  of $S^n$ into $
C^n$  which has a unique self-intersection
point $f(-1,0,\ldots,0)=f(1,0,\ldots,0)$. 

With respect to the canonical almost complex  
structure $J$ on $ C^n$,  the immersion $w$
is a Lagrangian immersion of $S^n$ into $
C^n$, which is called the Whitney immersion.
$S^n$ endowed with the  Riemannian metric
induced from the Whitney immersion is called
a Whitney $n$-sphere.

The example of the Whitney immersion
shows that the condition on the
$\delta$-invariants given above is sharp,
since $S^n\;(n\geq 2)$ has trivial
fundamental group and trivial first Betti
number; moreover, for each $k$-tuple
$(n_1,\ldots, n_k)\in\Cal S(n)$, the
Whitney $n$-sphere  satisfies
$\delta(n_1,\ldots, n_k)>0$ except at the
unique point of self-intersection.

Also, the assumptions on the
finiteness of $\pi_1(M)$ and vanishing of 
$b_1(M)$ given above are both necessary
for $n\geq 3$. This can be seen from the
following example: 

Let $F:S^1\to\hbox{\bf C}$ be the unit
circle in the complex plane given by
$F(s)=e^{is}$ and let $\iota:S^{n-1}\to E^n$
$(n\geq 3)$ be the unit hypersphere in $E^n$
centered at the origin. Denote by $f:S^1\times
S^{n-1}\to C^n\,$ the complex extensor
defined by $f(s,p)=F(s)\otimes \iota(p)$, $p\in
S^{n-1}$. Then $f$ is an isometric Lagrangian
immersion of $M=:S^1\times S^{n-1}$ into $C^n$
which carries each pair
$\{(u,p),\,(-u,-p)\}$  of points in $S^1\times
S^{n-1}$ to a point in $C^n$ (cf. [Chen
1997b]). Clearly, $\pi_1(M)=\hbox{\bf Z}$ and
$b_1(M)=1$, and moreover, for each $k$-tuple
$(n_1,\ldots,n_k)\in \Cal S(n)$, the
$\delta$-invariant  $\delta(n_1,\ldots,n_k)$
on $M$ is a positive constant. This example
shows that both the conditions on
$\pi_1(M)$ and $b_1(M)$ cannot be removed.

\subsection{ An inequality between scalar
curvature and mean curvature} 

Besides inequality (16.14), there is another
sharp inequality for Lagrangian submanifolds
in complex space forms.  

Let $\rho=\sum_{i\ne j} K(e_i\wedge e_j)$ denote
the scalar curvature of a Riemannian
$n$-manifold $M$, where
$e_1,\ldots,e_n$ is an orthonormal local frame.
The scalar curvature $\rho$ and the squared mean
curvature $H^2$ of a  Lagrangian submanifold  in 
complex space form $\tilde M^n(4c)$ satisfy the
following general sharp inequality:
$$H^2\geq {{n+2}\over{n^2(n-1)}}\rho
-\left( {{n+2}\over n}\right)c.\leqno (16.18)$$

Inequality (16.18) with $c=0$ and $n=2$ was
proved in [Castro-Urbano 1993]. Their proof
relies on complex analysis which is not
applicable to $n\geq 3$. The general
inequality was established in
[Borrelli-Chen-Morvan 1995] for
$c=0$ and arbitrary $n$; and in [Chen 1996b]
for
$c\ne 0$ and arbitrary $n$; and independently
by Castro and Urbano (1995), for
$c\ne 0$ with $n=2$, also using the method
of complex analysis. 

If $\tilde M^n(4c)= C^n$, the equality of
(16.18) holds identically if and only if
either the Lagrangian submanifold $M$ is an
open portion of a Lagrangian $n$-plane or,
up to dilations, $M$ is an open portion of
the Whitney immersion [Borrelli-Chen-Morvan
1995] (see also [Ros-Urbano 1998] for an
alternative proof).

Chen (1996b) proved that there
exists a one-parameter family of Riemannian
$n$-manifolds, denoted by
$P^n_a\, (a>1)$, which admit Lagrangian
isometric immersions into
$CP^n(4)$ satisfying the equality case of
the inequality (16.18) for $c=1$; and  there
are two one-parameter families of Riemannian
manifolds, denoted by
$\,C^n_a\, (a>1),D_a^n\,$ $ (0<a<1)$, and two
exceptional $n$-spaces, denoted by $F^n, L^n$,
which admit Lagrangian isometric immersion into
$ CH^n(-4)$,  satisfying the equality case of
the inequality for $c=-1$.  Besides the
totally geodesic ones, these are the only
Lagrangian submanifolds in $
CP^n(4)$ and in $CH^n(-4)$ which satisfy the
 equality case of (16.18) (see also
[Castro-Urbano 1995] for the case $n=2$). 

The explicit expressions of
those Lagrangian immersions of $P^n_a,
C^n_a, D^n_a, F^n$ and $L^n$ satisfying the
equality case of (16.18) were completely
determined in [Chen-Vrancken 1996].

I. Castro and F. Urbano (1995) showed that a
Lagrangian surface in $CP^2$ satisfies the
equality case of (16.14) for $n=2$ and $c=1$
if and only if the Lagrangian surface has
holomorphic twistor lift. 

\subsection{Characterizations
of parallel Lagrangian submanifolds } 

Compact Lagrangian submanifolds of
$CP^n(4c)$ with parallel second fundamental form
were completely classified by H. Naitoh in
[Naitoh 1981]. 

There are various pinching results for
totally real submanifolds in complex space
forms similar to K\"ahler submanifolds given
as follows. 

Y. Ohnita (1986b)
and F. Urbano (1986)  proved the following: Let
$M$ be a compact Lagrangian  submanifold
of a complex space form with  parallel mean
curvature vector. If $M$ has nonnegative
sectional curvature, then the second fundamental
form of $M$ is parallel. 

By applying the results of Ohnita and
Urbano, it follows that compact Lagrangian
submanifolds of
$CP^n(4)$ with parallel mean curvature vector
satisfying $K\geq 0$ must be the
products $T\times M_1\times\cdots\times M_k$,
where $T$ is a flat torus with $\dim M\geq k-1$
and each $M_i$ is one of the following:
$$ RP^r(1)\to CP^r(4)\;\;(r\geq
2)\;\;(\hbox{totally geodesic}),$$ $$
SU(r)/SO(r)\to CP^{(r-1)(r+2)/2}(4)\;\;(r\geq
3)\;\;(\hbox{minimal}),$$ $$
SU(2r)/Sp(r)\to CP^{(r-1)(2r+1)}(4)\;\;(r\geq
3)\;\;(\hbox{minimal}),\leqno (16.19)$$ $$SU(r)\to CP^{r^2-1}(4)\;\;(r\geq
3)\;\;(\hbox{minimal}),$$ $$
E_6/F_4\to CP^{26}(4)\;\;(\hbox{minimal}).
$$

S. Montiel, A.Ros and F. Urbano (1986) proved
that  if a compact Lagrangian minimal
submanifold of $CP^n(4)$ satisfies $Ric
\geq {3\over 4}(n-2)g$, then the second
fundamental form is parallel.

Combining this theorem   with the result
of H. Naitoh (1981), it follows that compact
Lagrangian minimal submanifold of $CP^n(4)$
satisfying   $Ric \geq {3\over 4}(n-2)g$ is
either one of the Lagrangian minimal
submanifolds given in (16.19)  or a minimal
Lagrangian flat torus in $CP^2(4)$.

Compact Lagrangian minimal submanifolds in
$CP^n(4)$ satisfying a pinching of scalar
curvature were studied in  [Chen-Ogiue 1974b, 
Shen-Dong-Guo 1995], among others. 

For scalar curvature pinching, we have the
following: If $M$ is an $n$-dimen\-sional
compact Lagrangian minimal submanifold in
$CP^n(4)$ whose  scalar curvature $\rho$ 
satisfies $\rho\geq 3(n-2)n/4$, then $M$ has
parallel second fundamental form.

F. Urbano (1989) proved that if $M^3$ is a
$3$-dimensional compact Lagrangian  submanifold
of a complex space form
$\tilde M^3(4c)$ with nonzero parallel mean
curvature vector, then $M^3$ is flat and has
parallel second fundamental form.

For complete Lagrangian submanifolds in $
C^n$ with parallel mean curvature vector,
F. Urbano (1989) and U. H. Ki and Y. H. Kim
(1996) proved the following.

Let $M$ be a complete Lagrangian submanifold
embedded in $ C^n$. If $M$ has parallel mean
curvature vector, then $M$ is either a
minimal submanifold or a product submanifold
$M_1\times\cdots\times M_k$, where each
$M_i$ is a Lagrangian submanifold embedded in some
$ C^{n_i}$ and each $M_i$ is also a minimal
submanifold of a hypersphere of $C^{n_i}$. 

J. S. Pak (1978) studied totally real
planar geodesic immersions and obtained that
if $f:M\to CP^m$ is a totally real isometric
immersion such that each geodesic $\gamma$
of $M$ is mapped into a 2-dimensional
totally real totally geodesic submanifold
of $CP^m$, then $M$ is a locally a compact
rank one symmetric space and the immersion
is rigid.  

\subsection{  Lagrangian $H$-umbilical
submanifolds and Lagrangian catenoid} 

Since there  do not exist  totally
umbilical Lagrangian submanifolds in complex
space forms except totally
geodesic ones, it is natural to look for the
``simplest'' Lagrangian submanifolds next to
totally geodesic ones. As a result, the
following notion of Lagrangian $H$-umbilical
submanifolds was introduced. 

A  non-totally geodesic Lagrangian submanifold
$M$ of a K\"ahler manifold is called a
Lagrangian $H$-umbilical submanifold if its
second fundamental form takes the following
simple  form:  $$
h(e_1,e_1)= \lambda Je_1,\quad
h(e_2,e_2)=\cdots = h(e_n,e_n)=\mu Je_1,$$ $$
h(e_1,e_j)=\mu Je_j,\quad h(e_j,e_k)=0,\; \;
j\not=k,
\;\;\;\; j,k=2,\ldots,n
\leqno(16.20)$$ for some suitable functions
$\lambda$ and
$\mu$ with respect to some suitable
orthonormal local frame field.  

Clearly,  a non-minimal  Lagrangian
$H$-umbilical  submanifold satisfies the
following two conditions:

(a) $JH$ is an eigenvector of the shape operator
$A_H$ and 

(b) the restriction of $A_H$ to
$(JH)^\perp$ is proportional to the identity map.

In fact, Lagrangian
$H$-umbilical  submanifolds are the
simplest Lagrangian submanifolds  satisfying
both conditions (a) and (b). In this way 
Lagrangian $H$-umbilical  submanifolds can
be considered as the ``simplest'' Lagrangian
submanifolds in  a complex space form next
to the totally geodesic ones. 

There exist ample examples of Lagrangian
$H$-umbilical  submanifolds in  a complex
space form. For instance, Chen (1997b)
showed that every complex extensor of the
unit hypersurface of $E^n$ via any unit
speed curve in the complex plane gives rise
to a Lagrangian $H$-umbilical submanifold
in $C^n$. Using this method one can
constructed ample examples of compact
Lagrangian submanifolds in a complex
Euclidean space.

 Lagrangian $H$-umbilical  submanifolds 
were classified in [Chen 1997a,1997b] in
connections with the notions of Legendre
curves and complex extensors. In
particular, Chen proved that, except the
flat ones, Lagrangian $H$-umbilical
submanifolds of dimension $\geq 3$ in a
complex Euclidean space are Lagrangian
pseudo-sphere ($\lambda=2\mu$) and complex
tensors of the unit hypersphere of a
Euclidean space. Moreover, except some
exceptional classes, Lagrangian
$H$-umbilical submanifolds in $CP^n$ and in
$CH^n$, $n\geq 3$, are obtained from Legendre
curves in $S^3$ or in $H^3_1$ via warped
products in some natural ways.
The intrinsic and extrinsic
structures of Lagrangian $H$-umbilical 
surfaces in complex space forms were
determined in [Chen 1997c]. In [Chen 1998e],
a representation formula for flat Lagrangian
$H$-umbilical submanifolds in complex
Euclidean spaces was discovered.

By a complex extensor of the unit hypersphere
$S^{n-1}$ of $ E^n$, we mean a Lagrangian
submanifold of $ C^n$ which is given by the
tensor product of the unit hypersphere and a
unit speed (real) curve in the complex line
$ C^1$. 
\smallskip
The  Lagrangian catenoid was constructed
by R. Harvey and H. B. Lawson (1982) which is
defined by
$$\aligned M_0=&\{(x,y)\in C^n=
E^n\times E^n: |x|y=|y|x,\,\\ &
\hbox{Im}\,(|x|+i|y|)^n=1,\, |y|<|x|\tan
{\pi\over n}\}.\endaligned\leqno(16.21)$$

Besides being a minimal Lagrangian
submanifold of $ C^n$, $M_0$ is
invariant under the diagonal action of
$SO(n)$ on $ C^n= E^n\times 
E^n$. 

I. Castro and F. Urbano (1997) characterized
Lagrangian catenoid as the only minimal
nonflat Lagrangian submanifold in $ C^n$
which is foliated by pieces of round
$(n-1)$-spheres of $ C^n$ (up to
dilations).

\subsection{Stability of Lagrangian submanifolds}

 The stability of minimal Lagrangian
submanifolds of a K\"ahler manifold was first
investigated by  B. Y. Chen, P. F. Leung
and T. Nagano in 1980 (cf. [Chen 1981a]).
In particular, they proved that the second
variational formula of a compact Lagrangian
submanifold $M$ in a
 K\"ahler manifold
$\tilde M$ is given by 
$$V''(\xi)=\int_M \left\{{1\over
2}||dX^\#||^2+(\delta X^\#)^2-\tilde
S(X,X)\right\}dV,\leqno(16.22)$$ where $JX=\xi$,
$X^\#$ is the dual 1-form of $X$ on $M$,
$\delta$ is the codifferential operator, and
$\tilde S$ is the Ricci tensor of $\tilde M$.

By applying (16.22), Chen, Leung and Nagano
proved the following (cf. [Chen 1981a]):

  Let $f:M\to\tilde M$ be a compact
Lagrangian minimal submanifold of a Kaeh\-ler
manifold $\tilde M$. 

(1) If $\tilde M$ has positive Ricci
curvature, then the index of $f$ satisfies
$i(f)\geq
\beta_1(M)$, where $\beta_1(M)$ denotes the
first Betti number of $M$.  In particular,
if the first cohomology group of $M$ is
nontrivial, that is, $H^1(M;\hbox{\bf R})\ne
0$, then $M$ is unstable;

(2) If $\tilde M$ has nonpositive Ricci
curvature, then $M$ is stable.

 Y. G. Oh (1990) introduced the notion
of Hamiltonian deformations in K\"ahler
manifolds. He considered normal variations
$V$ along a minimal Lagrangian submanifold
$M$ such that the 1-form
$\alpha_V=\left<JV,\,\cdot\,\right>$ is exact
and call such variations Hamiltonian variations.

A minimal Lagrangian submanifold is called 
Hamiltonian stable if the second variation is
nonnegative in the class of Hamiltonian
variations. 

Oh (1990) establishes the following Hamiltonian
stability criterion on Einstein-K\"ahler 
manifolds:  Let
$\tilde M$ be an Einstein-K\"ahler manifold with
$Ric=cg$, where $c$ is a constant. Then a minimal
Lagrangian submanifold $M$ is locally Hamiltonian
stable if and only if $\lambda_1(M)\geq c$, where
$\lambda_1(M)$ is the first nonzero eigenvalue of
the Laplacian acting on $C^\infty(M)$.  

According to [Lawson-Simons 1973] the Lagrangian
totally geodesic $RP^n(1)$ in $CP^n(4)$ is
unstable in the usual sense. In contrast, Oh's
result implies that the Lagrangian  totally
geodesic $RP^n(1)$ is Hamiltonian stable in
$CP^n(4)$.

I. Castro and F. Urbano (1998) constructed
examples of unstable Hamiltonian minimal
Lagrangian tori in $C^2$.

\subsection{Lagrangian immersions and Maslov class}

Let $\Omega$ denote the canonical
symplectic form on $ C^n$ defined by
$$\Omega(X,Y)=\left<JX,Y\right>.\leqno(16.23)$$
Consider the Grassmannian $\Cal L( C^n)$ of
all Lagrangian vector subspaces of $ C^n$. 
$\Cal L( C^n)$ can be identified with the
symmetric space $U(n)/O(n)$ in a natural
way. 

$U(n)/O(n)$ is a bundle over the circle
$S^1$ in $ C^1$ with the projection
$${\det}^2: U(n)/O(n)\to S^1,$$ where $\det^2$
is the square of the determinant.

For a Lagrangian submanifold $M$ in $
C^n$, the Gauss map takes the values in
$\Cal L( C^n)$ which yields the
following sequence:
$$M \to  \Cal L( C^n)\cong
U(n)/O(n) \to   S^1.\leqno(16.24)$$

If $ds$ denotes the volume form of $S^1$,
then $m_M:=(\det^2\circ G)^*(ds)$ is a
closed 1-form on $M$. The cohomology class  
$[m_M]\in H^1(M;\hbox{\bf Z})$ is called
the Maslov class of the Lagrangian
submanifold $M$. 

J. M. Morvan (1981) proved that the Maslov
form $m_M$ and the mean curvature vector
of a Lagrangian submanifold $M$ in $ C^n$
are related by
$$m_M(X)={1\over\pi}
\left<\right. J H,X\left
.\right>,\quad X\in TM.\leqno(16.25)$$ Hence, a
Lagrangian submanifold $M$ in $C^n$ is
minimal if and only if $\det^2\circ G$ is a
constant map.

Let $\xi$ be a normal vector field of a
Lagrangian submanifold $M$ of a K\"ahler
manifold
$\tilde M$. Denote by $\alpha_\xi$ the 1-form
on $M$ defined by 
$$\alpha_\xi(X)=\Omega(\xi,X)=
\left<J\xi,X\right>,\quad
X\in TM,\leqno(16.26)$$ where $\Omega$ is the
K\"ahler form of $\tilde M$. 

Chen and Morvan (1994) introduced the notion
of harmonic deformations in K\"ahler
manifolds: A normal vector field
$\xi$ of a Lagrangian submanifold $M$ is
called harmonic if the 1-form 
$\alpha_\xi$ associated with $\xi$ is a
harmonic 1-form. A normal variation of a
Lagrangian submanifold in a K\"ahler manifold
is called  harmonic if its variational
vector field is harmonic. 

A Lagrangian
submanifold $M$ of a K\"ahler manifold is
called harmonic minimal if it is a critical
point of the volume functional in the class of
harmonic variations.

  Chen and Morvan (1994) proved
that the  Maslov class of a Lagrangian
submanifold of an Einstein-K\"ahler manifold
 vanishes if and only if it
is harmonic minimal, thus  providing a
simple relationship between the calculus of
variations and Maslov class.
This result implies in particular that a
closed curve
$\gamma$ in a  K\"ahler manifold
$\tilde M$ with $\dim_{ R}\tilde M=2$ is
harmonic minimal if and only if it has zero
total curvature, that is,  $\int_\gamma
\kappa(s)ds=0$; an extension of a result
mentioned in [Vaisman 1987]. 

A Lagrangian submanifold $M$ is said to have
conformal Maslov form if $J H$ is
a conformal vector filed on $M$. The Whitney
immersion $w$ defined by (16.12) is known
to have conformal Maslov form.

 A. Ros and F. Urbano
(1998) proved that, up to dilations,
Whitney's  immersion is the only Lagrangian
immersion of a compact manifold with zero
first betti number and conformal Maslov form.

\vfill\eject

\section{$CR$-submanifolds of K\"ahler manifolds}

\subsection{Basic properties of
$CR$-submanifolds of K\"ahler
manifolds}

Let $M$ be a submanifold of a K\"ahler
manifold (or more generally, of an almost
Hermitian manifold) with almost complex structure
$J$ and metric $g$. At each point $x\in
M$, let
$\Cal D_x$ denote the maximal holomorphic
subspace of the tangent space $T_xM$, that is,
$\Cal D_x=T_xM\cap J(T_xM)$. If the dimension
of $\Cal D_x$ is the same for all $x\in M$, we
have a holomorphic distribution $\Cal D$ on
$M$. The submanifold $M$ is called a
$CR$-submanifold if there exists on $M$ a
holomorphic distribution
$\Cal D$ such that its orthogonal complement
$\Cal D^\perp$ is totally real, that is,
$J(\Cal D_x^\perp)\subset T^\perp_x M$, for
all $x\in M$ [Bejancu 1978].

Every real hypersurface of a
Hermitian manifold is a $CR$-submanifold.
 A $CR$-submanifold is called proper if it is
neither a K\"ahler submanifold ($\Cal D=TM$)
nor a totally real submanifold ($\Cal
D^\perp=TM$).

Blair and Chen (1979) proved that a submanifold
$M$ of a non-flat complex space form $\tilde M$
is a $CR$-submani\-fold  if and only
if the maximal holomorphic subspaces  define a
holomorphic distribution
$\Cal D$ on $M$ such that $\tilde R({\Cal
D},{\Cal D};{\Cal D}^\perp,{\Cal
D}^\perp)=\{0\},$ where
$\Cal D^\perp$ denotes the orthogonal
distribution of $\Cal D$ in $TM$, and
$\tilde R$ is the Riemann curvature tensor of
$\tilde M$. 

In this section, we denote by $h$
the complex rank of the holomorphic
distribution $\Cal D$ and by $p$ the real
rank of the totally real
distribution $\Cal D^\perp$ of a $CR$-submanifold
$M$ so that $\dim M=2h+p$.

Let $N$ be a differentiable manifold and $T_CN$
be the complexified tangent bundle of $N$. A
$CR$-structure on $N$ is a complex subbundle
$\Cal H$ of $T_CN$ such that 
${\mathcal H}\cap \bar{\mathcal H}=\{0\}$ and $\Cal H$ is involutive, that is,
for complex vector fields $U$ and $V$ in $\Cal H$,
$[U,V]$ is also in $\Cal H$. A manifold endowed
with a $CR$-structure is called a $CR$-manifold
[Greenfield 1968]. 

D. E. Blair and B. Y. Chen (1979)
proved that every $CR$-submanifold of a Hermitian
manifold is a $CR$-manifold.

In 1978  Chen proved the integrability
theorem for $CR$-submanifolds of a K\"ahler
manifold; namely, the totally real
distribution $\Cal D^\perp$ of a
$CR$-submani\-fold of a K\"ahler manifold
is always completely integrable. 
This theorem implies that every proper
$CR$-submanifold of a K\"ahler manifold is
foliated by totally real submanifolds. By
applying  this integrability theorem A.
Bejancu (1979) proved that a
$CR$-submani\-fold of a K\"ahler manifold is
mixed totally geodesic if and only if each
leaf of the totally real distribution is
totally geodesic in the $CR$-submanifold.

Chen's integrability theorem  was
extended to $CR$-submanifolds of various
families of Hermitian manifolds by various
geometers. For instance, this theorem was
extended to
$CR$-submanifolds of locally conformal symplectic
manifolds by Blair and Chen. Furthermore,
they  constructed $CR$-submanifolds in some
Hermitian manifolds with non-integrable 
totally real distributions [Blair-Chen 1979].

Let $M$ be a $CR$-submanifold with Riemannian
connection $\nabla$ and let $e_1,\ldots,e_{2h}$
be an orthonormal frame field of the holomorphic
distribution $\Cal D$. Put $\hat H=\,$trace $\hat
\sigma$, where
$\hat \sigma(X,Y)=(\nabla_XY)^\perp$ is the
component of
$\nabla_XY$ in the totally real distribution.
The holomorphic distribution  $\Cal D$ is called
 minimal if $\hat H=0$, identically. 

Although the holomorphic distribution is not
necessarily integrable in general, Chen
(1981c) proved that the holomorphic
distribution of a $CR$-submanifold is always a
minimal distribution.

Besides the minimality  Chen (1984a) also
proved   the following properties for the
holomorphic distributions: 

(1) If $M$ is a compact proper
$CR$-submanifold of a Hermitian symmetric space
of non-compact type, then  the holomorphic
distribution  is non-integrable.

(2) Let $M$ be a compact proper
$CR$-submanifold of the complex Euclidean space.
If the totally real distribution is a minimal
distribution, then the holomorphic distribution
is a non-integrable distribution.

A. Bejancu (1978)  obtained  a
necessary and sufficient condition for the
integrability of the holomorphic distribution:
Let $M$ be a $CR$-submanifold of a K\"ahler
manifold. Then the holomorphic distribution
$\Cal D$ is integrable if and only if the
second fundamental form of $M$ satisfies
$h(X,JY)=h(Y,JX)$ for any $X,Y$ tangent to $M$.

Chen (1981c) discovered a canonical cohomology
class $c(M)\in H^{2h}(M;\hbox{\bf R})$ for
every compact
$CR$-submanifold $M$ of a K\"ahler manifold. By
applying this cohomology class, he proved the
following: Let $M$ be a compact 
$CR$-submanifold of a K\"ahler manifold. If
the  cohomology group
$H^{2k}(M;\hbox{\bf R})=\{0\}$ for some
integer $k\leq h$, then either the
holomorphic distribution $\Cal D$ is not
integrable or the totally real distribution
$\Cal D^\perp$ is not  minimal. Chen's
cohomology class was used by S. Dragomir 
in his study concerning the minimality of
Levi distribution (cf. [Dragomir 1995,
Dragomir-Ornea 1998]).

A. Ros (1983) proved that if $M$ is an
$n$-dimensional compact minimal
$CR$-submani\-fold of $CP^m(4)$, then the
first nonzero eigenvalue of the Laplacian of
$M$ satisfies $\lambda_1(M)\leq
2(n^2+4h+p)/n.$

\subsection{Totally umbilical $CR$-submanifolds}

 Bejancu (1980) and  Chen
(1980c) investigated totally umbilical
$CR$-submanifolds of K\"ahler manifolds and
obtained the following: Let
$M$ be a $CR$-submanifold of a K\"ahler manifold.
If $M$ is totally umbilical, then either

(1) $M$ is totally geodesic, or

(2) $M$ is totally real, or

(3) the totally
real distribution is 1-dimensional, that is,
$p=1$.

Totally umbilical $CR$-submanifolds with
1-dimensional totally real distribution were
investigated in [Chen 1980c]. For instance, he
proved that  every totally umbilical
hypersurface of dimension $\geq 5$ in a K\"ahler
manifold has constant mean curvature. This result
is no longer true if the totally
umbilical hypersurface is 3-dimensional.

Totally geodesic $CR$-submanifolds of a K\"ahler
manifold are classified by Blair and Chen
(1979). In particular, they proved that if
$M$ is a totally geodesic $CR$-submanifold
of a K\"ahler manifold, then $M$ is a
$CR$-product. 

\subsection{Inequalities for $CR$-submanifolds}

For $CR$-submanifolds in complex space forms,
there is a sharp relationship between the
invariant $\delta_M={1\over 2}\rho-\inf K$ and the
squared mean curvature $H^2$   [Chen 1996a]:

 Let $M$ be an $n$-dimensional  
$CR$-submanifold in a complex space form $\tilde
M^m(4c)$. Then
\begin{align}  \delta_M \leq \begin{cases} 
\dfrac{n^2(n-2)}{2(n-1)} H^2 +\left\{
\dfrac{1}{2}(n+1)(n-2)+3h\right\}
c,\quad\hbox{if}\quad c>0;\\
\dfrac{n^2(n-2)}{2(n-1)}H^2,\quad \hbox{if}
\quad c=0;\\ \dfrac{n^2(n-2)}
{2(n-1)}H^2 + \dfrac{1}{2}(n+1)(n-2)c,\quad \hbox{if}
\quad c<0.\end{cases}\end{align}

There exist many $CR$-submanifolds in complex
space forms which satisfy the equality cases of
the above inequalities.

Proper $CR$-submanifolds of complex
hyperbolic spaces satisfying the
equality case were completely determined by
Chen and Vrancken (1997b) as follows:

 Let $U$ be a domain
of $\hbox{\bf C}$ and $\Psi:U\to
C^{m-1}$ be a nonconstant holomorphic curve
in $ C^{m-1}$. Define $z:E^2\times
U\to C^{m+1}_1$  by
$$z(u,t,w)=\Big(-1-{1\over
2}\Psi(w)\bar\Psi(w)+iu,-{1\over
2}\Psi(w)\bar\Psi(w)+iu,\Psi(w)\Big)e^{it}.$$
Then $\left<z,z\right>=-1$ and the image
$z(E^2\times U)$ in
$H^{2m+1}_1$ is invariant under the group action
of $H^1_1$. Moreover, away from points where
$\Psi'(w)=0$, the image
$\pi(E^2\times U)$, under the projection
$\pi:H^{2m+1}_1(-1)\to CH^m(-4)$, is a proper
$CR$-submanifold of $ CH^m(-4)$ which
satisfies $\delta_M = 
{{n^2(n-2)}\over{2(n-1)}}H^2 +
{1\over 2}(n+1)(n-2)c.$

Conversely, up to rigid motions of $
CH^m(-4)$, every proper $CR$-submanifold of
$CH^m(-4)$ satisfying the equality is
obtained in such way.
\smallskip
For a given submanifold $M$ of a K\"ahler
manifold,  let $P:TM\to TM$ denote the
tangential component of
$J:TM\to J(TM)$. 

M. Kon (1989)
proved the following: Suppose
$M$ is a compact orientable $n$-dimensional
minimal CR-submanifold of the complex projective
space $CP^m$. Suppose the Ricci tensor of
$M$ satisfies $$Ric(X, X)\ge (n-1)g(X,
X)+2g(PX, PX).\leqno (17.2)$$ Then $M$
is isometric to one of the following:

(1) a real projective space $
RP\sp n$;

(2) a complex projective space
$CP\sp {n/2}$; 

(3) a pseudo-Einstein real hypersurface
$\pi(S^{(n+1)/2}(\frac12)\times
S^{(n+1)/2}(\frac12)),$ of some

$CP^{(n+1)/2}$ in $CP\sp m$, where
$\pi:S^{2m+1}\to CP^m$ is the Hopf
fibration.

\subsection{$CR$-products}

The notion of $CR$-products was introduced
in [Chen 1981b]: A $CR$-submanifold
$M$ of a K\"ahler manifold
$\tilde M$ is called a $CR$-product if it is
locally a Riemannian product of a K\"ahler
submanifold $M^T$ and a totally real
submanifold $N^\perp$ of $\tilde M$.

Chen (1981b) showed that a submanifold $M$ of
a K\"ahler manifold is a $CR$-product if and
only if $\nabla P=0$, that is,  $P$ is
parallel with respect to the Levi-Civita
connection of $M$, where $P$ is the
endomorphism on the tangent bundle $TM$
induced from the almost complex structure
$J$ of $\tilde M$.

Let $f:M^\perp\to CP^p(4)$ be a Lagrangian
submanifold of
$CP^p(4)$. Then the composition $$CP^h\times
M^\perp @>i\times f >> CP^h\times CP^p 
@>S_{hp}>> CP^{h+p+hp}\leqno (17.3)$$ is a
$CR$-product in
$CP^{h+p+hp}$, where
$i:CP^h\to CP^h$ is the identity map and
$S_{hp}$ is the Segre embedding.

  A $CR$-product $M=M^T\times M^\perp$ in 
$CP^m$ is called a standard $CR$-product if
$m=h+p+hp$ and $M^T$ is a totally geodesic
K\"ahler submanifold of $CP^m$.

For $CR$-products in complex
space forms, Chen (1981b) proved the
following:

(1) A $CR$-product in a complex hyperbolic
space is non-proper, that is, it is either a
K\"ahler submanifold or a totally real
submanifold.

(2) A $CR$-product in  complex Euclidean
$m$-space $ C^m$ is a product submanifold
of a  complex linear subspace $ C^r$ of
$ C^m$ and a totally real submanifold in
a  complex linear subspace $
C^{m-r}$ of $ C^m$.

(3) If $M=M^T\times M^\perp$ is a 
$CR$-product of $CP^m(4)$, then 

(3.1) $m\geq h+p+hp$,

(3.2) the squared length $S$ of the second
fundamental form satisfies $S\geq 4hp$,

(3.3) if $m=h+p+hp$, then $M$ is a standard
$CR$-product, and

(3.4) if $S= 4hp$, then $M=M^T\times
M^\perp$ is a standard $CR$-product contained in a
totally geodesic K\"ahler submanifold
$CP^{h+p+hp}(4)$ of $ CP^m(4)$. 
Moreover, $M^T$ is an open portion of 
$CP^h(4)$ and $M^\perp$ is an open portion of 
$RP^p(1)$.

(4) If $M$ is a minimal $CR$-product in
$CP^m$, then the scalar curvature $\rho$ of
$M$ satisfies $$\rho\geq
4h^2+4h+p^2-p,\leqno(17.4)$$ with the equality
holding when and only when
$S=4hp$.

S. Maeda and N. Sato (1983) studied
$CR$-submanifolds
$M$ in a complex space form $\tilde M^m(4c)$
such that geodesics in $M$ are circles in
$\tilde M^m(4c)$ and obtained the following:
Let $M$ be a $CR$-submanifold in a complex space
form $\tilde M^m(4c)$. If  geodesics in
$M$ are circles in $\tilde M^m(4c)$, then  $M$ is
a $CR$-product.

\subsection{Cyclic parallel $CR$-submanifolds}

Concerning the covariant derivative of the second
fundamental form of 
$CR$-submanifolds of a complex space form $\tilde
M^m(4c)$, K. Yano and M. Kon (1980) (for
$c>0$) and Chen, G. D. Ludden and S. Montiel
(1984) (for
$c<0)$ proved  the following general inequality.

Let $M$ be a   $CR$-submanifold in a complex space form
$\tilde M^m(4c)$. Then the squared length of the
covariant derivative of the second fundamental
form satisfies
$$||\bar\nabla h||^2\geq 4c^2hp,\leqno(17.5)$$
with the equality holding if and only if
$M$ a cyclic-parallel
$CR$-submanifold, that is, $M$ satisfies
$$(\bar\nabla_X h)(Y,Z)+(\bar\nabla_Y
h)(Z,X)+(\bar\nabla_Z h)(X,Y)=0\leqno(17.6)$$
for $X,Y,Z$ tangent to $M$.

Let $H_1^{2m+1}(-1)$ denote the anti-de Sitter space time with constant sectional curvature
$-1$ and let $$\pi:H_1^{2m+1}(-1)\to CH^m(-4)\leqno(17.7)$$ denote the corresponding
Hopf fibration. For a submanifold $M$ of $CH^m(-4)$, let $\tilde M$ denote the pre-image of $M$. 

B. Y. Chen, G. D.  Ludden and S. Montiel
(1984) showed that a $CR$-submanifold
$M$ of $CH^m(-4)$ is cyclic-parallel if and only
if the preimage $\tilde M$ has parallel
second fundamental form in 
$H_1^{2m+1}(-1)$.  Similar result also holds for $CR$-submanifolds in $CP^m(4)$ [Yano-Kon
1983].

A submanifold of a real space
form is cyclic-parallel if and
only if it is a parallel submanifold.

A Riemannian manifold $M$ is called  a
two-point locally homogeneous space if it
is either flat or a rank one locally
symmetric space.  Chen and L. Vanhecke
(1981) proved that a Riemannian manifold is
a two-point locally homogeneous space if
and only if sufficiently small geodesic
hypersurfaces of $M$ are cyclic-parallel
hypersurfaces.

\subsection{ Homogeneous  and mixed foliate $CR$-submanifolds }

Y. Shimizu (1983) constructed homogeneous $CR$-submanifolds in $CP^n$ which are not
$CR$-products. Shimizu's results state as follows.

 Let $G/H$ be an irreducible Hermitian symmetric space of compact type. Denote by $\pi:S^{2n+1}(1)\to CP^n(4)$ the Hopf fibration. For a point $x\in S^{2n+1}$ denote by $N$ the $H$-orbit
of $x$ and $M=\pi(N)$. If the rank of $G/H$ is greater than one and if $N$ has the maximal dimension, then

(1) $M$ is a proper $CR$-submanifold
of $CP^n$ of codimension $rk(G/H)-1$,

(2)  $M$ is not a $CR$-product,

(3) $M$ has parallel mean curvature vector, and

(4)  $M$ has flat normal connection.

A $CR$-submanifold $M$ in a K\"ahler manifold is called mixed foliate if its holomorphic
distribution $\Cal D$ is integrable and its second fundamental form $h$ satisfies $h(X,Z)=0$
for $X$ in $\Cal D$ and $Z$ in $\Cal D^\perp$.

Mixed foliate $CR$-submanifolds in complex space forms were completely determined as follows.

(1) A complex projective space admits no mixed foliate proper $CR$-submani\-folds  [Bejancu-Kon-Yano 1981].

(2) A  $CR$-submanifold in $ C^m$ is mixed foliate if and only if it is a $CR$-product  [Chen 1981b].

(3) A  $CR$-submanifold in a complex hyperbolic space $CH^m$ is mixed foliate if and only if it is either a K\"ahler submanifold or a totally real submanifold  [Chen-Wu 1988].

\subsection{Nullity of $CR$-submanifolds}

T. Gotoh (1997) investigated the second variational formula of a compact minimal $CR$-submani\-fold in a complex projective space and estimated its nullity of the second variations to obtain the following.

Let $f:M\to CP^m$ be an $n$-dimensional compact minimal $CR$-submanifold of $CP^m$. 

(1) If $n$ is even, then the nullity of $M$ satisfies
$$n(f)\geq 2\left({n\over 2}+1\right)\left(
m-{n\over 2}\right),\leqno (17.8)$$ with equality
holding if and only if $M$ is a totally
geodesic K\"ahler submanifold;

(2) If $n$ is odd and equal to $m$, then the nullity of $M$ satisfies
$$n(f)\geq {{n(n+3)}\over 2},\leqno (17.9)$$ with
equality holding if and only if $M$ is a
totally real totally geodesic submanifold;

(3) If $n$ is odd and not equal to $m$, then
the nullity of $M$ satisfies
$$n(f)\geq n+1 +2\left({{n+1}\over 2}+1\right)\left( m-{{n+1}\over 2}\right),\leqno (17.10)$$ with equality holding if and only when
$$M=\pi\left( S^1\left(\sqrt{1\over {n+1}}\,\right)\times S^n\left(\sqrt{n\over
{n+1}}\,\right)\right)\subset CP^{(n+1)/2},$$ where $ CP^{(n+1)/2}$ is
embedded in $CP^m$ as a totally geodesic K\"ahler submanifold and
$\pi:S^{(n+1)/2+1}\to CP^{(n+1)/2}$ is the Hopf fibration.

\vfill\eject

\section{Slant submanifolds of K\"ahler manifolds}

Let $M$ be an $n$-dimensional Riemannian manifold
isometrically immersed in a K\"ahler manifold
$\tilde M$ with almost complex structure $J$ and
K\"ahler metric $g$.  For any vector $X$ tangent
to $M$ let $PX$ and $FX$ denote the tangential
and the normal components of $JX$,
respectively.  Then $P$ is an endomorphism of
the tangent bundle $TM$. For any nonzero vector
$X$ tangent to $M$ at a point $p\in M$, the
angle $\theta (X)$  between
$JX$ and the tangent space $T_pM$ is called the
 Wirtinger angle of $X$. 

A submanifold $M$ of $\tilde M$ is called slant if the Wirtinger angle $\theta (X)$ is constant
(which is independent of the choice of $x \in M$ and of $X \in T_{x}N$). The Wirtinger angle of a slant submanifold is called the  slant angle of the slant submanifold [Chen 1990]. 

K\"ahler submanifolds  and totally real submanifolds are nothing but slant
submanifolds with $\theta = 0$ and $\theta =
\pi /2$, respectively. A slant submanifold is
called proper if it is neither complex
nor totally real. 
In this sense, both $CR$-submanifolds
and slant submanifolds are generalizations of
both K\"ahler submanifolds and totally real
submanifolds.

Slant surfaces in almost Hermitian
manifolds do exist extensively. In fact,
Chen and Y. Tazawa (1990) proved the
following:

Let $f:M\to \tilde M$ be an embedding from
an oriented surface $M$ into an almost
Hermitian manifold $\tilde M$ endowed with
an almost complex structure $J$ and an
almost Hermitian metric $g$. If $f$ has no
complex tangent points, then for any
prescribed angle
$\theta\in (0,\pi)$, there exists an almost
complex structure $\tilde J$ on $\tilde M$
satisfying the following two conditions:

(i) $(\tilde M,g,\tilde J)$ is an almost
Hermitian manifold, and

(ii) $f$ is a $\theta$-slant surface with
respect to $(g,\tilde J)$.
 
By a complex tangent point of $f$, we mean
a point $x\in M$ such that the tangent
space $T_xM$ of $M$ at $x$ is invariant
under the action of the almost complex
structure $J$ on $M$.

\subsection{Basic properties of slant submanifolds} 

Proper  slant submanifolds are 
even-dimensional, such submanifolds do exist
extensively for any even dimension greater
than zero (cf. [Chen 1990, Tazawa 1994a,
1994b]).

Slant submanifolds of K\"ahler manifolds are
characterized by a simple condition; namely, $\,
P^{2}=\lambda I\,$ for a fixed real number
$\lambda \in [-1,0]$, where $I$ is the identity
map of the tangent bundle of the submanifold. 

A proper slant submanifold is called
K\"ahlerian slant  if the endomorphism $P$  
is parallel with
respect to the Riemannian connection, that is,
$\nabla P = 0$. A K\"ahlerian slant
submanifold is a K\"ahler manifold with
respect to the induced metric and  the
almost complex structure defined by
${\tilde J}= (\sec\theta)P$.

K\"ahler submanifolds, totally real submanifolds
and K\"ahlerian slant submanifolds satisfy
the condition: $\nabla P=0$. In general, let
$M$ be a submanifold of a K\"ahler manifold
$\tilde M$. Then $M$ satisfies $\nabla P=0$
if and only if $M$ is locally the Riemannian
product $M_1\times\cdots\times M_k$, where
each $M_i$ is a K\"ahler submanifold, a
totally real submanifold or a K\"ahlerian
slant submanifold of $\tilde M$.

Slant submanifolds have the following
topological properties:

(1)  If $M$  is a compact
$2k$-dimensional proper slant submanifold of a
K\"ahler manifold, then
$\,H^{2i}(M;\hbox{\bf R})\ne
\{0\}$ for $i=1,\ldots,k$ \ [Chen 1990].

(2) Let $M$ be a slant submanifold in a
complex Euclidean space. If $M$ is not
totally real, then $M$ is non-compact
[Chen-Tazawa 1991].

Although there do not exist  compact proper
slant submanifolds in complex Euclidean spaces,
there do exist compact proper slant submanifolds
in complex flat tori. 

The following result of Chen (1996c) provides
a Riemannian obstruction to the isometric
slant immersion in a flat K\"ahler manifold.

Let $M$ be a compact Riemannian
$n$-manifold with  finite fundamental group
$\pi_1(M)$. If there exists a $k$-tuple
$(n_1,\ldots,n_k)$ $\in \Cal S(n)$ such that 
$\delta(n_1,\ldots,n_k)>0$ on $M$, then $M$
admits no slant immersion into any flat
K\"ahlerian
$n$-manifold.

\subsection{Equivariant slant immersions}

S. Maeda, Y. Ohnita and S. Udagawa
(1993) investigated  slant immersions between
K\"ahler manifolds and obtained the following.

Let  $f:M\to N$ be an isometric immersion of an $m$-dimensional
compact K\"ahler manifold with K\"ahler form
$\omega_M$ into a K\"ahler manifold with K\"ahler form
$\omega_N$. Assume that the second Betti number
$b_2(M)=1$ and $f^*\omega_N$ is of type $(1,1)$.
Then the following three conditions are equivalent:

(1) $f$ is a slant immersion with slant angle
$\cos^{-1}({{|c|}\over m})$, for some nonnegative
constant $c$.

(2) $f^*\omega_N=({c\over m})\omega_M$,

(3) {\rm trace}$_g\,
f^*\omega_N=\sqrt{-1}c$ is a constant, where
$g$ is the K\"ahler metric on $M$.

A compact simply-connected homogeneous K\"ahler
manifold is simply called a K\"ahler $C$-space. Let
$f:M=G/H\to CP^m(4c)$ be an isometric immersion of
a compact homogeneous Riemannian manifold into
$CP^m(4c)$. The immersion $f$ is called
$G$-equivariant if there is a homomorphism
$\rho:G\to SU(m+1)$ such that $f(a\cdot
p)=\rho(a)f(p)$ for any $p\in M$ and $a\in G$. 

If $f:M\to CP^m(4c)$ is a map of a K\"ahler
manifold $M$ with $H_2(M;\hbox{\bf Z})\cong
\hbox{\bf Z}$. Denote by $\omega_M$ and
$\tilde
\omega$ the K\"ahler forms of $M$ and
$CP^m(4c)$, respectively. Let $S$ be a
positive generator of $H_2(M;\hbox{\bf Z})$.
Define the degree of $f$  by
deg$(f)={c\over \pi}[f^*\tilde\omega](S)$,
where
$[f^*\tilde\omega](S)$ is the evaluation of the
cohomology class $[f^*\tilde\omega]$ represented
by $f^*\tilde\omega$ at $S$.

The following theorem of Maeda,  Ohnita and
 Udagawa  provides some nice examples
of proper slant submanifolds in complex
projective spaces.

 Let  $f:M=G/H\to CP^m(4c)$
be a $G$-equivariant isometric immersion of an
$m$-dimensional K\"ahler $C$-space into a
complex projective space with K\"ahler form
$\tilde \omega$. Then  $f^*\tilde\omega$ is
of type
$(1,1)$ and {\rm trace}$_g\,\tilde \omega$ is
constant, where $g$ is the K\"ahler metric on
$M$. Moreover, if $b_2(M)=1$, then $f$ is a
slant immersion with slant angle given by
$$\cos^{-1}\left(|{deg}\,(f)|\cdot
{{\pi}\over{c\,vol(S)}}\right),$$ where
$S$ is a rational curve of $M$ which
represents the generator of
$H_2(M;\hbox{\bf Z})$.

\subsection{Slant surfaces in complex space forms}

Slant submanifolds of dimension two have some
special geometric properties. For instance, Chen
(1990) proved that a surface in a K\"ahler
manifold is a proper slant surface if and
only if it  is a K\"ahlerian slant surface. He
also showed that there do not exist flat
minimal proper slant surfaces in $C^2$. Also
Chen and Tazawa  proved in 1997 that there exist
no proper slant minimal surfaces in $CP^2$ and
in $CH^2$.

If the mean curvature of a complete
oriented proper slant surface in $ C^2$
is bounded below by a positive number, then
topologically it is either a circular
cylinder or a 2-plane [Chen-Morvan 1992].

Suppose $M$ is an immersed surface in a K\"ahler surface $\tilde M$ which is neither  K\"ahlerian  nor Lagrangian.  Then $M$ is a proper slant surface of
$\tilde M$ if and only if the shape operator of $M$ satisfies
$$A_{FX}Y=A_{FY}X\leqno (18.1)$$ for vectors
$X,Y$ tangent to $M$.

Applying this special property of the shape
operator for slant surfaces, Chen
(1995,1998b,1998d) proved that the squared
mean curvature and the Gaussian curvature
of a proper slant surface in a 2-dimensional
complex space form
$\tilde M^2(4c)$ satisfies $$H^2\geq
2K-2(1+3\cos^2\theta)c,\leqno(18.2)$$ where
$\theta$ denotes  the slant angle.

 There do not exist  proper
slant surfaces satisfying the equality case of
 inequality (18.2) for $c>0$. A proper slant
surface in a flat K\"ahler surface satisfies
the equality of inequality (18.2) if and
only if it is totally geodesic. Furthermore,
a proper slant surface in the complex
hyperbolic plane
$CH^2(-4)$ satisfying the equality case of 
inequality (18.2) is a surface of constant
Gaussian curvature $-{2\over 3}$ with slant
angle $\theta=\cos^{-1}({1\over 3})$.
Moreover, the immersion of such a slant
surface is rigid.

A submanifold $N$ of a pseudo-Riemannian
Sasakian manifold $(\tilde M,g,\phi,\xi)$
is called contact $\theta$-slant if the
structure vector field
$\xi$ of $\tilde M$ is tangent to
$N$ at each point of $N$ and, moreover, for
each unit vector $X$ tangent to $N$ and
orthogonal to $\xi$ at $p\in N$, the angle
$\theta(X)$ between $\phi(X)$ and $T_pN$ is
independent of the choice of $X$ and $p$.

Let $H^{2m+1}_1(-1)\subset C^{m+1}_1$ denote the
anti-de Sitter space-time and $\pi\colon\;
H^{2m+1}_1(-1)$ $\to  CH^m(-4)$ the corresponding
totally geodesic fibration (cf.
 section 16.3). Then every $n$-dimensional proper
$\theta$-slant submanifold $M$ in
$CH^{m}(-4)$ lifts to an $(n+1)$-dimensional
proper contact $\theta$-slant submanifold
$\pi^{-1}(M)$ in $H^{2m+1}_1(-1)$ via $\pi$.
Conversely, every proper contact
$\theta$-slant submanifold of
$H^{2m+1}_1(-1)$ projects to a proper
$\theta$-slant submanifold of $CH^{m}(-4)$ via
$\pi$. Similar correspondence also holds between
proper $\theta$-slant submanifolds of $CP^m(4)$
and proper contact $\theta$-slant submanifolds of
the Sasakian unit $(2m+1)$-sphere $S^{2m+1}(1)$.

The contact slant representation of the unique
proper slant surface in $CH^2(-4)$
which satisfies the equality case of (18.2)  in
$H^5_1(-1)\subset C^3_1$ has been determined by
Chen and Y. Tazawa in 1997. Up to rigid
motions of
$C^3_1$, this contact slant representation is
given by
$$z(u,v,t)=e^{it}\Big(1+\frac32\Big(\cosh \sqrt{2\over3}v-1\Big)+\frac{u^2}6e^{-\sqrt{2\over3}
v}-i{u\over\sqrt{6}} (1+e^{-\sqrt{2\over3}v}),$$ 
$$ \frac u3\Big( {1+2e^{-\sqrt{2\over3}v}}\Big)+\frac
i{6\sqrt{6}}e^{-\sqrt{2\over3}v}\Big(\Big(e^{\sqrt{2\over3}v}
-1\Big)\Big(9e^{\sqrt{2\over3}v}-3\Big)
+2u^2\Big),\leqno(18.3)$$ 
$$ \frac u{3\sqrt{2}}\Big(1-e^{-\sqrt{2\over3}v}\Big)+\frac i{12\sqrt{3}} \Big(6-15 e^{-\sqrt{2\over3}v}+
9e^{\sqrt{2\over3}v} +2e^{-\sqrt{2\over3}v}u^2\Big)\Big).$$
\vskip.1in

In 1990 Chen classified slant surfaces in
$C^2$ with parallel mean curvature vector:

Let $M$ be a  slant surface in $ C^2$ with
parallel mean curvature vector. Then
$M$ is one of the following surfaces:

(1) an open portion of the product surface of
two plane circles;

(2) an open portion of a circular cylinder
which is contained in a real hyperplane of
$ C^2$;

(3) a minimal slant surface.

Cases (1) and (2) occur only when $M$ is a
Lagrangian surface of $C^2$.

 J. Yang (1997) showed that a flat proper
slant surface with nonzero constant mean
curvature in
$ C^2$ is an open portion of a helical
cylinder and there do not exist proper slant
surfaces with nonzero constant mean curvature and
nonzero constant Gauss curvature in $ C^2$.

Y. Ohnita (1989) proved that totally geodesic
surfaces are the only  minimal slant surfaces
with constant Gauss curvature in complex
hyperbolic spaces. In contrast, Chen and
Vrancken (1997a) proved that for each 
constant $\theta,\,0<\theta<{\pi\over 2}$,
there exist complete $\theta$-slant surfaces
in the complex hyperbolic plane $CH^2$ with
nonzero constant mean  curvature and
constant negative Gaussian curvature.

Chen and Vrancken (1997a) also
proved the following:

(1) For a given  constant 
$\theta\,$ with $0<\theta\leq {\pi\over 2}$ and a given  function 
$\lambda$, there exist
infinitely many $\theta$-slant surfaces in 
$C^2$ with $\lambda$ as the prescribed mean
curvature function.

(2) For a given constant $\theta$ with
$0<\theta\leq {\pi\over 2}$ and  a given
function $K$, there exist infinitely many
$\theta$-slant surfaces in $C^2$ with $K$
as the prescribed Gaussian curvature function.

Slant surfaces in $ C^2$ were
completely classified by Chen and Y. Tazawa
(1991) for the following cases: 

(1) spherical slant surfaces;

(2) slant surfaces lying in a real
hyperplane of $ C^2$; or

(3) slant surfaces whose Gauss map has rank less than two.

For  case (1), they proved that a spherical
 surface in $ C^2$ is proper slant if
and only if it is locally a spherical
helical cylinder in a hypersphere $S^3$;
for case (2), the surfaces are doubly slant
and they are the unions of some open portions
of planes, circular cones and the tangent
developable surfaces obtained by
generalized helices; and for case (3) the
slant surfaces are unions of some special
flat ruled surfaces.

\subsection{Slant surfaces and almost complex structures}

Let ${ C}^{2}=( E^{4},J_{0})$ be the
complex Euclidean plane with the canonical
complex structure $J_0$. Then
$J_0$ is an orientation preserving isomorphism. 
Denote by $\Cal J$ the set of all almost complex
structures 
on $ E^4$ which are compatible with the
inner product $\left<\, ,\, \right>$, that
is, ${\Cal J}$ consists of all linear
endomorphisms $J$ of $
 E^{4}$ such that $ J^{2}= -I$ and
$\left<JX,JY\right>=\left<X,Y\right>$ for $X,Y$
$\in E^4$

An orthonormal basis $\{e_{1},e_{2},e_{3},e_{4}\}$
 on $ E^4$ is called a $J$-basis if
$Je_{1}=e_{2}, Je_{3}=e_{4}$. Any two $J$-bases
associated with the  same almost complex structure
have the same orientation. 

With respect to the
canonical orientation on
$ E^4$ one can divide
$\Cal J$ into two disjoint subsets $\Cal J^+$
and $\Cal J^-$ which consist of all positive and
all negative $J$-bases, respectively.

For an immersion $\phi$ of a Riemann surface $M$
into a K\"ahler manifold $N$, the K\"ahler angle
$\alpha$ of $\phi$ is defined to be the angle
between $J\phi_*(\partial/\partial x)$ and 
$\phi_*(\partial/\partial y)$, where
$z=x+\sqrt{-1}y$ is a local complex coordinate on
$M$ and $J$ the almost complex structure on $N$.

The relation between  $\theta$ and the
K\"ahler angle $\alpha$ for an immersion $\phi$
of a Riemann surface $M$ into a K\"ahler manifold
$N$ is $$\theta(X)=\min\,\{\alpha(T_pM),\pi-
\alpha(T_pM)\}$$
for any nonzero vector $X\in T_pM$.

The immersion $f$ of a Riemann surface in $N$ is
called holomorphic (respectively,
anti-holomorphic) if $\alpha\equiv 0$
(respectively, $\alpha\equiv\pi$). 

The following results of B. Y. Chen and Y.
Tazawa (1990) determine whether a surfaces
in $E^4$ is slant with respect to some
compatible almost complex structure on $E^4$:

(1)  Let $f : N\to E^4$ be a minimal
immersion. If there exists a compatible
complex structure ${\hat J}\in {\Cal J}^+$
(respectively, ${\hat J} \in {\Cal J}^-$)
such that the immersion is slant with
respect to $\,\,\hat J$, then 

(1-a) for any $\alpha \in [0,\pi]$, there
is a compatible complex structure $J_{\alpha}\in
{\Cal J}^+$ (respectively, $J_{\alpha}\in {\Cal
J}^-$) such that
$f$ is $\alpha$-slant with respect to the
complex structure $J_{\alpha}$, and

(1-b) the immersion $f$ is slant with
respect to any complex structure $J \in {\Cal J}^+$
(respectively, $J \in {\Cal J}^-$).

(2) If $f : N \to E^4$ is a non-minimal
immersion, then there exist at most two complex
structures $\,\,\pm J^{+} \in {\Cal J}^+$ and at
most two complex structures $\,\,\pm J^{-}\in
{\Cal J}^-$ such that the immersion $f$ is slant
with respect to them.

(3) If $f : N \rightarrow { 
C}^{2} = ( E^{4},J_{0})$ is  holomorphic,
then the immersion $f$ is slant with
respect to every complex structure $J \in
{\Cal J}^+$.

(4) If $f : N\rightarrow { C}^{2} =( E^{4},J_{0})$ is  anti-holomorphic, then the immersion $f$ is slant with respect to every complex structure $J \in {\Cal J}^-$. 

(5) If $f : N \to E^3$ is a non-totally geodesic minimal immersion, then $f :  N
\to E^{3} \subset E^4$ is not slant with respect to every
compatible complex structure on $ E^4$.

\subsection{Slant spheres in complex projective spaces}

For each $k=0,\ldots,m$, let $\psi_k:S^2\to
CP^m(4)$ be given by
$$\psi_k([z_0,z_1])=\left[g_{k,0}\left( {{z_0}\over{z_1}} \right),\ldots, g_{k,m}
\left({{z_0}\over{z_1}}\right)\right],\leqno
(18.4)$$ where $[z_0,z_1]\in CP^1=S^2$, and for
$j=0,\ldots,m$, $g_{k,j}(z)$ is given by
$$g_{k,j}(z)={{k!}\over{(1+z\bar
z)^k}}\sqrt{{m\choose j}}z^{j-k}\sum_p(-1)^p 
{j\choose {k-p}}{{m-j}\choose p}(z\bar
z)^p.\leqno(18.5)$$ It was proved by 
 Bolton, Jensen, Rigoli and Woodward (1988) that
each $\psi_k$ is a conformal minimal immersion
with constant Gaussian curvature
$4(m+2k(m-k))^{-1}$ and constant K\"ahler angle
$\alpha_k$ given by $$\tan^2\left(
{{\alpha_k}\over 2}\right)={{k(m-k+1)}\over
{(k+1)(m-k)}}.$$

Each $\psi_k$ is an embedding unless $m=2k$, in
which case $\psi_k$ is a totally real
immersion. 

The immersions $\psi_0,\ldots,\psi_m$ defined
above are called the Veronese sequence.
For Veronese sequence, J. Bolton, G. R. Jensen,
M. Rigoli and L. M. Woodward (1988) proved that

(1) Let $\psi:S^2\to CP^m(4)$ be a conformal
minimal immersion with constant Gaussian
curvature and assume that
$\psi(S^2)$ is not contained in any hyperplane of
$CP^m(4)$. Then, up to a holomorphic isometry of
$CP^m(4)$, the immersion $\psi$ is an
element of the Veronese sequence.

(2) Let  $\psi,\psi':S^2\to CP^m(4)$ be
conformal minimal immersions. Then $\psi,\psi'$
differ by a holomorphic isometry of $CP^m(4)$ if
and only if they have the same K\"ahler angle and
induced metrics at each point.

(3) If $\psi:S^2\to CP^m(4)$ is a totally
real minimal immersion, then $\psi$ is
totally geodesic. 

For each linearly full minimal immersion $\psi: S^2\to CP^n$, let 
$\psi_0,\psi_1,\cdots,\psi_n$ denote the corresponding Veronese sequence with
$\psi=\psi_k$ for some $k=0,1,\cdots,n$, where $\psi_0$ is holomorphic, called the
directrix of $\psi$. $\psi$ is called a
minimal immersion with position $k$. 

Z. Q. Li (1995) proved  the following: 

(1) Let $\psi\colon S\sp 2\to CP^n$
be a linearly full minimal immersion with position 2. Suppose the K\"ahler angle $\alpha$
is constant but the Gaussian curvature is not. If
$\alpha\not=0,\pi,\pi/2$ and the directrix $\psi_0$ of $\psi$ is unramified, then
$n\leq 10$ and $\tan^2(\alpha/2)=\frac 34$. 

(2) There are at least three
families of totally unramified minimal immersions $\psi: S^2\to CP^{10}$ such
that
$\psi$ is neither holomorphic, anti-holomorphic nor totally real, with constant K\"ahler
angle and nonconstant Gaussian curvature. Moreover, $\psi$ is homotopic to the Veronese
minimal immersion.

Y. Ohnita (1989) studied minimal surfaces
with constant curvature and constant K\"ahler
angle. He obtained the following.

 Let $M$ be a  minimal surface with constant
Gaussian curvature $K$ immersed fully in
$CP^m(4)$. If the K\"ahler angle $\alpha$ of $M$
is constant. 

(1) If $K>0$, then there exists some
constant $k$ with $0\leq k\leq m$ such that
$$K={4\over{2k(m-k)+m}},\quad\cos\alpha=\left({{
m-2k}\over 4}\right)\,K$$
and $M$ is an open submanifold of an element
of Veronese sequence.

(2) If $K=0$, then $M$ is totally real.

(3) $K<0$ is impossible.

A minimal surface in $CP^n$ is called superconformal if its harmonic sequence is
orthogonally periodic, and it is called pseudo-holomorphic if its harmonic sequence terminates at each end.

M. Sakaki (1996) proved the following:

 (1) Any superconformal minimal slant surface in $CP^3$ is
totally real. 

(2) Any pseudo-holomorphic minimal slant surface in $CP^4$ is either
holomorphic, anti-holomorphic, totally real or of constant curvature.

\vfill\eject

\section{Submanifolds of the nearly K\"ahler 6-sphere}

It was proved by E. Calabi (1958) 
that any oriented  submanifold
$M^6$ of the hyperplane  $\hbox{Im}\,\Cal O$ of
the imaginary octonions  carries a
$U(3)$-structure (that is, an almost Hermitian
structure). For instance, let
$S^6\subset \hbox{Im}\,\Cal O$ be the sphere of
unit imaginary vectors; then the right
multiplication by $u\in S^6$ induces a
linear transformation $J_u: {\Cal O} \to\Cal
O$ which is orthogonal and satisfies
$(J_u)^2= - I$. The operator $J_u$
preserves the 2-plane spanned by 1 and $u$ and
therefore preserves its orthogonal 6-plane
which may be identified with $T_uS^6$.
Thus $J_u$ induces an almost complex structure
on $T_uS^6$ which is compatible with the
inner product induced by the inner product of
$\Cal O$ and $S^6$ has an almost
complex structure.

The almost complex structure $J$ on $S^6$ is a
 nearly K\"ahler structure in
the sense that the (2,1)-tensor field $G$ on
$S^6$, defined by
$G(X,Y) = (\widetilde \nabla_XJ)(Y),$ is
skew-symmetric,
where $\widetilde \nabla$ denotes the Riemannian
connection on $S^6$. 

The group of automorphisms
of this nearly K\"ahler structure is the
exceptional simple Lie group $G_2$ which acts
transitively on $S^6$ as a group of isometries.

 A. Gray (1969) proved the following:

(1)  every almost complex
submanifold of the nearly K\"ahler $S^6$ is a
minimal submanifold, and 

(2) the nearly K\"ahler $S^6$ has no
4-dimensional almost complex submanifolds. 

\subsection{ Almost complex curves}

Almost complex curves, that is, real
2-dimensional almost complex submanifolds, in
the nearly K\"ahler $S^6$ have been studied by
various authors.

An almost complex curve in $S\sp 6$ is a
non-constant smooth map $f\colon S\to S\sp 6$,
from a Riemann surface $S$, whose differential
is complex linear. Such a map is necessarily a
weakly conformal harmonic map or, equivalently,
a weakly conformal branched minimal immersion.

Almost complex curves 
have ellipse of curvature a
circle, that is, the map
$v \mapsto h(v,v)$  describes a circle in
the normal space where $v \in UM_p$, $UM_p$
the unit hypersphere of $T_pM$, and
$h$ denotes the second  fundamental
form.  If the map
$v \mapsto (\bar\nabla_v h)(v,v)$,
$v \in UM_p$ also describes a circle, 
then $M^2$ is called superminimal.
This class of almost complex curves has been
investigated by Bryant (1982).  

In [Bryant 1982] a  Frenet formalism for
almost complex curves in $S^6$ was
developed; the first, second and third
fundamental forms are defined as holomorphic
sections of line bundles over $S$. In
particular, he showed that the third
fundamental form $III$, analogous to the
torsion of a space curve, plays a crucial
role. 

The hypothesis $III\ne 0$ is very restrictive
and Bryant proved that if $S= CP^1$,
then $III\ne 0$ is impossible. On the contrary,
he constructed almost complex curves $f :S\to
S^6$ with $III=0$ for any Riemann
surface $S$ such that the ramification
divisor of $f$ has arbitrarily large degree. He
also proved that if $S$ is compact and $f:S\to
S^6$ is an almost complex curve with $III=0$,
then $f$ is algebraic. In particular, it is
real analytic.

J. Bolton, L. Vrancken and L. M. Woodward (1994)
proved that there are four basic types of almost
complex curves in $S^6$;
namely,
\vskip.1in

(i)  linearly full in $S\sp 6$ and superminimal,

(ii) linearly full in $S^6$ but not superminimal, 

(iii) linearly full in some totally geodesic $S^5$ in $S^6$, and 

(iv) totally geodesic. 

\vskip.1in
They also provided metric criteria for
recognizing when a weakly conformal harmonic
map $f:S\to S^6$ is $O(7)$-congruent to an
almost complex curve of any one of the four
types. 

A surface $S$ in $S^{2m}$ is called
superminimal if $S$ is the image of a
horizontal holomorphic curve in the
Hermitian symmetric space
$SO(2m+1)/U(m)$ under the Riemannian submersion $\pi : 
SO(2m+1)/U(m) \to S^{2m}$ induced by the inclusions  $U(m)\subset
SO(2m)\subset SO(2m+1)$.
E. Calabi (1967) proved that
all minimal $2$-spheres in $S^{2m}$ are superminimal.

R. Bryant (1982) used a twistor
construction to obtain all almost complex 
curves in $S^6$ of type (I). Bryant's
construction involves consideration of the
twistor bundle $\pi:Q^5\to S^6$, where
$Q^5$ denotes the quadric in $CP^6$  given by
$$
Q^5 = \{[a+ib] \in C P^6:a,b \in E^7,\
|a|=|b|
\hbox{ and } a\perp b\}
$$
and $\pi:Q^5\to S^6$ is the map
$$
\pi([a+ib])=-{a\times b\over |a\times b|}.
$$

An alternative description of $\pi$ may be given in terms of a
projection between quotient spaces of $G_2$
as follows. Recall that the standard action
of $G_2$ on $E^7$ induces transitive actions
on $S^6$ and $Q^5$. If $\{e_1,\ldots ,e_7\}$
is the standard basis of $E^7$, then the
stabilizer of
$e_4 \in S^6$ may be identified with $SU(3)$ and
that of $[e_1+ie_5] \in Q^5$ with $U(2)$.
Using these identifications, if $g\in
U(2)$, then $$g(e_4)=-g(e_1\times
e_5)=-e_1\times e_5=e_4,$$ so that
$U(2)\subset SU(3)$. The above projection
$\pi$ may also be described as the standard
projection map
$\pi: Q^5= G_2/U(2) \to G_2/SU(3)=S^6.$

A holomorphic map $g:S\to Q^5$ is said to be
 superhorizontal if $g\times
g_z=0$, where $z=x+iy$ is a local complex
coordinate on $S$. Such maps are horizontal in the
sense that $g(S)$ intersects the fibres
orthogonally.  Although the map $\pi$ above is
not a  Riemannian submersion, if $g:S\to
Q^5$ is  holomorphic and superhorizontal,
the metrics induced on $S$  by $g$ and
$\pi\circ g$ are equal.
 
R. Bryant (1982) showed that there is a
one-to-one correspondence between linearly
full superhorizontal maps $g:S\to Q^5$ and
almost complex curves $f:S\to S^6$ of type
(I), where $g$  corresponds to the map
$f=\pi\circ g$. 

Bryant (1982) also gave a
``Weierstrass representation'' theorem for
almost complex curves in $S^6$ and proved 
that every compact Riemann surface admits an
infinite number of almost complex maps 
of type (I) into $S^6$.

Almost complex curves of
genus zero are necessarily of type (I) or (IV),
which have been studied by N. Ejiri (1986a)
who described all $S^1$-symmetric examples. 

A description of almost complex
curves of types (II) was given by  Bolton,
F. Pedit and Woodward  (1995). They showed 
 that almost complex curves of
type (II) all arise from solutions of the 
affine  2-dimensional $G_2$-Toda
equations.

 The method of construction of almost complex
curves of type (III) uses the Hopf
 fibration $\pi\colon
S^5\to CP^2$ to obtain 
these complex curves
horizontal lifts of suitable totally real branched minimal
immersions into $CP^2$.  In fact, if $f:S\to S^6$ is an almost complex
curve of type (III), by applying an
element of
$G_2$ if necessary we may assume that
$f(S)\subset S^5=S^6\cap P$ where $P$ is
 the hyperplane of $E^7$
given by 
$P=\{(x_1,\ldots ,x_7)\,:\,x_4=0\}.$
Then the map $\pi : S^5 \to CP^2$ given
by
$$\pi (x_1,x_2,x_3,0,x_5,x_6,x_7) = [x_1 +ix_5,x_2+ix_6,x_3+ix_7],
$$
is the Hopf fibration and it is shown in
[Bolton-Vrancken-Woodward 1994] that
$f$ is horizontal and $\pi\circ f: S
\to  CP^2$ is a totally real non-totally
geodesic weakly conformal harmonic map of
$S$ into $CP^2$. 

Conversely, if $S$ is simply-connected and if
$\psi:S\to CP^2$ is such a map, then among the
horizontal lifts of $\psi$ there are exactly three
almost complex curves (which are
all $G_2$-congruent). 

Using the theory of calibrations,
Palmer (1997) gave estimates for the
nullity and Morse index  of almost complex
curves in the nearly K\"ahler 6-sphere.

\subsection{Minimal surfaces of constant curvature in nearly K\"ahler 6-sphere}

The minimal surfaces of constant curvature in
the $6$-sphere have been known for some time.

R. L. Bryant (1985) proved that there are no
minimal surfaces in $S^n$ (in particular, in
$S^6$) of constant negative Gaussian
curvature. Flat minimal surfaces in $S^n$
were classified by K. Kenmotsu (1976). 

K. Sekigawa (1983) studied almost complex
curves of constant curvature in the nearly
K\"ahler 6-sphere and proved that, if the
Gaussian curvature
$K$ of an almost complex surface $M$ in the
nearly K\"ahler $S\sp 6$ is constant on $M$, then
$K=1$, $K=\tfrac16$ or $K=0$. Moreover, up to
$G_2$-congruence and conformal transformation
of the domain, the immersion $f:M\to S^6$ of
the almost complex curve is one of the
following:

\smallskip
(1) $K=1,\, M=S^2=\{(x,y,z)\in  E^3:x^2+y^2+z^2=1\}$ and
$$f(x,y,z)=(x,y,z,0,0,0,0),$$ for $(x,y,z)\in S^2$,
\smallskip

(2) $K={1\over 6},\, M=S^2$, and
$$ f(x,y,z)={1\over{2\sqrt{2}}}(\sqrt{3}x(-x^2
-y^2+4z^2),\sqrt{30}z(x^2
-y^2),$$ $$ \sqrt{5}y(3x^2 -y^2),
\sqrt{2}z(-3x^2-3y^2+2z^2),$$ $$
-\sqrt{3}y(-x^2-y^2+4z^2),
-2\sqrt{30}xyz,\sqrt{5}x(-x^2+3y^2)),$$
\smallskip
(3) $K=0,\, M= C^1,$ and
$$f(w)={1\over{\sqrt{6}}}\sum_{j=1}^3\left(
e^{\mu_jw-\overline{\mu_jw}}\nu_j+
e^{-\mu_jw+\overline{\mu_jw}}\bar\nu_j\right),$$
where $$\mu_1=1,\,\mu_2=\exp
({2\pi\over 3}),\mu_3=\exp ({{4\pi i}\over3}),
\nu_j={1\over\sqrt{2}}(e_j+ie_{j+4}),$$
for $j=1,2,3$; and 
$e_1,\ldots,e_7$ is the standard
basis of $ E^7$.

F. Dillen, B. Opozda, L. Verstraelen and L.
Vrancken (1987b) studied almost complex
surfaces $M$ in the nearly K\"ahler $S^6$ and
proved that

(1) If $0\leq K\leq 1/6$, then $K$ is
constant and $K=0$  or $K=1/6$.
 
 (2) If $1/6\leq K\leq 1$, then $K$ is
constant and $K=1/6$  or $K=1$.

The condition $K\leq 1$ in (2) is not
necessary, since it is always satisfied for
an almost complex  surface. (2) improves an
earlier result of Sekigawa (1983).

\subsection{Hopf hypersurfaces and almost complex curves}

Suppose that $M$ is a real
hypersurface in the nearly K\"ahler $S^6$.
Applying the almost complex structure 
$J$ on $S^6$ to the normal bundle of $M$, one
obtains a 1-dimensional distribution in $M$.
The 1-dimensional foliation induced by this
distribution is called the Hopf foliation, and
$M$ is said to be a Hopf hypersurface if this
foliation is totally geodesic. 

J. Berndt, J. Bolton and L. Woodward (1995)
proved that
a  hypersurface of the nearly
K\"ahler $S^6$ is a Hopf hypersurface if and
only if it is an open part of a tube around
an almost complex submanifold of $S\sp 6$. 

As the nearly K\"ahler six-sphere admits no
four-dimensional almost complex submanifolds,
this implies that $M$ is a Hopf hypersurface if
and only if it is an open part of either a
geodesic hypersphere or a tube around an almost
complex curve. As a
consequence, every Hopf hypersurface of $S\sp
6$ has exactly 1, 2, or 3 distinct
principal curvatures at each point. 

In the case
where $M$ is umbilical, it is an open subset of
a geodesic hypersphere. The Hopf hypersurface 
$M$ has exactly 2 distinct principal
curvatures if and only if $M$ is an open part
of a tube around a totally geodesic almost
complex curve in the nearly K\"ahler $S^6$. 

\subsection{Lagrangian submanifolds in nearly K\"ahler
6-sphere}.

\vskip.1in
\noindent{\bf  19.4.1. Ejiri's theorems for
Lagrangian submanifolds in $S^6$}

A 3-dimensional submanifold $M$ of the nearly
K\"ahler $S^6$ is called Lagrangian if the
almost complex structure $J$ on the nearly
K\"ahler 6-sphere carries each tangent space
$T_xM,\,x\in M$ onto the corresponding
normal space $T^\perp_xM$.

 N\. Ejiri (1981) proved that
a Lagrangian submanifold $M$ in $S^6$ is always
minimal and orientable.  He also proved that
if $M$ has constant sectional curvature, then
$M$ is either totally geodesic or has
constant curvature $1/16$. The first 
nonhomogeneous examples of Lagrangian
submanifolds in the nearly K\"ahler 6-sphere
were  described in [Ejiri 1986b].

F. Dillen, B. Opozda, L. Verstraelen and L.
Vrancken (1987a) proved that if $M$ is a
compact Lagrangian submanifold of
$S^6$ with $K>1/16$, then $M$ is a totally
geodesic submanifold. 

\vskip.1in
\noindent{\bf 19.4.2. Equivariant Lagrangian
submanifolds in $S^6$}

 K. Mashimo (1985) classified the
$G_2$-equivariant Lagrangian submanifolds
$M$ of the nearly K\"ahler 6-sphere. It
turns out that there are five models, and
every equivariant  Lagrangian submanifold in
the nearly K\"ahler 6-sphere is
$G_2$-congruent to one of the five models. 

These five models can be distinguished by the
following curvature properties:
\vskip.1in

 (1) $M^3$ is totally geodesic ($\delta(2)=2$),

(2)  $M^3$ has constant curvature $1/16$ ($\delta(2)=1/8$),

(3)  the curvature of $M^3$ satisfies  $1/16\leq K \leq 21/16$ ($\delta(2)=11/8$),
 
 (4)   the curvature of $M^3$ satisfies $-7/3\leq K \leq 1$ ($\delta(2)=2$), 

(5) the curvature of $M^3$ satisfies $-1\leq K \leq 1$ ($\delta(2)=2$).
\vskip.1in

F. Dillen, L. Verstraelen and L. Vrancken (1990)
characterized models (1), (2) and (3) as the
only compact Lagrangian submanifolds in $S^6$
whose sectional curvatures satisfy $K\geq
1/16$. They also obtained an explicit
expression for the Lagrangian submanifold of
constant curvature $1/16$ in terms of harmonic
homogeneous polynomials of degree 6. Using
these formulas, it follows that the immersion
has degree 24. Further, they also obtained an
explicit expression for model (3).

It follows from inequality (3.17) and Ejiri's result that the invariant
$$\delta(2):={\rho\over 2}-\inf K\leqno(19.1)$$
always satisfies $\delta(2)\leq 2$, for every Lagrangian submanifold of the nearly K\"ahler $S^6$. 

The models (1), (4) and (5) satisfy the equality $\delta(2)=2$ identically.  Chen, Dillen, Verstraelen and Vrancken (1995a) proved that these three models are the only Lagrangian
submanifolds of  the nearly K\"ahler $S^6$ with constant scalar  curvature that satisfy the
equality $\delta(2)=2$. 

Many further examples of Lagrangian
submanifolds in the nearly K\"ahler
$S^6$ satisfying the equality $\delta(2)= 2$
have been constructed in
[Chen-Dillen-Verstraelen-Vrancken 1995a,
1995b]. 

\vskip.1in
\noindent{\bf 19.4.3. Lagrangian submanifolds 
in $S^6$ satisfying $\delta(2)= 2$}

A Riemannian $n$-manifold $M$ whose Ricci
tensor has an eigenvalue of multiplicity at
least $n-1$ is called quasi-Einstein. R. 
Deszcz, F. Dillen, L. Verstraelen  and
 L. Vrancken (1997) proved that Lagrangian
submanifolds of the nearly K\"ahler 6-sphere
satisfying $\delta(2)=2$ are quasi-Einstein.

The complete classification of
Lagrangian submanifolds in the nearly
K\"ahler 6-sphere satisfying the equality 
$\delta(2)= 2$ was established by
Dillen and Vrancken (1996). More precisely,
they proved the following:

(1) Let $\phi: N_1 \to CP^2(4)$ be a
holomorphic curve in $CP^2(4)$, $PN_1$ 
the circle bundle over $N_1$ induced by the Hopf
fibration $\pi: S^5(1) \to
CP^2(4)$, and $\psi$ the isometric
immersion such that the following diagram
commutes:
\begin{align}\notag \CD PN_1 @>\psi>> S^5(1)\\ @V V\pi V @VV\pi V\\ N_1 @>\phi>> \Bbb CP^2(4) .\endCD \end{align}

Then, there exists a totally geodesic
embedding $i$ of $S^5$ into the nearly K\"ahler
6-sphere such that the
immersion $i \circ \psi : PN_1 \rightarrow S^6$
is a 3-dimensional Lagrangian immersion in $S^6$
satisfying equality $\delta(2)=2$.

(2) Let $\bar\phi: N_2 \to S^6$ be an
almost complex curve
(with second fundamental form $h$)  without
totally geodesic
points. Denote by $UN_2$ the
unit tangent bundle over $N_2$ and define a map
$$\bar\psi : UN_2 \rightarrow S^6 : v \mapsto
\bar\phi_\star(v)
\times {{h(v,v)}\over{\Vert h(v,v)
\Vert}}.\leqno{19.2}$$
Then $\bar\psi$ is a (possibly branched)
Lagrangian immersion into $S^6$ satisfying
 equality $\delta(2)=2$. Moreover, the
immersion is linearly full in $S^6$.

(3) Let $\bar\phi: N_2\to S^6$ be a
(branched) almost complex immersion.  Then, $SN_2$ is a
3-dimensional
(possibly branched) Lagrangian submanifold of
$S^6$ satisfying equality  $\delta(2)=2$.

(4) Let $f : M\to S^6$ be a Lagrangian
immersion which is not linearly full in $S^6$. 
Then $M$ automatically satisfies equality
$\delta(2)=2$ and there exists a totally
geodesic $S^5$, and a holomorphic immersion
$\phi:N_1\to  CP^2(4)$ such that $f$ is
congruent to
$\psi$, which is obtained from $\phi$ as in (1).

(5) Let $f: M\rightarrow S^6$
be a linearly full Lagrangian immersion of a
3-dimensional manifold satisfying equality 
$\delta(2)=2$.  Let $p$ be a non totally
geodesic point of $M$.
Then there exists a (possibly branched) almost complex curve
$\bar\phi : N_2 \rightarrow S^6$ such that
$f$ is locally around
$p$ congruent to $\bar \psi$, which is
obtained from $\bar \phi$ as in (3).

Let $f : S \rightarrow S^6$ be an almost
complex curve without totally geodesic
points.  Define
$$F : T_1S \rightarrow S^6(1) : v \mapsto
\frac{h(v,v)}{|| h(v,v) ||},\leqno (19.3)$$ where
$T_1S$ denotes the unit tangent bundle of
$S$.

 L. Vrancken (1997) showed that the
following:

(i) $F$ given by (19.2) defines a Lagrangian
immersion if and only if $f$ is superminimal,
and

(ii)  If $\psi : M \rightarrow S^6(1)$ be a
Lagrangian immersion which admits a unit
length Killing vector field whose integral
curves are great circles.  Then there exist
an open dense subset $U$ of $M$ such that
each point $p$ of $U$ has a neighborhood $V$
such that $\psi : V \rightarrow S^6$
satisfies $\delta(2)=2$, or $\psi: V
\rightarrow S^6$ is obtained as in (i).

\subsection{Further results}

L. Vrancken (1988) proved that a locally
symmetric Lagrangian submanifold of the
nearly K\"ahler $S^6$ has constant curvature 1 or
$1/16$. 

H. Li (1996) showed that if the Ricci
tensor of a compact Lagrangian submanifold in
the nearly K\"ahler $S^6$ satisfies $Ric\geq
(53/64) g$, then either $Ric=2g$ or the
submanifold is totally geodesic.

F. Dillen, B. Opozda, L. Verstraelen and L.
Vrancken (1988) proved that if a totally
real  surface $M$ in the nearly K\"ahler
$S^6$ is homeomorphic to a sphere, then
$M$ is totally geodesic.

K. Sekigawa studied $CR$-submanifolds in
the nearly K\"ahler $S^6$ and  proved that
there exists no proper $CR$-product in 
the nearly K\"ahler $S^6$, although there do
exist 3-dimensional $CR$-submanifolds in
the nearly K\"ahler $S^6$ whose totally
real and holomorphic distributions are
both integrable. 

Y. B. Shen (1998) studied slant minimal
surfaces in the nearly K\"ahler $S^6$. He
proved that if $f:M\to S^6$ is a minimal
slant isometric immersion of a complete 
surface of nonnagative Gauss curvature
$K$ in the nearly K\"ahler $S^6$ such
that $f$ is neither holomorphic nor
antiholomorphic, then either $K=1$ and
$f$ is totally geodesic, or $K=0$ and $f$
is either totally real or superminimal.
He also showed that if $f:S^2\to S^6$ is
a minimal slant immersion of a
topological 2-sphere which is neither
holomorphic nor antiholomorphic, then $f$
is totally geodesic.

\vfill\eject

\section{Axioms of submanifolds}

\subsection{Axiom of planes}

Historically the axiom of planes was
originally introduced by G. Riemann
in postulating the existence of a surface
$S$ passing through three given points
with the property that every straight line
having two points in $S$ is completely
contained in this surface. 

E. Beltrami (1835--1900) had shown in 1868
that a Riemannian manifold of constant
curvature satisfies the axiom of 2-planes
and F. Schur (1856--1932) proved in 1886 that
the converse is also true. The later result
was also obtained by L. Schl\"afli (1873) in
combination with the work of F. Klein
(1849--1925). 

In his 1928 book, \`E. Cartan defined the 
axiom of planes as follows: A Riemannian
$n$-manifold $M$, $n\geq 3$, is said to
satisfy the axiom of $k$-planes if, for
each point $x\in M$ and each
$k$-dimensional subspace $T'_x$ of the
tangent space $T_xM$, there exists a
$k$-dimensional totally geodesic submanifold
$N$ containing $x$ such that the tangent space
of $N$ at $x$ is $T'_x$, where $k$ is a fixed
integer $2\leq k<n$. 

\`E. Cartan's result states that real
space forms are the only Riemannian manifolds
of dimension $\geq 3$ which satisfy the axiom of
$k$-planes, for some $k$ with $2\leq k<n$.

\subsection{Axioms of spheres and  of totally umbilical submanifolds}

As a generalization of the axiom of $k$-planes,
D. S. Leung and K. Nomizu (1971) introduced the
axiom of $k$-spheres: for each point $x\in M$
and for each $k$-dimensional linear subspace
$T'_x$ of $T_xM$, there exists a
$k$-dimensional totally umbilical submanifold
$N$ of $M$ containing $x$ with parallel mean
curvature vector such that the tangent space
of $N$ at $x$ is $T'_x$. 

Leung and Nomizu proved that a Riemannian
manifold of dimension $n\geq 3$ satisfies the
axiom of $k$-spheres, $2\leq k<n$, if and only
if it is a real space form. 

The proof of Leung and Nomizu's result was
based on Codazzi's equation and the following
result of \'E. Cartan: A Riemannian manifold
$M$ of dimension $>2$ is a real space form
if and only if its curvature tensor $R$
satisfies
$R(X,Y,Z,X)=0$ for any orthonormal vector
fields $X,Y,Z$ in $M$.

 S. I. Goldberg and E.
M.  Moskal (1976) observed that the result
also holds if totally umbilical submanifolds
with parallel mean curvature vector are
replaced by submanifolds with parallel second
fundamental form. W. Str\"ubing (1979)
pointed out that totally umbilical
submanifolds with parallel mean curvature
vector can further be replaced by
submanifolds satisfying Codazzi equation for
real space forms.

J. A. Schouten in 1924 proved that a
Riemannian manifold of dimension $n\geq 4$ is
conformally flat if and only if it satisfies
the axiom of totally umbilical submanifolds of
dimension $k$, $3\leq k<n$. 

By the axiom of totally umbilical submanifolds
of dimension $k$, we mean that, for each point
$x\in M$ and for each
$k$-dimensional linear subspace
$T'_x$ of $T_xM$, there exists a
$k$-dimensional totally umbilical submanifold
$N$ of $M$ containing $x$ such that the
tangent space of $N$ at $x$ is $T'_x$.

K. L. Stellmacher (1951) showed that the same
result holds for Riemannian 3-manifold for
$k=2$.  

K. Yano and Y. Muto (1941) proved that
a Riemannian manifold of dimension $\geq 4$ 
is conformally flat if and only if it
satisfies the axiom of totally umbilical
surfaces with prescribed mean curvature vector.

 D. Van Lindt and L. Verstraelen (1981)
proved that a  
Riemannian manifold of dimension $n\geq 4$ is
conformally flat if and only if it satisfies
the axiom of conformally flat totally
quasiumbilical submanifolds of dimension $k$,
$3< k<n$. Here by the axiom of  conformally flat totally
quasiumbilical submanifolds of dimension $k$
we mean that, for each point $x\in M$ and for
each $k$-dimensional linear subspace
$T'_x$ of $T_xM$, there exists a
$k$-dimensional conformally flat totally
quasiumbilical submanifold
$N$ of $M$ containing $x$ such that the
tangent space of $N$ at $x$ is $T'_x$.

In 1975, Chen and Verstraelen proved that a  
Riemannian manifold $M$ of dimension $n\geq
4$ is conformally flat if and only if, for
each point $x\in M$ and for each
$k$-dimen\-sional ($2\leq k<n)$ linear
subspace
$T'_x$ of $T_xM$, there exists a
$k$-dimensional  submanifold
$N$ which passes through $x$ and which at $x$
tangent to $T'_x$ such that $N$ has flat
normal connection and commutative shape
operators.

\subsection{Axiom of  holomorphic $2k$-planes}

In 1955 K. Yano and Y. Mogi introduced the
axiom of holomorphic $2k$-planes on a K\"ahler
manifold $\tilde M$ as follows: for 
each point $x\in \tilde M$ and for each
holomorphic
$2k$-dimensional linear subspace
$T'_x$ of $T_xM$, there exists a
$2k$-dimensional totally
geodesic submanifold
$N$ of $\tilde M$ containing $x$ such that the
tangent space of $N$ at $x$ is $T'_x$. 

Yano and Mogi proved that a K\"ahler manifold of
real dimension $2n\geq 4$ is a complex space
form if and only if it satisfies the axiom of
holomorphic $2k$-planes for some $k$, $1\leq
k<n$. Goldberg and Moskal pointed out that the
same result holds if the
$2k$-dimensional totally geodesic submanifolds
are replaced by $2k$-dimensional totally
umbilical submanifolds with parallel mean
curvature vector. In 1982, O. Kassabov proved
that  the same result also holds if the $2k$-dimensional totally geodesic submanifolds
are replaced by $2k$-dimensional totally umbilical submanifolds.

\subsection{ Axiom of antiholomorphic $k$-planes}

The axiom of antiholomorphic $k$-planes was introduced in 1973 by  Chen and  Ogiue:  for 
each point $x\in \tilde M$ and for each totally real $k$-dimensional linear subspace
$T'_x$ of $T_x\tilde M$, there exists a $k$-dimensional totally
geodesic submanifold
$N$ of $\tilde M$ containing $x$ such that the
tangent space of $N$ at $x$ is $T'_x$. 

Chen and Ogiue proved that a K\"ahler manifold
of real dimension $2n\geq 4$ is a complex
space form if and only if it satisfies the
axiom of antiholomorphic $k$-planes for
$2\leq k\leq n$. The same result was also
obtained independently by K. Nomizu (1973b).

M. Harada (1974) pointed out that the same
result holds if the
$k$-dimensional totally geodesic submanifolds
were replaced by $k$-dimensional totally
umbilical submanifolds with parallel mean
curvature vector. 
Also, S. Yamaguchi and M. Kon
(1978) observed that the same result holds 
if the $k$-dimensional totally geodesic
submanifolds are replaced by $k$-dimensional
totally umbilical totally real submanifolds. 
 O. Kassabov proved that  the
same result also holds if the
$k$-dimensional totally geodesic submanifolds
are further replaced by $k$-dimensional
totally umbilical submanifolds.

D. Van Lindt and L. Verstraelen proved that a  
K\"ahler manifold of real dimension $2n> 4$
is a complex space form if and only if, for
each point $x\in M$ and for each
$k$-dimensional ($2\leq k<n)$ totally
real linear subspace
$T'_x$ of $T_xM$, there exists a totally real
$k$-dimensional  submanifold
$N$ which passes through $x$ and which at $x$
tangent to $T'_x$ such that $N$ has 
commutative shape operators and parallel
$f$-structure in the normal bundle. Here the
$f$-structure is the endomorphism on the
normal bundle induced from the almost complex
structure on the ambient space.

\subsection{Axioms of coholomorphic spheres}

Chen and Ogiue (1974a) introduced the axiom
of coholomorphic $(2k+\ell)$-spheres as
follows:   for  each point $x\in \tilde M$
and for each totally real
$(2k+\ell)$-dimensional
$CR$-plane section
$T'_x$ of $T_x\tilde M$, there exists a
$k$-dimensional totally
umbilical submanifold
$N$ of $\tilde M$ containing $x$ such that the
tangent space of $N$ at $x$ is $T'_x$. Chen
and Ogiue proved that  a K\"ahlerian manifold of
real dimension $2n\geq 4$ is locally flat
 if and only if it satisfies
the axiom of coholomorphic $(2k+\ell)$-spheres
for some integers $k$ and $\ell$ such that
$1\leq k,\ell<n$ and $2 k+\ell<2 n$.

An almost Hermitian manifold $M$ is called a
$RK$-manifold if its Riemann curvature tensor
$R$ and its almost complex structure $J$
satisfies $$R(X,Y,Z,W)=R(JX,JY,JZ,JW)$$ for
$X,Y,Z,W$ tangent to $M$. L. Vanhecke (1976)
studied $RK$-manifolds satisfying the axiom of
coholomorphic $(2k+1)$-spheres for some $k$ and
obtained characterization theorems for space
forms.

S. Tachibana and S. Kashiwada (1973) proved
that every geodesic hypersphere  $S$ with unit
normal vector field $\xi$ in a complex space
form is $J\xi$-quasiumbilical, that is, $S$ is
quasiumbilical with respect to $\xi$ and $J\xi$
is a principal direction with multiplicity
equal to either one or $n$, where $n=\dim S$.
The later case occurs only when the complex
space form is flat.

 A  hypersurface of a
K\"ahler manifold is said to be
$J\xi$-hypercylindric if $J\xi$ is a principal
direction and  the principal curvatures
other than the principal curvature associated
with $J\xi$ are zero.

 A K\"ahler manifold
$\tilde M$ is said to satisfy the axiom of
$J\xi$-quasiumbilical hypersurfaces if, for
each point $x\in \tilde M$ and for each
hyperplane $H$ of $T_x\tilde M$ with
hyperplane normal $\xi$, there exists a
$J\xi$-quasiumbilical hypersurface $N$
 containing $x$ such that the
tangent space of $N$ at $x$ is $H$. 

B. Y. Chen and L. Verstraelen (1980) proved
that a K\"ahler manifold of real dimension
$>4$ satisfies the axiom of
$J\xi$-quasiumbilical hypersurfaces if and
only if it is a complex space form. This
improves a result of L. Vanhecke and T. J.
Willmore (1977).

D.  van Lindt and L. Verstraelen showed that
a K\"ahler manifold of real dimension $>4$ is
locally flat if and only if,  for each point
$x\in\tilde M$ and for each hyperplane $L$
of $T_x\tilde M$ with hyperplane normal
$\xi$, there exists a $J\xi$-hypercylindric
hypersurface $N$ containing $x$ such that
the tangent space of $N$ at $x$ is $L$.

\subsection{Submanifolds contain many circles}

An ordinary torus contains exactly four
circles through each points. Since
each compact cyclide of Dupin in $E^3$
can be obtained from inversion of a
torus of revolution; thus it contains
four circles through each point. R.
Blum (1980) investigated the cyclide in
$E^3$ defined by
$$(x^2+y^2+z^2)^2-2ax^2-2by^2-2cz^2+d^2=0,\leqno(20.1)$$
where the real coefficients $a,b,c,d$ satisfy
the condition
$0<d<b\leq a,\, c<d.$ He proved that
(a) if $a\ne b$ and $c\ne -d$,  there
exist 6 circles through each point; (b)  if $a=b,\,c\ne -d$ or $a\ne
b,\, c=-d$, there exist 5 circles through
each point; and (c) if $a=b$ and $c= -d$,
there exist 4 circles through each point of
the cyclide. The case(c) represents a torus
of revolution. On the other hand, N. Takeuchi
(1987) showed that  a smooth compact
surface of genus one in
$E^3$ cannot contain seven circles
through each point. Obviously, there
exist infinitely many circles which
pass through each point of a round
sphere in $E^3$.

In 1984 K. Ogiue and R. Takagi proved
that a surface $M$ in $ E^3$ is locally a
plane or a sphere if, through each
point $p\in M$, there exist two
Euclidean circles such that (i)
they are contained in $M$ in some
neighborhood of $p$ and (ii) they
are tangent to each other at $p$.
Condition (i) alone is not
sufficient, as it is satisfied by
a torus of revolution. Ogiue and 
Takagi also generalized this to
obtain similar characterizations of
totally geodesic submanifolds and
extrinsic spheres of arbitrary
dimension in Riemannian manifolds.
In particular, they proved that a
2-dimensional surface $M$ in a
Riemannian manifold $N$ is
totally geodesic if through each
point $p\in M$ there exist three
geodesics of $N$ which lie in
$M$ in some neighborhood of $p$. 

R. Miyaoka and N. Takeuchi (1992) proved that 
a complete simply-connected surface in $E^3$
which contains two transversal circles through
each point must be a plane or a sphere.

K. Ogiue and N. Takeuchi (1992) proved 
that  a compact smooth surface of revolution
 which contains at least two circles through each
point is either a sphere or a hulahoop
surface, that is, a surface obtained by
revolving a circle around a suitable axis. A
hulahoop surface has 4, 5, or infinitely many
circles through each point. A
hulahoop surface, which is neither a
sphere nor an ordinary torus, contains
exactly 5 circles through each point.
Ogiue and Takeuchi also described
the concrete geometric construction of a
torus in Euclidean $3$-space containing five
circles through each point.

J. Arroyo, O. J. Garay and J. J. Menc\'ia
(1998) showed that if a compact surface of
revolution in $E^3$ contains at least two
ellipses through each point, then it is an
elliptic hulahoop surface, that is, a surface
obtained by revolving an ellipse around a suitable axis.

Without making distinction between real and nonreal circles, E. E. Kummer (1810--1893)
already observed in 1865 that a general cyclide has the property that there exist 10 circles through each generic point of the cyclide.

\vfill\eject

\section{Total absolute curvature}

\subsection{ Rotation index and total curvature of a curve}

Let $\gamma$ be closed smooth curve in the plane.
As a point moves along $\gamma$, the line
through a fixed point $O$ and parallel to the
tangent line of $\gamma$ rotates through an
angle $2n\pi$ or rotates $n$ times about $O$.
This integer $n$ is called the rotation index of
$\gamma$. If $\gamma$ is a simple closed curve,
$n=\pm 1$.

Two curves are said to be regularly homotopic if
one can be deformed to the other through a
family of closed smooth curves. Because the
rotation index is an integer and it varies
continuously through the deformation, it must be
constant. Therefore, two closed smooth curves
have the same rotation index if they are regular
homotopic. A theorem of W. C. Graustein
(1888--1941) and H. Whitney states that the
converse of this is also true; a result
suggested by Graustein whose proof was first
published in [Whitney 1937]. Hence, the only
invariant of a regular homotopy class is the
rotation index.

Let $\gamma(s)=(x(s),y(s))$ be  a
unit-speed smooth closed curve in $ E^2$.
Then
$$x''(x)=-\kappa(s) y'(s),\quad y''(x)=\kappa(s)
x'(s),$$ where $\kappa=\kappa(s)$ is the
curvature of the curve. If $\theta(s)$ denotes
the angle between the tangent line and the
$x$-axis, then  $d\theta=\kappa(s) ds$. Thus,
we have
$$\int_\gamma \kappa(s) ds=2n\pi,\leqno(21.1)$$
where $n$ is the rotation index of $\gamma$.

From (21.1), it follows that the total absolute
curvature of $\gamma$ satisfies
$$\int_\gamma |\kappa(s)| ds\geq
2\pi,\leqno(21.2)$$with the equality holding if
and only if $\gamma$ is a convex plane curve.

Inequality (21.2) was generalized  to closed
curves in $E^3$ by W. Fenchel (1905-- ) in
1929, and to closed curves in $ E^m,\, m>3,$
by K. Borsuk (1905--1982) in 1947. 

I. Fary (1922-- ) in 1949 and J. Milnor (1931-- ) in 1950 proved that if a closed curve $\gamma$ in $ E^m$ satisfies $$\int_\gamma |\kappa(s)| ds\leq 4\pi,$$ then $\gamma$ is unknotted.

\subsection{Total absolute curvature of Chern and Lashof}

Let $f:M\to  E^m$ be an isometric immersion
of a compact Riemannian $n$-manifold $M$ into
$ E^m$. Let $\nu_1(M)$ denote the unit
normal bundle of $f,$ and $S^{m-1}$ the unit
hypersphere centered at the origin of $
E^m$, and let $T_f:\nu_1(M)\to S^{m-1}$ be the
parallel translation. 

Denote by $\omega$ and $\Omega$ the volume
elements of $\nu_1(M)$ and of $S^{m-1}$,
respectively. For each $\xi\in\nu_1(M)$, we
have $$T^*_f\Omega=(\det A_\xi)\omega,$$
where $A_\xi$ is the shape operator of $f$ in
the direction $\xi$. 

As a generalization of the total absolute
curvature for a space curve, S. S. Chern and R.
K. Lashof (1957) defined the total absolute 
curvature of $f$ as follows:
$$\tau(f)={1\over {s_{m-1}}}\int_{\nu_1(M)}|
T^*_f\Omega|={1\over {s_{m-1}}}\int_{\nu_1(M)}|
\det A_\xi|\,\omega,\leqno(21.3)$$
where $s_{m-1}$ is the volume of the unit
$(m-1)$-sphere.

A function $\phi$ on $M$ which has only
nondegenerate critical point is called a Morse
function. Since $M$ is assumed to be compact,
each Morse function on $M$ has only a finite
number of critical points. The Morse number
$\gamma(M)$ of $M$ is defined as the least number
of critical points of  Morse functions on $M$.
The Morse inequalities imply that 
$$\gamma(M)\geq \beta(M; F)=\sum_{k=0}^n
\beta_k(M; F)\leqno(21.4)$$
for any field $ F$, where $\beta_k(M;
F)$ is the $k$-th betti number of $M$ over $
F$. For a (smooth) manifold $M$ of dimension
greater than five, the Morse number of $M$ is
equal to the number of cells in the smallest
CW-complex of the same simple homotopy type
as $M$ [Sharpe 1989a].

S. S. Chern and R. K. Lashof (1957,1958)
proved the following results.

 Let $f:M\to  E^m$ be an isometric immersion
of a compact Riemannian $n$-manifold
$M$ into $ E^m$. Then 

(1)  $\tau(f)\geq \gamma(M)\geq 2$;

(2) $\tau(f)=2$ if and only if $f$ is an
embedding and $f(M)$ is a convex hypersurface in
an affine $(n+1)$-subspace $E^{n+1}$ of $E^m$;

(3) if $\tau(f)<3$, then $M$ is homeomorphic
to $S^n$.

If $\tau(f)=3$, $M$ needs not be homeomorphic to
$S^n$. In fact, the Veronese embedding of the
real projective plane into $E^5$ satisfies
$\tau(f)=3$. 

J. Eells and N. H. Kuiper (1962) classified
compact manifolds which admit a Morse
function with three nondegenerate critical
points. They called such manifolds
``manifolds like projective spaces'', which
include the real, complex and quaternionic
projective planes and the Cayley plane. 

Applying Eells-Kuiper's result, it follows
that if  an immersion $f:M\to E^m$ of a
compact manifold $M$ into $ E^m$ satisfies
$\tau(f)<4$, then $M$ is homeomorphic either
to the sphere $S^n$ or else to one of the
manifolds like projective planes.

R. W. Sharpe (1989b) proved that for
manifolds $M$ of dimension greater than
five, the best possible lower bound for the
total absolute curvature $\tau(f)$ is the
Morse number $\gamma(M)$, as the immersion
$f:M\to E^m$ varies over all possibilities.

 Let $M$ be a compact $n$-manifold with $n>5$. Denote by
$\tau[i]$ the infimum of the total absolute 
curvature $\tau(j)$ as $j$ varies over all
immersions in the regular homotopy class of
the immersion $i: M\to E^m$. Sharpe (1989b)
also proved that, for $n>5$, if $m>n+1$ or if
$m=n+1$ is odd,  then $\tau[i]=\gamma(M)$. If
$m=n+1$ is even, then $\tau[i]=\max\{\gamma(M),2\vert d\vert \}$,
where $d$ is the normal degree of $i$ and
$\gamma(M)$ is the Morse number of $M$.
Examples are given of codimension-one
immersions of odd-dimensional spheres which
have arbitrary odd normal degree and which
attain the infimum of the total absolute
curvature in their regular homotopy class..

R. Langevin and H. Rosenberg (1976) proved that
 if the total absolute curvature of an
embedded surface of genus $g>1$ is 
$<2g+6$, then the surface is unknotted.
N. H. Kuiper and W. Meeks (1984) showed that
if the total absolute curvature
$\tau(f)$ of an embedding torus $f:T\to E^3$
is $\leq 8$, then $f(T)$ is unknotted.

Kuiper and Meeks (1984) also proved that,
for an embedding of a compact manifold
$f:M\sp n\to E^N$, the total absolute
curvature of $f$ satisfies
$$\tau(f)>\beta+4\sigma_1,\leqno (21.5)$$ where
$\beta$ is the sum of the mod 2 Betti
numbers and $\beta+\sigma_1$ is the minimal
number of generators of the fundamental group
of the complement of the image $f(M^n)$.

\subsection{Tight immersions}

An immersion $f:M\to E^m$ is called tight  (or
minimal total absolute curvature immersion)
if  $\tau(f)=b(M):=\min_\phi \beta(\phi)$,
where $\beta(\phi)$ denotes the number of
critical points of a Morse function $\phi$ on
$M$.  This condition is equivalent to
requiring that every height function that is
Morse have the minimum number of critical
points required by the Morse inequalities. Not
every compact manifold admits a tight
immersion. For instance, N. H. Kuiper (1958)
observed that the exotic $7$-sphere of J.
Milnor admits no immersion with minimal total
absolute curvature. This can be seems as
follows: Since a manifold $M$ homeomorphic to
a sphere admits a function with only two
critical points, $M$ satisfies $b(M)=2$.
Thus, if a manifold $M$ is homeomorphic to
$S^7$, then  a tight immersion  $f:M\to E^m$
would embed $M$ as a convex hypersurface in
an $ E^8\subset E^m$, and hence $M$ would be
diffeomorphic to the standard 7-sphere. 

D. Ferus (1967) proved that every embedding
$f$ of an exotic $n$-sphere $(n\geq 5)$ in $
E^{n+2}$ has total absolute curvature
$\tau(f)\geq 4$.

S. Kobayashi (1967b) showed that every
compact homogeneous K\"ahler manifold admits
a tight embedding. R. Bott and H.
Samelson (1958) proved that symmetric
$R$-spaces admit tight immersions. This
result was also proved independently by M.
Takeuchi and S. Kobayashi (1968). A tight
immersion of a symmetric
$R$-space is a minimal immersion into a
hypersphere.

In 1960's the theory of tight immersions
underwent substantial development and
reformulation. Since this notion is a
generalization of convexity, N. H. Kuiper 
called these ``convex immersions''. 
T. F. Banchoff (1965) first used
tight in conjunction with his study of the
two-piece property. 

Kuiper (1962) formulated
tightness in terms of intersections with
half-spaces and injectivity of induced maps on
homology and proved that his formulation is
equivalent to the minimal total absolute
curvature of manifolds which satisfy the
condition: the Morse number $\gamma(M)$ of
$M$ is equal to the sum $\beta(M; F)$ of the
Betti numbers for some field $F$.

 N. H. Kuiper (1961)
obtained smooth tight embeddings into $
E^3$ for all orientable surfaces, and smooth
tight immersions for all nonorientable surfaces
with Euler characteristic less than $-1$. He
also showed that smooth tight immersions of
the projective plane and the Klein bottle
into $E^3$ do not exist.  

The question of whether there is
a smooth tight immersion of the projective plane
with an attached handle into $ E^3$ has been
open for 30 years. In 1992, F. Haab proved that
no such  immersion exists.

 Kuiper and Meeks (1984) showed that if the genus
$g$ of a compact surface
$M$ is greater than 2, then there exists a
knotted tight embedding in $E^3$, whereas if
$g\leq 2$, there does not exist such an
embedding. Pinkall (1986a) showed that if
the Euler characteristic $\chi(M)$ of $M$
satisfies $\chi(M)<-9$, then every immersion $f$
of $M$ into $E^3$ is regularly homotopic to a
tight immersion. This is also true if $M$ is
orientable with genus $g\geq 4$. On the other
hand, there are immersions which are not
regularly homotopic to a tight immersion. This is
clearly true for any immersion of the projective
plane or Klein bottle.  Pinkall (1986a) also
showed that every tight immersion of the torus
$T^2$ is regularly homotopic to a standard
embedding, and thus there are no tight immersions
in the nonstandard regular homotopy class of
immersed tori. 

 Kuiper (1961) showed that smooth immersions
into $E^4$ of orientable surfaces exist for
every genus $\ge1$. He also showed that for
a substantial (that is, not contained in
any proper affine subspace) tight immersion
of a surface into
$ E^N$ one must have $N\leq 5$, with
equality only for surfaces projectively equivalent
to the Veronese surface and thus analytic.

Kuiper (1979) also proved the following:

(1) If $f:M^{2d}\to  E^N$ is a tight
substantial continuous embedding of a manifold
like a projective space, then $N\leq 3d+2$;

(2) Let $f:M^{2d}\to  E^{3d+2}$ be a tight
smooth substantial embedding of a compact
manifold with Morse number $\gamma(M)=3$. Then
$M^{2d}$ is algebraic. Moreover, it is the
union of its
$ E^{d+1}$-top-sets, smooth $d$-sphere $S^d$
that are quadratic $d$-manifolds.

T. F. Banchoff and N. H. Kuiper have produced tight
analytic immersions into $ E^3$ of all
orientable compact surfaces, while Kuiper has
produced tight analytic immersions into $
E^3$ of all nonorientable compact surfaces
with even Euler characteristic other than
zero.

G. Thorbergsson (1991) proved that an analytic
 tight immersion of a compact orientable surface
into $E^4$ which is substantial must be a
torus. He also showed that the surface is the
intersection of two developable ruled
hypersurfaces, possibly with singularities, with
two-dimensional rulings.

After the spheres, the $(k-1)$-connected 
$2k$-dimensional compact manifolds have the most
simple topology, for their homology groups vanish
in all dimensions except 0, $k$ and $2k$. Among
these so-called highly connected manifolds, the
only ones known to admit tight immersions into
some Euclidean space are the connected sums of
copies of $S\sp k\times S\sp k$, the projective
planes, and all surfaces except for the Klein
bottle and the projective plane with one handle
attached. Kuiper has
conjectured that the only $2k$-dimensional,
$(k-1)$-connected manifolds with trivial $k$-th
Stiefel-Whitney class that admit tight immersions
into $E^{2k+l}$ are homeomorphic to
$S^k\times S^k$. (The condition on the
$k$-th Stiefel-Whitney class follows from
$(k-1)$-connectedness for $k\neq 1,2,4,8$.) J.
J. Hebda (1984) and G. Thorbergsson (1986)
were able to construct counterexamples for
$l=1$ and $l=2$, respectively.  Assume the
immersion is analytic, R. Niebergall (1994) 
used top-set techniques to show that
Kuiper's conjecture is true for $l\ge 2$. 

C. S. Chen (1979) proved that if $f$ is a
substantial tight embedding of
$S^k{\times}S^{n-k}$ ($k/(n-k)\ne 2,{1\over
2}$) into $E^{n+2}$ whose image lies on an
ovoid, then $f$ is projectively equivalent
to a product embedding of two ovaloids of
dimensions $k$ and $n-k$, respectively.

M. van Gemmeren (1996)  generalized
tightness properties of immersed compact
manifolds to noncompact manifolds with a finite
number of ends. 

\subsection{Taut immersions}

T. Banchoff initiated the study of taut
immersions in 1970 by attempting to find
all tight surfaces which lie in a
hypersphere of a Euclidean $m$-space. Via
stereographic projection, this problem is
equivalent to the study of surfaces in
$E^m$ which have the spherical
two-piece property.

S. Carter and A. West (1972) extended the 
spherical two-piece property and defined an
immersion of a compact manifold $M$
to be taut if every nondegenerate Euclidean
distance function $L_p$ has the minimum
number of critical points. 

The property of tautness is preserved under 
Lie sphere transformations [Cecil-Chern 1987].
Furthermore, an embedding $f:M\to E^m$ is 
taut (or more precisely $F$-taut)
if and only if the embedding $\sigma\circ f:M\to
S^m$ has the property that every nondegenerate
spherical distance function has $\beta(M;F)$
critical points on $M$, where $\sigma:E^m\to
S^m-\{P\}$ is stereographic projection and
$\beta(M;F)$ is the sum of $F$-Betti numbers of
$M$ for any field $F$.

A spherical distance function
$d_p(q)=\cos^{-1}(\ell_p(q))$ is essentially
a Euclidean height function
$\ell_p(q)=p\cdot q$, for $p,q\in S^m$,
which has the same critical points as
$\ell_p$. Thus, the embedding $f$ is taut if
and only if the spherical embedding
$\sigma\circ f$ is tight, that is, every
nondegenerate height function $\ell_p$ has
$\beta(M;F)$ critical points on $M$. A tight
spherical embedding $F:M\to S^m\subset E^{m+1}$
is taut when regarded as an embedding of $M$ in
$E^{m+1}$ (cf. [Cecil-Ryan 1985]).

S. Carter and A. West (1972) observed that
if $f:M\to E^{n+1}$ is a compact
orientable embedded taut hypersurface, then
for sufficiently small $r>0$ the
hypersurface $f_r:M\to E^{n+1}$ defined by
$f_r(p)=f(p)+r\xi(p)$ if taut if and only if
$f$ is taut, where $\xi$ is a global unit
normal vector field of $f$. Carter and West
(1972) also pointed out their idea can be
generalized to taut embeddings of higher
codimension. Pinkall (1986b) proved that if
$M$ is a compact embedded submanifold of
dimension $n<m-1$ in $S^m$, then a tube
$T_r(M)$ of sufficiently small radius $r$
over $M$ is taut if and only if $M$ is {\bf
Z}$_2$-taut.

S. Carter and A. West (1972) also observed
that if
$M$ is noncompact but the immersion is
proper, then the Morse inequalities still
hold on compact subset of the form
$$M_r(L_p)=\{x\in M:L_p(x)\leq r\}$$ and the
notion of a taut immersion extends to this
case. 

Carter and West (1972) showed that a taut
immersion is tight; and it must be an
embedding, since if $p\in E^m$ were a double
point then the function $L_p$ would have two
absolute minima. Thus, tautness is a stronger
condition than tightness.  Via stereographic
projection and Chern-Lashof's theorem, Banchoff
and Carter-West observed that a taut embedding
of a sphere in a Euclidean space must be a
 round sphere, not merely convex.

S. Carter and A. West (1972) classified
 substantial taut embeddings as follows:

Let $f:M\to E^m$ be a substantial taut
embedding of an $n$-manifold. 

(a) If $M$ is compact, then $m\leq
\frac 12 n(n+3)$. In particular, if
$m=\frac12 n(n+3)$, then
$f$ is a spherical Veronese embedding of a
real projective space $RP^n$;

(b) If $M$ is noncompact, then $m\leq \frac
12 n(n+3)-1$. In particular, if $m= \frac 12
n(n+3)-1$, then $f(M)$ is the image under
stereographic projection of a Veronese
manifold, where the pole of the projection
is on the Veronese manifold.

Pinkall (1985a) proved that any taut
submanifold in a real space form is
Dupin. G. Thorbergsson (1983) proved
that a complete embedded proper Dupin
hypersurface in
$E^{n+1}$ is taut, and thus it must be
embedded. Pinkall (1985a) extended  this
result to compact submanifolds of higher
codimension for which the number of distinct
principal curvatures is constant on the unit
normal bundle. It remains as an open problem
whether Dupin implies taut without this
assumption.

 Taut surfaces in a Euclidean space
have been completely classified by Banchoff
(1970) and Cecil (1976) as follows:

Let $M$ be a taut surface substantially
embedded in a Euclidean $m$-space.

(c) If $M$ is a compact, then $M$ is a
round sphere or a ring cyclide in  $E^3$, a
spherical Veronese surface in $E^5$, or a
compact surface in $E^4$ related to one of
these by  stereographic projection;

(d) If $M$ is noncompact, then $M$ is a
plane, a circular cylinder, a parabolic
ring cyclide in $E^3$, or it is the image
in $E^4$ of a punctured spherical Veronese
surface under  stereographic projection.

Conversely, all of the surfaces listed in
(c) and (d) are taut.

Let $\psi:S^7\to S^4$ denote the Hopf
fibration defined by (10.4) in \S10.4.5. R.
Miyaoka and T. Ozawa (1988) proved that if
a compact submanifold of $S^4$ is taut,
then $\psi^{-1}(M)$ is also taut in $S^7$. 

Pinkall and Thorbergsson (1989b) classified,
up to diffeomorphism,
compact 3-manifolds which admit
taut embeddings into some Euclidean space.
They showed that there are seven such manifolds:
$S\sp 1\times S\sp 2$, $S^1\times RP^2$, the
twisted $S^2$-bundle over $S^1$, $S^3$, 
$RP^3$, the quaternion space
$S^3/\{\pm1,\pm i,\pm j,\pm k\}$, and the
torus $T^3$. Each has at least one known taut
embedding. The taut embeddings of $S\sp 1\times
S^2,S^ 3$ and the quaternion space have
been shown to be unique up to an appropriate
equivalence. $RP^3$ is known to
admit taut substantial embeddings of
codimension 2 and 5, but it is still unknown
whether it admits a taut substantial
embedding of codimensions 3 or 4. The
possible codimensions for taut substantial
embeddings of the other three manifolds have
been completely determined. Yet, the
complete geometric classification of taut
3-manifolds remains open.

For higher dimensional taut hypersurfaces,
Carter and West (1972) proved the following:

(e) Let $f:M\to E^{n+1}$ be a taut
embedding of a noncompact $n$-manifold with
$H_k(M;\hbox{\bf Z}_2)=\hbox{\bf Z}_2$ for
some $k,\,0<k<n,$ and $H_i(M;\hbox{\bf
Z}_2) =0$ for $i\ne 0,k$. Then $M$ is
diffeomorphic to $S^k\times E^{n-k}$ and
$f$ is a standard product embedding;

(f)  Let $f:M\to E^{n+q}$ be a
substantial taut embedding of a noncompact
$n$-manifold whose 
$\hbox{\bf Z}_2$-Betti numbers satisfy 
$\beta_k(M)=j>0$ for some 
$k,\,\frac n2<k<n$ and
$\beta_i =0$ for $i\ne 0,k$.
Then $q=1,\, j=1$ and $f$ embeds $M$ as a
 standard product $S^k\times E^{n-k}$ in
$E^{n+1}$.

Both the hypothesis $k>\frac n2$ in statement
(f) and the hypothesis that the codimension
is one in statement (e) are necessary
due to the example of the taut substantial
embedding of the M\"obius band
$M^2=RP^2-\{\hbox{a point}\}$ into
$E^4$ obtained from a Veronese surface $V^2$
in $S^4$ by stereographic projection with
respect to a pole on $V^2$.

Cecil and Ryan (1978a) studied taut
hypersurfaces of a Euclidean space and proved
the following:

(g) Let $M$ be a taut hypersurface in
$E^{2k+1}$ such that $H_i(M;\hbox{\bf Z}_2)=0$
for $i\ne 0,k,2k$. Then $M$ is a round
hypersphere,  a hyperplane,  a standard
product $S^k\times E^k$,  a ring cyclide, or a
parabolic cyclide;

(h) A taut compact hypersurface $M$ in
$E^{n+1}$ with the same {\bf Z}$_2$-homology
as $S^k\times S^{n-k}$ is a ring cyclide.

T. Ozawa (1986) showed that the
codimension of a taut substantial embedding of
the product of two spheres, $S^p\times S^
q, p<q$, is either 1 or $q-p+1$, and that a
connected sum of such products of two spheres
cannot be tautly embedded into any sphere. 

G. Thorbergsson (1983b) generalized statement
(h) to higher codimension. He obtained the
following:

(i) If $M^{2k}$ is a compact
$(k-1)$-connected, but not $k$-connected,
taut submanifold of $E^m$ which does not lie
in any  hypersphere of $E^m$, then either

(i-1) $m=2k+1$ and $M^{2k}$ is a cyclide of
Dupin diffeomorphic to $S^k\times S^k$, or

(i-2) $m=3k+1$ and $M^{2k}$ is diffeomorphic
to one of the projective planes $RP^2,
CP^2,$ $ HP^2$ or ${\Cal O}P^2$.

It is known that the compact focal submanifolds
of a taut hypersurface in $E^{n+1}$ need not be
taut. For instance, one focal submanifold of a
non-round cyclide of Dupin in $E^3$ is an
ellipse, which is tight but not taut. On the
other hand, S. Buyske (1989) used Lie
sphere-geometric techniques to show that if a
hypersurface
$M$ in $E^{n+1}$ is Lie equivalent to an
isoparametric hypersurface in $S^{n+1}$, then
each compact focal submanifold of $M$ is tight
in $E^{n+1}$.

R. Bott and H. Samelson (1958) proved that
the principal orbits of the isotropy
representations of symmetric space, known as
$R$-spaces, are taut. M. Takeuchi and S.
Kobayashi (1968) also proved the same result
independently. 

G. Thorbergsson (1988)
found some necessary conditions for the
existence of a taut embedding of a manifold
$M$. He proved the following:  

(j) Suppose that $M$ is tautly
embedded with respect to $F$ into $E^n$
and that $i>0$ is the smallest number such
that $H_i(M;F)$ is nontrivial. Then every
torsion element in $H_i(M;\bold Z)$ is of
order two. In particular, $H_i(M;\bold Z)$
is without torsion if the characteristic of
$F$ is not two. 

As a corollary,  Thorbergsson showed
that the only lens space admitting a taut
embedding is the real projective space. He
also showed that a coset space $M=G/H$
cannot be tautly embedded if one of the
following conditions is satisfied: 

(j-1) $G$ is simply connected with a torsion
element in the fundamental group of order
greater than two; 

(j-2) $G$ and $H$ are
simple and simply-connected and $H$ is a
subgroup with index greater than two. 

Thorbergsson gave several explicit examples
of homogeneous spaces which do not admit taut
embeddings because of these two conditions.

C. Olmos (1994) proved that if $M$ is a
substantial compact homogeneous submanifold
of a Euclidean space, then the following
five statements are equivalent: 

(0-1) $M$ is taut; 

(0-2) $M$ is Dupin; 

(0-3) $M$ is a submanifold with
constant principal curvatures; 

(0-4) $M$ is an orbit of the isotropy
representation of a symmetric space; 

(O-5) the first normal space of
$M$ coincides with the normal space. 
\vskip.1in

 T. E. Cecil and S. S. Chern (1987,1989)
established a relationship between a given
immersion $f$ and its induced Legendre
submanifold in $T_1S^{n+1}$ (cf. \S10.4.8).
They proved the following:  

Suppose  $L:M^n\to T_1S^{n+1}$ is a
Legendre submanifold whose M\"obius
projection is a taut immersion. If $\beta$
is a Lie transformation such that the
M\"obius  projection of $\beta\circ
L$ is an immersion, then the M\"obius
projection of $\beta\circ L$ is also taut.

Two hypersurfaces of $S^{n+1}$ are called Lie
equivalent if their induced Legendre
submanifolds are Lie equivalent, that is, 
there is a Lie sphere transformation which
carries one to the other. The result of
Cecil and Chern implies that  if two compact
hypersurfaces of $S^{n+1}$ are Lie equivalent,
then one of the two hypersurfaces is taut if
and only if the other is taut.

When $M$ is a compact embedded submanifold in
$S^m$ of codimension $>1$, it induces a
Legendre submanifold defined on the unit
normal bundle of $M$ in $S^m$. As with
hypersurfaces two submanifolds of arbitrary
codimension in $S^m$ are said to be Lie
equivalent if their induced Legendre
submanifolds are Lie equivalent. Cecil and
Chern's result also implies that  if two compact
embedded submanifolds of $S^m$ are Lie
equivalent, then one is {\bf Z}$_2$-taut if
and only if the other is {\bf Z}$_2$-taut.

Cecil and Chern (1987,1989) also
investigated the relationship between
Legendre submanifolds and Dupin submanifolds.

J. J. Hebda (1988) investigated possible
{\bf Z}$_2$-cohomology of some taut
submanifolds in spheres. E. Curtin (1994)
extended the notion of taut embeddings to
manifolds with boundary and obtained some
classification results.

Tight and taut submanifolds have also been
studied in hyperbolic space. In fact, Cecil
and Ryan (1979a,1979b) introduced three
classes of distance functions in hyperbolic
$m$-space $H^m$ whose level sets are spheres
centered at a point $p\in H^m$, equidistant
hypersurfaces from a hyperplane in $H^m$,
and horospheres equidistant from a fixed
horosphere. Suppose $M$ is a compact
embedded submanifold in $H^m$ and that
$\pi:H^m\to D^m$ is stereographic projection
of $H^m$ onto a disk $D^m\subset E^m$. They
proved that $\pi(M)$ is taut in $E^m$ if and
only if every nondegenerate hyperbolic
distance function of each of the three types
has $\beta(M;F)$ critical points on $M$.

G. Thorbergsson (1983b) obtained a result
similar to statement (i) for noncompact taut
submanifolds of $E^m$ and for taut compact
$(k-1)$-connected hypersurfaces in real
hyperbolic space $H^{2k+1}$.

\vfill\eject

\section{Total mean curvature}

\subsection{Total mean curvature of surfaces in Euclidean 3-space}

It is well-known that  the two most important
geometric invariants of a surface in a Euclidean
$m$--space  are the  Gaussian curvature $K$
and the squared mean  curvature  $H^2$.
According to Gauss' Theorema Egregium,
Gaussian curvature is an intrinsic
invariant and the integral of the Gaussian
curvature over a compact surface gives the
famous Gauss--Bonnet formula:
$$\int_M KdV=2\pi
\chi(M),\leqno(22.1)$$ where $\chi(M)$ denotes the
Euler number of $M$.

For a compact surface $M$ in $
E^3$, Chern-Lashof's inequality yields the
following inequality:
$$\int_M |K|dV\geq 4\pi(1+g),\leqno(22.2)$$
where $g$ denotes the genus of $M$. 

G. Thomsen initiated in 1923 the 
study of  total mean curvature: $$w(f)=\int_M
H^2dV\leqno(22.3)$$
of an immersion $f:M\to E^3$ of a surface
in $ E^3$.
Among others, G. Thomsen
studied the first variations of the
total mean curvature  and showed that the
Euler-Lagrange equation of (22.3) is given by
$$\Delta H+2H(H^2-K)=0.\leqno (22.4)$$

W. Blaschke proved in his 1923 book that the
total mean curvature of a compact surface in
$ E^3$  is a conformal invariant. Chen
(1973c,1974) extended Blaschke's result to
any submanifold of dimension $\geq 2$
in an arbitary Riemannian manifold. 

Using the inequality $H^2\geq K$ for a 
surface $M$ in $E^3$, T. J. Willmore
(1968) observed that combining
Gauss-Bonnet's formula and Chern-Lashof's 
inequality yields the following inequality: 
$$w(f)=\int_M H^2dV\geq4\pi,\leqno(22.5)$$ 
for an immersion $f:M\to E^3$ of a
compact surface $M$ in $ E^3$, with
equality holding if and only if
$M$ is a round sphere. 

The functional $w(f)$ defined by (22.3),
initially studied by G. Thomsen, is also
known as Willmore's functional and a surface
satisfying Thomsen's equation (22.4) is
called a stationary surface or a Willmore
surface. 

 An immersed surface $\bar M$ in a Euclidean
space is said to be conformally equivalent to
another immersed surface $M$ in a Euclidean
space if $\bar M$ can be obtained from $M$
via conformal mappings on Euclidean space.

\subsection{Willmore's conjecture}

Since the total mean curvature of an immersion
 $f:M\to E^3$ of a compact surface in
$ E^3$ is at least $4\pi$, it is a
natural question to determine the infinimum of
$w(f)$ among all immersions of a compact
surface $M_g$ of a given genus $g$, or among
all isometric immersions of a  compact
Riemannian surface. 

T. J. Willmore conjectured that if $f:M\to
E^3$ is an immersion of a torus, then
$w(f)\geq 2\pi^2$. 

Willmore's conjecture have been proved to be
true for various classes of immersed tori
in $ E^3$ (or more generally, in $
E^m,\, m\geq 3$). 

For instance, the following are known.

(1) If $f:M\to
E^3$ is a closed tube with fixed radius over a
closed curve in
$ E^3$, then $w(f)\geq 2\pi^2$, with the
equality holding if and only if it is a torus
of revolution whose generating circle has
radius $r$ and distance $(\sqrt{2}-1)\,r$ from
the axis of revolution [Shiohama-Tagaki
1970,  Willmore 1971,  Langer-Singer 1984].

(2)  If $f:M\to  E^3$ is a knotted torus,
then $w(f)\geq 8\pi$ [Chen 1971,1984b].

(3) If an immersed torus $M$ in
a Euclidean space $ E^m$ is conformally
equivalent to a flat torus in a
Euclidean space, then
$w(f)\geq 2\pi^2$, with the equality holding
if and only if $M$ is a conformal Clifford
torus, that is, $M$ is conformally equivalent
to a standard square torus in an affine $
E^4\subset E^m$ [Chen 1976b, 1984b].

(4) U. Hertrich-Jeromin and U.  Pinkall (1992)
proved that the conjecture is  true for 
elliptic tubular tori in $ E^3$. (A
special case of this was proved  by van de
Woestijne and Verstraelen in 1990, cf.
[Verstraelen 1990]).

(5) If an immersed tori in $ E^3$ has
self-intersections, then $w(f)\geq 8\pi$
[Li-Yau 1982].

(6) Li and Yau (1982) showed that the conjecture is true 
for certain bounded domain of the moduli
space of conformal structure on torus. 
Montiel and Ros (1985) proved that the
conjecture is true for a larger domain in
this moduli space.

\subsection{Further results on total mean curvature for surfaces in Euclidean space}

Chen (1973) proved that if $f:M\to  E^m$ is a compact surface in $ E^m$ which is
conformally equivalent to a compact surface in $ E^4\subset  E^m$ with nonnegative
Gaussian curvature and  $w(f)\leq (2+\pi)\pi$, then $M$ is homeomorphic to $S^2$.

Let $f:M\to E^4$ be an embedding of a compact surface $M$ into $ E^4$. The
fundamental group $\;\pi_1( E^4-f(M))\,$ of $\; E^4-f(M)\;$ is called the knot group of
$f$. The minimal number of generators of the knot group of $f$ is called the knot number
of $f$.

P. Wintgen (1978,1979) proved the following:

(a) If $f:M\to E^4$ is an embedding of a
compact surface $M$ into $ E^4$, then
$w(f)\geq 4\pi \rho$, where $\rho$ is the
knot number of $f$;

(b) If $f:M\to\ E^4$ is an immersion of a
compact oriented surface $M$ into $ E^4$,
then
$$w(f)\geq 4\pi(1+|I_f|-g),\leqno (22.6)$$ where
$I_f$ is the self-intersection number of $f$
and $g$ is the genus of $M$.

 Li and Yau (1982) showed that if $f:RP^2\to  E^3$ is an immersion of a real projective plane $RP^2$ into $ E^3$, then $w(f)\geq 12\pi$. R. Bryant (1987a) and, independently, R. Kusner (1987) found explicit immersions of $RP^2$ in $ E^3$ satisfying $w(f)=12\pi$.

Let $M$ be an oriented  surface immersed in
$ E^4$ and $\{X_1,X_2\}$ be an
orthonormal oriented frame field of $TM$. For
each point $x\in M$ and each unit tangent
vector of $M$ at $x$, we put $X=(\cos\theta)
X_1+(\sin\theta) X_2$. Then the second
fundamental form $h$ of $M$ satisfies
$$h(X,X)=
H+(\cos
2\theta)\left({{h(X_1,X_1)-h(X_2,X_2)}\over
2}\right)+(\sin 2\theta)h(X_1,X_2),$$  which
shows that $$E_x=\{h(X,X):X\in T_xM,\,
|X|=1\}$$  is an ellipse in the normal space
$T^\perp_xM$ centered at $ H$.
The ellipse $E_x$  is called the ellipse of
curvature at $x$. 

Let $\{e_3,e_4\}$ be an orthonormal oriented
frame field of the normal bundle of
$M$. The normal curvature $K^D$ of $M$ in
$ E^4$ is defined by
$$K^D=\left<R^D(X_1,X_2)e_4,e_3\right>.$$

 I. V. Guadalupe and L. Rodriguez
(1983) proved that if  $f:M\to E^4$ is an
immersion of a compact oriented surface $M$
into $ E^4$ and if the normal curvature
$K^D$ of $f$ is everywhere positive, then
$w(f)\geq 12\pi$, with the equality holding
if and only if the ellipse of curvature of
$f$ is always a circle.

For an immersion $f:M\to  E^3$ of a
compact surface into $ E^3$, W. K\"uhnel
and U. Pinkall (1986) showed the following:

(1) if $M$ is nonorientable with even Euler
number, then $w(f)\geq 8\pi$, and

(2) if $M$ has odd Euler number, then
$w(f)\geq 12\pi$.

 K\"uhnel and Pinkall (1986) also
proved that, for any genus $g$, there are
compact orientable surfaces of genus $g$
immersed in
$ E^3$ with
$\int H^2dV\leq 8\pi$.

\subsection{Total mean curvature for arbitrary submanifolds and applications} 

According to Nash's embedding theorem, every
compact  Riemannian $n$-manifold
can be isometrically embedded in
$E^{n(3n+11)/2}$. On the other hand,
most compact Riemannian $n$-manifolds
cannot be isometrically immersed in
$E^{n+1}$ as hypersurfaces. For instance,
every compact surface with nonpositive
Gaussian curvature cannot be isometrically
immersed in $E^3$.

Chen (1971) proved the following general
inequality for compact submanifolds in
Euclidean space for arbitrary dimension and 
arbitrary codimension:

 Let $M$ be a compact $n$-dimensional
submanifold of $ E^m$. Then
$$\int_M H^ndV\geq s_n,\leqno(22.7)$$ 
where $s_n$ is the volume of unit $n$-sphere.
The equality sign of
(22.7) holds if and only if $M$ is a 
convex planar curve when $n=1$; and $M$ is
embedded as a hypersphere in an affine
$(n+1)$--subspace of
$ E^m$ when $n>1$. 

If $n=1$, inequality (22.7) reduces to the
well-known Fenchel-Borsuk inequality for closed
curves in Euclidean space.

Some geometric applications of inequality
(22.7) are the following:

(1) If $M$ is a compact $n$-dimensional
minimal submanifold of the unit $m$-sphere
$S^m$, then the volume of $M$ satisfies
$\hbox{\rm vol}\,(M)\geq s_n,$
with the equality holding if and only if $M$ is a great
$n$-sphere of $S^m$;

(2) If $M$ is a compact $n$-dimensional
minimal submanifold of $ RP^m$ of constant
sectional curvature 1, then
$\hbox{\rm vol}\,(M)\geq {{s_n}/2},$  with the
equality holding if and only if $M$ is a
totally geodesic submanifold of $ RP^m$;

(3) if $M$ is a compact $n$-dimensional
minimal submanifold of  $\, CP^m$ of constant
holomorphic sectional curvature 4, then  
$\hbox{\rm vol}\,(M)\geq {{s_{n+1}}/{2\pi}}, $
with the equality holding  if and only if  
$M= CP^k$ which is embedded as a totally
geodesic complex submanifold of $CP^m(4)$;

(4) If $M$ is a compact $n$-dimensional
minimal submanifold of $\, HP^m$ of constant
quaternionic sectional curvature 4, then
$\hbox{\rm vol}\,(M)\geq {{s_{n+3}}/{2\pi^2}},$
with the equality holding  if and only if
$M= HP^k$ which is embedded as a totally
geodesic quaternionic submanifold of
$HP^m(4)$; and

(5) If $M$ is a compact $n$-dimensional
minimal submanifold of the Cayley plane 
$ {\Cal O}P^2$ of maximal sectional
curvature 4, then  $\hbox{\rm vol}\,(M)\geq
{{s_n}/{2^n}}.$

If $M$ is a compact $n$-dimensional
submanifold $M$ with nonnegative scalar
curvature in a Euclidean $m$-space, then the
total mean curvature of $M$ satisfying [Chen 1972]:
$$\int_M H^ndV\geq C(n)\beta(M),\leqno (22.8)$$
where $C(n)$ is a positive constant
 depending on $n$ and $\beta(M)$ is the
topological invariant: $\beta(M)=\max
\{\sum_{i=0}^n\beta_i(M;F):F\hbox{
fields}\}$, $\beta_i(M;F)$ the $i$-th
Betti number of $M$ over $F$.

Related to inequality (22.8), Chen
(1976c) proved that the total scalar
curvature of an arbitrary compact
$n$-dimensional submanifold $M$ in the
Euclidean $m$-space, regardless of
codimension, satisfies $$\int_M S^{n/2}dV \geq
\left(\left(\frac n2\right)^{n/2}
s_n\right)\beta(M).\leqno (22.9)$$
The equality sign of (22.9) holds if and
only if $M$ is embedded as a hypersphere in
an affine $(n+1)$-subspace of $E^m$.

\subsection{Some related results} 

In the 1973 AMS symposium held at
Standard University, Chen asked to find
the relationship between the total mean
curvature and  Riemannian invariants of a
compact submanifold in a Euclidean space
(cf. [Chen 1975]). In the late 1970s Chen
obtained a solution to this
problem; discovering a sharp relationship
between the total mean curvature of a
compact  submanifold and the order of the
immersion (cf. [Chen 1984b]). More
precisely, he proved the following:

 Let $M$ be an $n$-dimensional compact
submanifold of $ E^m$. Then
$$\left({{\lambda_p}\over n}\right)\hbox{\rm
vol}(M)\leq \int_M H^2dV\leq
\left({{\lambda_q}\over n}\right)\hbox{\rm
vol}(M), \leqno(22.10)$$
where $p$ and $q$ are the lower order and
the upper orders of
$f:M\to E^m$.  Either equality sign in
(22.10)  holds if and only if M is of
1-type.

Inequality (22.10) improves a result of
Reilly (1977) who proved that, for any
compact submanifold $M$ in $ E^m$, the
total mean curvature of $M$ satisfies
$$\int_M H^2dV\geq\left({{\lambda_1}\over
n}\right)
\hbox{\rm vol}(M).\leqno (22.11)$$

We also have the following sharp inequalities
for the total mean curvature [Chen 1987b,
Chen-Jiang 1995]:

 Let $f:M\to E^m$ be a compact
$n$-dimensional  submanifold of  $ E^m$
and let $c$ denote the distance from the
origin to the center of gravity of $M$ in
$ E^m$. Then

(i) if $M$ is contained in a closed ball
$\overline{B_0(R)}$ with radius $R$ centered 
at the origin, then
$M$ satisfies $$\int_M |H|^k dV\geq
{{\hbox{\rm vol}(M)}\over{(R^2-c^2)^{k/2}}},
\quad
k=2,3,\cdots,n,\leqno(22.12)$$ with 
equality  holding for some
$k\in\{2,3,\ldots,n\}$ if and only if $M$ is
a minimal submanifold of the hypersphere
$S_0^{m-1}(R)$ of radius $R$ centered at
the origin;  

(ii) if $M$ is contained in $E^m-B_0(r)$,
then $$\left( {{\lambda_p} \over n}\right)^2(r^2-c^2)\leq {1\over{\hbox{\rm
vol}(M)}}\int_M |H|^2 dV  \leq\left( {{\lambda_q}\over n}\right)^2(R^2-c^2)
\leqno(22.13)$$ where $p$ and $q$
denote the lower  and the upper orders of
$M$ in $ E^m$. Either equality of
$(22.13)$ holds if and only if $M$ is
a minimal submanifold of a hypersphere
centered at the origin; and 

(iii) if $f(M)$ is contained in a unit hypersphere  of $ E^m$, then the first
and the second nonzero eigenvalues of the Laplacian  of $M$ satisfy
$$\int_M  |H|^2dV\geq {1\over {n^2}}\{n(\lambda_1+\lambda_2)-\lambda_1
\lambda_2\}\hbox{\rm vol}(M),\leqno (22.14)$$
with equality sign of $(22.14)$ holding if
and only if either $f$ is of 1-type with 
order
$\{1\}$ or with order $\{2\}$, or $f$ is of
2-type  with order $\{1,2\}$.
\smallskip
Some easy applications of (iii) are the
following [Chen 1987b]: 

(1) If $M$ is an $n$-dimensional  compact
minimal submanifold of $RP^m(1)$, then the
first and the second nonzero eigenvalues of
the Laplacian of $M$ satisfy
$${m\over {2(m+1)}}\lambda_1\lambda_2\geq 
n(\lambda_1+\lambda_2-2n-2).\leqno(22.15)$$

(2) If $M$ is an $n$-dimensional compact
minimal submanifold of $CP^m(4)$, then
 
$${m\over
{2(m+1)}}\lambda_1\lambda_2\geq n(\lambda_1+
\lambda_2-2n-4),\leqno(22.16)$$
with equality holding if and
only if $M$ is one of the following compact
Hermitian symmetric spaces: 
$$  CP^k(4),\;\; CP^k(2),\;\;  SO(2+k)/SO(2)\times SO(k),$$ 
$$  CP^k(4)\times CP^k(4),\;\; U(2+k)/U(2)\times U(k),\; (k>2),$$ 
$$SO(10)/U(5),\;\;\hbox{\rm and}\;\; E_6/\hbox{\rm Spin}(10)\times T,$$ with an appropriate metric, where m is given respectively  by 
$$ k,\;{{k(k+3)}\over 2},\;
k+1,\; k(k+2),\; {{k(k+3)}\over
2},\;15,\;\;\hbox{\rm and}\;\;26.$$

(3) If $M$ is  an $n$-dimensional  compact
minimal submanifold of $HP^m(4)$, then
$${m\over {2(m+1)}}\lambda_1\lambda_2\geq n (\lambda_1+\lambda_2-2n-8),\leqno(22.17)$$ with
equality  holding if and only if  $M$ is a
totally geodesic quaternionic submanifold
and $n=4m$.

 For further applications of (i), (ii) and
(iii), see [Chen 1996d]. 
\smallskip
The following  conformal property of
$\,\lambda_1$vol$(M)\,$ as another
application of the notion of order was first
discovered in [Chen 1979b]:

 If a compact Riemannian
surface $M$  admits an order $\{1\}$ isometric embedding
into $ E^m$, then, for any compact surface
$\overline M\subset E^m$ which is
conformally equivalent to $M\subset  E^m$,
we have 
$$\lambda_1\hbox{\rm vol}(M)\geq \bar{\lambda}_1 {\rm vol}(\bar{ M}).\leqno (22.18)$$
Equality sign of $(22.18)$ holds if and only
if $\bar M$  also admits an isometric
embedding of order $\{1\}$. 

Some further applications of (22.10) and
(22.11) are the following [Chen 1983b,
1984b].

Let $M$ be an $n$-dimensional compact
submanifold of the unit hypersphere
$S^{m-1}$ of Euclidean $m$-space. Denote by $p$
and $q$  the lower order and the upper order
of $M$ in $ E^m$. Then

(a) If $M$ is mass-symmetric in $S^{m-1}$,
then $\lambda_1\leq \lambda_p\leq n.$
 In parti\-cular,  $\lambda_p=n$ if and only if $M$ is of 1-type
and of order $\{p\}$;

(b) If $M$ is of finite type, then
$\lambda_q \geq n$. In particular,
$\lambda_q=n$ if and only if $M$ is of 1-type and of
order $\{q\}$;

(c)  If $M$ is a compact
$n$-dimensional minimal submanifold of
$RP^m(1)$, then $\lambda_1$  of $M$
satisfies $\lambda_1\leq 2(n+1),$
with equality holding if and only if $M$ is a
totally geodesic $RP^n(1)$  in $RP^m(1)$;

(d)  If $M$ is a compact
$n$-dimensional minimal submanifold of
$CP^m(4)$, then 
$\lambda_1\leq 2(n+2),$ with 
equality holding if and only if 
$M$ is  holomorphically isometric to a
$CP^{n\over 2}(4)$, and $M$ is embedded as 
a complex totally geodesic  submanifold of
$CP^m(4)$.

(e)  If $M$ is a compact
$n$-dimensional minimal submanifold of $
HP^m(4)$, then $\lambda_1\leq
2(n+4),$ with  equality sign holding if and
only if  $M$ is a $HP^{n\over 4}(4)$ and 
$M$ is a quaternionic totally geodesic 
submanifold of 
$HP^m(4)$; and

(f) If $M$ is a compact
$n$-dimensional minimal submanifold of the Cayley plane $
{\Cal O}P^2(4)$ , then
 $\lambda_1\leq 4n$.

I. Dimitri\'c (1998) improved the result
of statement (f) to
$\lambda_1<676/15=45.0\bar 6$ for compact
minimal hypersurfaces in ${\Cal O}P^2(4)$. 

For compact K\"ahler submanifolds of $CP^m(4)$,
statement (d)  was due to A. Ros (1983), and
independently by N. Ejiri (1983).

The Euler-Lagrange equation of $w(f)=\int_M
H^ndV$ for a compact hypersurface $M$ in
$ E^{n+1}$ is given by
$$\Delta H^{n-1}+ H^{n-1}\left(nH^2-S\right)=0,\leqno (22.19)$$
where $S$ denotes the squared length of the second fundamental form.

Chen (1973d) proved that  hyperspheres in $ E^{n+1}$ are the only solutions of (22.19) when $n$ is odd (for compact submanifolds in Euclidean space
with higher codimension, see [Chen-Houh
1975]). On the other hand, there do exist
many solutions of (22.19) other than hyperspheres for even
$n$. For example, a torus of revolution in $
E^3$ whose generating circle has radius $r$ and distance
$(\sqrt{2}-1)\,r$ from the axis of revolution
is a solution of (22.19). Furthermore, the
stereographic projections of compact minimal
surfaces $M$ in $S^3$ also satisfy (22.19)
with $n=2$. 

U. Pinkall (1985c) found examples of compact embedded surfaces in $E^3$
satisfying (22.19) that are not
stereographic projections of compact
minimal surfaces in $S^3$ (see, also
[Weiner 1979, Pinkall-Sterling 1987]). Barros
and Garay [1998,1999] constructed many
submanifolds in spheres which satisfies (22.19).

\vfill\eject


\begin{thebibliography}{444 }

\bibitem{A1} Abe, K.,    
Applications of a Riccati type differential
equation to Riemannian manifolds with
totally geodesic distributions, T\^ohoku
Math. J., 25 ({1973}), 425--444.
\bibitem{A2}Abe, K., On a class of
hypersurfaces of $R\sp{2n+1}$, Duke
Math. J., 41 (1974), 865--874.

\bibitem{} Abe, K., Mori, H. and H. Takahashi,  
A parametrization of isometric immersions
between hyperbolic spaces, Geom.
Dedicata, 65 (1997), 31--46.

\bibitem{} Abe, T.,   Slant surfaces in Kaehlerian space forms, Math. J. Toyama Univ., 20 (1997), 37--48.

\bibitem{} Abi-Khuzam, F. F., A one-parameter family of minimal surfaces with four planar ends, Internat. J. Math., 6 ({1995}), 149--159. 

\bibitem{} Abresch, U., Isoparametric hypersurfaces with four or
six distinct principal curvatures. Necessary conditions on the
multiplicities, Math. Ann., 264 ({1983}), 283--302.,

\bibitem{} Abresch, U.,
Constant mean curvature tori in terms of
elliptic functions, J. Reine Angew Math., 374 ({1987}), 169--192.

\bibitem{} Aiyama, R. and K. Akutagawa, Kenmotsu-Bryant type
representation formulas for constant mean curvature surfaces in $H^3(-c^2)$ and
$S^3_1(c^2)$, Diff. Geom. Appl., 9 ({1998}), 251--272..

\bibitem{} Aiyama, R. and K. Akutagawa, Kenmotsu type  representation formula for surfaces
with prescribed mean curvature in $S^3(c^2)$ and adjusting it for CMC surfaces, {1997}, preprint.

\bibitem{} Alencar, H., do Carmo, M. and A. G. Colares, Stable
hypersurfaces with constant scalar curvature, Math. Z., 213 ({1993}),
117--131. 

\bibitem{} Alencar, H., do Carmo,  M. and H. Rosenberg, On first eigenvalue of the linearized
operator of the $r$th mean curvature of a hypersurface, Ann. Global Anal.
Geom., 11 ({1993}), 387--395.

\bibitem{} Alexandrov, A. D., Uniqueness theorems for surfaces in the large, Vestk
Leningrad Univ., 13 ({1958}), 5-8.

\bibitem{} Almgren, F. J. Jr.,     Some interior regularity theorems for
minimal surfaces and extension of
Berstein's theorem, Ann. Math.,
84 ({1966}), 277--292.

\bibitem{} Almgren, F. J. Jr.  and L. Simon, Existence of embedded
solutions of  Plateau's problem,  Ann. Scuola Norm. Sup. Pisa Cl. Sci., 6 ({1979}), 447--495. 

\bibitem{} Amaral, F. A., Toros minimos mergulhados em $S^3$, Doctoral
Thesis, Universidad Federal do Cear\'a, 
 (Fortaleza, Brazil),  December 1997.

\bibitem{} Anderson, M. T.,  Complete minimal varieties in hyperbolic space, Invent. Math., 69 ({1982}), 477--494. 

\bibitem{} Anderson, M. T.,  The compactification
of a minimal submanifold in Euclidean
space by the Gauss map, IHES preprint.

\bibitem{} Anderson, M. T.,  Curvature estimates for
minimal surfaces in $3$-manifolds, Ann.
Sci. \'Ecole Norm. Sup. , 18 (1985), 89--105.

\bibitem{} Aoust, L. S. V., Analyse infinit\'esimale des courbes
dans l'espace, Paris, {1876}.

\bibitem{} Arroyo, J., Garay, O. J. and J. J. Menc\'ia, On a family of surfaces of revolution of finite Chen-type, Kodai Math. J.,
21 ({1998}), 73--80.

\bibitem{} Arslan, K. and A. West, Nonspherical submanifolds with  pointwise $2$-planar normal sections,  Bull. London Math. Soc. 28 ({1996}), 88--92.

\bibitem{} Asperti, A. and M. Dajczer,   Conformally flat Riemannian manifolds as
hypersurfaces of the light cone, Canad. Math. Bull., 32 ({1989}), 281--285.

\bibitem{} Aubin, T.,   Equations du type Monge-Amp\`ere sur les vari\'et\'es k\"ahlerienne compactes, C. R. Acad. Sc. Paris, 283 (1976), 119--121.

\bibitem{} Aumann, G., Die Minimal hyperregelflaechen,
Manuscripta Math. 34 ({1981}), 293--304.

\bibitem{}  Baikoussis, C., Defever, F., Koufogiorgos, T. and L. Verstraelen,  Finite type
immersions of flat tori into Euclidean spaces, Proc. Edinburgh Math. Soc., 38 ({1995}), 413--420.

\bibitem{} Baikoussis, C., Blair, D. E., Chen, B. Y. and F. Defever, Hypersurfaces of
restricted type in Minkowski space, Geom. Dedicata, 62 ({1996}), 319--332.

\bibitem{}  Baikoussis, C. and L. Verstraelen, The Chen-type of the spiral surfaces,
Results Math. 28 ({1995}), 214--223.

\bibitem{} Baikousis, C. and T. Koufogiorgos, Isometric immersions of complete Riemannian
manifolds into Euclidean space, Proc. Amer. Math. Soc., 79 ({1980}), 87--88.

\bibitem{}   Bakelmeman, I. J. and B. E. Kantor, Existence of a hypersurface homeomorphic to the sphere in Euclidean space with a given mean curvature  (Russian), Geometry and
Topology, Leningrad. Gos. Ped. Inst. im. Gercena, Leningrad, 1 ({1974}), 3--10.

\bibitem{} Banchoff, T. F.,   Tightly embedded 2-dimensional
polyhedral manifolds, Amer. J. Math.
 87 ({1965}), 462--272.

\bibitem{} Banchoff, T. F.,   The spherical two-piece property and tight surfaces in
spheres, J. Differential Geometry, 4 (1970), 193--205.

\bibitem{} Barbosa, J. L.,{1975}
On minimal immersions of $S^{2}$ into
$S^{2m}$, Trans. Amer. Math. Soc., 210
, 75--106.

\bibitem{} Barbosa, J. L. and M. do Carmo
,  {1974}  Stable minimal surfaces
, Bull. Amer. Math. Soc., 80,
581--583.

\bibitem{} Barbosa, J. L. and M. do Carmo ,  On the size of stable
minimal surfaces in $R^3$
, Amer. J. Math., 98 (1976),
515--528.

\bibitem{} Barbosa, J. L. and M. do Carmo A necessary condition
for a metric in $M^n$ to be minimally
immersed in $R^{n+1}$, An. Acad. Brasil.
Ci\^enc., 50 (1978), 451--454.

\bibitem{} Barbosa, J. L. and M. do Carmo, Stability of minimal
surfaces and eigenvalues of the Laplacian,
Math. Z., 173 (1980), 13--28.

\bibitem{} Barbosa, J. L. and M. do Carmo, Stability of minimal
surfaces in spaces of constant
curvature, Bol. Soc. Brasil. Mat., 11 (1980), 1--10.

\bibitem{} Barbosa, J. L. and M. do Carmo, Stability of minimal surfaces in a
3-dimensional hyperbolic space, Arch. Math., 36 ({1981}), 554--557.

\bibitem{} Barbosa, J. L. and M. do Carmo Stability of hypersurfaces with constant mean
curvature, Math. Z., 185 (1984), 339--353.

\bibitem{} Barbosa, J. L., do Carmo, M. and J. Eschenburg, Stability of
hypersurfaces of constant mean curvature in Riemannian manifolds, Math. Z., 197 ({1988}), 123--138.

\bibitem{} Barbosa, J. L. and A. G. Colares, Minimal surfaces in $R^3$, Lecture Notes in Math., Springer-Verlag, 1195. ({1986}).


\bibitem{} Barbosa, J. L. and A. G. Colares,  Stability of hypersurfaces with constant $r$-mean curvature, Ann. Global Anal. Geom., 15
(1997), 277--297.

\bibitem{} Barbosa, J. L.,  Dajczer, M. and L. P. Jorge, Minimal ruled submanifolds in spaces of constant curvature, Indiana Univ. Math. J. , 33
({1984}), 531--547.

\bibitem{}  Barros, M. and B. Y. Chen,   Stationary 2-type surfaces in a hypersphere, J. Math. Soc. Japan, 39 (1987), 627-648. 

\bibitem{}  Barros, M. and B. Y. Chen,   Spherical submanifolds which are of 2-type via the
second standard immersion of the sphere,
Nagoya Math. J., 108 (1987), 77--91.  

\bibitem{} Barros, M. and O. J. Garay,  2-type surfaces in $S^3$, Geom. Dedicata, 24 ({1987}), 329--336.

\bibitem{} Barros, M. and O. J. Garay,  Hopf submanifolds in $S^7$ which are
Willmore-Chen submanifolds, Math. Z.,
228 (1988), 121--129.

\bibitem{} Barros, M. and O. J. Garay,  Willmore-Chen tubes in homogeneous spaces in
warped product spaces, Pacific J. Math. , 188 (1987), 201--207.

\bibitem{} Beckenbach, E. F. and T. Rad\'o, Subharmonic functions and surfaces of negative curvature, Trans. Amer. Math. Soc., 35 (1933), 662--674.
 
\bibitem{} Beez, R.,     Zur Theorie der Kr\"ummungmasses
von Mannigfaltigkeiten h\"oherer Ordunung,
Z. Math. Physik, 21({1876}), 373--401.

\bibitem{} Bejancu, A.,     $CR$-submanifolds of a Kaehler manifold I, Proc. Amer. Math.
Soc., 69 ({1978}), 134--142.

\bibitem{} Bejancu, A.,  On the geometry of leaves on a
$CR$-submanifold, Ann. St. Univ. Al. I. Cuza Iasi, 25 (1979), 393--398.

\bibitem{} Bejancu, A.,   Umbilical $CR$-submanifolds of a Kaehler manifold,
Rend. Mat., 12 (1980), 439--445.

\bibitem{} Bejancu, A.,    Geometry of $CR$-submanifolds, D. Reidel Publishing
Co., Dordrecht, 1986.

\bibitem{}  Bejancu, A., Kon, M. and K. Yano, CR-submanifolds of a
complex space form, J. Differential
Geometry, 16 (1981), 137--145.

\bibitem{} Beltrami, E.,    Teoria fondamentale degli spazii di
curvatura constante, Ann. di Mat., 2 (1868), 232--255.

\bibitem{} B\'erard, P.,    Remarques sur l'\'equation de J. Simons, Pitmon
Monog. and Surveys in Pure and Appl. Math., 52 (1991), 47--57.

\bibitem{} B\'erard, P. and G. Besson, On the number of bound
states and estimates on some geometric
invariants, Partial differential
equations (Rio de Janeiro, 1986),
Lecture Notes in Math., 1324, 30--40 ({1988}).


\bibitem{} B\'erard, P. and G. Besson, Number of bound
states and estimates on some geometric
invariants, J. Funct. Anal., 94 (1990), 375--396.

\bibitem{} B\'erard, P., do Carmo, M. and W. Santos, The  index
of constant mean curvature surfaces in
hyperbolic $3$-space, Math. Z., 224 (1997), 313--326.

\bibitem{} Berger, M.,   Les espaces symmetriques non compacts, Ann. Sci. \'Ecole Norm.
Sup., 74 ({1957}), 85--177.

\bibitem{} Berger, E.,  Bryant, R. and P. Griffiths, The Gauss equations and rigidity of isometric embeddings,  Duke Math. J., 50 ({1983}), 803--892.

\bibitem{} Bernatzki, F., The Plateau-Douglas problem fo nonorientable minimal surfaces, 
Manusc. Math., 79 ({1993}), 73--80.

\bibitem{} Berndt, J., Bolton, J. and L. M. Woodward, Almost complex
curves and Hopf hypersurfaces in the nearly Kaehler $6$-sphere,  Geom. Dedicata ,
56 ({1995}), 237--247. 

\bibitem{} Besse, A., Manifolds all of whose geodesics are
closed, Springer-Verlag 1978.

\bibitem{} Besson, G., Courtois, G. and S. Gallot, Entropies et rigidites des espaces localement symetriques de courbure strictement negative, Geom. Funct. Anal., 5 ({1995}),
731--799. 

\bibitem{} Bianchi, L., Lezioni di geometria differenziale II, Pisa, {1903}.

\bibitem{} Blair, D. E.,  A generalization of the catenoid, Canad. J. Math. 27 ({1975}), 231--236.

\bibitem{} Blair, D. E.,  Contact Manifolds in Riemannian Geometry, Lecture Notes in
Math., 509, 1976.

\bibitem{}  Blair, D. E. and B. Y. Chen, On $CR$-submanifolds of
Hermitian manifolds, Israel J. Math., 34 ({1979}), 353--363.

\bibitem{} Blair, D. E.,  Dillen, F., Verstraelen, L. and L. Vrancken,
 Calabi curves as holomorphic Legendre
curves and Chen's inequality, Kyunpook Math. J. 35 ({1996}), 407--416.

\bibitem{} Blair, D. E.  and J. R. Vanstone, A generalization of
the helicoid, Minimal Submanifolds and
Geodesics, Kaigai Publ., Tokyo, {1978},13--16.

\bibitem{} Blanusa, D., \"Uber die Einbettung hyperbolischer R\"aume
in euklidische R\"aume, Monatsch. Math., 59 (1955), 217--229.

\bibitem{}  Blaschke, W., Vorlesungen \"uber Differentialgeometrie, vol. III, Springer
Verlag, Berlin,  {1929}.

bibitem{}  Blaschke, W., Sulla geometria
differenziale delle superficie $S_2$ nello
spazio euclideo $S_4$, Ann. Mat. Pura
Appl., 28 (1949), 205--209.

\bibitem{} Bleecker, D. D. and J. L. Weiner, Extrinsic bounds on $\lambda_1$ of $\Delta$ on a compact manifold , Comment Math. Helv. 51 ({1976}), 601--609.

\bibitem{} Bloss, D., Elliptische funktionen und vollst\"andige
Minimalfl\"achen, Doctoral Thesis, Freien Universit\"at Berlin, {1989}. 

\bibitem{} Blum, R.,   Circles on surfaces in the Euclidean 3-space ,  Lecture Notes in Math.,
792 ({1980}), 213--221. 

\bibitem{} B\"ohme, R. and A. J. Tromba,  The index theorem for classical minimal surfaces, 
Ann. Math., 113 ({1981}), 447--499. 

\bibitem{} Bolton, J., Jensen, G. R., Rigoli, M. and L. M. Woodward,   On conformal minimal
immersions of $S^2$ into $CP^n$, Math. Ann., 279 ({1988}), 599--620.

\bibitem{}  Bolton, J.,  Pedit F. and L. M.
Woodward,{1995} Minimal surfaces
and the affine Toda field model, J. Reine
Angew. Math., 459, 119--150.

\bibitem{}   Bolton, J.,  Scharlach, C., Vrancken L. and L. M. Woodward
, Lagrangian submanifolds
of $CP^3$ satisfying Chen's equality, {1998}.

\bibitem{}  Bolton, J., Vrancken, L. and L. M. Woodward, On almost
complex curves in the nearly Kaehler
$6$-sphere, Quart. J. Math. Oxford
Ser., 45 ({1994}), 407--427.

\bibitem{}  Bolton, J., Vrancken, L. and L. M. Woodward,  Totally
real minimal surfaces with non-circular ellipse of curvature in the 
nearly Kaeh\-ler 6-sphere,  J. London Math. Soc., 56 (1997), 625--644.

\bibitem{}  Bolton, J. and L. M. Woodward, Special submanifolds of $S^6$ with its $G_2$
geometry, Geometry, topology and physics (Campinas, 1996), de Gruyter,
Berlin,  59--68, {1997}.

\bibitem{} Bombieri, E., de Giorgi, E. and E. Giusti, Minimal cones
and the Bernstein problem, Invent. Math. 7 ({1969}), 243--268.

\bibitem{} Bonnet, O., Deuxi\`eme note sur les surfaces \`a lignes de courbure
sph\'eriques, C. R. Acad. Sc. Paris, 36 (1853), 389--391.

\bibitem{} Bonnet, O.,  Troisi\`eme note sur les surfaces \`a lignes des courbure
planes et sph\'eriques, C. R. Acad. Sc. Paris, 36 (1853), 585--587.

\bibitem{} Bonnet, O.,  M\'emoire sur l'emploi d'un nouveau syst\`eme de variables dans
l'\'etude des surfaces courbes, J. Math. Pures Appl., 2 (1860), 153--266.

\bibitem{} Bonnet, O.,  M\'emoire sur la the\'eorie des surfaces applicables sur
une surface donn\'ee, J. l'\'Ecole Polytech., 42 (1867), 31--151.

\bibitem{} Borel, A. and J. P. Serre,    Sur certains sousgroupes des groupes de Lie compacts, Comm. Math. Helv., 27 ({1953}), 128--139.

\bibitem{}  Borisenko, A. A.and Yu. A.
Nikolaevskii, Grassmann manifolds
and Grassmann image of submanifolds,
Russian Math. Survey, 46 ({1991}), 
45--94.

\bibitem{}  Borrelli, V., Chen, B. Y. and J. M. Morvan, Une caract\'erisation g\'om\'etrique de la sph\`ere de Whitney, C. R. Acad. Sci. Paris S\'er. I Math., 321 ({1995}), 1485--1490.

\bibitem{} Borsuk, K.,    Sur la courbure totale des courbes, Ann. de la Soc.
Polonaise, 20 ({1947}), 251--265.

\bibitem{} Bott, R. and H. Samelson,  Applications of the
theory of Morse to symmetric spaces, Amer. J. Math, 80 ({1958}), 964-1029.

\bibitem{} de Brito, F. F., Power series with Hadamard gaps and
 hyperbolic complete minimal surfaces, Duke Math. J., 68 ({1992}), 297--300. 

\bibitem{} Brinkmann, W. H.,  On Riemannian
spaces conformal to Euclidean space, Proc. Nat. Acad. Sci. U.S.A., 9 ({1923}), 1--3.

\bibitem{} Brothers, J. E.,  Stability of minimal
orbits, Trans. Amer. Math. Soc., 294 ({1986}), 537--552.

\bibitem{} Bryant, R. L., Submanifolds and special structures on the
octonians, J. Differential Geom.  17 ({1982}), 185--232.

\bibitem{} Bryant, R. L.,  A duality theorem for Willmore surfaces, J. Differential Geom., 20 (1984), 23--53.

\bibitem{} Bryant, R. L.,   Minimal surfaces of constant curvature in
$S^n$, Trans. Amer. Math. Soc., 290 (1985), 259--271.

\bibitem{} Bryant, R. L.,  Surfaces in conformal
 geometry, The mathematical heritage of
Hermann Weyl, Proc. Sympos. Pure Math., 48, Amer. Math.
Soc., 227--240, 1987.

\bibitem{} Bryant, R. L.,  Surfaces of mean curvature one in hyperbolic space, Ast\'erisque  1987, 154-155

\bibitem{} Bryant, R. L.,  Minimal Lagrangian  submanifolds of Kaehler-Einstein
manifolds, Differential geometry and
differential equations (Shanghai, 1985),  Lecture Notes in Math., 1255, 1--12. 

\bibitem{} Bryant, R. L., Griffiths, P. A. and D. Yang, Characteristics and existence of isometric embeddings, Duke Math. J., 50 ({1983}), 893--994.

\bibitem{}  Burago, Yu. D. and V. A. Zalgaller, Geometric Inequalities,
Springer-Verlag, Berlin, {1988} .

\bibitem{} Burns, J. M.,  Homotopy of compact symmetric
spaces, Glasgow Math. J., 34 ({1992}),
221--228. 

\bibitem{} Burns, J. M.,  Conjugate loci of totally geodesic submanifolds of symmetric
spaces, Trans. Amer. Math. Soc., 337 (1993), 411--425.

\bibitem{} Burns, J. M. and M. J. Clancy,  Polar sets as nondegenerate critical submanifolds in symmetric spaces, Osaka J. Math., 31({1994}), 533--559.

\bibitem{} Burstin, C., Ein Betrag zum Problem der Einbettung der
Riemannschen R\"aume in eukildische
R\"au\-me, Math Sb., 38 ({1931}), 74--85.

\bibitem{} Buyske, S., Lie sphere transformations and the focal sets of certain taut immersions, Trans. Amer. Math. Soc., 311 ({1989}),117--133.

\bibitem{}  Caffarelli, L., Nirenberg, L. and J. Spruck, Nonlinear second order
elliptic equations. IV. Starshaped compact
Weingarten hypersurfaces, Current Topics 
in Partial Differential Equations, 
Kinokuniya Company Ltd., Tokyo,
1--26 (1986).

\bibitem{} Calabi, E.,    Isometric imbedding of
complex manifolds, Ann. Math., 58 ({1953}), 1--23.

\bibitem{} Calabi, E.,  Construction and
properties of some $6$-dimensional almost
complex manifolds, Trans. Amer. Math. Soc., 87 (1958), 407--438.


\bibitem{} Calabi, E.,  Minimal immersions of surfaces in Euclidean spheres, J. Differential
Geometry, 1 (1967), 111--126.

\bibitem{} Calabi, E.,  Quelques applications de l'analyse complexe aux
surfaces d'aire minima, Topics in Complex Manifolds (ed. H. Rossi),
University of Montreal, 59--81, 1968.

\bibitem{} Callahan, M.,  Hoffman, D. and H. Karcher,  A family of singly
periodic minimal surfaces invariant under a
screw motion, Experiment. Math., 2 ({1993}), 157--182. 

\bibitem{} Callahan, M., Hoffman, D. and W. H.
Meeks, III, The structure of singly-periodic minimal surfaces,
 Invent. Math., 99 ({1990}), 455--481.

\bibitem{} do Carmo, M. and M. Dajczer, Riemannian metrics induced by two
immersions, Proc.. Amer. Math. Soc., 86 ({1982} ), 115--119.

\bibitem{} do Carmo, M. and M. Dajczer,  Rotation hypersurfaces in
spaces of constant curvature, Trans.
Amer. Math. Soc. , 277 (1983), 685--709.

\bibitem{} do Carmo, M. and M. Dajczer,  A rigidity theorem for
higher codimensions, Math. Ann. , 274 (1986), 577--583.

\bibitem{} do Carmo, M., Dajczer, M. and F.
Mercuri, Compact conformally flat hypersurfaces,
Trans. Amer. Math. Soc. 288 ({1985}), 189--203.

\bibitem{} do Carmo, M. and A. M. Da Silveira, Index and
total curvature of surfaces with constant
mean curvature,  Proc. Amer. Math. Soc., 110 ({1990}), 1009--1015.

\bibitem{} do Carmo, M. P. and H. B. Lawson, On Alexandrov-Bernstein theorems in hyperbolic space,  Duke Math. J., 50 ({1983}),
995--1003.

\bibitem{}  do Carmo, M. and C. K. Peng,   Stable complete minimal surfaces in $R^3$ are
planes, Bull. Amer. Math. Soc. 1 ({1979}), 903--905.

\bibitem{}  do Carmo, M. and C. K. Peng,  Stable complete minimal hypersurfaces,
Proc. Beijing Symp. Diff. Equa.  Diff. Geom., 3 (1980), 1349--1358.

\bibitem{} M. do Carmo and N. R. Wallach, 
Minimal immersions of spheres into spheres, Ann. Math., 93 ({1971}), 43--62.

\bibitem{}  Cartan, \'E.,  
Sur certains hypersurfaces de l'espace
conforme r\'eel a cinq dimensions, Bull.
Soc. Math. France, 46 ({1918}), 84--105.

\bibitem{}  Cartan, \'E.,   Sur la vari\'et\'es de courboure
constante d'un espace euclidien ou non euclidien, Bull. Soc. Math. France, 47 (1919), 125--160.

\bibitem{}  Cartan, \'E.,   Sur une classe remarquable d'espaces de
Riemann, Bull. Soc. Math. France,
54 ,(1926), 214--264; , ibid, 55 (1927), 114--134.

\bibitem{}  Cartan, \'E., Sur la possibilit\'e de plonger un espace
riemannien donn\'e dans un espace euclidien,
Ann. Soc. Math. Polon., 6 (1927), 1--17.

\bibitem{}  Cartan, \'E., Lecons sur la G\'eom\'etrie des Espaces de
Riemann, Gauthier-Villars, Paris, 1928.

\bibitem{}  Cartan, \'E.,  Familles de surfaces
isoparam\'etroques dans les espaces \`a
courbure constante, Ann. Mat. , 17 (1938), 177--191.

\bibitem{}  Cartan, \'E.,  Sur des familles remarquables d'hypersurfaces
isoparametriques dans les espaces sphe\-riques,  Math. Z., 45 (1939), 335--367.

\bibitem{}  Cartan, \'E.,  Sur quelques familles remarquables d'hypersurfaces, C. R.
Congr\`es Math. Li\`ege, 30--41, 1939.

\bibitem{}  Cartan, \'E.,  Sur des familles  d'hypersurfaces isoparam\'etriques des
espaces sph\'eriques \`a 5 et \`a 9 dimensions, Univ. Nac. Tucum\'an. Revista A., 1 (1940), 5--22.

\bibitem{} Carter, S. and U. Dursun, Isoparametric and
Chen submanifolds, Geometry and topology
of submanifolds, VIII, 41--45, {1995}.

\bibitem{} Carter, S. and U. Dursun,   On generalised Chen
and $k$-minimal immersions, Beitr\"age
Algebra Geom., 38 (1997), 125--134.

\bibitem{}  Carter, S. and A. West, Tight and taut immersions, Proc. London Math. Soc.,
25 ({1972}), 701--720.

\bibitem{}  Carter, S. and A. West,  Totally focal
embeddings, J. Differential Geometry, 13 (1978), 251--261.

\bibitem{}  Carter, S. and A. West,  A characterization of isoparametric
hypersurfaces in spheres, J. London Math. Soc.  26 (1082), 183--192.

\bibitem{}  Carter, S. and A. West,  Isoparametric systems
and transnormality,  Proc. London Math. Soc. ,51 (1985), 520--542.

\bibitem{}  Carter, S. and A. West,  Isoparametric and totally focal submanifolds, Proc.
London Math. Soc., 60 (1990), 609--624.
 
\bibitem{} Castro, I. and F. Urbano,  Lagrangian surfaces in
the complex Euclidean plane with conformal
Maslov form, T\^o\-hoku Math. J.,
45, ({1993}), 656--582.

\bibitem{} Castro, I. and F. Urbano,  Twistor holomorphic
Lagrangian surface in the complex projective
and hyperbolic pl\-anes, Ann. Global Anal. Geom., 13 (1995), 59--67.

\bibitem{} Castro, I. and F. Urbano,  Examples of unstable
Hamiltionian-minimal Lagrangian tori in $C^2$,  Compos. Math.,. 111 (1998), 1--14.

\bibitem{}  Cayley, A., On the cyclide, Quart. J.  Pure  Appl. Math., 12 ({1873}),148-165.
 
\bibitem{} Cecil, T., Geometric applications of critical point theory to
submanifolds of complex projective space, Nagoya Math. J., 55 ({1974}), 5--31.

\bibitem{} Cecil, T.,  Taut immersions of  noncompact surfaces into a Euclidean
$3$-space, J. Differential Geometry, 11 (1976), ,451--459.

\bibitem{} Cecil, T.,  On the Lie curvatures of Dupin hypersurfaces, Kodai Math. J., 13 (1990), 143--153.

\bibitem{} Cecil, T.,  Lie sphere geometry, Springer-Verlag, New York-Berlin, 1992.

\bibitem{} Cecil, T.,  Taut and Dupin submanifolds, in Tight and Taut
Submanifolds,  (T. Cecil and S.-S. Chern,
eds.),  MSRI Publications, Cambridge Univ. Press,  32 (1997), 135-180.

\bibitem{}  Cecil, T. E. and S. S. Chern, Tautness and Lie sphere
geometry, Math. Ann., 278 ({1987}), 381--399.

\bibitem{}  Cecil, T. E. and S. S. Chern, Dupin submanifolds in Lie sphere geometry,  Lecture Notes in Math., 1369, 1--48 ({1989}).

\bibitem{} Cecil, T. E. and G. E. Jensen, Dupin hypersurfaces
with four principal curvatures,
preprint, {1997}.

\bibitem{} Cecil, T. E. and G. E. Jensen,  Dupin hypersurfaces
with three principal curvatures,
Invent. Math., 132 (1998), 121--178.

\bibitem{} Cecil, T. E.  and P. Ryan,  Focal sets, taut embeddings and the
 cyclides of Dupin, Math. Ann., 236 ({1978}), 177--190.

\bibitem{} Cecil, T. E.  and P. Ryan, Focal sets of
submanifolds, Pacific J. Math., 78 (1978), 27--39.

\bibitem{} Cecil, T. E.  and P. Ryan, Tight and taut 
immersions into hyperbolic space, J.
London Math. Soc. 19 (1979), 561--572.

\bibitem{} Cecil, T. E.  and P. Ryan, Distance functions
and umbilic submanifolds of hyperbolic
space, Nagoya Math. J., 74 (1979), 67--75.

\bibitem{} Cecil, T. E.  and P. Ryan,
Conformal geometry and the cyclides of
Dupin, Canad. J. Math., 32 (1980),
767--782. 

\bibitem{} Cecil, T. E.  and P. Ryan, Tight spherical
 embeddings, Global differential
geometry and global analysis, Lecture Notes
in Math., 838 (1981), 94--104.

\bibitem{} Cecil, T. E.  and P. Ryan, Focal sets and real 
hypersurfaces in complex projective
space, Trans. Amer. Math. Soc., 269 (1982), 481--499.

\bibitem{} Cecil, T. E.  and P. Ryan,Tight and Taut
Immersions of Manifolds, Pitman Researh
Notes in Math., 107, 1985.

\bibitem{}  Chang, S., A closed hypersurface with constant scalar
 and mean curvatures in $S^4$ is
isoparametric,  Comm. Anal. Geom. 1 ({1993}), 71--100. 

\bibitem{} Chen, B. Y.,  Minimal hypersurfaces of
 an $m$-sphere, Proc. Amer. Math. Soc., 29 (1971), 375-380.

 \bibitem{} Chen, B. Y.,  On total curvature of immersed
 manifolds, I,  Amer. J. Math.  93 (1971), 148--162.

 \bibitem{} Chen, B. Y.,  On total curvature of immersed
 manifolds, II,  Amer. J. Math.  94 (1972), 899--907.

 \bibitem{} Chen, B. Y.,  On total curvature of immersed
 manifolds, III,  Amer. J. Math.  95 (1973), 636--642.
 
 \bibitem{} Chen, B. Y.,  Minimal surfaces with
constant Gauss curvature, Proc. Amer. Math.
Soc., 34 (1971), 504--508.

 \bibitem{} Chen, B. Y.,  On the surface with
parallel mean curvature vector, Indiana Univ.
Math. J.  22 (1973), 655--666.

 \bibitem{} Chen, B. Y.,  Geometry of
submanifolds, Mercel Dekker, New York, 1973.
.
 \bibitem{} Chen, B. Y.,  An invariant of conformal mappings, Proc. Amer. Math.
Soc., 40  (1973), 563--564.

 \bibitem{} Chen, B. Y.,  On a variational
problem on hypersurfaces, J. London Math. Soc., 6 (1973), 321--325.

 \bibitem{} Chen, B. Y.,  Some conformal
invariants of submanifolds and their
applications, Boll. Un. Mat. Ital  10 (1974), 380--385.

 \bibitem{} Chen, B. Y.,  Mean curvature
vector of a submanifold, Proc.
Symp. Pure Math. Amer. Math. Soc.,
27-I, 119--123, 1975.

 \bibitem{} Chen, B. Y.,  Extrinsic spheres in
K\"ahler manifolds, Michigan Math. J., 23 (1976), 327--330.

 \bibitem{} Chen, B. Y.,  Total mean curvature
of immersed surfaces in $E^m$, Trans. Amer. Math. Soc., 218 (1976), 333--341.

 \bibitem{} Chen, B. Y.,  Some relations between differential geometric invariants
and topological invariants of submanifolds, Nagoya Math. J., 60 (1976), 1--6.

 \bibitem{} Chen, B. Y.,  Extrinsic
spheres in compact symmetric spaces are
intrinsic spheres, Michigan Math. J., 24 (1977), 265--271.

 \bibitem{} Chen, B. Y.,  Totally umbilical
submanifolds of Cayley plane, Soochow J. Math., 3 (1977), 1--7.

 \bibitem{} Chen, B. Y.,  Totally
umbilical submanifolds of
quaternion-space-forms, J. Austral. Math. Soc.
Ser. A, 26 (1978), 154--162.

 \bibitem{} Chen, B. Y.,    On the total curvature of immersed manifolds, IV,
Bull. Inst. Math. Acad. Sinica, 7 (1979), 301--311.

 \bibitem{} Chen, B. Y., 
 Conformal mappings and first eigenvalue of Laplacian on surfaces, Bull. Inst. Math. Acad. Sinica, 7 (1979), 395--400.

 \bibitem{} Chen, B. Y.,  Classification of
totally umbilical submanifolds in symmetric
space, J. Austral. Math. Soc. Ser. A, 30 (1980), 129--136.

 \bibitem{} Chen, B. Y.,  Surfaces with parallel
normalized mean curvature vector, Monat. Math. 90 (1980), 185--194.

 \bibitem{} Chen, B. Y., Totally umbilical
submanifolds of Kaehler manifolds, Archiv Math.
 36 (1980), 83-91.
 
 \bibitem{} Chen, B. Y.,  Geometry of 
submanifolds and its applications, Science
University of Tokyo, Tokyo, Japan, 1981.

 \bibitem{} Chen, B. Y.,   $CR$-submanifolds of a Kaehler manifold, I, II, J. Differential
Geometry, 16 (1981), 305--322; 493--509

 \bibitem{} Chen, B. Y.,  Cohomology of $CR$-submanifolds, Ann. Fac. Sc. Toulouse Math. Ser. V,, 3 (1981), 167--172.

 \bibitem{} Chen, B. Y.,  Differential  geometry of submanifolds with planar
normal sections,  Ann. Mat. Pura Appl., 130 (1982), 59--66.

 \bibitem{} Chen, B. Y.,  Classification of  surfaces with planar normal sections, J. Geometry, 20 (1983), 122--127. 

 \bibitem{} Chen, B. Y.,  On the first eigenvalue of Laplacian of compact
minimal submanifolds of rank one
symmetric spaces, Chinese J. Math., 11 (1983), 259--273.

 \bibitem{} Chen, B. Y.,   Some non-integrability
theorems of holomorphic distributions, Alg.
Diff. Topology/Glo\-bal Diff. Geom.,
Teubner-Verlag, 36--48, 1984.

 \bibitem{} Chen, B. Y.,  Total mean curvature and submanifolds of finite type, World
Scientific Publ., Singapore-New Jersey-London, 1984.

 \bibitem{} Chen, B. Y., A new approach
to compact symmetric spaces and applications,
Katholieke Universiteit Leuven, Leuven, Belgium, 1987.

 \bibitem{} Chen, B. Y.,  Some estimates of
total tension and their applications, Kodai Math.
J., 10 (1987), 93-101.

 \bibitem{} Chen, B. Y.,  Surfaces of finite
type in Euclidean 3-space, Bull. Soc.
Math. Belg. Ser. B, 39 (1987), 243--254.

 \bibitem{} Chen, B. Y.,  Null 2-type surfaces
in Euclidean space, Algebra, Analysis
and Geometry, World Scientific, 1--18, 1988.

 \bibitem{} Chen, B. Y.,  Geometry of slant submanifolds, Katholieke Universiteit
Leuven, Leuven, Belgium, 1990 .

 \bibitem{} Chen, B. Y.,  Linearly independent,
orthogonal, and equivariant immersions, Kodai Math. J. 14 (1991), 341--349.

 \bibitem{} Chen, B. Y.,  Local rigidity
theorems of 2-type hypersurfaces in a
hypersphere, Nagoya Math. J., 122 (1991), 139--148. 

 \bibitem{} Chen, B. Y.,   Some pinching and
classification theorems for minimal
submanifolds, Archiv Math., 60 (1993), 568--578.

 \bibitem{} Chen, B. Y.,   A Riemannian invariant
and its applications to submanifold theory,
Results in Math., 27 (1995), 17--26

 \bibitem{} Chen, B. Y.,   A general inequality for
submanifolds in complex space forms and its
applications, Arc\-hiv Math., 67 (1996), 519--528.

 \bibitem{} Chen, B. Y.,  Jacobi's elliptic
functions and Lagrangian immersions,
Proc. Royal Soc. Edinburgh, 126 (1996), 687--704.

 \bibitem{} Chen, B. Y.,   Some geometric 
inequalities and their applications, Proc. 1st Intern. Workshop
Diff. Geom. (Taegu), 31--67, 1996.

 \bibitem{} Chen, B. Y.,   A report of
submanifolds of finite type, Soochow J. Math., 22 (1996), 117--337.

 \bibitem{} Chen, B. Y.,  Mean curvature and
shape operator of isometric immersions in
real space forms, Glasgow Math. J.  38 (1996), 87--97.

 \bibitem{} Chen, B. Y.,  Strings of Riemannian
invariants, inequalities, ideal immersions
and their  applications, Proc. Third Pacific Rim Geom. 
Conf., Intern. Press, Cambridge, MA, 7--60 (1998).

 \bibitem{} Chen, B. Y.,  Interaction of
Legendre curves and Lagrangian
submanifolds, Israel J. Math., 99 (1997), 69--108.

 \bibitem{} Chen, B. Y.,   Complex extensors
and Lagrangian submanifolds in complex
space forms, T\^ohoku Math. J., 49 (1997), 277--297.

 \bibitem{} Chen, B. Y.,  Some new
obstructions to minimal and Lagrangian
isometric immersions, Japan. J.
Math., 26 (2000), 105-127.

 \bibitem{} Chen, B. Y.,    A vanishing theorem
for Lagrangian immersions into Einstein-Kaehler manifolds and its
applications, Soochow J. Math., 24 (1998), 157--164.

 \bibitem{} Chen, B. Y.,  On slant surfaces, Taiwanese J. Math., 2 (1998).

 \bibitem{} Chen, B. Y.,  Intrinsic and extrinsic structures of Lagrangian
surfaces in complex space forms,
Tsu\-kuba J. Math., 22 (1998), 657--680.

 \bibitem{} Chen, B. Y.,  Special slant surfaces and a basic inequality,
Results in Math., 33 (1998), 65--78.

 \bibitem{} Chen, B. Y.,   Relations between
Ricci curvature and shape operator for
submanifolds with arbitrary codimensions,
Glasgow Math. J., 41 (1999), 33--41.

 \bibitem{} Chen, B. Y.,  Representation of
flat Lagrangian $H$-umbilical submanifolds
in complex Euclidean spa\-ces, T\^ohoku Math. J., 50 (1999),13--20.

\bibitem{} Chen, B. Y., Deprez, J. and P. Verheyen,   Immersions, dans un espace
euclidien, d'un espace sym\'etrique
compact de rang un \'a g\'eod\'esiques simples, C. R.
Acad. Sc. Paris, 304 ({1987}), 567--570. 

\bibitem{} Chen, B. Y., Deprez, J., Dillen, F.,  Verstraelen, L. and L. Vrancken,   
 Finite type curves, Geometry and Topology of Submanifolds, 2, 76--110, {1990}. 

\bibitem{} Chen, B. Y. and F. Dillen,   Surfaces of finite
type and constant curvature in the 3-sphere, C. R.
Math. Rep. Acad. Sci. Canada, 12 (1990), 47-49.

\bibitem{} Chen, B. Y. and F. Dillen,   Quadrics of
finite type, J. Geometry, 38 (1990), 16--22. 

\bibitem{} Chen, B. Y.,  Dillen, F. and H. Song, Quadric hypersurfaces
of finite type, Colloq. Math., 63 (1992), 145--152. 

\bibitem{} Chen, B. Y., Dillen, F., Verstraelen, L. and L. Vrancken,  Ruled surfaces of
finite type, Bull. Austral. Math. Soc., 42 (1990), 447--453. 

\bibitem{} Chen, B. Y., Dillen, F., Verstraelen, L. and L. Vrancken,  A variational minimal principle characterizes submanifolds of finite type, C.R. Acad. Sc. Paris, 317 (1993),
961--965 .

\bibitem{} Chen, B. Y., Dillen, F., Verstraelen, L. and L. Vrancken,   Totally real
submanifolds of $CP^n$ satisfying a basic
inequality, Archiv Math., 63 (1994), 553--564.

\bibitem{} Chen, B. Y., Dillen, F., Verstraelen, L. and L. Vrancken,  Two equivariant
totally real immersions into the nearly K\"ahler $6$-sphere and their
characterization, Japan. J. Math. (N.S.) 21 (1995), 207--222.

\bibitem{} Chen, B. Y., Dillen, F., Verstraelen, L. and L. Vrancken,   Characterizing a
class of totally real submanifolds of $S\sp
6$ by their sectional curvatures, T\^ohoku Math. J., 47 (1995), 185--198.

\bibitem{} Chen, B. Y., Dillen, F., Verstraelen, L. and L. Vrancken,   An exotic totally real
minimal immersion of $S^3$ in $CP^3$ and its
characterization, Proc. Royal Soc. Edinburgh Ser. A, Math., 126 (1996), 153--165

\bibitem{} Chen, B. Y., Dillen, F., Verstraelen, L. and L. Vrancken,   Lagrangian isometric
immersions of a real space form $M^n(c)$ into a complex space form $\widetilde
M^n(4c)$, Math. Proc. Cambridge Philo. Soc., 124 (1998), 107--125.

\bibitem{} Chen, B. Y. and C. S. Houh, On stable
submanifolds with parallel mean curvature, Quart. J. Math. Oxford, 26 ({1975}), 229--236.

\bibitem{} Chen, B. Y. and C. S. Houh,  Totally real
submanifolds of a quaternion projective space, Ann. Mat. Pura Appl., 70 (1979), 185--199.

\bibitem{} Chen, B. Y. and S. Ishikawa, Biharmonic
surfaces in pseudo-Euclidean spaces, Memoirs Fac. Sci. Kyushu
Univ. Ser. A, Math. , 45 (1991), 323--347.

\bibitem{} Chen, B. Y. and S. Ishikawa, On classification of some surfaces of
revolution of finite type , Tsukuba J. Math. 17 (1993), 287--298.

\bibitem{} Chen, B. Y. and S. Jiang, Inequalities between
volume, center of mass, circumscribed radius,
order, and mean curvature, Bull. Soc.
Math. Belg. (New Series), 2 ({1995}), 75--85.

\bibitem{} Chen, B. Y.  and W. E. Kuan,  The Segre imbedding and
its converse, Ann. Fac. Sc. Toulouse Math.
Ser. V, 7 ({1985}), 1--28.

\bibitem{} Chen, B. Y. and S. J. Li,  Classification of
surfaces with pointwise planar normal sections and its application to
Fomenko's conjecture,  J. Geom., 26 ({1986}), 21--34.

\bibitem{} Chen, B. Y. and S. J. Li,  3-type
hypersurfaces in a hypersphere ,
Bull. Soc. Math. Belg. S\'er. B , 43 (1991), 135--141.

\bibitem{} Chen, B. Y.  and G. D. Ludden,  Surfaces with mean
curvature vector parallel in the normal bundle,
Nagoya Math. J., 47 (1972),161--168.

\bibitem{} Chen, B. Y.  and H. S. Lue, Differential geometry of $SO(n+2)/SO(2)\times SO(n)$, I, Geom. Dedicata, 4 (1975), 253--261.

\bibitem{} Chen, B. Y.  and H. S. Lue,  On normal connection
of Kaehler submanifolds, J. Math. Soc.
Japan, 27 (1975), 550--556.

\bibitem{} Chen, B. Y.  and H. S. Lue, Some 2-type submanifolds and applications, Ann.
Fac. Sc. Toulouse Math. Ser. V, 9 (1988), 121--131.

\bibitem{} Chen, B. Y. Ludden, G. D. and S. Montiel,    Real submanifolds of a Kaehler manifold, Algebra, Groups and Geom., 1 ({1984}), 176--212

\bibitem{} Chen, B. Y. and J. M. Morvan, Geometrie des surfaces
lagrangiennes de $C^2$, J. Math. Pures Appl, 66 ({1987}), 321--325.

\bibitem{} Chen, B. Y. and J. M. Morvan,  Cohomologie des
sous-vari\'et\'es $\alpha$-obliques, C. R. Acad. Sci. Paris, 314 (1992), 931--934.

\bibitem{} Chen, B. Y. and J. M. Morvan,  Deformations of isotropic submanifolds in K\"ahler
manifolds, J. Geom. Phys., 13 (1994), 79--104.

\bibitem{} Chen, B. Y. and T. Nagano, Totally geodesic
submanifolds of symmetric spaces, I, Duke Math. J., 44 (1977), 745--755.

\bibitem{} Chen, B. Y. and T. Nagano, Totally geodesic
submanifolds of symmetric spaces, II, Duke
Math. J., 45 (1978), 405--425.

\bibitem{} Chen, B. Y. and T. Nagano, Un invariant
g\'eometriques Riemannien, C. R. Acad.
Sc. Paris, 295 (1982), 389--391.

\bibitem{} Chen, B. Y. and T. Nagano, Harmonic metrics,
harmonic tensors and Gauss map, J. Math. Soc. Japan, 36 (1984), 295-313.

\bibitem{} Chen, B. Y. and T. Nagano, A Riemannian invariant
and its applications to a problem of Borel and
Serre, Trans. Amer. Math. Soc. , 308 (1988), 273--297.

\bibitem{} Chen, B. Y. and K. Ogiue,  Some
extrinsic results for Kaehler submanifolds,
Tamkang J. Math., 4 (1973), 207--213

\bibitem{} Chen, B. Y. and K. Ogiue,   Some
characterization of complex space forms, Duke Math. J., 40 (1973), 797--799

\bibitem{} Chen, B. Y. and K. Ogiue,  A characterization of
the complex sphere, Michigan Math. J., 21 (1974), 231--232
.
\bibitem{} Chen, B. Y. and K. Ogiue,  On totally real
submanifolds, Trans. Amer. Math. Soc., 193 (1974), 257--266.

\bibitem{} Chen, B. Y. and K. Ogiue,  Two theorems on Kaehler
manifolds, Michigan Math. J., 21 (1974), 225--229.
\bibitem{} Chen, B. Y. and K. Ogiue, Some characterizations
of complex space forms in terms of Chern
classes, Quart. J. Math. Oxford, 26 (1975), 459--464.

\bibitem{} Chen, B. Y. and M. Okumura, Scalar curvature,
inequality and submanifold, Proc. Amer.
Math. Soc. 38 (1973), 605--608.

\bibitem{} Chen, B. Y. and M. Petrovic,  On spectral
decomposition of immersions of finite type, Bull.
Austral. Math. Soc., 44 (1991),117--129.

\bibitem{} Chen, B. Y. and Y. Tazawa, Slant surfaces of
codimension two, Ann. Fac. Sc. Toulouse,
Math. Ser. V., 11 (1990), 29--43 .

\bibitem{} Chen, B. Y. and Y. Tazawa, Slant submanifolds in
complex Euclidean spaces, Tokyo J. Math.14  (1991), 101--120.

\bibitem{} Chen, B. Y. and Y. Tazawa, Slant submanifolds of complex projective and
complex hyperbolic spaces, Glasgow Math. J. , 42  (2000), 439--454.

\bibitem{} Chen, B. Y. and L. Vanhecke, Differential geometry of
geodesic spheres, J.  Reine  Angew. Math., 325 (1981), 28--67.

\bibitem{} Chen, B. Y. and P. Verheyen, Sous-vari\'et\'es dont
les sections normales sont des
g\'eod\'e\-siques, C.R. Acad. Sc.
Paris, 293 ({1981}), 611-613 .
\bibitem{} Chen, B. Y. and P. Verheyen,  Totally umbilical
submanifolds of Kaehler manifolds, Bull.
Math. Soc. Belg., 35 (1983), 27--44.

\bibitem{} Chen, B. Y. and P. Verheyen, Submanifolds with
geodesic normal sections, Math. Ann., 269 (1984), 417--429.

\bibitem{} Chen, B. Y. and L. Verstraelen, A characterization of
quasiumbilical submanifolds and its applications, Boll. Un. Mat. Ital., 14 (1977), 49--57. 

\bibitem{} Chen, B. Y. and L. Verstraelen, Errata to ``A characterization of
quasiumbilical submanifolds and its applications'', Boll. Un. Mat. Ital., 14 (1977), 634.

\bibitem{} Chen, B. Y. and L. Verstraelen,  Hypersurfaces of
symmetric spaces, Bull. Inst. Math. Acad. Sinica, 8 (1980), 201--236.

\bibitem{} Chen, B. Y.  and L. Vrancken,  Lagrangian submanifolds 
satisfying a basic inequality, Math. Proc.
Cambridge Phil. Soc., 120 (1996), 291--307.

\bibitem{} Chen, B. Y.  and L. Vrancken,  Existence and
uniqueness theorem for slant immersions and its
applications, Results Math.  31 (1997), 28--39.

\bibitem{} Chen, B. Y.  and L. Vrancken,   $CR$-submanifolds of
complex hyperbolic spaces satisfying a basic equality, Israel J.
Math. 110 (1999), 341--358.

\bibitem{} Chen, B. Y.  and L. Vrancken,   Lagrangian
submanifolds of the complex hyperbolic
space,  Tsukuba J. Math. 26 (2002), no. 1, 95--118.

\bibitem{} Chen, B. Y.  and L. Vrancken,  Addendum to ``Existence and
uniqueness theorem for slant immersions and its
applications'', Results Math.  39 (2001), 18--22.

\bibitem{} Chen, B. Y. and B. Wu,  Mixed foliate
CR-submanifolds in a complex hyperbolic space
are nonproper, Internat. J. Math. \& Math. Sci. 11 ({1988}), 507--515.

\bibitem{} Chen, B. Y. and S. Yamaguchi,  Classification of
surfaces with totally geodesic Gauss
image, Indiana Univ. Math. J., 32 ({1983}), 143--154.

\bibitem{} Chen, B. Y. and S. Yamaguchi,   Submanifolds with
totally geodesic Gauss images, Geom. Dedicata, 15 (1984), 313--322.

\bibitem{} Chen, B. Y.  and J. Yang,
 Elliptic functions, Theta function and hypersurfaces
satisfying a basic equality,  Math. Proc. Cambridge Phil. Soc, 125 (1999), 463--509.

\bibitem{} Chen, B. Y. and K. Yano, Sous-vari\'et\'e localement
conformes \`a un espace euclidien,
C. R. Acad. Sc. Paris, 275(1972), 123--126.

\bibitem{} Chen, B. Y. and K. Yano, Umbilical submanifolds with respect to a
nonparallel normal direction, J. Differential Geometry, 8 (1973), 589--597.

\bibitem{} Chen, C. C., Complete minimal surfaces with
total curvature $-2\pi$, Boll. Soc. Mat. Brasil. 10 ({1979}), 71--76.

\bibitem{} Chen, C. C., Elliptic functions and non-existence of complete minimal
surfaces of certain type, Proc. Amer.
Math. Soc.  79 (1980), 289--293.

\bibitem{} Chen, C. C. and F. Gackstatter,
 Elliptische und hyperelliptische
Funktionen und vollst\"andige
Minimalfl\"achen vom Enneperschen Type,
Math. Ann., 259 ({1982}), 359--369.

\bibitem{} Chen, C. S., Tight embedding and projective
transformation, Amer. J. Math., 101 ({1979}), 1083--1102.

\bibitem{} Chen, G. and X. Zou,  Rigidity of compact
submanifolds in a unit sphere, Kodai Math. J., 18 ({1995}), 75--85.

\bibitem{} Chen, Q., On the total curvature and area growth of
minimal surfaces in $R^n$, Manus. Math., 92 ({1997}), 135--142.

\bibitem{} Cheng, J. H.,   An integral formula on the scalar curvature of algebraic manifolds, Proc. Amer. Math. Soc., 81 ({1981}), 451--454

\bibitem{} Cheng, Q. M, Complete minimal hypersurfaces in $S^4(1)$ with constant
scalar curvature, Osaka J. Math., 27 ({1990}), 885--892. 

\bibitem{} Cheng, Q. M,   The rigidity of Clifford torus $S^1(\sqrt{1/n})\times
S^{n-1}(\sqrt{(n-1)/n})$ , Comm. Math. Helv., 71 (1996), 60--69.

\bibitem{} Cheng, S. Y. and J. Tysk, An index characterization of the catenoid
 and index bounds for minimal surfaces in
$R^4$, Pacific J. Math., 134 ({1988}), 251--260.

\bibitem{} Cheng, S. Y. and J. Tysk, Schrodinger operators and index bounds for
minimal submanifolds,  Rocky Mountain J.
Math., 24 (1994), 977--996. 

\bibitem{} Chern, S. S.,  Minimal surfaces in an Euclidean space of $N$ dimensions,
Diff. and Comb. Topology, Princeton, 187--198, {1965}.

\bibitem{} Chern, S. S.,   On Einstein hypersurfaces in a Kaehlerian manifold of
constant holomorphic sectional curvature, J.
Differential Geometry, 1, 21--31, 1967.

\bibitem{} Chern, S. S.,  On minimal spheres in
the four sphere, Studies and Essays
Presented to Y. W. Chen, Taiwan, 137--150, 1970.

\bibitem{} Chern, S. S., do Carmo, M. and S. Kobayashi,  Minimal
submanifolds of a sphere with second
fundamental form of constant length,
Functional Analysis and Related Fields, Springer-Verlag, , 59--75, 1970.

\bibitem{} Chern, S. S. and N. H. Kuiper,  Some theorems on the
isometric imbedding of compact Riemann
manifolds in Euclidean space, Ann. Math., 156 (1952), 422--430.

\bibitem{} Chern, S. S. and R. K. Lashof,  On the total curvature
of immersed manifolds, Amer. J. Math., 79 (1957), 306--318. 

\bibitem{} Chern, S. S. and R. K. Lashof,   On the total curvature
of immersed manifolds, II, Michigan Math. J., 5 (1958), 5--12. 

\bibitem{}  Chern, S. S. and R. Osserman,  Complete minimal
surfaces in euclidean $n$-space, J. Analyse Math. ,19 ({1967}), 15--34.

\bibitem{} Choe, J., Index, vision number and stability of  complete minimal surfaces, Arch. Rational Mech. Anal., 109 ({1990}), 195--212.

\bibitem{} Choi, H. I., Meeks, W. H., III and B. White, A
rigidity theorem for properly embedded minimal
surfaces in $R^3$,  J. Differential Geometry, 32 (1990), 65--76. 

\bibitem{}  Choi, H. I. and R. Schoen, The space of minimal
embeddings of a surface into a
three-dimensional manifold of positive
Ricci curvature,  Invent. Math., 81 (1985), 387--394.

\bibitem{} Christoffel, E. B., \"Uber einige allgemeine
Eigenschaften der Minimumsfl\"achen, J.
Reine Angew. Math., 67 (1867), 218-228.

\bibitem{}  Codazzi, D., Sulle coordinate curvilinee, 
Ann. di Mat. Pura Appl., 2 (1868), 1101--119.

\bibitem{} Collin,  P., Topologie et courbure
des surfaces minimales proprement plong\'ees
de $R^3$, Ann. of Math.  145 (1997), 1--31.

\bibitem{} Console, S. and A. Fino, Homogeneous
structures on Kaehler submanifolds of complex
projective spaces,  Proc. Edinburgh Math. Soc., 39 (1996), 381--395.

\bibitem{} Copson, E. C. and H. S. Ruse,  Harmonic Riemannian
spaces,  Proc. Royal Soc. Edinburgh, 60 (1939), 117--133. 

\bibitem{} Costa, C. J., Example of a complete minimal immersion in
$R^3$ of genus one and three embedded
ends,  Bol. Soc. Brasil. Mat. , 15 (1984), 47--54. 

\bibitem{} Costa, C. J.,Classification of complete minimal surfaces in
 $R^3$ with total curvature $-12\pi$, Invent. Math., 105 (1991), 273--303. 

\bibitem{} Costa, C. J., Complete minimal surfaces in $R^3$ of finite
topology and infinite total curvature,  Mat. Contemp., 4 (1993), 79--93.

\bibitem{} Costa, C. J. and P. A. Q. Simoes, Complete minimal
surfaces of arbitrary genus  in a slab of
$R^3$, Ann. Inst. Fourier (Grenoble), 46 (1996), 535--546.

\bibitem{} Courant, R., The existence of
minimal surfaces of least area bounded a
prescribed Jordan arcs and prescribed
surfaces, Proc. Nat. Acad. Sci. USA, 24 (1938), 97--102.

\bibitem{} Courant, R.,  The existence of minimal surfaces of given topological
structure under prescribed boundary condition, Acta Math., 72 (1940), 51--98.

\bibitem{} Curtin, E., Tautness for manifolds with boundary, Houston J. Math., 20 ({1994}),
409--424.

\bibitem{} Dajczer, M. and L. A. Florit, On conformally flat
submanifolds, Comm. Anal. Geom., 4 ({1996}), 261--284. 

\bibitem{} Dajczer, M. and L. A. Florit,  A class of austere
submanifolds, Illinois J. Math., 45 (2001), 735--755.

\bibitem{} Dajczer, M. and L. A. Florit,  On Chen's basic equality, Illinois J. Math., 42 (1998),
97--106.

\bibitem{} Dajczer, M. and D. Gromoll,  Real Kaehler
submanifolds and uniqueness of the Gauss
map, J. Differential Geometry, 22 (1985), 13--28.

\bibitem{} Dajczer, M. and D. Gromoll, The Weierstrass representation for
complete minimal real Kaehler
submanifolds of codimension two, Invent. Math., 119 (1985), 235--242.

\bibitem{} Dajczer, M. and L. Rodriquez, Rigidity of real
Kaehler submanifolds, Duke Math. J., 53 ({1986}), 211--220.

\bibitem{} Dajczer, M. and L. Rodriquez,  Complete real Kaehler
submanifolds, J. Reine Angew. Math., 419 (1991), 1--8.

\bibitem{} Dajczer, M. and R. Tojeiro, On flat surfaces in
space forms, Houston J. Math., 21 (1995), 319--338. 

\bibitem{} Dajczer, M. and R. Tojeiro, Flat totally real submanifolds of $CP^n$
and the symmetric generalized wave equation,
T\^ohoku Math. J. 47 (1995), 117--123.
 
\bibitem{} Dal Lago, W., Garc\'ia, A. and C. S\'anchez, Planar
normal sections on the natural imbedding of a flag manifold, Geom. Dedicata,
53 ({1994}), 223--235.

\bibitem{} Damek, E. and E. Ricci, A class of nonsymmetric harmonic Riemannian
spaces,  Bull. Amer. Math. Soc. (N.S.), 27 (1992),139--142. 

\bibitem{} Dao, C. T., Minimal real currents on compact Riemannian
 manifolds, Math. USSR Izv., 11 ({1977}), 807-820.

\bibitem{} Darboux, G., Lecons sur la th\'eorie g\'en\'erale des
surfaces, Premi\`ere partie, Gauthier-Villars, Paris, {1887}.

\bibitem{} Da Silveira, A. M., Stability of complete noncompact
surfaces with constant mean curvature, 
Math. Ann., 277 (1987), 629--638.

\bibitem{}  Defever, F., Hypersurfaces of $E^4$
satisfying $\Delta {H} = \lambda{H}$ ,
Michigan Math. J., 44 (1997), 355-363.

\bibitem{}  Defever, F., Hypersurfaces of $ E^4$ with harmonic mean curvature vector , 
Math. Nachr., 196 (1998), 61--69.

 \bibitem{} Defever, F. Deszcz, R. and L. Verstraelen, The compact cylcides of Dupin and a conjecture of B.-Y. Chen, J. Geometry, , 46 ({1993}), 33--38.

\bibitem{} Defever, F. Deszcz, R. and L. Verstraelen, The Chen--type of
noncompact cylcides of Dupin , Glasgow
Math. J., 36 (1994), 71--75.

\bibitem{} Defever, F.,  Mihai, I. and L. Verstraelen,  B. Y. Chen's
inequality for $C$-totally real submanifolds in Sasakian space
forms, Boll. Un. Mat. Ital. (B), 11 ({1997}) 365--374.

\bibitem{} Delaunay, C., Sur la surface de r\'evolution dont la
courbure moyenne est constante, J. Math. Pures Appl., 6 (1841), 309--320.

\bibitem{}  Deprez, J., Immersions of finite type of compact
homogeneous Riemannian manifolds  Doctoral Thesis,
Katholieke Universiteit Leuven, Belgium 1988.

\bibitem{} Deszcz, R.,  Dillen, F.,  Verstraelen, L. and L. Vrancken,
  Quasi-Einstein totally real
 submanifolds of $S^6(1)$, Tohoku Math. J. (2) 51 (1999), 461--478. 

\bibitem{} DeTurck, D. and W. Ziller,  Minimal isometric
immersions of spherical space forms in
spheres,  Comment. Math. Helv., 67 (1992), 428--458. 

\bibitem{} Dillen, F., Ruled submanifolds of finite type,
Proc. Amer. Math. Soc., 114 (1992), 795--798.

\bibitem{} Dillen, F., Opozda, B., Verstraelen, L., and L. Vrancken, On totally real
$3$-dimensional submanifolds of the nearly Kaehler $6$-sphere, Proc. Amer. Math.
Soc., 99 (1987), 741--749.
\bibitem{} Dillen, F., Opozda, B., Verstraelen, L., and L. Vrancken,   On almost complex
surfaces of the nearly Kaehler $6$-sphere, Zb. Rad. (Kragujevac), 8 (1987), 5--13.

\bibitem{} Dillen, F., Opozda, B., Verstraelen, L., and L. Vrancken,  On totally
real surfaces of the nearly Kaehler
$6$-sphere, Geom. Dedicata, 7 (1988), 325--334.

\bibitem{} Dillen, F., Pas, J. and L.  Verstraelen, On surfaces of finite type in  Euclidean
3-space, Kodai Math. J., 13 (1990), 10--21.

\bibitem{} Dillen, F., Petrovic, M. and L. Verstraelen, Einstein,
conformally flat and semi-symmetric 
submanifolds satisfying Chen's equality, Israel J. Math., 100 (1997), 163--169


\bibitem{}  Dillen, F., Verstraelen, L. and L. Vrancken, On almost
complex surfaces of the nearly Kaehler $6$-sphere. II, Kodai Math. J., 10 (1987), 261--271.

\bibitem{}  Dillen, F., Verstraelen, L. and L. Vrancken,   Classification of
totally real $3$-dimensional submanifolds of
$S^6(1)$ with $K\ge 1/16$, J. Math. Soc. Japan, 42 (1990), 565--584.

\bibitem{}  Dillen, F., Verstraelen, L. and L. Vrancken,   On problems of U. Simon
concerning minimal submanifolds of the nearly Kaehler $6$-sphere, Bull. Amer.
Math. Soc. (N.S.) , 19 (1988), 433--438.

\bibitem{} Dillen, F., Verstraelen, L., Vrancken, and G. Zafindratafa,  Classification of
polynomial translation hypersurfaces of
finite type, Results Math., 27 (1995), 244--249. 

\bibitem{}  Dillen, F. and L. Vrancken, Totally real submanifolds in $S\sp 6(1)$ satisfying Chen's equality, Trans. Amer. Math. Soc., 348 (1996), 1633--1646.

\bibitem{}  Dimitri\'c, I., Quadric representation
and submanifolds of finite type,  Doctoral
Thesis, Michigan State University, 1989.

\bibitem{}  Dimitri\'c, I.,  Spherical
hypersurfaces with low type quadric
representation, Tokyo J. Math., 13 (1990), 469--492.

\bibitem{}  Dimitri\'c, I., 1-type submanifolds of
the complex projective space, Kodai
Math. J., 14 (1991), 281--295.

\bibitem{}  Dimitri\'c, I.,  Quadric representation
of a submanifold, Proc. Amer. Math.
Soc., 114 (1992), 201--210.

\bibitem{}  Dimitri\'c, I.,  {1992b} Submanifolds of $E^m$
with harmonic mean curvature vector,
Bull. Inst. Math. Acad. Sinica, 20 (1992), 53--65.

\bibitem{}  Dimitri\'c, I., Quadric representation of a  submanifold and
spectral geometry, Proc. Symp. Pure Math. , 54-3 (1993), 155--167.

\bibitem{}  Dimitri\'c, I.,  $CR$ submanifolds of $HP^m$ and hypersurfaces of the Cayley
plane whose Chen-type is 1, preprint, 1998.

\bibitem{} Dombrowski, P., Jacobi fields, totally geodesic foliations and
geodesic differential forms, Results Math., 1 (1978), 156--194.

\bibitem{} Dombrowski, P.,  150 years after Gauss'
``Disquisitiones generales circa superficies curvas'', Ast\'erisque, 62 (1979).

\bibitem{} Dombrowski, P.,  Differentialgeometrie. Ein Jahrhundert
Mathematik 1890-1990, Dokum. Gesch.
Math. Vieweg, Braunschweig, 6, 323--360, 1990.

\bibitem{} Dorfmeister, J. and G. Haak, Meromorphic potentials
and smooth surfaces of constant mean
curvature,  Math. Z., 224 (1997), 603--640. 

\bibitem{} Dorfmeister, J. and E. Neher, An algebraic approach
to isoparametric hypersurfaces in spheres. I, II, T\^ohoku Math. J.  35 (1983),
187--224;  , ibid , 35 (1983), 225--247.,

\bibitem{} Dorfmeister, J. and E. Neher,  Isoparametric hypersurfaces, case
$g=6,\;m=1$,  Comm. Algebra , 13 (1985), 2299--2368. 

\bibitem{} Dorfmeister, J., Pedit F. and H. Wu, Weierstrass type
representations of harmonic maps into
symmetric spaces, Comm. Anal. Geom., 6 (1998), 633--668.

\bibitem{} Douady, A. and R. Douady, Changements de
cadres \`a des surfaces minimales, preprint, 1995.

\bibitem{}  Douglas, J., Solution of the problem of Plateau,
Trans. Amer. Math. Soc., 33 (1931), 263--321.

\bibitem{}  Douglas, J.,  The problem of Plateau for two contours, J.
Math. Phys., 10 (1931), 35--359.

\bibitem{}  Dragomir, S., On pseudo-hermitian
immersions between strictly pseudoconvex
$CR$ manifolds, Amer. J. Math., 117 (1995), 169--202.

\bibitem{} Dragomir, S. and L. Ornea, Locally conformal K\"ahler
geometry, Birkh\"auser, Boston-Basel-Berlin, 1998.

\bibitem{}  Dupin, C., Applications de g\'eom\'etrie et de
m\'echanique, Paris,  1822.

\bibitem{}  Eells, J. and N. H. Kuiper,  Manifolds which are
like projective planes, Inst. Hautes Etudes Sc. Publ. Math., 14 (1962), 5--46.

\bibitem{}  Eells, J. and L. Lemaire, A report on harmonic maps, Bull. London Math.
Soc., 10 (1978),1--68.

\bibitem{}  Eells, J. and L. Lemaire,  Another report on
harmonic maps, Bull. London Math. Soc., 20 (1988), 385--524.

\bibitem{} Efimov, N. V., The appearance of singularities on a
surface of negative curvature
(Russian), Mat. Sb., 64 (1964), 286--320.

\bibitem{}  Ejiri, N.,  Compact minimal
submanifolds of a sphere with positive Ricci
curvature, J. Math. Soc. Japan, 31 (1979), 251--256.

\bibitem{}  Ejiri, N.,  Minimal immersions of
Riemannian products into real space forms,
Tokyo J. Math., 2 (1979), 63--70.

\bibitem{}  Ejiri, N., Totally real submanifolds in a 6-sphere,
Proc. Amer. Math. Soc., 83 (1981), 759--763.

\bibitem{}  Ejiri, N.,  Totally real minimal immersions of
$n$-dimensional real space forms into
$n$-dimensional complex space forms, Proc. Amer. Math. Soc., 84 (1982),
243--246.

\bibitem{}  Ejiri, N., The index of minimal immersions of $S^{2}$ into
$S^{2n}$,  Math. Z., 184 (198), 127--132. 

\bibitem{}  Ejiri, N., Equivariant minimal immersions of $S^2$
into $S^{2m}$, Trans. Amer. Math. Soc., 297 (1986), 105--124.

\bibitem{}  Ejiri, N.,  Calabi lifting and
surface geometry in $S^4$, Tokyo J. Math., 9 (1986), 297--324.

\bibitem{}  Ejiri, N.,  Two applications of 
the unit normal bundle of a minimal surface
in $R^N$, Pacific J. Math., 147 (1991), 291--300.

\bibitem{} Ejiri, N. and M. Kotani, Index and flat
ends of minimal surfaces, Tokyo J. Math., 16 (1993), 37--48.

\bibitem{} Enneper, A.,
Die cyklischen Fl\"achen, Z. Math. u.
Phys., 10 (1869), 315--359.

\bibitem{} Enomoto, K.,   Umbilical points on surfaces in
$R^N$, Nagoya Math. J., 100 (1985), 135--143.

\bibitem{} Erbacher, J.,  Isometric immersions of constant
mean curvature and triviality of the normal
connection, Nagoya Math. J., 45 (1971), 139--165.

\bibitem{}  Eschenburg, J. H. and V.
Schroeder, Tits distance
of Hadamard manifolds and isoparametric
hypersurfaces, Geom. Dedicata, 40 (1991), 97-101.

\bibitem{}  Escher, C. M., Minimal isometric immersions of inhomogeneous
spherical space forms into spheres---a
necessary condition for existence, Trans.
Amer. Math. Soc., 348 (1996), 3713--3732.

\bibitem{}  Do Espirito-Santo, N.,  Complete minimal
surfaces in $R^3$ with type Enneper end,
Ann. Inst. Fourier (Grenoble), 44 (1994), 525--557. 

\bibitem{} Euler, L., Methodus inviendi lineas curvas maximi
minimive proprietate gaudentes, sive,
solutio problematis isoperimetrici
latissimo sensu accepti, Lausann \&
Genevae, Bousquet, 1744.
\bibitem{} Euler, L., Euler \`a Lagrange,
Berolini, die 6 Sept. 1755, Oeuvres de
Lagrange, XIV, 144--146.

\bibitem{}  Fang, F., Topology of
Dupin hypersurfaces with six principal 
curvatures, preprint, 1995, Bielefeld.

\bibitem{}  Fang, F.,  Multiplicities of principal curvatures of isoparametric
 hypersurfaces, preprint, Max Planck Institut f\"{u}r
Mathematik.

\bibitem{}  Fang, F., On the topology of isoparametric hypersurfaces
with four distinct principal curvatures,  Proc.
Amer. Math. Soc.,  127 (1996), 259--264.

\bibitem{} Fang, Y., On the Gauss map of complete minimal
 surfaces with finite total curvature,
Indiana Univ. Math. J., 42 (1943), 1389--1411. 

\bibitem{} Fang, Y., On minimal annuli in
a slab, Comm. Math. Helv. 69 (1994), 417--430.

\bibitem{} Fang, Y.,  Total curvature of
branched minimal  surfaces, Proc. Amer. Math. Soc., 124 (1996), 1895--1898.

\bibitem{} Fang, Y., Lectures on minimal surfaces in $R^3$, Australian Nat.
Univ, Australia, 1996.

\bibitem{} Fang, Y. and W. H. Meeks, III, Some
global properties of complete minimal
surfaces of finite topology in $R^3$, Topology, 30 (1991), 9--20. 

\bibitem{}  Fang, Y. and F. Wei, On uniqueness of 
Riemann's examples, Proc. Amer. Math. Soc., 126 (1998), 1531--1539.

\bibitem{} Fary, I., Sur la courbure totale d'une
courve gauche faisant un noeud,
Bull. Soc. Math. France, 77 (1949), 128--138.

\bibitem{} Fenchel, W.,   \"Uber die Kr\"ummung und Wind\"ung
geschlossener Raumkurven, Math. Ann., 101 (1929), 238--252.

\bibitem{}  Ferapontov, E. V.,  Dupin hypersurfaces and integrable hamiltonian
systems of hydrodynamic type, which do not possess Riemann invariants,  Diff.
Geom. and its Appl. 5 (1995), 121-152.

\bibitem{}  Ferapontov, E. V.,   Isoparametric hypersurfaces in spheres, integrable
nondiagonalizable systems of hydrodynamic
type, and $N$-wave systems,  Diff. Geom. and its Appl., 5 (1995), 335-369.

\bibitem{} Ferus, D.,  \"Uber die absolute Totalkru\"ummung
h\"oher-dimensionaler Konten, Math. Ann., 171 (1967), 81--86.

\bibitem{} Ferus, D.,   On the type number
of hypersurfaces in spaces of constant
curvature, Math. Ann., 187 (1970), 310--316.

\bibitem{} Ferus, D.,  On the completeness of
nullity foliations, Michigan Math. J., 18 (1971), 61--64.

\bibitem{} Ferus, D.,   The torsion form of
submanifolds in $E^n$, Math. Ann., 193 (1971), 114--120.

\bibitem{} Ferus, D.,  On isometric
immersions between hyperbolic spaces, Math. Ann., 205 (1973), 193--200.

\bibitem{} Ferus, D., Immersions with parallel second fundamental 
form, Math. Z., 140 (1974), 87--93.

\bibitem{} Ferus, D.,  Isometric immersions of
constant curvature manifolds,  Math. Ann., 217 (1975), 155--156.

\bibitem{} Ferus, D.,  Symmetric submanifolds of Euclidean space,  Math.
Ann., 247 (1980), 81--93.

\bibitem{} Ferus, D., Karcher, H. and H. F.
M\"unzner, Cliffordalgebren und neue
isopara\-metrische Hyperfl\"achen , Math. Z., 177 ({1981}), 479--502.

\bibitem{} Ferus, D. and F. Pedit,  Isometric
immersions of space forms and soliton theory,  Math. Ann., 305 (1996), 329--342. 

\bibitem{} Ferus, D. and S. Schirrmacher, Submanifolds in
Euclidean space with simple geodesics, Math. Ann., 260 (1982), 57--62.

\bibitem{} Fialkow, A.,  Hypersurfaces
of spaces of constant curvature, Ann. Math. , 39 (1938), 762--785.

\bibitem{} Fischer, W. and E. Koch, Spanning minimal
surfaces, Philos. Trans. Roy. Soc.
London Ser. A, 354 (1996), 2105--2142.

\bibitem{} Fischer-Colbrie, D., On complete minimal
surface with finite Morse index in
three-manifolds, Invent. Math., 82 (1985), 121--132.

\bibitem{} Fischer-Colbrie, D. and R. Schoen,  The structure of
complete stable minimal surfaces in 3-manifolds of nonnegative scalar
curvature, Comm. Pure Appl. Math., 33 (1980), 199--211.

\bibitem{}  Fladt, K. and A.Baur, Analytische Geometrie
spezieller Fl\"achen und Raumkurven,
Friedr. Vieweg, Braunsch\-wieg, 1975.

\bibitem{} Fomenko, A. T.,  Minimal compacta in Riemannian
manifolds and Reifenberg's conjecture,
Math. USSR-Izv., 6 (1972),1037--1066.

\bibitem{} Foulon, P. and F.
Labourie,{1992} Sur les varietes
compactes asymptotiquement harmoniques,
Invent. Math., 109, 97--111.

\bibitem{} Frankel, T. J.  Manifolds with positive
curvature, Pacific J. Math., 11 (1961), 165--171.

\bibitem{} Frankel, T. J.  On the fundamental
group of a compact minimal submanifold,
Ann. Math., 83 (1966), 68--73.

\bibitem{} Frohman, C.,  The topological uniqueness of triply periodic minimal surfaces in
$R^3$, J. Differential Geometry, 31 (1990), 277--283. 

\bibitem{}  Fujimoto, H., On the number of exceptional values of the
Gauss maps of minimal surfaces, J. Math.
Soc. Japan, 40 (1988), 235--247.

\bibitem{}  Fujimoto, H.,  Modified defect
relations for the Gauss map of minimal
surface, II, J. Differential Geometry, 31 (1990), 365--385.

\bibitem{}  Fujimoto, H.,  Value distribution
theory of the Gauss map of minimal surfaces
in $R^m$, Aspects of Mathematics,
Friedr. Vieweg \& Sohn, Braunschweig, 1993.

\bibitem{} Furuhata, H., Construction and classification of 
isometric minimal immersions of Kaehler
manifolds into Euclidean spaces, Bull. London Math. Soc., 26 (1994), 487--496. 

\bibitem{} Fwu, C. C., Total absolute curvature of
submanifolds in compact  symmetric spaces of rank one, Math. Z., 172 (1980) , 245--254.

\bibitem{} Garay, O. J., Spherical Chen surfaces which are mass-symmetric and of 2-type, J.
Geometry, 33 (1988), 39--52.

\bibitem{} Garay, O. J., Finite type cones shaped on spherical submanifolds, Proc.
Amer. Math. Soc., 104 (1988), 868--870.

\bibitem{} Garay, O. J., On a certain class of finite type surfaces of
revolution, Kodai Math. J. , 11 (1988), 25--31.

\bibitem{} Garay, O. J., An extension of
Takahashi's theorem , Geom. Dedicata, 34 (1990), 105--112.

\bibitem{} Garay, O. J.,Orthogonal surfaces with constant mean
curvature in the Euclidean $4$-space,
Ann. Glo\-bal Anal. Geom. 12 (1994), 79--86.

\bibitem{} Gauchman, H., Minimal submanifolds of a sphere with
bounded second fundamental form, Trans.
Amer. Math. Soc., 298 (1986), 778--791.

\bibitem{} Gauss, C. F., Disquisitiones generales circa
superficies curvas, Comment. Soc. Sci. Gotting. Recent. Classis Math., 6, 1827.

\bibitem{} Gheysens, L., Verheyen, P. and L. Verstraelen, Sur les
surfaces $\Cal A$ ou les surfaces de Chen, C. R. Acad. Sc. Paris Ser. I.
Math.  292 (1981), 913--916.

\bibitem{} de Giorgi, E.,  Una estensione del teorema de
Bernstein, Ann. Scoula Norm. Sup. Piza, Sci. Fis. Mat., 19 (1965), 79--85.

\bibitem{} Goldberg, S. I. and S. Kobayashi, Holomorphic
bisectional curvature, J. Differential Geometry, 1 (1967), 225--233.

\bibitem{} Goldberg, S. I. and E. M. Moskal,  The axioms of
spheres in Kaehler geometry, Kodai Math. Sem. Rep., 27 (1976), 188--192

\bibitem{} Gotoh, T.,  Compact minimal $CR$-submanifolds with the
least nullity in a complex projective space, Osaka J. Math., 34 (1997), 175--197
.

\bibitem{}  Gray, A., Almost complex submanifolds of the six
sphere, Proc. Amer. Math. Soc., 20 (1969), 277--279.

\bibitem{} Greene, R. E.,  Isometric embeddings of Riemannian
and pseudo-Riemannian manifolds, Mem. Amer. Math. Soc., 97, 1970.

\bibitem{} Greenfield, S., Cauchy-Riemann equations in several
variables, Ann. Scoula Norm. Sup. Pisa, 22 (1968), 275--314.

\bibitem{} Greub, W. and D. Socolescu, The local isometric
imbedding in $R^3$ of two-dimensional
Riemannian manifolds with Gaussian curvature changing sign arbitrarily,
Libertas Math., 14 (1994), 27--39.

\bibitem{} Grove, K. and S. Halperin,  Dupin
hypersurfaces, group actions and the
double mapping cylinder, J. Differential Geometry, 26 (1987), 429--459.

\bibitem{} Gromov, M. L.,   A topological technique for
construction of solutions of differential equations and inequalities, Proc.
Intern. Congr. Math., 2 (1971), 221--225.

\bibitem{} Gromov, M. L., Pseudoholomorphic curves in symplectic
 manifolds, Invent. Math. , 82 (1985), 307--347.

\bibitem{} Gromov, M. L.,  Partial differential relations, 
Springer-Verlag, Berlin-New York, 1986.

\bibitem{} Gromov, M. L. and V. A. Rokhlin, Embeddingsand immersions in
Riemannian geometry, Russian Math. Surveys, 25 (1970), no. 5, 1--57.

\bibitem{} Gr\"uter, M. and J. Jost, On embedded minimal
disks in convex bodies, Ann. Inst. H. Poincar\'e Anal. Non Lin\'eaire, 3 (1986), 345--390. 

\bibitem{} Guadalupe, I. V. and L. Rodriguez,  Normal curvature of
surfaces in space forms, Pacific J. Math., 106 (1983), 95--103.

\bibitem{} Gulliver, R., Existence of surfaces with prescribed mean
curvature vector, Math. Z., 131 (1973), 117--140.

\bibitem{} Gulliver, R., Necessary conditions for submanifolds and currents with
prescribed mean curvature vector, Seminar on Minimal Submanifolds, Ann. of
Math. Stud., 103, Princeton Univ. Press, Princeton, N.J.,  225--242, 1983.

\bibitem{} Guilliver, R., Osserman, R. and H. L. Roydon, A theory
of branched immersions of surfaces, Amer. J. Math., 95 (1973), 750--812.


\bibitem{} Gulliver, R. D. and F. D. Lesley, On boundary
branch points if minimizing surface, Arch. Rat. Mech. Anal., 2 (1973), 20--25. 

\bibitem{} Haab, F., Immersions tendues de surfaces dans 
$E^3$, Comment. Math. Helv. 67 (1992), 182--202. 

\bibitem{} Hall, P.,   A Picard theorem with an application to
minimal surfaces, Trans. Amer. Math.
Soc. , 314 (1989), 597--603.

\bibitem{} Hall, P.,  A Picard theorem with an
application to minimal surfaces. II,
Trans. Amer. Math. Soc., 325 (1991), 895--902.

\bibitem{} Hadamard, J., Sur les \'el\'ements lin\'eaires
\`a plusieurs dimensions, Bull. Sci. Math., 25 (1901), 57--60.

\bibitem{} Hano, J., Einstein complete intersections in complex
projective space, Math. Ann.,
216 (1975), 197--208.

\bibitem{} Harada, M., On Kaehler manifolds satisfying the axiom of
anti-holomorphic 2-spheres, Proc. Amer. Math. Soc., 43 (1974), 186--189.

\bibitem{} Hardt, R. and L. Simon, Boundary regularity and embedded solutions for the
oriented Plateau's problem, Ann. Math., 110 (1979), 439--486.

\bibitem{} Harle, C. E., Rigidity of hypersurfaces of constant
scalar curvature, J. Differential Geometry, 5 (1971), 85--111.
 
\bibitem{} Hartman, P. and L. Nirenberg, On spherical image
maps whose Jacobians do not change sign,
Amer. J. Math., 81 (1959), 901--920.

\bibitem{} Harvey, R. and  H. B. Lawson, Calibrated geometries, Acta
Math., 148 (1982), 47--157.

\bibitem{} Hasanis, T. and T. Vlachos,  Spherical  2-type
hypersurfaces, J. Geometry, 40 (1991), 82--94.

\bibitem{} Hasanis, T. and T. Vlachos,  Surfaces of finite
type with constant mean curvature, Kodai
Math. J., 16 (1993), 244--252.

\bibitem{} Hasanis, T. and T. Vlachos,  Hypersurfaces with
constant scalar curvature and constant mean
curvature, Ann. Global Anal. Geom.  13 (1995), 69--77.

\bibitem{} Hasanis, T. and T. Vlachos,  Hypersurfaces in $E^4$
with harmonic mean curvature vector field, Math. Nachr.  172 (1995), 145--169.

\bibitem{} Hashimoto, H. and K. Sekigawa, Minimal surfaces in a
$4$-dimensional sphere,  Houston J.
Math., 21 (1995), 449--464. 

\bibitem{} Hass, J., Pitts, J. T. and J. H. Rubinstein, Existence of
unstable minimal surfaces in manifolds with homology and applications to triply periodic
minimal surfaces, Differential geometry: partial differential equations on manifolds
(Los Angeles, CA, 1990) Proc. Sympos. Pure Math., 54 (1993), I, 147--162.

\bibitem{} Hebda, J.,  Some new tight embeddings which cannot be  made
taut, Geom. Dedicata, 17 (1984), 49--60.

\bibitem{} Hebda, J.,The possible cohomology of certain types
 of taut submanifolds, Nagoya Math. J. , 111 (1988), 85--97. 

\bibitem{} Heinz, E., On the nonexistence of a surface of
constant mean curvature with finite area
and prescribed rectifiable boundary,
Arch. Rat. Mech. Anal., 35 (1969), 249--252.
 
\bibitem{} Helgason, S., Totally geodesic spheres in compact
symmetric spaces, Math. Ann., 165 (1966), 309--317.

\bibitem{} Helgason, S., Differential geometry,
Lie groups and symmetric spaces, Academic Press, New York, 1978.

\bibitem{}  Henke, W.,  \"Uber die Existenz isometrischer
Immersionen der Kodimension 2 von sph\"arischen Rau\-mformen
in Standard-R\"aume , Math. Ann., 222 (1976), 89--95.

\bibitem{} Hertrich-Jeromin, U. and U. Pinkall, Ein Beweis der
Willmoreschen Vermutung fur Kanaltori.,  J. Reine Angew. Math., 430 (1992), 21--34.

\bibitem{}  Hilbert, D.,  \"Uber Fl\"achen von konstanter
Gausscher Kr\"ummung, Trans. Amer. Math. Soc.  2 (1901), 87--99.

\bibitem{} Hildebrandt, S., Boundary behavior of minimal surfaces,
Arch. Rat. Mech. Anal., 35 (1969), 47--82.

\bibitem{} Hoffman, D.,  Surfaces of constant mean curvature
in constant curvature manifolds, J.
Differential Geometry, 8 (1973), 161--176

\bibitem{} Hoffman, D., New examples of singly-periodic minimal surfaces and their
qualitative behavior, Contemp. Math., 97--10, 101, 1989.

\bibitem{} Hoffman, D., Karcher, H. and H. Rosenberg, Embedded
minimal annuli in $R^3$ bounded by a pair of straight lines, Comment.
Math. Helv., 66 (1991), 599--617.

\bibitem{} Hoffman, D. and W. H.  Meeks, III, A complete embedded
minimal surface in $E^3$ with genus one and
three ends, J. Differential
Geometry, 21 (1985), 109--127.

\bibitem{} Hoffman, D. and W. H.  Meeks, III,  A variational approach
to the existence of complete embedded
minimal surfaces, Duke Math. J., 57 (1988), 877--893.

\bibitem{} Hoffman, D. and W. H.  Meeks, III,  The asymptotic
behavior of properly embedded minimal
surfaces of finite topology, J. Amer.
Math. Soc., 2 (1989), 667--682.

\bibitem{} Hoffman, D. and W. H.  Meeks, III,   Embedded minimal
surfaces of finite topology, Ann. of
Math., 131 (1990), 1--34.

\bibitem{} Hoffman, D. and W. H.  Meeks, III,  The strong halfspace
theorem for minimal surfaces, Invent.
Math., 101 (1990), 373--377.

\bibitem{}  Hoffman, D. and R. Osserman, The geometry of
the generalized Gauss map, Mem. Amer. Math. Soc.   236, 1980.

\bibitem{} Hoffman, D. and W. H.  Meeks, III,  The area of the
generalized Gaussian image and the stability of minimal surfaces in $S^n$
and $R^n$, Math. Ann., 260 (1982), 437--452.

\bibitem{} Hoffman, D., F.  Wei and H. Karcher, Adding
handles to the helicoid,  Bull. Amer. Math. Soc. (N.S.), 29 (1993), 77--84. 

\bibitem{} Holmgren, E., Sur les surfaces \`a courbure constant
n\'egative, C. R. Acad. Sci. Paris, 134 (1902), 740--743.

\bibitem{} Hopf, H., \"Uber Fl\"achen mit einer Relation zwischen
den Hauptkr\"ummungen, Math. Nachr., 4 (1951), 232--249.

\bibitem{}  Hong, J., Isometric embedding in $R^ 3$
complete noncompact nonnegatively curved
surfaces,. Ma\-nusc. Math., 94 (1997), 271--286.

\bibitem{} Hong, J. and C. Zuily, Isometric embedding of
the 2-sphere with nonnegative curvature in
$R^3$,  Math. Z., 219 (1995), 323--334. 

\bibitem{}  Hong, S. L., Isometric immersions of manifolds with planar
geodesics into Euclidean space, J. Differential Geometry, 7 (1973), 259--278.

\bibitem{} Hong, Y., Helical immersions in Euclidean spaces,
Indiana Univ. Math. J., 35 (1986), 29-43 .

\bibitem{} Hong, Y. and C. S. Houh, Helical immersions and normal
sections, Kodai Math. J., 8 (1985), 171-192 . 

\bibitem{}  Hou, Z. H.,  Hypersurfaces in a sphere with constant mean curvature, Proc. Amer. Math. Soc., 125 (1997), 1193-1196.

\bibitem{} Houh, C. S.,  On totally real bisectional curvature, Proc. Amer.
Math. Soc., 56 (1976), 261--263.

\bibitem{} Houh, C. S., Lagrangian submanifolds of quaternionic Kaehlerian manifolds satisfying Chen's equality,  Beitr\"age zur
Alg. Geom., 39 (1998), 413--421.

\bibitem{} Howard, R., The kinematic formula in Riemannian 
homogeneous spaces, Mem. Amer. Math. Soc., 106 (1993), no.509.

\bibitem{} Hsiang, W. T., Hsiang, W. Y. and I. Sterling,
On the construction of codimension two
minimal immersions of exotic spheres into
Euclidean spheres, Invent. Math., 82 (1985), 447--460. 

\bibitem{} Hsiang, W. Y.,  Remarks on closed
minimal submanifolds in the standard
Riemannian $m$-sphere, J. Differential
Geometry, 1 (1967), 257--267.

\bibitem{} Hsiang, W. Y.,  Generalized 
rotational hypersurfaces of constant mean
curvature in  Euclidean spaces, I, J. Differential Geom., 17 (1982), 337--356.

\bibitem{} Hsiang, W. Y.,  New examples of minimal imbeddings of
$S^{n-1}$ into $S^{n}(1)$--the spherical Bernstein problem for $n=4$,
$5$, $6$,  Bull. Amer. Math. Soc. (N.S.), 7 (1982), 377--379.

\bibitem{} Hsiang, W. Y., On generalization of theorems of A. D.
 Alexandrov and C. Delaunay on hypersurfaces
of constant mean curvature,  Duke Math. J. , 49 (199=82), 485--496.

\bibitem{} Hsiang, W. Y. and H. B. Lawson,  Minimal
submanifolds of low cohomogeneity, J.
Differential Geometry, 5 (1971), 1--38.

\bibitem{} Hsiang, W. Y.,  Teng, Z. H. and W. C. Yu, New
examples of constant mean curvature
immersions of $(2k-1)$-spheres into
Euclidean $2k$-space,  Ann. of Math., 117 (1983), 609--625.

\bibitem{} Hsiang, W. Y. and W. C. Yu,  A generalization of a
theorem of Delaunay,  J. Differential
Geometry, 16 (1981), 161--177. 

\bibitem{} Hsiung, C. C., Some integral formulas for
hypersurfaces, Math. Scand., 2 (1954), 286--294. 

\bibitem{} Hu, Z. J.  and G. S. Zhao, Isometric immersions from the hyperbolic space
$H^2(-1)$ into $H^3(-1)$ , Proc. Amer. Math. Soc., 125 (1977), 2693--2697. 

\bibitem{} Hu, Z. J.  and G. S. Zhao,  Classification of isometric  immersions of the hyperbolic space  $H^2$ into $H^3$, Geom. Dedicata, 65 (1977), 47--57.

\bibitem{} Huber, A., On subharmonic functions and differential
geometry in the large, Comment. Math. Helv., 32 (1957), 13--72.

\bibitem{} Iseri, H., On the existence of minimal surfaces with
 singular boundaries, Proc. Amer. Math. Soc.  124 (1996), 3493--3500.

\bibitem{} Itoh, T., Addendum to: "On Veronese manifolds", J.
Math. Soc. Japan, 30 (1978), 73--74.

\bibitem{} Jacobowitz, H. and J. D. Moore,  The Cartan-Janet
theorem for conformal embeddings,
Indiana Univ. Math. J., 23 (1973), 187--203.

\bibitem{}  Janet, M.,  Sur la possibilit\'e de plonger un
espace riemannien donn\'e dans un espace
euclidien, Ann. Soc. Math. Polon., 5 (1926), 38--43.

\bibitem{} Ji, M. and G. Y. Wang,  Minimal surfaces in
Riemannian manifolds, Mem. Amer. Math.
Soc., 104 (1993), no. 495. 

\bibitem{}  Jorge, L. P. and W. Meeks III, The topology of complete minimal surfaces
of finite total Gaussian curvature, Topology, 22 (1983), 203--221.

\bibitem{} Jorge, L. P. and F. Mercuri,  Minimal immersions
into space forms with two principal
curvatures,  Math. Z., 187 (1984), 325--333. 

\bibitem{} Jorge, L. P. and F. Xavier, On the existence of
complete bounded minimal surfaces in
$R^n$, Bol. Soc. Brasil. Mat., 10 (1979), 171--173.,

\bibitem{} Jorge, L. P. and F. Xavier, A complete minimal surface in $R^3$ between
two parallel planes, Ann.  Math., 112 (1980), 203--206.

\bibitem{} Jost, J., Conformal mappings and the Plateau-Douglas
 problem in Riemannian manifolds, 
J. Reine Angew. Math., 359 (1985), 37--54.

\bibitem{} Jost, J. and M. Struwe, Morse-Conley theory
for minimal surfaces of varying topological
type, Invent. Math., 102 (1990), 465--499.

\bibitem{}  Kagan, V. F.,  The fundamentals of the theory of surfaces
in tensor presentation. Part two. Surfaces in space, Transformations and
Deformations of Surfaces. Special Questions   (Russian). Gosudarstv. Izdat.
Tehn.-Teor. Lit., Moscow-Leningrad, 1948.

\bibitem{} K\"ahler, E., \"Uber eine bemerkenswerte Hermitiesche
Metrik, Abh. Math. Sem. Univ. Hamburg, 9 (1933), 173--186.

\bibitem{} Kapouleas, N.,  Complete constant mean curvature
surfaces in Euclidean three-space,
Ann. Math., 131 (1990), 239--330.

\bibitem{} Kapouleas, N.,  Compact constant mean curvature
surfaces in Euclidean three-space,
J. Differential Geometry, 33 (1991), 683--715.

\bibitem{} Karcher, H., The triply periodic minimal surfaces of Alan
Schoen and their constant mean
curvature companions,  Manusc. Math., 64 (1989), 291--357. 

\bibitem{} Karcher, H.,  Construction of minimal surfaces, Surveys in
Geometry, Univ. of Tokyo, 1--96, 1989. 

\bibitem{} Karcher, H., Pinkall, U. and I. Sterling, New minimal
surfaces in $S^3$, J. Differential
Geometry, 28 (1988), 169--185.

\bibitem{} Karcher, H. and K. Polthier, Construction of
triply periodic minimal  surfaces,  Philos. Trans. Roy. Soc. London Ser. A, 354 (1996), no. 1715, 2077--2104. 

\bibitem{} Kassabov, O.,  On the axiom of spheres in Kaehler
geometry, C. R. Acad. Bulgare Sc., 35 (1982), 303--306.

\bibitem{} Kenmotsu, K., On minimal immersions
of $R^{2}$ into $S^{N}$, J. Math. Soc.
Japan, 28 (1976), 182--191.

\bibitem{} Kenmotsu, K.,  Weierstrass formula
for surfaces of prescribed mean curvature 
, Math\. Ann. , 245 (1979), 89--99

\bibitem{} Kenmotsu, K.,  On minimal immersions of $R^2$ into
$P^n(C)$ , J. Math. Soc. Japan, 37 (1985), 665--682.

\bibitem{} Kenmotsu, K. and C. Xia,  Hadamard-Frankel type theorems
for manifolds with partially positive
curvature, Pacific J. Math., 176 (1996), 129--139. 

\bibitem{} Ki, U. H. and Y. H. Kim, Totally real
submanifolds of a complex space form,
Internt. J. Math. \& Math. Sci. , 19 (1996), 39--44.

\bibitem{} Kim Y. H. and E. K. Lee,  Surfaces of
Euclidean 4-space whose geodesics are
$W$-curves, Nihonkai Math. J., 4 (1993), 221--232.

\bibitem{} Kitagawa, Y., Embedded flat tori in the unit $3$-sphere, J. Math. Soc. Japan, 47 (1995), 275--296. 

\bibitem{} Klein, F., Vorlesungen \"uber h\"ohere Geometrie,
third ed. Spinger-Verlag, Berlin, 1926.

\bibitem{} Klotz, T. and R. Osserman,  Complete surfaces
in $E\sp{3}$ with constant mean curvature, Comment. Math. Helv. , 41 (1966), 313--318.

\bibitem{} Klotz, T. and L. Sario,  Existence of complete minimal surfaces of arbitrary
connectivity and genus, Proc. Nat. Acad. Sci. USA , 54 (1965), 42--44.

\bibitem{} Kobayashi, S., Compact homogeneous hypersurfaces, Trans.
Amer. Math. Soc., 88 (1958), 137--143.

\bibitem{} Kobayashi, S., Hypersurfaces of
complex projective space with constant scalar
curvature, J. Differential Geometry, 1 (1967), 369--370.

\bibitem{} Kobayashi, S.,  Imbeddings of
homogeneous spaces with minimum total
curvature, T\^ohoku Math. J., 19 (1967), 63--74.

\bibitem{} Kobayashi, S., Isometric imbeddings of compact symmetric
spaces  , T\^ohoku Math. J. , 20 (1968), 21--25.

\bibitem{} Kobayashi, S.,  Transformation groups in
differential geometry, Springer Verlag, Berlin, 1972.

\bibitem{} Kobayashi, S.,  Differential geometry of
complex  vector bundles, Princeton University Press, Princeton, 1987.

\bibitem{} Kobayashi, S. and K. Nomizu, Foundations in
Differential Geometry, I, 1963, Wiley-Interscience, New York; 

\bibitem{} Kobayashi, S. and K. Nomizu,  Foundations in
Differential Geometry, II, 1969
Wiley-Interscience, New York.

\bibitem{} Kokubu, M., Weierstrass representation for minimal surfaces
in hyperbolic space, T\^ohoku Math. J., 49 (1997), 367--377.

\bibitem{} Kon, M., On some complex submanifolds in Kaehler
manifolds, Canad. J. Math. , 26 (1974),
1442--1449. 

\bibitem{} Kon, M.,  Minimal CR submanifolds immersed in a complex
 projective space, Geom. Dedicata, 31 (1989), 357--368.

\bibitem{} Korevaar, N. J.,  Sphere theorems via Alenxandrov
for constant Weingarten curvature
hypersurfaces, J. Differential Geometry, 30 (1988), 221-223.

\bibitem{} Korevaar, N. J., Kusner, R.,
Meeks, W. H. III and B. Solomon, Constant mean curvature
surfaces in hyperbolic space, Amer. J. Math., 114 (1992), 1--43.

\bibitem{} Korevaar, N. J., Kusner, R. and B. Solomon,  The structure
of complete embedded surfaces with constant
mean curvature, J. Differential Geometry, 30 (1989), 465--503.

\bibitem{} Kotani, M., The first eigenvalue of homogeneous minimal
hypersurfaces in a unit sphere
$S^{n+1}(1)$, T\^ohoku Math. J., 37 (1985), 523--532.

\bibitem{} Kotani, M.,   An immersion of an
$n$-dimensional real space form into an
$n$-dimensional complex space form, Tokyo J. Math., 9 (1986), 103--113

\bibitem{} Kotani, M.,   A decomposition
theorem of 2-type immersions, Nagoya Math. J., 118 (1990), 55--66.

\bibitem{} Kowalski, O., Properties of hypersurfaces which are characteristic for
spaces of constant curvature, Ann. Scoula Norm. Sup. Pisa, 26 (1972), 233-245.

\bibitem{} Kronheimer, P. B.and T. S. Mrowka, The genus of
embedded surfaces in the projective plane, Math. Res. Lett., 1 (1994), 797--808. 

\bibitem{} K\"uhnel, W. and U. Pinkall, On total mean
curvature, Quart. J. Math. Oxford, 37 (1986), 437--447.

\bibitem{} Kuiper, N. H., On conformally flat
spaces in the large, Ann. Math., 50 (1949), 916--924.

\bibitem{} Kuiper, N. H.,  On $C^1$-isometric
imbeddings, Indag. Math. 17 (1955), 683--689.

\bibitem{} Kuiper, N. H.,  Immersions with minimal total absolute curvature,
Coll. de g\'eom. diff., Centre Belge de Recherches Math., Bruxelles, 1958, 75--88.

\bibitem{} Kuiper, N. H.,  Convex immersions of closed surfaces in $E^3$,
Comm. Math. Helv., 35 (1961), 85--92.

\bibitem{} Kuiper, N. H.,  On convex maps,
Nieuw Archief voor Wisk., 10 (1962), 147--164.

\bibitem{} Kuiper, N. H.,  Tight embeddings and maps.
Submanifolds of  geometrical class three in
$E^{N}$, The Chern Symposium 1979, Springer-Verlag, 97--145.

\bibitem{} Kuiper, N. H. and W. Meeks III, Total curvature for knotted surfaces,
Invent. Math., 77 (1984), 25--69. 

\bibitem{} Kummer, E. E.,  \"Uber die Fl\"achen vierten
Grades auf welchen Schaaren von Kegelschnitten liegen, J. Reine
Angew. Math., 64 (1865), 66--76.

\bibitem{} Kusner, R., Conformal geometry and complete minimal surfaces, Bull. Amer. Math.
Soc. (N.S.), 17 (1987), 291--295.

\bibitem{} Lagrange, J. L.,  Essai d'une nouvelle m\'ethode pour d\'eterminer les maxima et les minima des formules int\'egrales ind\'efinies, Miscellanea Taurinensia, 2 (1760/1),
173--195.

\bibitem{} Lagrange, J. L.,  Oeuvres
I-XIV, (ed. J. A. Serret), Gauthier-Villars, (1867--1892).

\bibitem{} Lancaster, G. M., Canonical metrics for
certain conformally Euclidean spaces of
dimension three and codimension one, Duke Math. J. , 40 (1973), 1--8.

\bibitem{} Lancret, M. A., M\'emoire sur les courbes \`a double
courbure , M\'em. des Sav. E\'trangers, 1 (1806), 416--454.

\bibitem{} Langer, J. and D. A. Singer,  Curves in the hyperbolic
plane and mean curvature of tori in 3-space,  Bull. London Math. Soc. 16 (1984), 531--534.

\bibitem{} Langevin, R. and H. Rosenberg, On curvature integrals and
knots, Topology, 15 (1976), 405--416.,

\bibitem{} Langevin, R. and H. Rosenberg, A maximum principle at infinity for minimal
surfaces and applications, Duke Math. J., 57 (1988), 819--828.

\bibitem{} Laura, E., Sopra la propagazione di onde in un mezzo
indefinito, in Scritti mat. offerti ad E. D'Ovidio, 253--278, 1918.

\bibitem{} Lawson, H. B.,  Complete minimal
surfaces in $S^3$, Ann. Math., 92 (1970), 335--374.

\bibitem{} Lawson, H. B.,  Some intrinsic
characterizations of minimal surfaces, J.
D'Analyse Math., 24 (1971), 151--161.

\bibitem{} Lawson, H. B.,   Minimal varieties in
real and complex geometry, S\'em. Math. Sup.,
No. 57, l'Univ. de Montr\'eal, Canada, 1973.

\bibitem{} Lawson, H. B.  and J. Simons, On stable currents and
their applications to global problems, Ann. Math., 98 (1973), 427--450.

\bibitem{} Le Khong Van, Curvature estimate for the volume growth of 
globally minimal submanifolds, Math. Ann., 296 (1993), 103--118. 

\bibitem{} Le Khong Van and A. T. Fomenko , Lagrangian
manifolds and the Maslov index in the theory of minimal surfaces, Sov.
Math. Dokl., 37 (1988), 330--333.

\bibitem{} Lee, Y. I., Lagrangian minimal surfaces in  Kaehler-Einstein surfaces of negative scalar curvature, Comm. Anal. Geom., 2 (1994), 579--592. 

\bibitem{} Leung, D. S. and K. Nomizu,  The axiom of spheres in
Riemannian geometry, J. Differential Geometry, 5 (1971), 487--489.

\bibitem{} Leung, P. F.,  On the topology of a
compact submanifold of a sphere with bounded second fundamental form,
Manusc. Math., 79 (1993), 183--185.

\bibitem{} Leung, P. F.,  On the curvature of minimal submanifolds  in a sphere,
Geom. Dedicata, 56 (1995), 5--6. 

\bibitem{} Levy, H., Forma canonica dei $ds^2$ per i quali si annulano i simboli di Riemann a cinque indici, Rend. Acad. Lincei, 3 (1926), 65--69.

\bibitem{} Lewy, H., On the existence of a closed convex surface realizing a given
Riemannian metric, Proc. Nat. Acad. Sci. USA, 24 (1938), 104--106.

\bibitem{} Lewy, H., On the boundary behavior of minimal
surfaces, Proc. Nat. Acad. Sci. USA, 37 (1951), 103--110.

\bibitem{} Li, A, M. and J. M. Li,  An intrinsic rigidity
theorem for minimal submanifolds in a
sphere,  Arch. Math., 58 (1992), 582--594. 

\bibitem{} Li, H., The Ricci curvature of totally real
$3$-dimensional submanifolds of the nearly
Kaehler $6$-sphere,  Bull. Belg. Math. Soc. , 3 (1996), 193--199. 

\bibitem{} Li, P. and S. T. Yau, A new conformal invariant and
its applications to the Willmore conjecture and the first eigenvalue of compact
surface, Inven. Math. , 69 (1982), 269--291.

\bibitem{} Li, S.-J., Null 2-type surfaces in $E^m$ with parallel
normalized mean curvature vector, Math.
J. Toyama Univ., 17 (1994), 23--30.

\bibitem{} Li, S.-J., Null 2-type Chen surfaces, Glasgow Math. J., 37 (1995), 233--242.

\bibitem{} Li, X. X., Minimal immersions of $S^2$ into $S^{2m}(1)$
with degree $2m+2$,  Kodai Math. J., 18 (1995), 351--364.

\bibitem{} Li, Z. Q., Counterexamples to the conjecture on minimal
 $S^2$ in $CP^n$ with
constant Kaehler angle, Manuscripta
Math., 88 (1995), 417--431.

\bibitem{} Liao, R. J.,  Scalar curvature of
Kaehler submanifolds in a complex projective
space, Kexue Tongbao, 32 (1987), 865--869.

\bibitem{} Liao, R. J.,   Sectional curvature of
Kaehler submanifolds in a complex projective
space, Chinese Ann. Math. Ser. B, 9 (1988),  292--296.

\bibitem{}  Lichnerowicz, A.,  Sur les espaces
riemanniens compl\`etement harmoniques,
Bull. Soc. Math. France, 72 (1944), 146--168.

\bibitem{}  Lichnerowicz, A.,  Sur
une in\'egalit\'e relative aux espaces
riemanniens compl\`etement harmoniques,
C. R. Acad. Sci. Paris, 218 (1944), 436--437.

\bibitem{}  Lichnerowicz, A.,  Sur les espaces
riemanniens compl\`etement harmoniques, C.
R. Acad. Sci. Paris, 218 (1944), 493--495.

\bibitem{} Liebmann, H., \"Uber die Verbiegung
der geschlossenen Fl\"achen positiver
Kr\"ummung, Math. Ann., 53 (1900), 81--112.

\bibitem{} Li-Jost, X., Uniqueness of minimal surfaces in Euclidean 
and hyperbolic $3$-space, Math. Z., 217 (1994), 275--285. 

\bibitem{} Lin, C. S., The local isometric embedding in $R^3$ of
two-dimensional Riemannian manifolds with
nonnegative curvature, J. Differential
Geometry, 21 (1985), 213--230.

\bibitem{} Lin, C. S.,  The local isometric embedding in $R^3$ of
two-dimensional Riemannian manifolds with
Gaussian curvature changing sign
cleanly, Comm. Pure Appl. Math., 39 (1986), 867--887.

\bibitem{}   Liouville, J.,  Note au sujet de l'article
pr\'{e}ced\'{e}nt, J.  Math. Pure  Appl.  12 (1847), 265-290.

\bibitem{} Little, J. A., Manifolds with planar geodesics, J.
Differential Geometry , 11 (1976), 265--285.

\bibitem{} L\'opez, F. J., The classification of complete minimal 
surfaces with total curvature greater than
$-12\pi$, Trans. Amer. Math. Soc., 334 (1992), 49--74. 

\bibitem{} L\'opez, F. J.,  A complete minimal Klein
 bottle in $E^3$, Duke Math. J., 71 (1993), 23--30.

\bibitem{} L\'opez, F. J.,  On complete nonorientable
 minimal surfaces with low total curvature,
Trans. Amer. Math. Soc., 348 (1996), 2737--2758. 

\bibitem{} L\'opez, F. J. and F. Martin, Complete nonorientable
minimal surfaces and symmetries , Duke Math. J., 79 (1995), 667--686.

\bibitem{} L\'opez, F. J. and F. Martin,  Complete
nonorientable minimal surfaces with the
highest symmetry group, Amer. J. Math., 119 (1997), 55--81.

\bibitem{} L\'opez, F. J., Martin, F. and D. Rodriguez, Complete
minimal surfaces derived from the
Chen-Gackstatter genus two example,
pre\-print, 1997. 

\bibitem{} L\'opez, F. J., Ritor\'e, M. and F. Wei, A
characterization of Riemann's minimal
surfaces, J. Differential Geoemtry, 47(1997), 376--397. 

 \bibitem{}  L\'opez, F. J. and A. Ros, Complete minimal
surfaces with index one and stable constant
mean curvature surfaces, Comment. Math. Helv.  64 (1989), 34--43.

 \bibitem{}  L\'opez, F. J. and A. Ros, On embedded complete
minimal surfaces of genus zero,  J. Differential Geometry, 33 (1991), 293--300.

\bibitem{} L\'opez, J., Bilinear relations on minimal surfaces, preprint, 1995. 

\bibitem{} P\`orez, J. and A. Ros, Some uniqueness and
nonexistence theorems for embedded minimal
surfaces, Math. Ann., 295 (1993), 513--525.

\bibitem{} P\`orez, J. and A. Ros, The space of properly
embedded minimal surfaces with finite
total curvature, Indiana Univ. Math. J., 45 (1996), 177--204.

\bibitem{}  L\"ubke, M., Chernklassen von  Hermite-Einstein-Vektorbundeln ,
Math. Ann., 260 (1982), 133--141.

\bibitem{} Ludden, G. D., Okumura, M. and K. Yano, A totally real surface
in $CP^{2}$ that is not totally 
geodesic, Proc. Amer. Math. Soc., 53 (1975), 186--190.

\bibitem{} Lumiste, \"U., Die $n$-dimensionalen Minimalfl\"achen mit
einer $(n-1)$-dimen\=sionalen asymptotischen
Richtung im jedem Punkte, Tartu Riikl.
\"Ul. Toimetised, 62 (1958), 117--141.

\bibitem{} Lumiste, \"U. and M. V\"aljas, On geometry of totally
quasiumbilical submanifolds, Tartu
Riikl. \"Ul. Toimetised , 836 (1989), 172--185. 

\bibitem{} Luo, H. P., Stable hypersurfaces with constant mean
curvature in $R^n$, Hiroshima Math. J., 26 (1996), 587--591. 

\bibitem{} Maeda, S. and Sato, N.,  On submanifolds all of
whose geodesics are circles in a complex space
form, Kodai Math. J., 6 (1983), 157--166.

\bibitem{} Maeda, S., Ohnita, Y. and S. Udagawa, On slant
immersions into K\"ahler manifolds, Kodai
Math. J., 16 (1993), 205--219.

\bibitem{} Mainardi, G.,  Su la teria generale delle
superficie, Giornale dell' Istituteo Lombardo di Sci. Lett. Art., 9 (1956), 385--398.

\bibitem{} Martin, F. and D. Rodriguez, A new approach to the
construction of complete minimal surfaces
derived from the genus two Chen-Gackstatter
example, Illinois J. Math., 41 (1997), 171--192. 

\bibitem{} Mashimo, K., Order of standard isometric immersions of
CROSS as helical geodesic immersions, Tsukuba J. Math., 7 (1983), 257--263.

\bibitem{} Mashimo, K.,  Homogeneous totally real submanifolds of $S^6$,
Tsukuba J. Math.  9 (1985), 185--202.

\bibitem{} Mashimo, K.,  On the stability of Cartan embeddings of compact symmetric
spaces, Archiv Math., 58 (1992), 500--508. 

\bibitem{} Mashimo, K., and H. Tasaki, Stability of closed Lie
subgroups in compact Lie groups, Kodai Math. J.  13 (1990), 181--203.


\bibitem{} Mashimo, K., and H. Tasaki,  Stability of maximal tori in compact Lie groups, Algebras, Groups and Geometries  7 (1990), 114--126.

\bibitem{} Matsuyama, Y., Rigidity of hypersurfaces with
constant mean curvature, T\^ohoku Math. J., 28 (1976), 199--213.

\bibitem{}  Maxwell, J. C., On the cyclide,
Quart. J. Pure Appl. Math., 34, 1867.

\bibitem{}  Meeks, W. H., III,  The geometry and
conformal structure of triply periodic
minimal surfaces in $R^3$, Ph. D. Thesis,
Univ. of California, Berkeley, 1975.

\bibitem{}  Meeks, W. H., III,  The classification of complete minimal surfaces with total
curvature greater than $-8\pi$, Duke Math. J., 48 (1981), 523--535.

\bibitem{}  Meeks, W. H., III,  The theory of triply periodic minimal surfaces, Indiana Univ.
Math. J., 39 (1990), 877--936.

\bibitem{}  Meeks, W. H., III,  The geometry of periodic minimal surfaces, Proc. Symp.
Pure Math. AMS, 54-I (1993), 327--374.

\bibitem{}  Meeks, W. H., III and H. Rosenberg, The geometry, topology,
and existence of doubly periodic minimal
surfaces, C. R. Acad. Sci. Paris S\'er. I Math., 306 (1988), 605--609.

\bibitem{}  Meeks, W. H., III and H. Rosenberg,  The global theory of doubly periodic minimal
surfaces,  Invent. Math., 97 (1989), 351--379.

\bibitem{}  Meeks, W. H., III and H. Rosenberg, The maximum principle
at infinity for minimal surfaces in flat three manifolds, Comment. Math. Helv.,
65 (1990), 255--270.

\bibitem{}  Meeks, W. H., III and H. Rosenberg, The geometry of
periodic minimal surfaces,  Comment. Math. Helv., 68 (1993), 538--578.

\bibitem{}  Meeks, W. H., III and H. Rosenberg, The geometry and
conformal structure of properly embedded
minimal surfaces of finite topology in
$R^3$, Invent. Math., 114 (1993), 625--639.

\bibitem{}  Meeks, W. H., III and B. White, Minimal
surfaces bounded by convex curves in
parallel planes, Comment. Math. Helv., 66 (1991), 263--278.

\bibitem{}  Meeks, W. H., III and B. White, The space of minimal 
annuli bounded by an extremal pair of
planar curves,  Comm. Anal. Geom., 1 (1993), 415--437. 

\bibitem{} Meeks, W. H., III and S. T. Yau, The existence of
embedded minimal surfaces and the problem of
uniqueness , Math. Z., 179 (1982), 151--168.

\bibitem{} J. Meusnier, M\'emoire sur la
courbure des surfaces, M\'emoires  Math. Phys., 10 (1785), 477--510.

\bibitem{} Milnor, J. W.,  On the total curvature of knots,
Ann. Math., 52 (1950), 248--257.

\bibitem{} Milnor, J. W.,   On manifolds
homeomorphic to the 7-sphere, Ann. Math., 64 (1956), 399--405.

\bibitem{} Minding, F., Wie sich entscheiden l\"asst, ob zwei
gegebene krumme Fl\"achen auf einander abwickelbar sind oder nicht: nebst
Bemerkungen \"uber die F\"achen von unver\"anderlichen Kr\"ummung\-smasse,
J. Reine Angew. Math., 19 (1839), 370--387.

\bibitem{} Minding, F.,  Beitr\"age zur Theorie der k\"urzesten Linien und krummen
Fl\"achen, J. Reine Angew. Math., 20 (1940), 323--327.

\bibitem{} Minkowski, H.,   Volume und Oberfl\"ache, Math. Ann., 57 (1903), 447--495.

\bibitem{} Miyaoka, R., Complete hypersurfaces  in the space form with three principal
curvatures, Math. Z., 179 (1982), 345--354.

\bibitem{} Miyaoka, R.,  Compact Dupin hypersurfaces with three
principal curvatures , Math. Z., 187 (1984), 433-452.

\bibitem{} Miyaoka, R., Taut embeddings and Dupin hypersurfaces, 
Lecture Notes in Math., Springer, Berlin-New York, 1090 (1984), 15--23. 

\bibitem{} Miyaoka, R., Dupin hypersurfaces and a Lie invariant, Kodai
Math. J., 12 (1989), 228--256. 

\bibitem{} Miyaoka, R., Dupin hypersurfaces with six principal curvatures, Kodai
Math. J., 12 (1989), 308--315.

\bibitem{} Miyaoka, R., The linear isotropy group of $G_ 2/SO(4)$, the Hopf fibering and isoparametric hypersurfaces,  Osaka J. Math., 30 (1993), 179--202. 

\bibitem{} Miyaoka, R., The homogeneity of isoparametric hypersurfaces with
six principal curvature,  preprint, 1998. 

\bibitem{} Miyaoka R. and T. Ozawa, Construction of taut embeddings and Cecil-Ryan
conjecture,  Geometry of manifolds, (ed. K. Shiohama)  Academic Press, New
York, 181--189, 1988.

\bibitem{} Miyaoka R. and K. Sato, On complete minimal
surfaces whose Gauss map misses two
directions,  Arch. Math., 63 (1994), 565--576.

\bibitem{} Miyaoka, R. and N. Takeuchi, A note on Ogiue-Takagi
conjecture on a characterization of Euclidean
$2$-spheres, Mem. Fac. Sci. Kyushu Univ.
Ser. A , 46 (1992), 129--135. 

\bibitem{} Miyata, Y., 2-type
surfaces of constant curvature in $S^n$,
Tokyo J. Math., 11 (1988), 157--204.

\bibitem{} Miyazawa, T. and G. Chuman, On certain subspaces of
Riemannian recurrent spaces, Tensor, N. S., 23 (1972), 253--260.

\bibitem{} Mo, X. and R. Osserman, On the Gauss map and
total curvature of complete minimal
surfaces and an extension of Fujimoto's
theorem, J. Differential Geometry, 31 (1990), 343--355.

\bibitem{} Montiel, S. and A.  Ros, Minimal immersions of
surfaces by the first eigenfunctions and
conformal area, Invent. Math. ,83 (1985), 153--166.

\bibitem{} Montiel, S. and A.  Ros, Schr\"odinger operators
 associated to a holomorphic map,
Global differential geometry and global 
analysis, 147--174, Lecture Notes in Math., 1481 (1990), 147--174.

\bibitem{} Montiel, S. and A.  Ros, 
Compact hypersurfaces: the Alexandrov theorem
for higher order mean curvature , 
Differential geometry, 279--296, Pitman
Monographs Surveys Pure Appl. Math., 52, 1991,

\bibitem{} Montiel, S.,  Ros, A. and F. Urbano,  Curvature
pinching and eigenvalue rigidity for minimal
submanifolds, Math. Z., 191 (1986), 537--548.

\bibitem{} Moore, H., Minimal submanifolds with finite total
scalar curvature, Indiana Univ. Math.
J., 45 (1996), 1021--1043.

\bibitem{} Moore, J. D.,  Isometric immersions of riemannian
products, J. Differential Geometry ,
5 (1971), 159--168.

\bibitem{} Moore, J. D.,  Isometric immersions
of space forms in space forms, Pacific
J. Math., 40 (1972), 157--166.

\bibitem{} Moore, J. D.,   
Equivariant embeddings of Riemannian
homogeneous spaces, Indiana Univ. Math. J., 25 (1976), 271--279.

\bibitem{} Moore, J. D.,    Conformally flat
submanifolds of Euclidean space, Math.
Ann., 225 (1977), 89--97.

\bibitem{} Moore, J. D.,  On extendability of
 isometric immersions of spheres,
Duke Math. J., 85 (1976), 685--699. 

\bibitem{} Moore, J. D. and J. M. Morvan, Sous-vari\'et\'es conform\'ment plates
de codimension quatre, C. R. Acad. Sc.
Paris, 287 (1978), 655--657.

\bibitem{} Morgan, F., Almost every curve in $R^3$ bounds a unique
area minimizing surface, Invent. Math., 45 (1978), 253--297.

\bibitem{} Mori, H.,  A note on the
stability of minimal surfaces in the
3-dimensional unit sphere, Indiana
Univ. Math. J., 26 (1977), 977--980.
,
\bibitem{} Mori, H.,  Stable complete
constant mean curvature  surfaces in $R^3$
and $H^3$, Trans. Amer. Math. Soc., 278 (1983), 671--687.

\bibitem{} Morrey, C. B., The problem of Plateau on a Riemannian
manifold, Ann. Math. 49 ({1948}), 807--851.

\bibitem{} Morrey, C. B., Multiple integrals in the calculus of
variations, Springer-Verlar, New York-Berlin, 1966.

\bibitem{}  Morvan, J. M., Classe de Maslov d'une sous-vari\'et\'e Lagrangienne et
minimalit\'e , C. R. Acad. Sc. Paris, 292 (1981), 633--636.

\bibitem{} Morvan, J. M. and G. Zafindratafa, Conformally
flat submanifolds, Ann. Fac. Sci. Toulouse Math., 8 (1986), 331--348.

\bibitem{}  M\"unzner, H. F., Isoparametrische Hyperfl\"achen in
Sph\"aren, I, Math. Ann. , 251 (1980), 57--71.

\bibitem{}  M\"unzner, H. F.,   Isoparametrische Hyperfl\"achen in Sph\"aren, II , Math.
Ann., 256 (1981), 215--232.

\bibitem{} Myers, S. B., Riemannian manifolds with positive mean
curature, Duke Math. J., 8 (1941), 401--404.

\bibitem{} N. Nadirashvili, Hadamard's and Calabi-Yau's conjectures on negatively 
curved and minimal surfaces, Invent. Math., 126 (1996), 457--465.

\bibitem{}  Nagano, T., On the minimal eigenvalues of the Laplacian in
Riemannian manifolds, Sci. Papers Coll. Gen. ed. Univ. Tokyo , 11 (1961), 177-182.

\bibitem{}  Nagano, T.,  Transformation groups
on compact symmetric space, Trans. Amer. Math. Soc.  118 (1965), 428--453.

\bibitem{}  Nagano, T., The involutions
of compact symmetric spaces, I, Tokyo J. Math., 11 (1988), 57--79.

\bibitem{}  Nagano, T.,  The involutions of
compact symmetric spaces, II, Tokyo J. Math., 15 (1992), 39--82.

\bibitem{}  Nagano, T., and B. Smyth, Sur les surfaces
minimales hyperelliptiques dans un tore, C.
R. Acad. Sci. Paris S\'er. , 280 (1975), 1527--1529.

\bibitem{}  Nagano, T., and B. Smyth,  Minimal varieties and
harmonic maps in tori, Comment. Math. Helv., 50 (1975), 249--265.

\bibitem{}  Nagano, T., and B. Smyth,  Minimal varieties in
tori,  Differential geometry, Proc. Sympos., Pure Math., 27-1 (1975), 189--190.

\bibitem{}  Nagano, T., and B. Smyth,  Minimal surfaces in tori
by Weyl groups, Proc. Amer. Math. Soc., 61 (1976), 102--104 .

\bibitem{}  Nagano, T., and B. Smyth,  Periodic minimal surfaces, Comment. Math. Helv., 53 (1978), 29--55.

\bibitem{}  Nagano, T., and B. Smyth,  Periodic minimal surfaces
and Weyl groups, Acta Math., 145 (1980), 1--27.

\bibitem{}  Nagano, T. and M. Sumi, The structure of the symmetric spaces with applications , Geometry of Manifolds (K. Shiohama, ed.)
Academic Press, 111-128, 1989.

\bibitem{}  Nagano, T. and M. Sumi,  The spheres in symmetric
spaces, Hokkaido Math. J., 20 (1991), 331--352.

\bibitem{}  Nagano, T., and T. Takahashi, Homogeneous hypersurfaces in Euclidean space, J. Math. Soc. Japan, 12 (1960), 1--7.

\bibitem{} Naitoh, H., Totally real parallel submanifolds, Tokyo
J. Math., 4 (1981), 279--306.

\bibitem{} Naitoh, H.,   Parallel submanifolds of complex space forms I,
Nagoya Math. J., 90 (1983), 85--117;  II , ibid,
91 (1984), 119--149.

\bibitem{} Naitoh, H. and M. Takeuchi,  Totally real  submanifolds and symmetric bounded domain, Osaka J. Math., 19 (1982), 717--731.

\bibitem{} Nakagawa, H., On a certain minimal immersion of a Riemannian
manifold into a sphere, Kodai Math. J. , 3 (1980), 321--340.

\bibitem{} Nakagawa, H. and K. Ogiue,  Complex space forms immersed
in complex space forms, Trans. Amer. Math. Soc., 219 (1976), 289--297.

\bibitem{} Nakagawa, H. and R. Takagi,  On locally symmetric
Kaehler submanifolds in complex projective
space, J. Math. Soc. Japan, 28 (1976), 638--667.

\bibitem{} Nakamura, G. and Y. Maeda, Local isometric
embedding problem of Riemannian $3$-manifold into $R^6$,
Proc. Japan Acad. Ser. A Math. Sci., 62 (1986), 257--259.

\bibitem{}  Nash, J. F.,  $C^1$-isometric imbeddings, Ann. Math. , 60 (1954),
383--396.

\bibitem{}  Nash, J. F.,  The imbedding problem for Riemannian
manifolds, Ann.  Math.  63 (1956), 20--63.

\bibitem{} Nayatani, S., On the Morse index of complete minimal
surfaces in Euclidean space, Osaka J. Math., 27 (1990), 441--451.

\bibitem{} Nayatani, S.,   Lower bounds for the
Morse index of complete minimal surfaces in
Euclidean $3$-space, Osaka J. Math., 27 (1990), 453--464.

\bibitem{} Nayatani, S.,  Morse index and Gauss
maps of complete minimal surfaces in Euclidean
$3$-space, Comment. Math. Helv., 68 (1993), 511--537. 

\bibitem{} Niebergall, R., Dupin hypersurfaces in $R^5$, I, Geom.
Dedicata, 40 (1991), 1--22.

\bibitem{} Niebergall, R.,  Dupin hypersurfaces in
$R^5$, II, Geom. Dedicata, 41 (1992), 5--38.

\bibitem{} Niebergall, R.,  Tight analytic immersions of highly connected manifolds,
Proc. Amer. Math. Soc., 120 (1994), 907--916. 

\bibitem{} Nikolaevskij, Yu. A., Totally umbilical submanifolds
in $G(2,n)$. I ,  Ukrain. Geom. Sb. No. 34 (1991), 83--98.

\bibitem{} Nikolaevskij, Yu. A.,  Classification of multidimensional submanifolds in a
Euclidean space with a totally geodesic
Grassmann image, Russian Acad. Sci. Sb. Math., 76 (1993), 225--246.

\bibitem{} Nikolaevskij, Yu. A.,  Totally umbilical submanifolds of symmetric
spaces, Mat. Fiz. Anal. Geom., 1 (1994), 314--357.

\bibitem{} Nirenberg, L., The Weyl and Minkowski problems in differential
geometry in the large, Comm. pure Appl. Math.  6 (1953), 337--394.

\bibitem{} Nitsche, J. C. C., Lectures on minimal surfaces, Cambridge
Univ. Press, Oxford, 1989.

\bibitem{} N\"olker, S, Isometric immersions of warped products, Differential Geom. Appl., 6 (1996), 1--30.

\bibitem{} Nomizu, K.,  Some results in \'E. Cartan's theory on isoparametric
families of hypersurfaces , Bull. Amer. Math. Soc., 79 (1973), 1184--1188.

\bibitem{} Nomizu, K., Conditions for constancy of the holomorphic sectional curvature,
J. Differential Geometry, 8 (1973),
335--339. 

\bibitem{} Nomizu, K., Isometric immersions of the hyperbolic plane
into the hyperbolic space, Math. Ann.,
205 (1974), 181--192.

\bibitem{} Nomizu, K.,  \'Elie Cartan's work on isoparametric families of hypersurfaces, 
Differential geometry (Proc. Sympos. Pure Math.,  27-1 (1975), 191--200.

\bibitem{} Nomizu, K. and B. Smyth, Differential geometry
of complex hypersurfaces, II, J. Math. Soc. Japan, 20 (1968), 498--521.

\bibitem{} Nomizu, K. and B. Smyth,  A formula of Simons'
type and hypersurfaces with constant mean
curvature, J. Differential Geometry, 3 (1969), 367--377.

\bibitem{} Nomizu, K., and K. Yano, On circles and spheres in
Riemannian geometry, Math. Ann., 210 (1974), 163--170.

\bibitem{} Ogiue, K., Positively curved complex submanifolds immersed in a complex
projective space, II, Hokkaido Math.
J., 1 (1972), 16--20.

\bibitem{} Ogiue, K.,  On Kaehler immersions, Canad. J. Math., 24 (1972),
1178--1182.

\bibitem{} Ogiue, K., Differential geometry of Kaehler
submanifold, Adv. in Math., 13 (1974), 73--114.

\bibitem{} Ogiue, K., Positively curved  totally real minimal submanifolds
immersed in a complex projective space,  Proc. Amer. Math. Soc., 56 (1976), 264--266.

\bibitem{} Ogiue, K.,   Some recent topics in the theory of
submanifolds, Sugaku Expositions, 4 (1991), 21--41.

\bibitem{} Ogiue, K. and R. Takagi, A submanifold which
contains many extrinsic circles, Tsukuba J. Math., 8 (1984), 171--182. 

\bibitem{} Ogiue, K. and N. Takeuchi, A geometric construction
of a torus which contains five circles
through each point, Bull. Tokyo Gakugei
Univ.,  1992, 15--18.

\bibitem{} Ogiue, K. and N. Takeuchi Hulahoop surfaces, J.
Geometry, 46 (1993), 127--132.


\bibitem{}  Oh, Y. G., Second variation and stability of minimal
Lagrangian submanifolds in Kaehler
manifolds, Invent. Math., 101 (1990), 501--519.

\bibitem{} Ohnita, Y., Stable minimal submanifolds
in compact rank one symmetric spaces
, T\^ohoku Math. J., 38 (1986), 199--217.

\bibitem{} Ohnita, Y., Totally real submanifolds with nonnegative
sectional curvature, Proc. Amer. Math. Soc., 97 (1986), 474--478.

\bibitem{} Ohnita, Y.,  On stability of minimal
submanifolds in compact symmetric space , Compos. Math., 64 (1987), 157--189.

\bibitem{} Ohnita, Y.,  Minimal surfaces with
constant curvature and K\"ahler angle in complex
space forms, Tsukuba J. Math., 13 (1989), 191--207.

\bibitem{} Okumura, M., Submanifolds and a pinching problem on the
second fundamental tensor, Trans. Amer.
Math. Soc., 178 (1973), 285--291.

\bibitem{} Oliker, V. I., Hypersurfaces in $R^{n+1}$ with
prescribed Gaussian curvature and related equations of Monge-Ampere type, 
Comm. Partial Differential Equations, 9 (1984), 807--838.

\bibitem{} Oliker, V. I.,The problem of embedding $S^n$ into
$R^{n+1}$ with prescribed Gauss
curvature and its solution by variational
methods, Trans. Amer. Math. Soc., 295  (1986), 291--303.

\bibitem{} de Oliveira Filho, G.,  Compactification of minimal submanifolds
of hyperbolic space,  Comm. Anal. Geom., 1 (1993), 1--29.

\bibitem{} de Oliveira Filho, M. E. G. G., Superf\'icies M\'imimas
nao-orient\'aveis no $R^n$, Ph. D. Thesis,
IMEUSP, 1984.

\bibitem{} de Oliveira, M. E. G. G. and \'E. Toubiana,
Surfaces non-orientables de genre deux,
Bol. Soc. Brasil. Mat. (N.S.), 24 (1993), 63--88.

 \bibitem{} Olmos, C., Homogeneous submanifolds of higher rank and
parallel mean curvature,  J. Differential Geometry 39 (1994), 605--627.

 \bibitem{} Olmos, C.,   Orbits of rank one and
parallel mean curvature, Trans. Amer. Math. Soc., 347 (1995), 2927--2939. 

\bibitem{} Olmos, C. and C. S\'anchez, A geometric
characterization of the orbits of
$s$-representations,  J. Reine Angew. Math., 420 (1991), 195--202.
 
\bibitem{} O'Neill, B.,  Immersion of manifolds of
non-positive curvature, Proc. Amer. Math. Soc., 11  (1960), 132--134.

\bibitem{} O'Neill, B. and E. Stiel,   Isometric immersions of
constant curvature manifolds, Michigan Math. J., 10 (1963), 335--339.

\bibitem{} Osserman, R.,  Minimal surfaces in the large,
Comment. Math. Helv. , 35 (1961), 65--76.

\bibitem{} Osserman, R.,  Global properties of minimal surfaces in $E^3$ and $E^n$,
Ann. Math., 80 (1964), 340--364.

\bibitem{} Osserman, R.,   Minimal varieties, Bull. Amer.
Math. Soc., 75 (1969), 1092--1120.

\bibitem{} Osserman, R.,   A survey of minimal surfaces, Van Nostrand, 1969.

\bibitem{} Osserman, R.,  A proof of the regularity everywhere of the classical
solution to Plateau's problem, Ann. Math., 91 (1970), 550--569.

\bibitem{} Osserman, R.,   On the convex hull
property of immersed maniflods, J.
Differential Geometry, 6 (1971), 267--270.

\bibitem{} Osserman, R.,  Isoperimetric and
related inequalities, Proc. Symp. Pure Math., 27 (1975), 207--215.

\bibitem{} Osserman, R.,  The isoperimetric
inequality, Bull. Amer. Math. Soc., 84 (1978), 1182--1238.

\bibitem{} Osserman, R.,  Minimal surfaces, Gauss maps, total curvature, eigenvalue
estimates, and stability, The Chern
Symposium 1979, Springer-Verlag, 199--228.

\bibitem{} Osserman, R., Curvature in the eighties, Amer. Math.
Monthly, 97 (1990), 731--756.

\bibitem{} Osserman, R. and M. Ru, An estimate for the
Gauss  curvature of minimal surfaces in
$R^m$ whose Gauss map omits a set of
hyperplanes, J. Differential Geometry, 46 (1997), 578--593.

\bibitem{} Osserman, R. and M. Schiffer, Doubly-connected
minimal surfaces,  Arch. Rational Mech. Anal., 58 (1974), 285--307. 

\bibitem{} Otsuki, T., On the existence of
solutions of a system of quadratic equations
and its geometrical application, Proc. Japan Acad., 29 (1953), 99--100.

\bibitem{} Otsuki, T.,  Isometric imbedding of Riemannian
manifolds in a Riemannian manifold, J. Math. Soc. Japan, 6 (1954), 221--234.

\bibitem{} Otsuki, T.,  Note on the isometric
imbedding of compact Riemannian manifolds
in Euclidean spaces, Math. J. Okayama Univ., 5 (1956), 95--102.

\bibitem{} Otsuki, T.,  Minimal hypersurfaces in a Riemannian
manifold of constant curvature, Amer. J. Math., 92 (1970), 145--173.

\bibitem{} Otsuki, T.,  On principal normal vector fields of submanifolds in a
Riemannian manifold of constant curvature, J. Math. Soc. Japan, 22 (1970), 35--46.

\bibitem{} Ozawa, T., On critical sets of distance functions to a
 taut submanifold, Math. Ann., 276 (1986), 91--96.

\bibitem{}  Ozeki, H. and M. Takeuchi,  On some types of isoparametric hypersurfaces in spheres, I, T\^ohoku Math J.  27 (1975), 515--559.

\bibitem{}  Ozeki, H. and M. Takeuchi,  On some
types of isoparametric hypersurfaces in
spheres, I, II, T\^ohoku Math J.  28 (1976), 7--55.

\bibitem{} Pak, J. S., Planar geodesic submanifolds in complex 
space forms, Kodai Math. J., 1 (1978), 187--196.

\bibitem{} Pak, J. S. and K. Sakamoto, Submanifolds with
proper $d$-planar geodesics immersed in
complex projective spaces, T\^hoku Math. J., 38 (1986),297--311.

\bibitem{} Pak, J. S. and K. Sakamoto,  $4$-planar geodesic
Kaehler immersions into a complex
projective space,  Proc. Amer. Math. Soc., 102 (1988), 995--999. 

\bibitem{}  Palmer, B., Surfaces of constant mean curvature in space forms, 
Thesis, Stanford University, 1986.

\bibitem{}  Palmer, B.,   Stability of minimal
hypersurfaces, Comment. Math. Helv., 66 (1991), 185--188. 

\bibitem{}  Palmer, B.,  Calibrations and Lagrangian submanifolds, T\^ohoku
Math. J. , 50 (1998), 303--315.

\bibitem{} Peetre, J., A generalization of Courant's nodal line
theorem, Math. Scand., 5 (1959), 15--20.

\bibitem{} Peng, C. K. and Z. Hou, A remark on the
isoparametric polynomials of degree 6, 
 Lecture Notes in Math., 1369 (1989), 222-224.

\bibitem{} Peng, C. K. and C. L. Terng,  The scalar curvature of minimal
hypersurfaces in sphere, Math. Ann., 266 (1983), 105--113.

\bibitem{} Peng, C. K. and Z. Z. Tang, Brouwer degrees of gradient maps of
isoparametric functions. , Sci. China Ser. A , 39 (1996), 1131--1139.

\bibitem{}  P\'erez, J., On  singly-periodic minimal surfaces with
planar ends, Trans. Amer. Math. Soc., 349 (1997), 2371--2389.

\bibitem{} P\'erez, J. and A. Ros, The space of properly embedded minimal surfaces with
finite total curvature, Indiana Univ. Math. J., 45 (1996), 177--204.

\bibitem{} Peterson, S. P., Arithmetic distance on compact symmetric
spaces, Geom. Dedicata, 23 (1987), 1--14.

\bibitem{}  Pinkall, U., Dupin hypersurfaces,  Math. Ann. , 270 (1985), 427--440. 

\bibitem{}  Pinkall, U., Dupin hypersurfaces,  Dupinsche
Hyperflachen in $E^4$, Manuscripta Math., 51 (1985), 89--119. 

\bibitem{}  Pinkall, U., Dupin hypersurfaces,  Hopf tori in $S^3$, 
Invent. Math., 81 (1985), 379--386. 

\bibitem{}  Pinkall, U., Dupin hypersurfaces, Tight surfaces and
regular homotopy,  Topology, 25 (1986), 475--481.

\bibitem{}  Pinkall, U., Dupin hypersurfaces, Curvature
properties of taut submanifolds, Geom. Dedicata, 20 (1986), 79--83. 

\bibitem{}  Pinkall, U., Dupin hypersurfaces, Compact conformally flat
hypersurfaces, Conformal Geometry, Max
Planck Institu f\"ur Math., Bonn, 217--236, 1988.

\bibitem{}  Pinkall, U. and I. Sterling,  Willmore surfaces, Math. Intellignercer, 9 (1987),
38--43. 

\bibitem{}  Pinkall, U. and I. Sterling,  On the classification of constant mean
curvature tori, Ann. Math., 130 (1989), 407--451. 

\bibitem{}  Pinkall, U. and G. Thorbergsson, Deformations of Dupin
hypersurfaces, Proc. Amer. Math. Soc., 107 (1989), 1037--1043.

\bibitem{}  Pinkall, U. and G. Thorbergsson,  Taut $3$-manifolds,
Topology, 28 (1989), 389--401. 


\bibitem{} Pinl, M., B-Kugelbilder reeler Minimalfl\"achen in 
$R^4$, Math. Z., 59 (1953), 290--295.

\bibitem{} Pitts, J. T. and J. H. Rubinstein, Equivariant minimax and
minimal surfaces in geometric three-manifolds,  Bull. Amer. Math. Soc.
(N.S.), 19 (1988), 303--309. 

\bibitem{} Pogorelov, A. W., Continuous maps of bounded
variations (Russian), Dokl. Acad. Nauk SSSR, 111 (1956), 757--759.

\bibitem{} Pogorelov, A. W., Extensions of the theorem of
Gauss on spherical representation to the case of surfaces of bounded extrinsic
curvature (Russian), Dokl. Acad. Nauk SSSR, 111 (1956), 945--947.

\bibitem{} Pogorelov, A. W., The Minkowski
multidimensional problem (Russian), Nauka, Moscow, 1975.

\bibitem{} Poisson, S. D., Trait\'e de Mecanique, 2 volumes, 2nd ed., Paris, 1833.

\bibitem{} Polthier, K., New periodic minimal surfaces in $H\sp 3$, Workshop on Theoretical and Numerical Aspects of Geometric Variational Problems, Canberra, 1990,  Austral. Nat.
Univ., 201--210. 

\bibitem{} Rademacher, H., Conformal and isometric immersions of conformally flat Riemannian manifolds into spheres and Euclidean spaces,  Conformal Geometry, Max
Planck Institu f\"ur Math., Bonn, 191--216, 1988.

\bibitem{}  Rad\'o, T., On Plateau's problem, Ann. Math., 31 (1930),
457--469.

\bibitem{}  Rad\'o, T.,  The problem of the least area and the problem of
Plateau, Math. Z., 32 (1930), 763--796.

\bibitem{} Reckziegel, H., Submanifolds with prescribed mean curvature vector
field, Manusc. Math., 13 (1974), 69--71.

\bibitem{} Reckziegel, H., Kr\"ummungsfl\"achen von isometrischen
Immersionen in R\"aume konstanter Kr\"ummung, Math. Ann., 223 (1976), 169--181.

\bibitem{} Reckziegel, H.,  Completeness of curvature
surfaces of an isometric immersion, J. Differential Geometry, 14 (1979), 7--20.

\bibitem{} Reckziegel, H.,  Horizontal lifts of isometric immersions into the bundle space of
a pseudo-Riemannian submersion, Lecture Notes in Math.,
Springer Verlag, 1156 (1985), 264--279.

\bibitem{}  Reilly, R. C., Variational properties of functions of the
mean curvatures for hypersurfaces in space
forms,  J. Differential Geometry, 8 (1973), 465--477. 

\bibitem{}  Reilly, R. C.,  On the first eigenvalue of the Laplacian for compact
submanifolds of Euclidean space , Comm. Math. Helv., 52 (1977), 525--533.

\bibitem{} Ricci, G., Sulla classificazione delle forme differenziali quadratiche, Rend. di
Lincei, 41 (1888), 203--207.

\bibitem{} Ricci, G., Sulla teoria intrinseca delle superficie ed in ispecie di quelle di secondo
grado, Atti R. Ist. Ven. di Lett. ed Arti, 6 (1894), 445--488.

\bibitem{} Riemann, B., Gesammelte mathematische Werke, second
edition, Teubner, Leipzig, 1892.

\bibitem{} Ripoll, J. B., Helicoidal minimal surfaces in hyperbolic space, Nagoya Math. J., 114 (1989), 65--75. 

\bibitem{} Ritor\'e, M. and A. Ros, The spaces of index one minimal surfaces and stable constant mean curvature surfaces embedded in flat three
manifolds, Trans. Amer. Math. Soc., 348 (1996), 391--410.

\bibitem{} Romon, P., A rigidity theorem for Riemann's minimal surface, Ann. Inst. Fourier
(Grenoble), 43 (1993), 485--502. 

\bibitem{} Ros, A.,  Spectral geometry of $CR$-minimal submanifolds
in the complex projective space, Kodai Math. J., 6 (1983), 88--99.

\bibitem{} Ros, A.,  Positively curved Kaehler submanifolds, Proc. Amer. Math. Soc.,
93 (1985), 329--331.

\bibitem{} Ros, A.,  A characterization of seven compact Kaehler submanifolds by
holomorphic pinching, Ann. Math., 121 (1985), 377--382.

\bibitem{} Ros, A.,  Kaehler submanifolds in the complex projective space, Lecture Notes in
Math., Springer-Verlag, 1209 (1986), 259--274.

\bibitem{} Ros, A., Compact hypersurfaces with constant scalar curvature and a congruence theorem, J. Differential Geometry, 27 (1987),
215--220. 

\bibitem{} Ros, A., Compact hypersurfaces with constant higher order mean curvature, Rev. Mat. Ibero\-amer., 3 (1988), 447--453. 

\bibitem{} Ros, A., A two-piece property for compact minimal surfaces in a
three-sphere, Indiana Univ. Math. J., 44 (1995), 841--849.

\bibitem{} Ros, A., Embedded minimal surfaces: forces, topology and
symmetries, Calc. Var. Partial Diff. Equat., 4 (1996), 469--496. 

\bibitem{}  Ros, A. and F. Urbano, Lagrangian submanifolds of
$ C^n$ with conformal Maslov form and the
Whitney sphere, J. Math. Soc. Japan, 50 (1998), 203--226.

\bibitem{} Ros, A. and L. Verstraelen,  On a conjecture of K.
Ogiue, J. Differential Geometry, 19 (1984), 561--566.

\bibitem{} Rosenberg, H. and \'E. Toubiana, A cylindrical type
complete minimal surface in a slab of $E^3$,  Bull. Sci. Math.  111 (1987),
241--245. 

\bibitem{} Ross, M., Complete nonorientable minimal surfaces in $R^3$, Comment.
Math. Helv., 67 (1992), 64--76. 

\bibitem{} Ross, M.,The second variation of nonorientable 
minimal submanifolds,  Trans. Amer. Math. Soc., 349 (1997), 3093--3104. 

\bibitem{} Ross, M. and C. Schoen,  Stable quotients of
periodic minimal  surfaces,  Comm. Anal. Geom., 2 (1994), 451--459. 

\bibitem{} Rossman, W., Minimal surfaces in $R^3$ with
 dihedral symmetry, T\^ohoku Math. J., 47 (1995), 31--54.

\bibitem{} Rouxel, B., Chen submanifolds, Geometry and Topology of
Submanifolds, VI, 185--198, 1994.

\bibitem{} Ru, M.,  On the Gauss map of minimal surfaces immersed in
$R^n$, J. Differential Geometry, 34 (1991), 411--423.

\bibitem{} Ruh, E. A., Minimal immersions of 2-spheres in $S^4$,
Proc. Amer. Math. Soc., 28 (1971), 219--222.

\bibitem{} Ruh, E. A. and J. Vilms,  The tension field of a Gauss
map, Trans. Amer. Math. Soc., 149 (1970), 569--573.

\bibitem{}  Ryan, P. J., Homogeneity and some curvature conditions for
hypersurfaces, T\^ohoku Math. J., 21 (1969), 363--388.

\bibitem{}  Ryan, P. J., Hypersurfaces with
parallel Ricci tensor, Osaka J. Math., 8 (1971), 251--259.

\bibitem{}  Ryan, P. J.,  K\"ahler manifolds as
real hypersurfaces, Duke Math. J., 40  (1973), 207--213.

\bibitem{} Sacksteder, R.,  The rigidity of hypersurfaces , J.
Math. Mech., 11 (1962), 929--939.

\bibitem{} Sakaki, M., On minimal surfaces with constant Kaehler 
angle in $CP^3$ and $CP^4$, Tsukuba J. Math., 20 (1996), 33--44.

\bibitem{} Sakamoto, K., Planar geodesic immersions, T\^ohoku Math. J., 29 (1977), 25--56.

\bibitem{} Sakamoto, K.,  Helical immersions into a unit sphere, Math.
Ann., 261 (1982), 63-80.

\bibitem{} Sakamoto, K.,  Helical immersions of compact Riemannian manifolds into a unit
sphere, Trans. Amer. Math. Soc., 288 (1985), 765--790.,

\bibitem{} Sakamoto, K.,  Helical immersions , Lecture Notes in
Math., 1201 (1986), 230--241.

\bibitem{} S\'anchez, C. U., $k$-symmetric submanifolds of $R^N$ ,
Math. Ann., 270 (1985), 297--316.

\bibitem{} S\'anchez, C. U.,  A  characterization of extrinsic $k$-symmetric
submanifolds on $R^N$, Rev. Un. Mat. Argentina, 38 (1992), 1--15. 

\bibitem{} S\'anchez, C. U.,  The invariant of
Chen-Nagano on flag manifolds, Proc.
Amer. Math. Soc., 118 (1993), 1237--1242.

\bibitem{}  Sasaki, S., On complete flat surfaces in hyperbolic
$3$-space, Kodai Math. Sem. Rep., 25 (1973), 449--457.

\bibitem{} Sato, K., Construction of higher genus minimal
 surfaces with one end and finite total
curvature, T\^ohoku Math. J.  48 (1996), 229--246.

\bibitem{} Schlafli, L.,  Nota alla Memoria del sig. Beltrami, Ann.  Mat. Pura
Appl., 5 (1873), 178-193.

\bibitem{} Schoen, A. H.,  Infinite periodic minimal surfaces
without self-intersections, Technical Note
D-5541, NASA, Cambridge, Mass, 1970 .

\bibitem{} Schoen, R., Uniquenss, symmetry, and embeddedness of
minimal surfaces, J. Differential Geometry, 18 (1983), 791--809.

\bibitem{} Schouten, J. A., \"Uber die konforme Abblildung $n$-dimensionaler
Mannigfaltigkeiten mit quadratischer Massbestimmung auf eine Mannigfaltigkeit
mit eukidischer Massbestimmung, Math. Z., 11 (1924), 38--88.

\bibitem{} Schouten, J. A., Ricci calculus,
Springer-Verlag, Berlin, 1954.

\bibitem{} Schouten, J. A. and D. van Dantzig,  \" Uber unit\" are
 Geometrie, Math. Ann., 103 (1930), 319-346.

\bibitem{} Schouten, J. A. and D. van Dantzig, \" Uber  unit\" are
 Geometrie  konstanter Kr\" ummung,
Proc. Kon. Nederl. Akad. Amsterdam, 34 (1931), 1293-1314.

\bibitem{} Schur, F.,  \"Uber den Zusammenhang der R\"aume konstanten
Kr\"ummungs\-masses mit den projektiven
R\"aumen, Math. Ann., 27 (1863), 537--576
.

\bibitem{} Schwarz, H. A.,  Gesammelte mathematische
Abhandlungen, Springer-Verlag, Berlin, 1890.

\bibitem{} Segre, B., Una propriet\'a caratteristica di tre sistemi $\infty^1$
de superficie, Att. Acc. Sci. Torino, 59 (1924), 666--671.

\bibitem{} Segre, B.,Famiglie di ipersuperficie isoparametrische
negli spazi euclidei ad un qualunque numero
di dimensoni, Atti Accad. Naz. Lincei
Rend. Cl. Sc. Fis Mat. Natur., 27 (1938), 203--207. 

\bibitem{} Segre, C., Sulle vaiet\`a che rappresentano le coppie
di punti di due piani o spazi, Rend.
Cir. Mat. Palermo, 5 (1891), 192--204. 

\bibitem{} Sekigawa, K., Almost complex submanifolds of a
$6$-dimensional sphere, Kodai Math. J., 6 (1983), 174--185.

\bibitem{} Sekigawa, K.,  Some $CR$-submanifolds in
a 6-dimensional sphere, Tensor (N.S.), 41 (1984), 13--20.

\bibitem{} Sharpe, R. W., Total absolute curvature and embedded Morse numbers, J.
Differential Geometry, 28 (1989), 59--92. 

\bibitem{} Sharpe, R. W., A proof of the Chern-Lashof conjecture in
 dimensions greater than five, Comment. Math. Helv., 64 (1989), 221--235. 

\bibitem{} Shen, Y. B., Complete submanifolds in $E^{n+p}$ with parallel mean curvature,
Chinese Ann. Math. Ser. B , 6 (1985),
345--350.  

\bibitem{} Shen, Y. B., On compact Kaehler submanifolds in $CP^{n+p}$ with
nonnegative sectional curvature, Proc. Amer. Math. Soc., 123 (1995), 3507--3512. 

\bibitem{} Shen, Y. B.,   The nearly  K\"ahler structure and minimal surfaces
in $S^6$ , Chinese Ann. Math., 19B (1998), 87--96.

\bibitem{} Shen, Y. B., Dong, Y. X. and X. Y. Guo, On scalar
curvature for totally real minimal
submanifolds in $CP^n$, Chinese Sci. Bull. , 40 (1995), 621--626.

\bibitem{} Shen, Y. B. and X. H. Zhu, Stable complete
minimal hypersurfaces in $R^{n+1}$, Amer. J. Math.,120 (1998), 103--116.

\bibitem{} Shiffman, M., On surfaces of stationary area bounded
by two circles, or convex curves, in
parallel planes, Ann. Math., 63 (1956), 77--90. 

\bibitem{} Shimizu, Y.,  On a construction of
homogeneous $CR$-submanifolds in a complex
projective space, Comm. Math. Univ. San. Pauli, 328 (1983), 203--207.

 \bibitem{} Shiohama, K. and R. Tagaki, A characterization
of a standard torus in $E^3$, J. Differential Geometry, 4 (1970), 477--485.

\bibitem{} Shirokov, P. A. C., Constant
fields of vectors and tensors of
econd order in Riemannian manifolds (Russian), Izv. Fiz.-Mat.
Obshchestva Kazan. Univ., 25 (1925),  86--114.

\bibitem{} Sikorav, J. C., Non-existence de sous-vari\'et\'e lagrangienne exacte dans
{\bf C}$^n$ (d'apr\`es Gromov), Aspects Dynamiques et Topologiques des Groupes
Infinis de Transformation de la M\'ecanique Lyon, Travaux en Cours, Hermann, Paris,
25 (1986),  95--110.

\bibitem{} Simon, U. and A. Weinstein, Anwendungen der de Rhamschen Zerlegung auf Probleme der lokalen Fl\"achentheorie, Manuscr. Math., 1  (1969),139--146.

\bibitem{} Simons, J., Minimal varieties in riemannian manifolds, Ann. Math., 88 (1968), 62--105.

\bibitem{} Singley, D., Smoothness theorems for the principal
curvatures and principal vectors of a
hypersurface, Rocky Mountain J. Math., 5 (1975), 135--144.

\bibitem{} Smyth, B., Differential geometry of complex
hypersurfaces, Ann. Math., 85 (1967), 246--266.

\bibitem{} Smyth, B., Homogeneous complex
hypersurfaces, J. Math. Soc. Japan,
20 (1968), 643--647.

\bibitem{} Smyth, B., Submanifolds of constant
mean curvature, Math. Ann., 205 (1973), 265--280.

\bibitem{} Smyth, B. and F. Xavier, Efimov's theorem in dimension
greater than two, Invent. Math., 90 (1987), 443--450.

\bibitem{} Solomon, B., The harmonic analysis of cubic isoparametric
minimal hypersurfaces. I. Dimensions $3$ and
$6$,  Amer. J. Math., 112 (1990),
157--203. 

\bibitem{} Solomon, B., The harmonic analysis of cubic isoparametric
 minimal hypersurfaces. II. Dimensions
$12$ and $24$, Amer. J. Math., 112 (1990), 205--241.

\bibitem{} Solomon, B.,  Quartic isoparametric
 hypersurfaces and quadratic forms, Math. Ann., 293 (1992), 387--398e.

\bibitem{} Somigliana, C., Sulle relazione fra il pricipio di Huygens e l'ottica
geometrica, Atti Acc. Sci. Torino, 54 (1918), 974--979.

\bibitem{} Soret, M., Maximum principle at infinity for complete 
minimal surfaces in flat $3$-manifolds, 
Ann. Global Anal. Geom., 13 (1995), 101--116. 

\bibitem{} Spivak, M.,  A comprehensive introduction to differential
geometry, I--V, Publish or Perish, Inc. Wilmington, Delaware, 1970.

\bibitem{}  Spruck, J., The elliptic sinh-Gordon equation and the
construction of toroidal soap bubbles, 
Lecture Notes in Math., 1340  (1986), 275--301. 

\bibitem{} Stellmacher, J. L.,  Geometrische Deutung konform
invarianter Eigenschafter des Riemannschen
Raumes, Math. Ann., 123 (1951), 34--52.

\bibitem{} Stolz, S., Multiplicities of Dupin hypersurfaces, 
Preprint, Max-Planck-Institut f\"ur Mathematik, Bonn, 1997. 

\bibitem{} Str\"ubing, W.,  Symmetric submanifolds of Riemannian
manifolds, Math. Ann., 245 (1979), 37--44.

\bibitem{} S\"uss, W.,  Zur relativen Differentialgeometrie. V: \"Uber
Eihyperfl\"achen im $R^{m+1}$, T\^ohoku Math. J., 31 (1929), 202--209.

\bibitem{}  Szab\"o, Z. I.,  The Lichnerowicz conjecture on harmonic manifolds, J.
Differential Geometry, 31 (1990), 1--28. 

\bibitem{} Tachibana, S. and T. Kashiwada,  On a characterization of spaces of constant
holomorphic curvature in terms of geodesic hyperspheres, Kyungpook Math. J.,
13 (1973), 109--119.

\bibitem{} Takagi, R., On the principal curvatures of homogeneous
hyperspheres in a  sphere, Differential Geometry, in honor of K. Yano, Kinokuniya,
Tokyo, 469--481, 1971.

\bibitem{} Takagi, R.,  A class of hypersurfaces
with constant principal curvatures in a
sphere, J. Differential Geometry, 11 (1976), 225--233.

\bibitem{}  Takagi, R. and M. Takeuchi,  Degree of symmetric
Kaehlerian submanifolds of a complex
projective space, Osaka J. Math., 14 (1977), 501--518.

\bibitem{}   Takahashi, T.,  Minimal immersions of Riemannian manifolds, J. Math. Soc.
Japan, 18 (1966), 380--385.

\bibitem{}   Takahashi, T.,   Hypersurface with
parallel Ricci tensor in a space of constant
holomorphic sectional curvature, J. Math. Soc. Japan, 19 (1967), 199--204.

\bibitem{}   Takahashi, T.,  An isometric immersion
of a homogeneous Riemannian manifold of
dimension $3$ in the hyperbolic space, J. Math. Soc. Japan , 23 (1971), 649--661.

\bibitem{}   Takahashi, T.,  A note on Kaehlerian
hypersurfaces of spaces of constant
curvature, Kumanoto J. Sci. (Math.), 9 (1972), 21--24.

\bibitem{}   Takahashi, T.,  Isometric immersions of Riemannian homogeneous
manifolds, Tsukuba J. Math., 12 (1988), 231-233.

\bibitem{}  Takagi, R. and T. Takahashi, On the principal curvatures of homogeneous
hypersurfaces in a sphere, Differential
geometry (in honor of Kentaro Yano),  Kinokuniya, Tokyo, 469--481, 1972.

\bibitem{}   Takeuchi, M.,  Homogeneous K\"ahler
submanifolds in complex projective spaces,
Japan. J. Math., 4 (1978), 171--219.

\bibitem{}   Takeuchi, M.,   Parallel submanifolds of space forms, Manifolds
and Lie Groups, in honor of Y. Matsushima,
(J. Hano et al. eds.), Birkh\"auser,
Boston,  429-447, 1981.

\bibitem{}   Takeuchi, M.,   Stability of certain minimal submanifolds of compact Hermitian
symmetric spaces, T\^o\-hoku Math. J., 36 (1984), 293--314.

\bibitem{}   Takeuchi, M.,  Two-number of symmetric
 $R$-spaces, Nagoya Math. J. , 115 (1989), 43--46.
 
 \bibitem{}   Takeuchi, M., 
Proper Dupin hypersurfaces generated by
symmetric submanifolds,  Osaka J. Math.
28 (1991), 153--161. 

\bibitem{} Takeuchi, M. and S.
Kobayashi,{1968} Minimal
imbeddings of $R$-spaces, J. Differential
Geometry, 2, 203--215.

\bibitem{} Takeuchi, N.,{1987}
A closed surface of genus one in
 $E^3$ cannot contain seven circles through
each point,  Proc. Amer. Math.
Soc., 100 , 145--147.

\bibitem{} Tanaka, M. S.
,{1995} Stability of minimal
submanifolds in symmetric spaces, Tsukuba J.
Math., 19, 27--56.

\bibitem{} Tasaki, H.
,{1985} Certain minimal or homologically volume minimizing
submanifolds in compact symmetric spaces,
Tsukuba J. Math., 9, 117--131.

\bibitem{}   Takeuchi, M., Integral geometry under cut loci in compact
symmetric spaces,  Nagoya Math. J. 137 (1995), 33-53.

\bibitem{} Tazawa, Y., Construction of slant submanifolds, I ,
Bull. Inst. Math. Acad. Sinica, 22 (1994), 153--166.

\bibitem{} Tazawa, Y.,   Construction of slant
submanifolds, II, Bull. Soc. Math. Belg. (New Series), 1 (1994), 569--576.

\bibitem{} Thayer, E. C., Higher genus Chen-Gackstatter surfaces and
the Weierstrass representation for
infinite genus surfaces, Experimental Math., 4 (1995), 11--31. 

\bibitem{} Thi, D. C.,  Minimal real currents on compact Riemannian
manifolds, Math. USSR-Izv., 11 (1977), 807--820.

\bibitem{} Thorbergsson, G., Dupin hypersurfaces, Bull. London Math. Soc., 15 (1983),
493--498. 

\bibitem{} Thorbergsson, G., Highly connected taut
submanifolds, Math. Ann., 265 (1983), 399--405.

\bibitem{} Thorbergsson, G., Tight immersions of highly connected 
manifolds, Comment. Math. Helv., 61 (1986), 102--121. 

\bibitem{} Thorbergsson, G., Homogeneous spaces
without taut embeddings,  Duke Math. J., 57 (1988), 347--355.

\bibitem{} Thorbergsson, G.,  Tight analytic surfaces, Topology, 30 (1991), 423--428. 

\bibitem{} Thomsen, G., \"Uber konforme Geometrie I: Grundlagen der konformen
F\"achentheorie, Abh. Math. Sem. Univ. Hamburg, 1923, 31--56.

\bibitem{} Tojo, K., Normal homogeneous spaces admitting totally geodesic
hypersurfaces, J. Math. Soc. Japan, 49 (1997), 781--815.

\bibitem{} Tojo, K.,  Extrinsic hyperspheres of naturally reductive
homogeneous spaces, Tokyo J. Math., 20 (1997), 35--43.

\bibitem{} Toth, G., New construction for spherical minimal
immersions, Geom. Dedicata, 67 (1997), 187--196.

\bibitem{} Tomi, F., A finiteness results in the free boundary
value problem for minimal surfaces, Ann. Inst. Henri Poincar\'e, 3 (1986), 331--343.

\bibitem{}  Tomi, F. and A. J. Tromba, The index theorem for
minimal surfaces of higher genus, Mem. Amer. Math. Soc., 117 (1995), no. 560.

\bibitem{} Toubiana, \'E.,  On the minimal surfaces
of Riemann,  Comment. Math. Helv., 67 (1992), 546--570. 

\bibitem{} Toubiana, \'E.,Surfaces minimales non orientables de 
genre quelconque,  Bull. Soc. Math. France, 121 (1993), 183--195.

\bibitem{} Traizet, M.,  Construction de surfaces minimales en
recollant des surfaces de Scherk, Ann.
Inst. Fourier (Grenoble), 46 (1996), 1385--1442.
 
\bibitem{} Treibergs, A., Existence and convexity for hyperspheres of
 prescribed mean curvature, Ann. Scuola Norm. Sup. Pisa Cl. Sci., 12 (1985), 225--241. 

\bibitem{} Treibergs, A. and S. W. Wei, Embeddied hyperspheres
with prescribed mean curvature, J.
Differential Geometry, 18 (1983), 513--521.

\bibitem{} Tromba, A. J., On the number of simply connected minimal
surfaces spanning a curve, Mem. Amer. Math. Soc., 194, 1977.

\bibitem{} Tso, K., Convex hypersurfaces with prescribed
 Gauss-Kronecker curvature, J. Differential Geometry, 34 (1991), 389--410. 

\bibitem{} Tsukada, K., Helical geodesic immersions of compact rank one
symmetric spaces into spheres, Tokyo J. Math.,
6 (1983), 267--285.

\bibitem{} Tsukada, K., Parallel Kaehler submanifolds of Hermitian
 symmetric spaces, Math. Z. 190 (1985), 129--150. 

\bibitem{} Tsukada, K., Parallel submanifolds in a quaternion projective space,
Osaka J. Math., 22 (1985), 187--241.

\bibitem{} Tsukada, K.,  Parallel submanifolds of
Cayley plane, Sci. Rep. Niigata Univ. Ser. A, 21 (1985), 19--32.

\bibitem{} Tsukada, K.,  Totally geodesic submanifolds of Riemannian manifolds and
curvature-invariant subspaces, Kodai
Math. J.  19 (1996), 395--437.

\bibitem{} Tuzhilin, A. A.,  Morse-type indices for two-dimensional minimal surfaces in
$R^3$ and $H^3$, Math. USSR-Izv., 38 (1992), 575--598.

\bibitem{} Tysk, J., Eigenvalue estimates with applications to minimal surfaces,
Pacific J. Math., 128 (1987), 361--366.

\bibitem{} Tysk, J.,Finiteness of index and total scalar  curvature for minimal hypersurfaces, 
Proc. Amer. Math. Soc., 105 (1989), 429--435. 

\bibitem{} Umehara, M., Kaehler submanifolds of complex space forms ,
Tokyo J. Math., 10 (1987), 203--214.

\bibitem{} Umehara, M.,  Einstein Kaehler
submanifolds of a complex linear or hyperbolic
space, T\^ohoku Math. J., 39 (1987), 385--389.

\bibitem{} Umehara, M and K. Yamada,  A parametrization of the
Weierstrass formulae and perturbation of
complete minimal surfaces in $R^3$
into the hyperbolic $3$-space, J. Reine Angew. Math., 432 (1992), 93--116.


\bibitem{} Umehara, M and K. Yamada,  Complete surfaces of
constant mean curvature $1$ in the hyperbolic
$3$-space, Ann. Math., 137 (1993),
611--638.

\bibitem{} Umehara, M and K. Yamada,  Surfaces of constant
mean curvature $c$ in
$H^3(-c^2)$ with prescribed hyperbolic
Gauss map, Math. Ann., 304 (1996), 203--224.

\bibitem{} Umehara, M and K. Yamada,  A duality  on
CMC-$1$ surfaces in hyperbolic space, and a
hyperbolic analogue of the Osserman
inequality, Tsukuba J. Math., 21 (1997), 229--237 .

\bibitem{} Umehara, M and K. Yamada,  Geometry of surfaces of
constant mean curvature $1$ in the hyperbolic
$3$-space, Sugaku Expositions, 10 (1997), 41--55.

\bibitem{} Urakawa, H., Spectrum  of the Jacobi operator on a minimal
hypersurface of the Euclidean space,
Geometry and its applications (Yokohama, 1991), 257--268.

\bibitem{}  Urbano, F., Totally real minimal submanifolds of a complex projective space, Proc. Amer. Math. Soc., 93 (1985), 332--334.

\bibitem{}  Urbano, F.,  Negatively curved totally
real submanifolds, Math. Ann. , 273 (1986), 345--348.

\bibitem{}  Urbano, F.,  Totally real  submanifolds of a complex projective space, Geometry and Topology of Submanifolds, I, 198--208, 1989.

\bibitem{}  Urbano, F.,  Minimal surfaces with low index in the three-dimensional
sphere,  Proc. Amer. Math. Soc., 108 (1990), 989--992.

\bibitem{} Vaisman, I., Sympletic Geometry and Secondary
Characteristic Classes, Progress in  Math., Birkh\"u\-ser, 72, 1987.

\bibitem{} van Gemmeren, M., Total absolute curvature and tightness of
noncompact manifolds, Trans. Amer. Math. Soc., 348 (1986), 2413--2426.

\bibitem{} van Lindt, D. and L. Verstraelen,  Some axioms of
Einsteinian and conformally flat hypersurfaces, J. Differential Geometry,
16 (1981), 205--212.

\bibitem{} van Lindt, D. and L. Verstraelen,   A survey on axioms
of submanifolds in Riemannian and
Kaehlerian geometry, Colloq. Math.,
54 (1987), 193--213.

\bibitem{} Vanhecke, L., The axiom  of coholomorphic $(2p+1)$-spheres for
some almost Hermitian manifolds, Tensor (N.S.) , 30 (1976), 275--28.

\bibitem{} Vanhecke, L. and T. J. Willmore, Umbilical
hypersurfaces of Riemannian, Kaehler and
nearly Kaehler manifolds, J. Univ. Kuwait, 4 (1977), 1--8.

\bibitem{}  Verheyen, P.,   Submanifolds with geodesic normal sections are
helical , Rend. Sem. Mat. Univ. Politecn. Torino,
43 (1985), 511-527.

\bibitem{} Verstraelen, L., Curves and surfaces of finite Chen type,
Geometry and Topology of Submanifolds, III (1990), 304--311.

\bibitem{} Volkov, Ju. A. and S. M. Vladimirova, Isometric
immersions of the Euclidean plane in
Lobaceskii space, Mat. Zametki, 10 (1971), 327--332.
 
\bibitem{} Voss, A., Zur Theorie der Transformation 
quadratischer Differentialausdr\"ucke
und der Kr\"um\-mung h\"oherer
Mannigfaltigekeiten, Math. Ann., 16 (1880), 129--178.

\bibitem{} Voss, K.,  Einige differentialgeometrische Kongruenzs\"atze
f\"ur geschlos\-sene Fl\"achen und 
Hyperfl\"a\-chen , Math. Ann., 131 (1956), 180--218.

\bibitem{} Vrancken, L.,  Locally symmetric
submanifolds of the nearly Kaehler $S^6$,
Algebras Groups Geom., 5 (1988), 369--394. 

\bibitem{} Vrancken, L.,   Killing vector fields and Lagrangian
submanifolds of the nearly Kaehler $S^6$,  J. Math. Pure Appl. , 77 (1998), 631--645.

\bibitem{} Wallach, N. R.,  Minimal immersions of symmetric
spaces into spheres, in Symmetric Spaces, M. Dekker, 1--40, 1972.

\bibitem{} Wang, C. P., Surfaces in M\"obius geometry, Nagoya
Math. J., 125 (1992), 53-72.

\bibitem{} Wang, Q. M., On the topology of Clifford isoparametric
hypersurfaces, J. Differential Geometry, 27 (1988), 55--66.

\bibitem{} Wang, S. P. and S. W. Wei, Bernstein conjecture
in hyperbolic geometry, Seminar on
Minimal Submanifolds (ed. E. Bombieri),
Ann. Math. Studies, 103 (1983), 339--358.

\bibitem{} Wang, X. J., Existence of convex hypersurfaces with
 prescribed Gauss-Kronecker curvature,
Trans. Amer. Math. Soc., 348 (1996), 4501--4524. 

\bibitem{} Wei, F., Some existence and uniqueness theorems for doubly periodic minimal
surfaces,  Invent. Math., 109 (1992), 113--136.

\bibitem{} Wei, F.,  Adding handles to the
Riemann examples, preprint, 1995.

\bibitem{}  Wei, S. W., Plateau's problem in symmetric spaces,
Nonlinear Anal., 12 (1988), 749--760.

\bibitem{}  Weiner, J. L., On a problem of Chen, Willmore, et al., Indiana Univ. 
Math. J., 288 (1979), 19--35.

\bibitem{}  Weingarten, J., \"Uber eine Klasse aufeinander
abwickelbarer Fl\"achen, J. Reine Angew. Math., 59 (1861), 296--310.

\bibitem{} Weitsman, A. and F. Xavier, Some function theoretic
properties of the Gauss map for hyperbolic
complete minimal surfaces, Michigan Math. J., 34 (1987), 275--283.

\bibitem{}  Wente, H., Counter-example to the Hopf
conjecture,  Pacific J. Math., 121 (1986), 193--244.

\bibitem{} Weyl, H., \"Uber die Bestimmung einer geschlossen
konvexen Fl\"ache durch ihr Linielement, Z\"urich Naturf. Ges., 61 (1916), 40--72.

\bibitem{} White, B., Complete surfaces of finite total curvature, J. Differential Geometry, 26 (1987), 315--326; Correction ,
ibid, 28, 359--360.

\bibitem{} Whitney, H., Differentiable manifolds, Ann.
Math., 37 (1936), 645--680.

\bibitem{} Whitney, H., On regular closed
vurves in the plane, Compos. Math., 4 (1937), 276--284.
 
\bibitem{}  White, J. H.,  A global invariant of conformal mappings in space, Proc.
Amer. Math. Soc. , 38 (1973), 162--164.

\bibitem{} Willmore, T. L.,  Mean curvature of immersed surfaces,
An. Sti. Univ. ``Al. I. Cuza'' Iasi, Sec. I. a Mat. (N.S.),
14 (1968), 99--103.

\bibitem{} Willmore, T. L.,  Mean curvature of Riemannian immersions,
J. London Math. Soc., 3 (1971), 307--310.

\bibitem{} Willmore, T. L.,  Tight immersions and total absolute curvature,
Bull. London Math. Soc., 38 (1973), 129-151.

\bibitem{} Willmore, T. L., Total curvature in Riemannian geometry, Ellis
Horwood Limited, London, 1982.

\bibitem{}  Wintgen, P.,  On the total curvature of
surfaces in $E^4$, Colloq. Math., 39 (1978), 289--296.

\bibitem{}  Wintgen, P.,   Sur l'in\'egalit\'e de Chen-Willmore, C. R. Acad. Sc. Paris, 27 (1979), 19--35.

\bibitem{} Wirtinger, W., Eine Determinantenidentit\"at und ihre Anwendung
auf analytische Gebilde in Euclidischer und Hermitischer Massbestimmung, Monatsh.
Math. Phys., 44 (1936), 343--365.

\bibitem{}  Wohlgemuth, M., Minimal surfaces of higher genus with
finite total curvature, Arch. Rational Mech. Anal., 137 (1997), 1--25.

\bibitem{} Wolf, J. A.,Elliptic spaces in Grassmann manifold,
Illinois J. Math., 7 (1963), 447--462.

\bibitem{} Wolfson, J. G., Minimal surfaces in Kaehler surfaces and
Ricci curvature, J. Differential Geometry, 29 (1989), 281--294.

\bibitem{} Wu, B., A finiteness theorem for isoparametric
hypersurfaces, Geom. Dedicata,
50 (1994), 247--250.

\bibitem{} Xavier, F., The Gauss map of a complete nonflat minimal
surface cannot omit $7$ points of the
sphere, Ann. of Math., 113 (1981), 211--214.

\bibitem{}  Xia, C., A pinching theorem for minimal submanifolds 
in a sphere , Arch. Math., 57 (1991), 307--312.

\bibitem{}  Xia, C., A sphere theorem for
submanifolds in a manifold with pinched
positive curvature, Monats. Math. ,
124 (1997), 365--368.

\bibitem{} Xu, Y., Symmetric minimal surfaces in $R^3$, Pacific J.
Math., 171 (1995), 275--296.

\bibitem{} Yamaguchi, S. and M. Kon, Kaehler manifolds
satisfying the axiom of anti-invariant
2-spheres, Geom. Dedicata, 7 (1978), 403--406.

\bibitem{}  Yamaguchi, S., Nemoto, H. and N. Kawabata, Extrinsic
spheres in a K\"ahler manifold, Michigan Math. J., 31 (1984), 15--19.

\bibitem{} Yang, H. C. and Q. M. Cheng, Chern's Conjecture on
minimal hypersurfaces ,  Math.
Z., 227 (1998), 377--390.

\bibitem{} Yang, J., On slant surfaces with constant mean curvature
in $C^2$, J. Geometry, 59 (1997), 184--201.

\bibitem{} Yang, K., Complete minimal surfaces of finite total
curvature, Kluwer Academic Publ., Dordrecht, 1994.

\bibitem{} Yano, K. and B. Y. Chen, Minimal submanifolds
of a higher dimensional sphere, Tensor, (N.S.), 22 (1971), 369-373.

\bibitem{} Yano K. and S. Ishihara, Submanifolds with parallel mean
 curvature vector, J. Differential Geometry , 6 (1971), 95--118.

\bibitem{} Yano K. and M. Kon,  Generic submanifolds, Ann. Mat.
Pura Appl., 123 (1980), 59--92.

\bibitem{} Yano K. and M. Kon, $CR$-submanifolds of
Kaehlerian and Sasakian manifolds,
Birkh\"auser, Boston-Basel, 1983.

\bibitem{} Yano K. and I. Mogi, On real representations of
Kaehlerian manifolds, Ann. Math., 61 (1955), 170--189.

\bibitem{} Yano K. and Y. Muto,  Sur la th\'eorie des espaces \`a
connection conforme normale et la
g\'eom\'etrie conforme des espaces de
Riemann, J. Fac. Sci. Imp. Univ. Tokyo, 4 (1941), 117--169.

\bibitem{} Yau, S. T.,  Submanifolds with constant mean curvature,
I, Amer. J. Math., 96 (1974), 346--366.

\bibitem{} Yau, S. T., Submanifolds with constant mean curvature,
II, Amer. J. Math., 97 (1975), 76--100.

\bibitem{} Yau, S. T.,   Calabi's conjecture and some new results in
algebraic geometry, Proc. Nat. Acad. Sci.
U.S.A., 74 (1977), 1798--1799.

\bibitem{} Zhao, Q., The indices, the nullities and the stability of
totally geodesic submanifolds in the complex quadratic hypersurfaces:
$Q_m=SO(m+2)/SO(m)\times SO(2)$, Proc.
Amer. Math. Soc. , 124 (1996), 2501--2512.

\bibitem{} Zheng, F, Isometric embedding of Kaehler manifolds
 with nonpositive sectional curvature,
Math. Ann., 304 (1996), 769--784.

\bibitem{} Zhang, S. P., On complete minimal immersion $\chi:
RP^2-\{a,b\}\to E^3$ with total curvature $-10\pi$,  Lecture Notes in Math.,
Springer-Verlag, 1369 (1989), 339--350.

\bibitem{} Zheng, Y., Submanifolds with flat normal bundle,
Geom. Dedicata, 67 (1997), 295--300.

\end{thebibliography}
 \end{document}